Centro de Investigación y de Estudios Avanzados del I.P.N.

Unidad Zacatenco

Departamento de Matemáticas

# Álgebras de Rees, Subanillos Monomiales y Problemas de Optimización Lineal

Tesis que presenta

**Luis Alfredo Dupont García**

para obtener el Grado de
Doctor en Ciencias
en la Especialidad de
Matemáticas

Director de Tesis: Dr. Rafael Heraclio Villarreal Rodríguez

México, D.F. Junio del 2010.

CENTER FOR RESEARCH AND ADVANCED
STUDIES OF THE
NATIONAL POLYTHECHNIC INSTITUTE

Campus Zacatenco

Department of Mathematics

# Rees algebras, Monomial Subrings and Linear Optimization Problems

A DISSERTATION PRESENTED BY

**Luis Alfredo Dupont García**

TO OBTAIN THE DEGREE OF
DOCTOR IN SCIENCE
IN THE SPECIALITY OF
MATHEMATICS

Thesis Advisor: Dr. Rafael Heraclio Villarreal Rodríguez

# Resumen

En esta tesis estamos interesados en estudiar propiedades algebraicas de álgebras monomiales, que pueden vincularse a estructuras combinatorias, por ejemplo a gráficas y conglomerados, y a problemas de optimización. Una meta aquí es establecer puentes entre el álgebra conmutativa, combinatoria y optimización. Estudiamos la normalidad y la propiedad Gorenstein—así como el módulo canónico y el $a$-invariante—de álgebras de Rees y subanillos que surgen de problemas de optimización lineal. Particularmente, estudiamos propiedades algebraicas de ideales de aristas y álgebras asociadas a conglomerados uniformes con la propiedad de máximo-flujo mínimo-corte o con la propiedad del empacamiento. También estudiamos propiedades algebraicas de las álgebras de Rees simbólicas de ideales de aristas de gráficas, los ideales de aristas de conglomerados de clanes de gráficas de comparabilidad, y anillos de Stanley-Reisner.

A continuación describimos el contenido de esta tesis. El objetivo del Capítulo 1 es estudiar propiedades de redondeo entero de sistemas lineales para obtener información sobre las propiedades algebraicas de las álgebras de Rees y de los subanillos monomiales, y viceversa. Estudiamos los módulos canónicos, $a$-invariantes, y la propiedades de normalidad y Gorenstein, de los subanillos monomiales que surgen de sistemas con la propiedad del redondeo entero. Relacionamos las propiedades algebraicas de las álgebras de Rees y de los subanillos monomiales con la propiedad del redondeo entero presentando un teorema de dualidad. La normalidad de un ideal monomial es expresada en términos del poliedro bloqueador y de la propiedad de descomposición entera. Para los ideales de aristas de conglomerados esta propiedad determina totalmente su normalidad. Para los sistemas que se surgen de los clanes de gráficas perfectas se obtienen expresiones explícitas para el módulo canónico y el $a$-invariante.

En el Capítulo 2 estudiamos propiedades algebraicas y combinatorias de ideales y álgebra de conglomerados uniformes con la propiedad de máximo-flujo mínimo-corte. Sean $\mathcal{C}$ un conglomerado uniforme y $A$ su matriz de incidencia. Bajo ciertas condiciones probamos que $\mathcal{C}$ es crítica por vértices. Si $\mathcal{C}$ satisface la propiedad de máximo-flujo mínimo-corte, probamos que $A$ se diagonaliza sobre los números enteros en la matriz identidad. Suponiendo que $\mathcal{C}$ tiene un emparejamiento perfecto y su número de cubierta es igual a 2, probamos que la propiedad de empacamiento implica que $A$ se diagonaliza sobre los números enteros en la matriz identidad. Si la matriz $A$ es balanceada, demostramos que cualquier triangulación regular del cono generado por las columnas de $A$ es unimodular. Presentamos algunos ejemplos que demuestran que nuestros resultados solamente son válidos para conglomerados uniformes. Estos resultados se relacionan con las propiedades de ser libre de torsión, con la normalidad de ciertos anillos inducidos y con grupos abelianos. También se relacionan con la teoría de bases de Gröbner de ideales tóricos y anillos de Ehrhart.

El Capítulo 3 se dedica a la construcción de conglomerados Cohen-Macaulay con propiedades de optimización combinatoria y a examinar paralelizaciones de conglomerados y sus ideales de aristas. Sea $\mathcal{C}$ un conglomerado uniforme. Demostramos que si $\mathcal{C}$ satisface la propiedad de empacamiento (resp. la propiedad de máximo-flujo mínimo-corte), entonces existe un conglomerado Cohen-Macaulay uniforme $\mathcal{C}_1$ que satisface la propiedad de empacamiento (resp. la propiedad de máximo-flujo mínimo-corte) tal que $\mathcal{C}$ es un menor de $\mathcal{C}_1$. Para ideales de aristas de conglomerados arbitrarios probamos que la normalidad es cerrada bajo paralelizaciones. A continuación mostramos algunas aplicaciones a ideales de aristas y conglomerados que están relacionados con una conjetura de Conforti y Cornuéjols así como a problemas de máximo-flujo mínimo-corte.

Luego, en el Capítulo 4 centramos nuestra atención en los ideales de aristas asociados a conglomerados de clanes en gráficas de comparabilidad. Sean $(P, \prec)$ un conjunto finito parcialmente ordenado y $G$ su gráfica de comparabilidad. Si $\mathrm{cl}(G)$ es el conglomerado de clanes maximales de $G$, probamos que $\mathrm{cl}(G)$ cumple la propiedad de máximo-flujo mínimo-corte y que su ideal de aristas es normalmente libre de torsión. Demostramos que los ideales de aristas de conglomerados uniformes admisibles completos son normalmente libres de torsión.

Finalmente, en el Capítulo 5 consideramos álgebras de Rees simbólicas, cubiertas por vértice y representaciones irreducibles de conos de Rees. La relación entre las caretas del cono de Rees de conglomerados y las $b$-cubiertas (por vértices) irreducibles es examinada. Sean $G$ una gráfica simple y $I_c(G)$ su ideal de cubiertas por vértices. Damos una descripción gráfica teórica de las $b$-cubiertas irreducibles de $G$, es decir, describimos los generadores minimales del álgebra de Rees simbólica de $I_c(G)$. Como una aplicación, recuperamos una descripción explícita de el cono de Rees de una gráfica. Luego, estudiamos las $b$-cubiertas irreducibles del bloqueador de $G$, es decir, estudiamos los generadores minimales del álgebra de Rees simbólica de el ideal de aristas de $G$. Damos una descripción gráfica teórica de las $b$-cubiertas irreducibles del bloqueador de $G$. Se demuestra que están en correspondencia uno a uno con las subgráficas inducidas irreducibles de $G$. Como un subproducto obtenemos un método, utilizando bases de Hilbert, para obtener todas las subgráficas inducidas irreducibles de $G$. Las gráficas irreducibles son estudiadas. Mostramos cómo construir gráficas irreducibles y damos un método para construir $b$-cubiertas irreducibles del bloqueador de $G$.

# Abstract

In this thesis we are interested in studying algebraic properties of monomial algebras, that can be linked to combinatorial structures, such as graphs and clutters, and to optimization problems. A goal here is to establish bridges between commutative algebra, combinatorics and optimization. We study the normality and the Gorenstein property—as well as the canonical module and the $a$-invariant—of Rees algebras and subrings arising from linear optimization problems. In particular, we study algebraic properties of edge ideals and algebras associated to uniform clutters with the max-flow min-cut property or the packing property. We also study algebraic properties of symbolic Rees algebras of edge ideals of graphs, edge ideals of clique clutters of comparability graphs, and Stanley-Reisner rings.

The contents of this thesis are as follows. The aim of Chapter 1 is to study integer rounding properties of linear systems to gain insight about the algebraic properties of Rees algebras and monomial subrings, and viceversa. We study canonical modules, $a$-invariants, and the normality and Gorenstein properties, of monomial subrings arising from systems with the integer rounding property. We relate the algebraic properties of Rees algebras and monomial subrings with integer rounding properties by presenting a duality theorem. The normality of a monomial ideal is expressed in terms of blocking polyhedra and the integer decomposition property. For edge ideals of clutters these properties completely determine their normality. For systems arising from cliques of perfect graphs explicit expressions for the canonical module and the $a$-invariant are given.

In Chapter 2 we study algebraic and combinatorial properties of ideals and algebras of uniform clutters with the max-flow min-cut property. Let $\mathcal{C}$ be a uniform clutter and let $A$ be its incidence matrix. Under certain conditions we prove that $\mathcal{C}$ is vertex critical. If $\mathcal{C}$ satisfies the max-flow min-cut property, we prove that $A$ diagonalizes over the integers to an identity matrix. Assuming that $\mathcal{C}$ has a perfect matching and its covering number is equal to 2, we prove that the packing property implies that $A$ diagonalizes over the integers to an identity matrix. If the matrix $A$ is balanced, we show that any regular triangulation of the cone generated by the columns of $A$ is unimodular. Some examples are presented to show that our results only hold for uniform clutters. These results are related to the normality and torsion freeness of blowup rings and to abelian groups. They are also related to the theory of Gröbner bases of toric ideals and to Ehrhart rings.

Chapter 3 is devoted to the construction of Cohen-Macaulay clutters with combinatorial optimization properties and to the examination of parallelizations of clutters and their edge ideals. Let $\mathcal{C}$ be a uniform clutter. We prove that if $\mathcal{C}$ satisfies the packing property (resp. max-flow min-cut property), then there is a uniform Cohen-Macaulay clutter $\mathcal{C}_1$ satisfying the packing property (resp. max-flow min-cut property) such that $\mathcal{C}$ is a minor of $\mathcal{C}_1$. For arbitrary edge ideals of clutters we prove that the normality property is closed under parallelizations. Then we show some applications to edge ideals and clutters which are related to a conjecture of Conforti and Cornuéjols and to max-flow min-cut problems.

Then in Chapter 4 we focus our attention on edge ideals of clique clutters of comparability graphs. Let $(P, \prec)$ be a finite poset and let $G$ be its comparability graph. If $\mathrm{cl}(G)$ is the clutter of maximal cliques of $G$, we prove that $\mathrm{cl}(G)$ satisfies the max-flow min-cut property and that its edge ideal is normally torsion free. We prove that edge ideals of complete admissible uniform clutters are normally torsion free.

Finally in Chapter 5 we consider symbolic Rees algebras, vertex covers and irreducible representations of Rees cones. The relation between facets of Rees cones of clutters and irreducible $b$-vertex covers is examined. Let $G$ be a simple graph and let $I_c(G)$ be its ideal of vertex covers. We give a graph theoretical description of the irreducible $b$-vertex covers of $G$, i.e., we describe the minimal generators of the symbolic Rees algebra of $I_c(G)$. As an application we recover an explicit description of the edge cone of a graph. Then we study the irreducible $b$-vertex covers of the blocker of $G$, i.e., we study the minimal generators of the symbolic Rees algebra of the edge ideal of $G$. We give a graph theoretical description of the irreducible binary $b$-vertex covers of the blocker of $G$. It is shown that they are in one to one correspondence with the irreducible induced subgraphs of $G$. As a byproduct we obtain a method, using Hilbert bases, to obtain all irreducible induced subgraphs of $G$. Irreducible graphs are studied. We show how to build irreducible graphs and give a method to construct irreducible $b$-vertex covers of the blocker of $G$.

# Reconocimientos

Este trabajo de tesis se realizó dentro del proyecto de investigación:

*Ideales Tóricos, Normalizaciones y Funciones de Hilbert*. Fondo Sectorial SEP-CONACyT. Clave: 49251-F.

Agradezco al *Consejo Nacional de Ciencia y Tecnología* la beca que me otorgó para la realización de mis estudios de doctorado.

Junio del 2010.
Ciudad de México — Luis Alfredo Dupont García

# Contents

# Introduction

Let $R = K[x_1, \ldots, x_n]$ be a polynomial ring over a field $K$ and let $v_1, \ldots, v_q$ be the column vectors of a matrix $A = (a_{ij})$ whose entries are non-negative integers. For technical reasons, we shall always assume that the rows and columns of the matrix $A$ are different from zero. As usual we use the notation $x^a := x_1^{a_1} \cdots x_n^{a_n}$, where $a = (a_1, \ldots, a_n) \in \mathbb{N}^n$.

The *ideals* and *monomial algebras* considered here are: (a) the *monomial ideal*:

$$I = (x^{v_1}, \ldots, x^{v_q}) \subset R$$

generated by $F = \{x^{v_1}, \ldots, x^{v_q}\}$, (b) the *integral closure* of the ideal $I$:

$$\overline{I} = (\{x^a \in R |\, \exists\, p \in \mathbb{N} \setminus \{0\}; (x^a)^p \in I^p\});$$

(c) the *Rees algebra*:

$$R[It] := R \oplus It \oplus \cdots \oplus I^i t^i \oplus \cdots \subset R[t],$$

where $t$ is a new variable, (d) the *symbolic Rees algebra*:

$$R_s(I) = R \oplus I^{(1)} t \oplus \cdots \oplus I^{(i)} t^i \oplus \cdots \subset R[t],$$

where $I$ is a square-free monomial ideal and $I^{(i)}$ is its *ith* symbolic power, (e) the *extended Rees algebra*

$$R[It, t^{-1}] := R[It][t^{-1}] \subset R[t, t^{-1}],$$

(f) the *monomial subring*

$$K[F] = K[x^{v_1}, \ldots, x^{v_q}] \subset R$$

spanned by $F = \{x^{v_1}, \ldots, x^{v_q}\}$, (g) the *integral closure* of $K[F]$ in its field of fractions:

$$\overline{K[F]} = K[\{x^a |\, a \in \mathbb{Z}\mathcal{A} \cap \mathbb{R}_+\mathcal{A}\}], \tag{1}$$

where $\mathcal{A} = \{v_1, \ldots, v_q\}$, $\mathbb{Z}\mathcal{A}$ is the subgroup of $\mathbb{Z}^n$ spanned by $\mathcal{A}$, and $\mathbb{R}_+\mathcal{A}$ is the cone generated by $\mathcal{A}$, (h) the *homogeneous monomial subring*

$$K[Ft] = K[x^{v_1}t, \ldots, x^{v_q}t] \subset R[t]$$

spanned by $Ft$, (i) the *homogeneous monomial subring*

$$K[Ft \cup \{t\}] = K[x^{v_1}t, \ldots, x^{v_q}t, t] \subset R[t]$$

spanned by $Ft \cup \{t\}$, (j) the *homogeneous monomial subring*

$$S = K[x^{w_1}t, \ldots, x^{w_r}t] \subset R[t],$$

where $w_1, \ldots, w_r$ is the set of vectors $\alpha \in \mathbb{N}^n$ such that $0 \leq \alpha \leq v_i$ for some $i$, and (k) the *Ehrhart ring*

$$A(P) = K[\{x^a t^i \,|\, a \in \mathbb{Z}^n \cap iP; i \in \mathbb{N}\}] \subset R[t]$$

of a lattice polytope $P$.

Recall that $K[F]$ is called *integrally closed* or *normal* if $K[F] = \overline{K[F]}$. The subring $K[F]$ equals $K[\mathbb{N}\mathcal{A}]$, the semigroup ring of the subsemigroup $\mathbb{N}\mathcal{A}$ of $\mathbb{N}^n$ generated by $\mathcal{A}$. Thus, from Eq. (1), $K[F]$ is normal if and only if

$$\mathbb{N}\mathcal{A} = \mathbb{Z}\mathcal{A} \cap \mathbb{R}_+\mathcal{A}.$$

The description of the integral closure given in Eq. (1) can of course be applied to any of the monomial subrings considered here. In particular if $\mathcal{A}'$ is the set

$$\mathcal{A}' = \{e_1, \ldots, e_n, (v_1, 1), \ldots, (v_q, 1)\},$$

where $e_i$ is the *ith* unit vector, then

$$R[It] = K[x_1, \ldots, x_n, x^{v_1}t, \ldots, x^{v_q}t] = K[\mathbb{N}\mathcal{A}'],$$

$\mathbb{Z}\mathcal{A}' = \mathbb{Z}^{n+1}$, and $R[It]$ is normal if and only if $\mathbb{N}\mathcal{A}' = \mathbb{Z}^{n+1} \cap \mathbb{R}_+\mathcal{A}'$. A dual characterization of the normality of $R[It]$ will be given in Proposition 1.2.8. Recall that the *integral closure* of $I^i$, denoted by $\overline{I^i}$, is the ideal of $R$ given by

$$\overline{I^i} = (\{x^a \in R \,|\, \exists\, p \in \mathbb{N} \setminus \{0\}; (x^a)^p \in I^{pi}\}),$$

see for instance [87, Proposition 7.3.3]. The ideal $I$ is called *normal* if $I^i = \overline{I^i}$ for all $i$. If $\overline{I} = I$, the ideal $I$ is called *integrally closed*. Note that the ideal $I$ is normal if and only if $R[It]$ is normal [87].

Recall that the Ehrhart ring $A(P)$ is always normal [10]. A finite set $\mathcal{A} \subset \mathbb{R}^n$ is called an *Hilbert basis* if $\mathbb{N}\mathcal{A} = \mathbb{R}_+\mathcal{A} \cap \mathbb{Z}^n$. An *integral Hilbert basis* is a Hilbert basis consisting of integral vectors only. Note that if $\mathcal{A} = \{v_1, ..., v_q\}$ is an integral Hilbert basis, then the ring $K[x^{v_1}, ..., x^{v_q}]$ is normal.

Next we introduce a combinatorial stucture which generalizes the notion of a graph [50], and introduce its algebraic counterpart which generalizes the notion of an edge ideal of a graph [85]. These two notions will play a major role.

A *clutter* $\mathcal{C}$ with finite vertex set $X = \{x_1, ..., x_n\}$ is a family of subsets of $X$, called edges, none of which is included in another. Let $\mathcal{C}$ be a clutter. The set of vertices and edges of $\mathcal{C}$ are denoted by $V(\mathcal{C})$ and $E(\mathcal{C})$ respectively. The *edge ideal* of $\mathcal{C}$, denoted by $I(\mathcal{C})$, is the ideal of $R$ generated by all monomials $x_e = \prod_{x_i \in e} x_i$ such that $e \in E(\mathcal{C})$. The map

$$\mathcal{C} \longmapsto I(\mathcal{C})$$

gives a one to one correspondence between the family of clutters and the family of square-free monomial ideals. Edge ideals of clutters also correspond to simplicial complexes via the Stanley-Reisner correspondence [75]. The *Alexander dual* of $I(\mathcal{C})$ is the ideal of $R$ given by

$$I(\mathcal{C})^\vee = \cap_{e \in E}(e),$$

where $E = E(\mathcal{C})$ is the edge set of $\mathcal{C}$ and $(e)$ is the prime ideal generated by the vertices in $e$. The Alexander dual of $I(\mathcal{C})$ is also denoted by $I_c(\mathcal{C})$ and is also called the *ideal of covers* of $\mathcal{C}$. The *dual* $I^*$, of an edge ideal $I$, is the ideal of $R$ generated by all monomials $x_1 \cdots x_n / x_e$ such that $e$ is an edge of $\mathcal{C}$. The *incidence matrix* of $\mathcal{C}$ is the vertex-edge matrix whose columns are the characteristic vectors of the edges of $\mathcal{C}$. Notice that if $I = I(\mathcal{C})$ is minimally generated by $F = \{x^{v_1}, \ldots, x^{v_q}\}$, then $A$ turns out to be the incidence matrix of $\mathcal{C}$.

Below we introduce the main results of this thesis. The techniques used in this dissertation come from commutative algebra [10, 84], combinatorial optimization [70, 71] and graph theory [20, 50]

**Contents of Chapter 1** The aim of this chapter is to study max-flow min-cut properties of clutters and integer rounding properties of various systems of linear inequalities—and their underlying polyhedra—to gain insight about the algebraic properties of these algebras and viceversa. Systems with integer rounding properties and clutters with the max-flow min-cut property come from linear optimization problems [70, 71].

Let us recall some of the linear optimization properties that are considered here. The linear system $x \geq 0; xA \geq \mathbf{1}$ has the *integer rounding property* if

$$\max\{\langle y, \mathbf{1}\rangle \,|\, y \geq 0; Ay \leq w; y \in \mathbb{N}^q\} = \lfloor \max\{\langle y, \mathbf{1}\rangle \,|\, y \geq 0; Ay \leq w\}\rfloor$$

for each integral vector $w$ for which the right hand side is finite. Here $\langle\, , \rangle$ is the standard inner product and $\mathbf{1}$ denotes the vector $(1, \ldots, 1)$. The linear system $x \geq 0; xA \leq \mathbf{1}$ has the *integer rounding property* if

$$\lceil \min\{\langle y, \mathbf{1}\rangle \,|\, y \geq 0;\ Ay \geq a\}\rceil = \min\{\langle y, \mathbf{1}\rangle \,|\, Ay \geq a;\ y \in \mathbb{N}^q\}$$

for each integral vector $a$ for which $\min\{\langle y, \mathbf{1}\rangle \,|\, y \geq 0;\ Ay \geq a\}$ is finite. The system $xA \leq \mathbf{1}$ is said to have the *integer rounding property* if

$$\lceil \min\{\langle y, \mathbf{1}\rangle \,|\, y \geq 0;\ Ay = a\}\rceil = \min\{\langle y, \mathbf{1}\rangle \,|\, Ay = a;\ y \in \mathbb{N}^q\}$$

for each integral vector $a$ for which $\min\{\langle y, \mathbf{1}\rangle \,|\, y \geq 0;\ Ay = a\}$ is finite.

A rational polyhedron $Q$, i.e., with rational vertices, is said to have the *integer decomposition property* if for each natural number $k$ and for each integer vector $a$ in $kQ$, $a$ is the sum of $k$ integer vectors in $Q$. Recall that $kQ$ is equal to $\{ka \,|\, a \in Q\}$.

In Section 1.1 we study the normality of general monomial ideals. We are able to characterize this property in terms of blocking polyhedra and the integer decomposition property (see Theorem 1.1.2). As an interesting consequence using a result of Baum and Trotter [4] we describe the normality of a monomial ideal in terms of the integer rounding property:

**Corollary 1.1.8** *The ideal $I$ is normal if and only if the system $xA \geq \mathbf{1}; x \geq 0$ has the integer rounding property.*

There are already some characterizations of the normality of monomial ideals [29, Theorem 4.4]. There are also useful membership tests based on linear optimization to decide whether a given monomial lies in the integral closure of a monomial ideal [19, Proposition 3.5].

We use the theory of blocking and antiblocking polyhedra [4, 37, 70] to describe when the systems

$$x \geq 0;\ xA \leq \mathbf{1}, \quad x \geq 0; xA \geq \mathbf{1}, \quad xA \leq \mathbf{1},$$

have the *integer rounding property* in terms of the normality of the monomial algebras considered here.

One of the main results of Section 1.2 is:

**Theorem 1.2.5** *The system $x \geq 0$; $xA \leq \mathbf{1}$ has the integer rounding property if and only if the subring $S = K[x^{w_1}t, \ldots, x^{w_r}t]$ is normal.*

We present a duality between the integer rounding property of the systems $x \geq 0; xA \geq \mathbf{1}$ and $x \geq 0; xA^* \leq \mathbf{1}$ valid for matrices with entries in $\{0,1\}$, where $a_{ij}^* = 1 - a_{ij}$ is the $ij$-entry of $A^*$. This duality is extended to a duality between the monomial subrings associated to $A$ and $A^*$.

The duality theorem—which is a major result of Section 1.2—is:

**Theorem 1.2.11** *Let $A$ be the incidence matrix of a clutter. If $v_i^* = \mathbf{1} - v_i$ and $A^*$ is the matrix with column vectors $v_1^*, \ldots, v_q^*$, then the following are equivalent:*

(a) *$R[I^*t]$ is normal, where $I^* = (x^{v_1^*}, \ldots, x^{v_q^*})$.*

(b) *$S = K[x^{w_1}t, \ldots, x^{w_r}t]$ is normal.*

(c) *$\{-e_1, \ldots, -e_n, (v_1, 1), \ldots, (v_q, 1)\}$ is a Hilbert basis.*

(d) *$x \geq 0; xA^* \geq \mathbf{1}$ has the integer rounding property.*

(e) *$x \geq 0; xA \leq \mathbf{1}$ has the integer rounding property.*

Then we present some interesting consequences of this duality. First of all we recover one of the main results of [91] showing that if

$$P = \{x | x \geq 0; xA \leq \mathbf{1}\}$$

is an integral polytope, i.e., $P$ has only integer vertices, and $A$ is a $\{0,1\}$-matrix, then the Rees algebra $R[I^*t]$ is normal (see Corollary 1.2.14). This result is related to perfect graphs. Indeed if $P$ is integral, then $v_1^*, \ldots, v_q^*$ correspond to the maximal cliques (maximal complete subgraphs) of a perfect graph [14, 62]. Second we show that if $A$ is the incidence matrix of the collection of basis of a matroid, then all systems

$$x \geq 0; xA \geq \mathbf{1}, \quad x \geq 0; xA^* \geq \mathbf{1}, \quad x \geq 0; xA \leq \mathbf{1}, \quad x \geq 0; xA^* \leq \mathbf{1}$$

have the integer rounding property (see Corollary 1.2.15). Third we show that if $A$ is the incidence matrix of a graph, then $R[It]$ is normal if and only if $R[I^*t]$ is normal (see Corollary 1.3.4). We give an example to show that this result does not extend to arbitrary uniform clutters (see Example 1.3.5), i.e., to clutters having all edges of the same size. If $A$ is the incidence matrix of a graph $G$, we characterize when $I^*$ is the Alexander dual of the edge ideal of the complement of $G$ (see Proposition 1.3.6). If $G$ is a triangle-free graph, we show a duality between the normality of $I = I(G)$ and that of the Alexander dual of the edge ideal of the complement of $G$ (see Corollary 1.3.7). We show an example of an edge ideal of a graph whose Alexander dual is not normal (see Example 1.3.8). In [91] it is shown that this is never the case if the graph is perfect, i.e., the Alexander dual of the edge ideal of a perfect graph is always normal.

A clutter $\mathcal{C}$ satisfies the *max-flow min-cut* (MFMC) property if both sides of the LP-duality equation

$$\min\{\langle a, x\rangle \,|\, x \geq 0; xA \geq \mathbf{1}\} = \max\{\langle y, \mathbf{1}\rangle \,|\, y \geq 0; Ay \leq a\}$$

have integer optimum solutions $x$, $y$ for each non-negative integral vector $a$.

We recover one of the main results of [45] showing that if $A$ is the incidence matrix of a clutter $\mathcal{C}$, then $\mathcal{C}$ satisfies the max-flow min-cut property if and only if the set covering polyhedron

$$Q(A) = \{x \,|\, x \geq 0; xA \geq \mathbf{1}\}$$

is integral and $R[It]$ is normal (see Corollary 1.2.17).

Another main result of Section 1.2 is:

**Theorem 1.2.20** *If the system $xA \leq \mathbf{1}$ has the integer rounding property, then $K[F]$ is normal and $\mathbb{Z}^n/\mathbb{Z}\mathcal{A}$ is a torsion-free group. The converse holds if $|v_i| = d$ for all $i$. Here $|v_i| = \langle v_i, \mathbf{1}\rangle$.*

As a consequence of this result we prove: (i) If $A$ is the incidence matrix of a connected graph $G$, then the system $xA \leq \mathbf{1}$ has the integer rounding property if and only if $G$ is a bipartite graph (see Corollary 1.2.21), and (ii) Let $A$ be the incidence matrix of a clutter $\mathcal{C}$. If $\mathcal{C}$ is uniform, i.e., all its edges have the same size, and $\mathcal{C}$ has the max-flow min-cut property, then the system $xA \leq \mathbf{1}$ has the integer rounding property (see Corollary 1.2.22).

If $A$ is the incidence matrix of a bipartite graph, a nice result of Section 1.3 shows that the system $x \geq 0; xA \leq \mathbf{1}$ has the integer rounding property if and only if the extended Rees algebra $R[It, t^{-1}]$ is normal.

Let $G$ be a graph with incidence matrix $A$. If $G$ is a connected graph, in Theorem 1.3.2 we are able to prove that the system $x \geq 0; xA \leq \mathbf{1}$ has the integer rounding property if and only if the induced subgraph of the vertices of any two vertex disjoint odd cycles of $G$ is connected. Other equivalent descriptions of this property are also presented (see Theorem 1.3.3).

Before stating the main results of Section 1.5, we need to introduce the canonical module and the $a$-invariant. The subring $S$ is a standard $K$-algebra because $\langle (w_i, 1), e_{n+1}\rangle = 1$ for all $i$. If $S$ is normal, then by a formula of Danilov and Stanley [18] the canonical module of $S$ is the ideal of $S$ given by

$$\omega_S = (\{x^a t^b | \, (a, b) \in \mathbb{N}\mathcal{B} \cap (\mathbb{R}_+\mathcal{B})^{\rm o}\}), \tag{2}$$

where $\mathcal{B} = \{(w_1, 1), \ldots, (w_r, 1)\}$ and $(\mathbb{R}_+\mathcal{B})^{\rm o}$ is the relative interior of $\mathbb{R}_+\mathcal{B}$. This expression for the canonical module of $S$ is central for our purposes. There are some methods, based on combinatorial optimization, that have been used to study canonical modules of edge subrings of bipartite graphs [79]. Our approach to study canonical modules is inspired by these methods.

Recall that the *$a$-invariant* of $S$, denoted by $a(S)$, is the degree as a rational function of the Hilbert series of $S$ [87, p. 99]. Thus we may compute $a$-invariants using the program *Normaliz* [11]. This invariant plays a key role in the theory of Hilbert functions [10]. If $S$ is normal, then $S$ is Cohen-Macaulay [56] and its $a$-invariant is given by

$$a(S) = -\min\{\, i \,|\, (\omega_S)_i \neq 0\}, \tag{3}$$

see [10, p. 141] and [87, Proposition 4.2.3].

In Section 1.5 we study the canonical module and the $a$-invariant of $S$, when $S$ arises from an integer rounding problem.

For use below let $P = \{x | \, x \geq 0; xA \leq \mathbf{1}\}$; let $\mathrm{vert}(P)$ be the set of vertices of $P$ and let $\ell_1, \ldots, \ell_p$ be the set of all maximal elements of $\mathrm{vert}(P)$ (maximal with respect to $<$). For each $1 \leq i \leq p$ there is a unique positive integer $d_i$ such that the non-zero entries of $(-d_i\ell_i, d_i)$ are relatively prime.

The main results of Section 1.5 are as follows.

**Theorem 1.5.2** *If the system $x \geq 0; xA \leq \mathbf{1}$ has the integer rounding property, then the canonical module of $S = K[x^{w_1}t, \ldots, x^{w_r}t]$ is given by*

$$\omega_S = \left(\left\{ x^a t^b \Big| \, (a, b) \begin{pmatrix} -d_1\ell_1 & \cdots & -d_p\ell_p & e_1 & \cdots & e_n \\ d_1 & \cdots & d_p & 0 & \cdots & 0 \end{pmatrix} \geq \mathbf{1} \right\}\right), \tag{4}$$

*and the a-invariant of S is equal to* $-\max_i\{\lceil 1/d_i + |\ell_i|\rceil\}$. *Here* $|\ell_i| = \langle \ell_i, \mathbf{1}\rangle$.

This result complements and generalizes a result of [23] valid only for incidence matrices of clutters.

Using the description above for $\omega_S$, we then prove the following two results.

**Theorem 1.5.3** *Assume that the system* $x \geq 0$*;* $xA \leq \mathbf{1}$ *has the integer rounding property. If S is Gorenstein and* $c_0 = \max\{|\ell_i|\colon\, 1 \leq i \leq p\}$ *is an integer, then* $|\ell_k| = c_0$ *for each* $1 \leq k \leq p$ *such that* $\ell_k$ *has integer entries.*

**Theorem 1.5.4** *Assume that the system* $x \geq 0; xA \leq \mathbf{1}$ *has the integer rounding property. If* $-a(S) = 1/d_i + |\ell_i|$ *for* $i = 1, \ldots, p$*, then S is Gorenstein.*

As a consequence of Theorems 1.5.3 and 1.5.4 we obtain that if $P$ is an integral polytope, i.e., it has only integral vertices, then $S$ is Gorenstein if and only if $a(S) = -(|\ell_i| + 1)$ for $i = 1, \ldots, p$ (see Corollary 1.5.5). Based on a computer analysis, using the program *Normaliz* [11], we conjecture a possible description of all Gorenstein subrings $S$ in terms of the vertices of $P$ (see Problem 1.5.6).

We are able to give an explicit description of the canonical module of $S$ and its $a$-invariant when $\mathcal{C}$ is the clutter of maximal cliques of a perfect graph (Theorem 1.5.7). The $a$-invariant of general subrings arising from clutters seems to be closely related to the combinatorial structure of the clutter (see [79, Proposition 4.2] and Theorem 1.5.7).

We also examine the Gorenstein and complete intersection properties of subrings arising from systems with the integer rounding property of incidence matrices of graphs (see Proposition 1.5.11).

**Contents of Chapter 2** Let $\mathcal{C}$ be a clutter with vertex set $X = \{x_1, \ldots, x_n\}$ and let $I = I(\mathcal{C}) = (x^{v_1}, \ldots, x^{v_q})$ be its edge ideal. The set of vertices and edges of $\mathcal{C}$ are denoted by $V(\mathcal{C})$ and $E(\mathcal{C})$ respectively. We shall always assume that $\mathcal{C}$ has no isolated vertices, i.e., each vertex $x_i$ occurs in at least one edge of $\mathcal{C}$.

Let $A$ be the *incidence matrix* of $\mathcal{C}$ whose column vectors are $v_1, \ldots, v_q$. The *set covering polyhedron* of $\mathcal{C}$ is given by:

$$Q(A) = \{x \in \mathbb{R}^n |\, x \geq 0;\, xA \geq \mathbf{1}\},$$

A subset $C \subset X$ is called a *minimal vertex cover* of $\mathcal{C}$ if: (i) every edge of $\mathcal{C}$ contains at least one vertex of $C$, and (ii) there is no proper subset of $C$ with the first property. The map $C \mapsto \sum_{x_i \in C} e_i$ gives a bijection between the minimal vertex covers of $\mathcal{C}$ and the integral vectors of $Q(A)$. A polyhedron is called an *integral polyhedron* if it has only integral vertices. A clutter is called *d-uniform* or *uniform* if all its edges have exactly $d$ vertices.

We begin in Section 2.1 by introducing various combinatorial properties of clutters. We then give a simple combinatorial proof of the following result of [42]:

**Proposition 2.1.2** *If $\mathcal{C}$ is a d-uniform clutter whose set covering polyhedron $Q(A)$ is integral, then there are $X_1, \ldots, X_d$ mutually disjoint minimal vertex covers of $\mathcal{C}$ such that $X = \cup_{i=1}^{d} X_i$.*

Example 2.1.3 shows that this result fails if we drop the uniformity hypothesis. For use below we denote the smallest number of vertices in any minimal vertex cover of $\mathcal{C}$ by $\alpha_0(\mathcal{C})$. The clutter obtained from $\mathcal{C}$ by deleting a vertex $x_i$ and removing all edges containing $x_i$ is denoted by $\mathcal{C} \setminus \{x_i\}$. A set of pairwise disjoint edges of $\mathcal{C}$ is called *independent* or a *matching* and a set of independent edges of $\mathcal{C}$ whose union is $X$ is called a *perfect matching*. We then prove:

**Proposition 2.1.6** *Let $\mathcal{C}$ be a d-uniform clutter with a perfect matching such that $Q(A)$ is integral. Then $\mathcal{C}$ is vertex critical, i.e., $\alpha_0(\mathcal{C} \setminus \{x_i\}) < \alpha_0(\mathcal{C})$ for all $i$.*

A simple example is shown to see that this result fails for non uniform clutters with integral set covering polyhedron (Remark 2.1.7).

The main result of Section 2.2 is:

**Theorem 2.2.6** *If $\mathcal{C}$ is a uniform clutter with the max-flow min-cut property, then*

(a) $\Delta_r(A) = 1$, *where* $r = \mathrm{rank}(A)$.

(b) $\mathbb{N}\mathcal{A} = \mathbb{R}_+\mathcal{A} \cap \mathbb{Z}^n$, *where* $\mathcal{A} = \{v_1, \ldots, v_q\}$.

Here $\Delta_r(A)$ denotes the greatest common divisor of all nonzero $r \times r$ subdeterminants of $A$, $\mathbb{N}\mathcal{A}$ denotes the semigroup generated by $\mathcal{A}$, and $\mathbb{R}_+\mathcal{A}$ denotes the cone generated by $\mathcal{A}$. As interesting consequences we obtain that if $\mathcal{C}$ is a $d$-uniform clutter with the max-flow min-cut property, then $A$ diagonalizes over $\mathbb{Z}$—using row and column operations—to an identity matrix (see Corollary 2.2.8) and $\mathcal{C}$ has a perfect matching if and only if $n = d\alpha_0(\mathcal{C})$ (see Corollary 2.2.10). In Example 2.2.7 we show that the uniformity hypothesis is essential in the two statements of Theorem 2.2.6.

Section 2.3 deals with the diagonalization problem (see Conjecture 2.2.16) for clutters with the packing property (see Definition 2.2.11). The following is one of the main result of this section. It gives some support to Conjecture 2.2.16.

**Theorem 2.3.1** *Let $\mathcal{C}$ be a $d$-uniform clutter with a perfect matching such that $\mathcal{C}$ has the packing property and $\alpha_0(\mathcal{C}) = 2$. If $A$ has rank $r$, then*

$$\Delta_r \begin{pmatrix} A \\ \mathbf{1} \end{pmatrix} = 1.$$

As an application we obtain the next result which gives some support to a Conjecture of Conforti and Cornuéjols [15] (see Conjecture 2.2.14).

**Corollary 2.3.3** *Let $\mathcal{C}$ be a $d$-uniform clutter with a perfect matching such that $\mathcal{C}$ has the packing property and $\alpha_0(\mathcal{C}) = 2$. If $v_1, \ldots, v_q$ are linearly independent, then $\mathcal{C}$ has the max-flow min-cut property.*

The other main result of Section 2.3 is:

**Theorem 2.3.4** *Let $\mathcal{C}$ be a $d$-uniform clutter with a partition $X_1, \ldots, X_d$ of $X$ such that $X_i = \{x_{2i-1}, x_{2i}\}$ is a minimal vertex cover of $\mathcal{C}$ for all $i$. If $I = I(\mathcal{C})$ is minimally non-normal and $\mathcal{C}$ satisfies the packing property, then* $\mathrm{rank}(A) = d + 1$.

Regular and unimodular triangulations are introduced in Section 2.4. There is a relationship between the Gröbner bases of the toric ideal of $K[F]$ and the triangulations of $\mathcal{A}$ [76], which has many interesting applications. We make use of the theory of Gröbner bases and convex polytopes, which was created and developed by Sturmfels [77], to prove the following main result of Section 2.4:

**Theorem 2.4.6** *Let $A$ be a balanced matrix with distinct column vectors $v_1, \ldots, v_q$. If $|v_i| = d$ for all $i$, then any regular triangulation of the cone $\mathbb{R}_+\{v_1, \ldots, v_q\}$ is unimodular.*

Here $|v_i|$ denotes the sum of the entries of the vector $v_i$. Recall that a matrix $A$ with entries in $\{0, 1\}$ is called *balanced* if $A$ has no square submatrix of odd order with exactly two 1's in each row and column. If we do not require the uniformity condition $|v_i| = d$ for all $i$ this result is false, as is seen in Example 2.4.7. What makes this result surprising is the fact that not all balanced matrices are unimodular (see Example 2.4.7). This result gives some support to Conjecture 2.4.3: If $\mathcal{C}$ is a uniform clutter that satisfies the max-flow min-cut property, then the rational polyhedral cone $\mathbb{R}_+\{v_1,\ldots,v_q\}$ has a unimodular regular triangulation.

**Contents of Chapter 3** Let $\mathcal{C}$ be a clutter with vertex set $X=\{x_1,\ldots,x_n\}$ and let $I=I(\mathcal{C})=(x^{v_1},\ldots,x^{v_q})$ be its edge ideal.

Edge ideals of graphs were introduced and studied in [73, 85]. Edge ideals of clutters correspond to simplicial complexes via the Stanley-Reisner correspondence [75] and to facet ideals [31, 94]. The Cohen-Macaulay property of edge ideals has been recently studied in [12, 30, 32, 49, 66, 81] using a combinatorial approach based on the notions of shellability, linear quotients, unmixedness, acyclicity and transitivity of digraphs, and the König property.

The aim of this chapter is to study the behavior, under certain operations, of various algebraic and combinatorial optimization properties of edge ideals and clutters such as the Cohen-Macaulay property, the normality, the torsion freeness, the packing and the max-flow min-cut properties. The study of edge ideals from the combinatorial optimization point of view was initiated in [7, 79] and continued in [19, 29, 40, 42, 45, 91], see also [53] . The Cohen-Macaulay and normality properties are two of the most interesting properties an edge ideal can have, see [10, 32, 75, 87] and [58, 84] respectively.

Let $\mathcal{C}$ be a clutter and let $I=I(\mathcal{C})$ be its edge ideal. Recall that a *deletion* of $I$ is any ideal $I'$ obtained from $I$ by making a variable equal to 0. A *deletion* of $\mathcal{C}$ is a clutter $\mathcal{C}'$ that corresponds to a deletion $I'$ of $I$. Notice that $\mathcal{C}'$ is obtained from $I'$ by considering the unique set of square-free monomials that minimally generate $I'$. A *contraction* of $I$ is any ideal $I'$ obtained from $I$ by making a variable equal to 1. A *contraction* of $\mathcal{C}$ is a clutter $\mathcal{C}'$ that corresponds to a contraction $I'$ of $I$. This terminology is consistent with that of [16, p. 23]. The *duplication* of a vertex $x_i$ means extending $X$ by a new vertex $x_i'$ and replacing $E(\mathcal{C})$ by

$$E(\mathcal{C})\cup\{(e\setminus\{x_i\})\cup\{x_i'\}\,|\,x_i\in e\in E(\mathcal{C})\}.$$

A clutter obtained from $\mathcal{C}$ by a sequence of deletions and duplications of vertices is called a *parallelization* of $\mathcal{C}$ and a clutter obtained from $\mathcal{C}$ by a sequence of deletions and contractions of vertices is called a *minor* of $\mathcal{C}$, see Section 3.1. It is known that the normality of $I(\mathcal{C})$ is closed under minors [29]. One of our main results shows that the normality of $I(\mathcal{C})$ is closed under parallelizations:

**Theorem 3.1.3** *Let $\mathcal{C}$ be a clutter and let $\mathcal{C}'$ be a parallelization of $\mathcal{C}$. If $I(\mathcal{C})$ is normal, then $I(\mathcal{C}')$ is normal.*

The ideal $I = I(\mathcal{C})$ is called *normally torsion free* if $I^i = I^{(i)}$ for all $i$, where $I^{(i)}$ is the *ith* symbolic power of $I$. As an application we prove that if $I(\mathcal{C})$ is normally torsion free and $\mathcal{C}'$ is a parallelization of $\mathcal{C}$, then $I(\mathcal{C}')$ is normally torsion free (Corollary 3.1.12).

Let $A$ be the incidence matrix of a clutter $\mathcal{C}$, i.e., $A$ is the matrix with column vectors $v_1, \ldots, v_q$. A remarkable result of [45] (cf. [42, Theorem 4.6]) shows that $I(\mathcal{C})$ is normally torsion free if and only if $\mathcal{C}$ has the max-flow min-cut property. This fact makes a strong connection between commutative algebra and combinatorial optimization. Some other interesting links between these two areas can be found in [29, 42, 47, 48, 53, 81, 91] and in the references there. It is known [71, Chapter 79] that a clutter $\mathcal{C}$ satisfies the max-flow min-cut property if and only if all parallelizations of the clutter $\mathcal{C}$ satisfy the König property (see Definition 3.1.7). As another application we give a proof of this fact using that the integrality of the polyhedron $\{x \mid x \geq 0; xA \geq \mathbf{1}\}$ is closed under parallelizations and minors and using that the normality of $I(\mathcal{C})$ is preserved under parallelizations and minors (Corollary 4.1.3).

A clutter $\mathcal{C}$ satisfies the *packing property* (PP for short) if all minors of $\mathcal{C}$ satisfy the König property. We say that a clutter $\mathcal{C}$ is *Cohen-Macaulay* if $R/I(\mathcal{C})$ is a Cohen-Macaulay ring, see [65]. The other main result of this chapter is:

**Theorem 3.2.3** *Let $\mathcal{C}$ be a $d$-uniform clutter on the vertex set $X$. Let*

$$Y = \{y_{ij} \mid 1 \leq i \leq n;\ 1 \leq j \leq d-1\}$$

*be a set of new variables, and let $\mathcal{C}'$ be the clutter with vertex set $V(\mathcal{C}') = X \cup Y$ and edge set*

$$E(\mathcal{C}') = E(\mathcal{C}) \cup \{\{x_1, y_{11}, \ldots, y_{1(d-1)}\}, \ldots, \{x_n, y_{n1}, \ldots, y_{n(d-1)}\}\}.$$

*Then the edge ideal $I(\mathcal{C}')$ is Cohen-Macaulay. If $\mathcal{C}$ satisfies* PP (*resp. max-flow min-cut*), *then $\mathcal{C}'$ satisfies* PP (*resp. max-flow min-cut*).

It is well known that if $\mathcal{C}$ satisfies the max-flow min-cut property, then $\mathcal{C}$ satisfies the packing property [16] (see Corollary 2.2.13). Conforti and Cornuéjols [15] conjecture that the converse is also true. Theorem 3.2.3 is interesting because it says that for uniform clutters it suffices to prove the conjecture for Cohen-Macaulay clutters, which have a rich structure. The Conforti-Cornuéjols conjecture has been studied in [22, 42, 45] using an algebraic approach based on certain algebraic properties of blowup algebras.

**Contents of Chapter 4** Let $(P, \prec)$ be a *partially ordered set* (*poset* for short) on the finite vertex set $X = \{x_1, \ldots, x_n\}$ and let $G$ be its *comparability graph*. Recall that the vertex set of $G$ is $X$ and the edge set of $G$ is the set of all unordered pairs $\{x_i, x_j\}$ such that $x_i$ and $x_j$ are comparable. A *clique* of $G$ is a subset of the set of vertices that induces a complete subgraph.

The *clique clutter* of $G$, denoted by $\mathrm{cl}(G)$, is the clutter with vertex set $X$ whose edges are exactly the maximal cliques of $G$ with respect to inclusion. One of our main algebraic results shows that the edge ideal $I = I(\mathrm{cl}(G))$ of $\mathrm{cl}(G)$ is normally torsion free (see Theorem 4.2.2), i.e., $I^i = I^{(i)}$ for $i \geq 1$, where $I^{(i)}$ is the *ith* symbolic power of $I$ (see Section 4.2). To prove this result we first show that the clique clutter of $G$ has the max-flow min-cut property (see Theorem 4.1.7). Then we use a result of [45] showing that an edge ideal $I(\mathcal{C})$, of a clutter $\mathcal{C}$, is normally torsion free if and only if $\mathcal{C}$ has the max-flow min-cut property. There are some other nice characterizations of the normally torsion free property that can be found in [57].

Let $f_1, \ldots, f_q$ be the maximal cliques of $G$ and let $v_k = \sum_{x_i \in f_k} e_i$ be the characteristic vector of $f_k$ for $k = 1, \ldots, q$, where $e_i$ is the *ith* unit vector of $\mathbb{R}^n$. The matrix $A$ with column vectors $v_1, \ldots, v_q$ is called the *vertex-clique matrix* of $G$ or the *incidence matrix* of $\mathrm{cl}(G)$. A *colouring* of the vertices of a graph is an assignment of colours to the vertices of the graph in such a way that adjacent vertices have distinct colours. The *chromatic number* of a graph is the minimal number of colours in a colouring. A graph is called *perfect* if for every induced subgraph $H$, the chromatic number of $H$ equals the size of the largest complete subgraph of $H$. It is well known that comparability graphs are perfect [71]. Thus $G$ is a perfect graph or equivalently the polytope

$$P(A) = \{x \mid x \geq 0;\ xA \leq \mathbf{1}\}$$

is integral, i.e., it has only integral vertices. We complement this fact by observing that the *set covering polyhedron*

$$Q(A) = \{x \mid x \geq 0;\ xA \geq \mathbf{1}\}$$

is also integral (see Corollary 4.1.9). Comparability graphs are interesting objects of study [3]. They have been nicely characterized in [39]. By [2] the clique clutter of any induced subgraph of $G$ satisfies the König property (see

Definition 4.1.2). In the terminology of [8] comparability graphs are clique-perfect.

Let $D_1, \ldots, D_s$ be the minimal vertex covers of $\overline{G}$, where $\overline{G}$ is the complement of $G$. The ideal of *vertex covers* of $\overline{G}$ is the square-free monomial ideal

$$I_c(\overline{G}) = (x^{u_1}, \ldots, x^{u_s}) \subset R,$$

where $x^{u_k} = \prod_{x_i \in D_k} x_i$. Note that $I_c(\overline{G})$ and $I(\mathrm{cl}(G))$ are dual of each other in the sense that $u_i + v_i = \mathbf{1}$ for $i = 1, \ldots, q$. In [91] it is shown that $I_c(\overline{G})$ is a normal ideal. We complement this fact by showing that the edge ideal $I(\mathrm{cl}(G))$ of $\mathrm{cl}(G)$ is a normal ideal (see Theorem 4.2.2). This is surprising because in general this duality does not preserve normality (see Example 1.3.8). As an application we prove that edge ideals of complete admissible uniform clutters are normally torsion free (see Theorem 4.2.3).

**Contents of Chapter 5** Let $\mathcal{C}$ be a clutter with vertex set $X = \{x_1, \ldots, x_n\}$ and let $I = I(\mathcal{C}) = (x^{v_1}, \ldots, x^{v_q})$ be its edge ideal.

The *blowup algebra* studied in this chapter is the *symbolic Rees algebra*:

$$R_s(I) = R \oplus I^{(1)}t \oplus \cdots \oplus I^{(i)}t^i \oplus \cdots \subset R[t],$$

where $t$ is a new variable and $I^{(i)}$ is the *ith* symbolic power of $I$. Closely related to $R_s(I)$ is the *Rees algebra* of $I$:

$$R[It] := R \oplus It \oplus \cdots \oplus I^i t^i \oplus \cdots \subset R[t].$$

The study of symbolic powers of edge ideals was initiated in [73] and further elaborated on in [1, 29, 42, 45, 51, 78, 91]. By a result of Lyubeznik [63], $R_s(I)$ is a $K$-algebra of finite type. In general the minimal set of generators of $R_s(I)$ as a $K$-algebra is very hard to describe in terms of $\mathcal{C}$ (see [1]). There are two exceptional cases. If the clutter $\mathcal{C}$ has the max-flow min-cut property, then by a result of [45] we have $I^i = I^{(i)}$ for all $i \geq 1$, i.e., $R_s(I) = R[It]$. If $G$ is a perfect graph, then the minimal generators of $R_s(I(G))$ are in one to one correspondence with the cliques (complete subgraphs) of $G$ [91], see also [78]. We shall be interested in studying the minimal set of generators of $R_s(I)$ using polyhedral geometry. Let $G$ be a graph and let $I_c(G)$ be the Alexander dual of $I(G)$, see definition below. Some of the main results of this chapter are graph theoretical descriptions of the minimal generators of $R_s(I(G))$ and $R_s(I_c(G))$. In Sections 5.1 and 5.2 we show that both algebras encode combinatorial information of the graph which can be decoded using integral Hilbert bases.

The *Rees cone* of $I$, denoted by $\mathbb{R}_+(I)$, is the polyhedral cone consisting of the non-negative linear combinations of the set

$$\mathcal{A}' = \{e_1, \ldots, e_n, (v_1, 1), \ldots, (v_q, 1)\} \subset \mathbb{R}^{n+1},$$

where $e_i$ is the *ith* unit vector.

Let $\mathfrak{p}_1, \ldots, \mathfrak{p}_s$ be the minimal primes of the edge ideal $I = I(\mathcal{C})$ and let

$$C_k = \{x_i | \, x_i \in \mathfrak{p}_k\} \quad (k = 1, \ldots, s)$$

be the corresponding minimal vertex covers of $\mathcal{C}$, see [87, Proposition 6.1.16]. Recall that the primary decomposition of the edge ideal of $\mathcal{C}$ is given by

$$I(\mathcal{C}) = (C_1) \cap (C_2) \cap \cdots \cap (C_s).$$

In particular observe that the height of $I(\mathcal{C})$ equals the number of vertices in a minimum vertex cover of $\mathcal{C}$. The *ith symbolic power* of $I$ is given by

$$I^{(i)} = S^{-1}I^i \cap R \;\text{ for } i \geq 1,$$

where $S = R \setminus \cup_{k=1}^{s} \mathfrak{p}_i$ and $S^{-1}I^i$ is the localization of $I^i$ at $S$. In our situation the *ith* symbolic power of $I$ has a simple expression: $I^{(i)} = \mathfrak{p}_1^i \cap \cdots \cap \mathfrak{p}_s^i$, see [87]. The Rees cone of $I$ is a finitely generated rational cone of dimension $n + 1$. Hence by the finite basis theorem [93, Theorem 4.11] there is a unique irreducible representation

$$\mathbb{R}_+(I) = H_{e_1}^+ \cap H_{e_2}^+ \cap \cdots \cap H_{e_{n+1}}^+ \cap H_{\ell_1}^+ \cap H_{\ell_2}^+ \cap \cdots \cap H_{\ell_r}^+ \qquad (5)$$

such that each $\ell_k$ is in $\mathbb{Z}^{n+1}$, the non-zero entries of each $\ell_k$ are relatively prime, and none of the closed halfspaces $H_{e_1}^+, \ldots, H_{e_{n+1}}^+, H_{\ell_1}^+, \ldots, H_{\ell_r}^+$ can be omitted from the intersection. Here $H_a^+$ denotes the closed halfspace $H_a^+ = \{x | \, \langle x, a \rangle \geq 0\}$ and $H_a$ stands for the hyperplane through the origin with normal vector $a$, where $\langle \; , \; \rangle$ denotes the standard inner product. The *facets* (i.e., the proper faces of maximum dimension or equivalently the faces of dimension $n$) of the Rees cone are exactly:

$$F_1 = H_{e_1} \cap \mathbb{R}_+(I), \ldots, F_{n+1} = H_{e_{n+1}} \cap \mathbb{R}_+(I), H_{\ell_1} \cap \mathbb{R}_+(I), \ldots, H_{\ell_r} \cap \mathbb{R}_+(I).$$

According to [29, Lemma 3.1] we may always assume that $\ell_k = -e_{n+1} + \sum_{x_i \in C_k} e_i$ for $1 \leq k \leq s$, i.e., each minimal vertex cover of $\mathcal{C}$ determines a facet of the Rees cone and every facet of the Rees cone satisfying $\langle \ell_k, e_{n+1}\rangle = -1$ must be of the form $\ell_k = -e_{n+1} + \sum_{x_i\in C_k} e_i$ for some minimal vertex cover $C_k$ of $\mathcal{C}$. This is quite interesting because this is saying that the Rees cone of $I(\mathcal{C})$ is a carrier of combinatorial information of the clutter $\mathcal{C}$. Thus we can extract the primary decomposition of $I(\mathcal{C})$ from the irreducible representation of $\mathbb{R}_+(I(\mathcal{C}))$.

Rees cones have been used to study algebraic and combinatorial properties of blowup algebras of square-free monomial ideals and clutters [29, 42]. Blowup algebras are interesting objects of study in algebra and geometry [82].

The ideal of *vertex covers* of $\mathcal{C}$ is the square-free monomial ideal

$$I_c(\mathcal{C}) = (x^{u_1}, \ldots, x^{u_s}) \subset R,$$

where $x^{u_k} = \prod_{x_i\in C_k} x_i$. Often the ideal $I_c(\mathcal{C})$ is called the *Alexander dual* of $I(\mathcal{C})$. The clutter $\Upsilon(\mathcal{C})$ associated to $I_c(\mathcal{C})$ is called the *blocker* of $\mathcal{C}$, see [16]. Notice that the edges of $\Upsilon(\mathcal{C})$ are precisely the minimal vertex covers of $\mathcal{C}$. If $G$ is a graph, then $R_s(I_c(G))$ is generated as a $K$-algebra by elements of degree in $t$ at most two [51, Theorem 5.1]. One of the main result of Section 5.1 is a graph theoretical description of the minimal generators of $R_s(I_c(G))$ (see Theorem 5.1.9). As an application we recover an explicit description given in [80] of the edge cone of a graph (Corollary 5.1.10).

The symbolic Rees algebra of the ideal $I_c(\mathcal{C})$ can be interpreted in terms of "$k$-vertex covers" [51] as we now explain. Let $a = (a_1, \ldots, a_n) \neq 0$ be a vector in $\mathbb{N}^n$ and let $b \in \mathbb{N}$. We say that $a$ is a *$b$-vertex cover* of $I$ (or $\mathcal{C}$) if $\langle v_i, a\rangle \geq b$ for $i = 1, \ldots, q$. Often we will call a $b$-vertex cover simply a *$b$-cover*. This notion plays a role in combinatorial optimization [71, Chapter 77, p. 1378] and algebraic combinatorics [51, 52].

The *algebra of covers* of $I$ (or $\mathcal{C}$), denoted by $R_c(I)$, is the $K$-subalgebra of $K[t]$ generated by all monomials $x^a t^b$ such that $a$ is a $b$-cover of $I$. We say that a $b$-cover $a$ of $I$ is *reducible* if there exists an $i$-cover $c$ and a $j$-cover $d$ of $I$ such that $a = c + d$ and $b = i + j$. If $a$ is not reducible, we call $a$ *irreducible*. The irreducible 0 and 1 covers of $\mathcal{C}$ are the unit vectors $e_1, \ldots, e_n$ and the incidence vectors $u_1, \ldots, u_s$ of the minimal vertex covers of $\mathcal{C}$, respectively. The minimal generators of $R_c(I)$ as a $K$-algebra correspond to the irreducible covers of $I$. Notice the following dual descriptions:

$$\begin{aligned}
I^{(b)} &= (\{x^a \mid \langle (a, b), \ell_i\rangle \geq 0 \text{ for } i = 1, \ldots, s\}) \\
&= (\{x^a \mid \langle a, u_i\rangle \geq b \text{ for } i = 1, \ldots, s\}), \\
J^{(b)} &= (\{x^a \mid \langle a, v_i\rangle \geq b \text{ for } i = 1, \ldots, q\}),
\end{aligned}$$

where $J = I_c(\mathcal{C})$. Hence $R_c(I) = R_s(J)$ and $R_c(J) = R_s(I)$.

In general each $\ell_i$ occurring in Eq. (5) determines a minimal generator of $R_s(I_c(\mathcal{C}))$. Indeed if we write $\ell_i = (a_i, -d_i)$, where $a_i \in \mathbb{N}^n$, $d_i \in \mathbb{N}$, then $a_i$ is an irreducible $d_i$-cover of $I$ (Lemma 5.1.3). Let $F_{n+1}$ be the facet of $\mathbb{R}_+(I)$ determined by the hyperplane $H_{e_{n+1}}$. Thus we have a map $\psi$:

$$\begin{array}{rcl} \{\text{Facets of } \mathbb{R}_+(I(\mathcal{C}))\} \setminus \{F_{n+1}\} & \stackrel{\psi}{\longrightarrow} & R_s(I_c(\mathcal{C})) \\ H_{\ell_k} \cap \mathbb{R}_+(I) & \stackrel{\psi}{\longrightarrow} & x^{a_k} t^{d_k}, \text{ where } \ell_k = (a_k, -d_k) \\ H_{e_i} \cap \mathbb{R}_+(I) & \stackrel{\psi}{\longrightarrow} & x_i \end{array}$$

whose image provides a good approximation for the minimal set of generators of $R_s(I_c(\mathcal{C}))$ as a $K$-algebra. Likewise the facets of $\mathbb{R}_+(I_c(\mathcal{C}))$ give an approximation for the minimal set of generators of $R_s(I(\mathcal{C}))$. In Example 5.1.6 we show a connected graph $G$ for which the image of the map $\psi$ does not generate $R_s(I_c(G))$. For balanced clutters, i.e., for clutters without odd cycles, the image of the map $\psi$ generates $R_s(I_c(\mathcal{C}))$. This follows from [42, Propositions 4.10 and 4.11]. In particular the image of the map $\psi$ generate $R_s(I_c(\mathcal{C}))$ when $\mathcal{C}$ is a bipartite graph. It would be interesting to characterize when the irreducible representation of the Rees cone determine the irreducible covers.

The *Simis cone* of $I$ is the rational polyhedral cone:

$$\mathrm{Cn}(I) = H^+_{e_1} \cap \cdots \cap H^+_{e_{n+1}} \cap H^+_{(u_1,-1)} \cap \cdots \cap H^+_{(u_s,-1)},$$

Simis cones were introduced in [29] to study symbolic Rees algebras of square-free monomial ideals. If $\mathcal{H}$ is an integral Hilbert basis of $\mathrm{Cn}(I)$, then $R_s(I(\mathcal{C}))$ equals $K[\mathbb{N}\mathcal{H}]$, the semigroup ring of $\mathbb{N}\mathcal{H}$ (see [29, Theorem 3.5]). This result is interesting because it allows us to compute the minimal generators of $R_s(I(\mathcal{C}))$ using Hilbert bases. The program *Normaliz* [11] is suitable for computing Hilbert bases. There is a description of $\mathcal{H}$ valid for perfect graphs [91].

If $G$ is a perfect graph, the irreducible $b$-covers of $\Upsilon(G)$ correspond to cliques of $G$ [91] (cf. Corollary 5.2.5). In this case, setting $\mathcal{C} = \Upsilon(G)$, it turns out that the image of $\psi$ generates $R_s(I_c(\Upsilon(G)))$. Notice that $I_c(\Upsilon(G))$ is equal to $I(G)$.

In Section 5.2 we introduce and study the concept of an irreducible graph. A $b$-cover $a = (a_1, \ldots, a_n)$ is called *binary* if $a_i \in \{0, 1\}$ for all $i$. We present a graph theoretical description of the irreducible binary $b$-vertex covers of the blocker of $G$ (see Theorem 5.2.7). It is shown that they are in one to one correspondence with the irreducible induced subgraphs of $G$. As a byproduct we obtain a method, using Hilbert bases, to obtain all irreducible induced subgraphs of $G$ (see Corollary 5.2.10). We give a simple procedure to build irreducible graphs (Proposition 5.2.17) and give a method to construct irreducible $b$-vertex covers of the blocker of $G$ with high degree relative to the number of vertices of $G$ (see Corollaries 5.2.22 and 5.2.23).

# Chapter 1

# Duality, a-invariants and Canonical Modules of Rings Arising from Linear Optimization Problems

The aim of this chapter is to study integer rounding properties of various systems of linear inequalities to gain insight about the algebraic properties of Rees algebras of monomial ideals and monomial subrings. We study the normality and the Gorenstein property—as well as the canonical module and the $a$-invariant—of Rees algebras and subrings arising from systems with the integer rounding property. We relate the algebraic properties of Rees algebras and monomial subrings with integer rounding properties and present a duality theorem. The normality of a monomial ideal is expressed in terms of blocking polyhedra and the integer decomposition property. For edge ideals of clutters this property completely determine their normality. For systems arising from cliques of perfect graphs explicit expressions for the canonical module and the $a$-invariant are given. The combinatorial notions considered here come from linear optimization problems.

Let $R = K[x_1, \ldots, x_n]$ be a polynomial ring over a field $K$ and let $v_1, \ldots, v_q$ be the column vectors of a matrix $A = (a_{ij})$ whose entries are non-negative integers. We shall always assume that the rows and columns of $A$ are different from zero. As usual we use the notation $x^a := x_1^{a_1} \cdots x_n^{a_n}$, where $a = (a_1, \ldots, a_n) \in \mathbb{N}^n$.

The *monomial algebras* considered here are: (a) the *Rees algebra*

$$R[It] := R \oplus It \oplus \cdots \oplus I^i t^i \oplus \cdots \subset R[t],$$

where $I = (x^{v_1}, \ldots, x^{v_q}) \subset R$ and $t$ is a new variable, (b) the *extended Rees algebra*

$$R[It, t^{-1}] := R[It][t^{-1}] \subset R[t, t^{-1}],$$

(c) the *monomial subring*

$$K[F] = K[x^{v_1}, \ldots, x^{v_q}] \subset R$$

spanned by $F = \{x^{v_1}, \ldots, x^{v_q}\}$, (d) the *homogeneous monomial subring*

$$K[Ft] = K[x^{v_1}t, \ldots, x^{v_q}t] \subset R[t]$$

spanned by $Ft$, (e) the *homogeneous monomial subring*

$$K[Ft \cup \{t\}] = K[x^{v_1}t, \ldots, x^{v_q}t, t] \subset R[t]$$

spanned by $Ft \cup \{t\}$, (f) the *homogeneous monomial subring*

$$S = K[x^{w_1}t, \ldots, x^{w_r}t] \subset R[t],$$

where $w_1, \ldots, w_r$ is the set of all vectors $\alpha \in \mathbb{N}^n$ such that $0 \leq \alpha \leq v_i$ for some $i$, and (g) the *Ehrhart ring*

$$A(P) = K[\{x^a t^i \,|\, a \in \mathbb{Z}^n \cap iP; i \in \mathbb{N}\}] \subset R[t]$$

of a lattice polytope $P$.

The aim of this chapter is to study max-flow min-cut properties of clutters and integer rounding properties of various systems of linear inequalities—and their underlying polyhedra—to gain insight about the algebraic properties of these algebras and viceversa. Systems with integer rounding properties and clutters with the max-flow min-cut property come from linear optimization problems [70, 71]. The precise definitions will be given in Section 1.2.

Before stating our main results, we recall a few basic facts about the normality of monomial subrings. According to [87] the *integral closure* of $K[F]$ in its field of fractions can be expressed as

$$\overline{K[F]} = K[\{x^a \,|\, a \in \mathbb{Z}\mathcal{A} \cap \mathbb{R}_+\mathcal{A}\}], \tag{1.1}$$

where $\mathcal{A} = \{v_1, \ldots, v_q\}$, $\mathbb{Z}\mathcal{A}$ is the subgroup of $\mathbb{Z}^n$ spanned by $\mathcal{A}$, and $\mathbb{R}_+\mathcal{A}$ is the cone generated by $\mathcal{A}$. The subring $K[F]$ equals $K[\mathbb{N}\mathcal{A}]$, the semigroup ring of $\mathbb{N}\mathcal{A}$. Recall that $K[F]$ is called *integrally closed* or *normal* if $K[F] = \overline{K[F]}$. Thus $K[F]$ is normal if and only if

$$\mathbb{N}\mathcal{A} = \mathbb{Z}\mathcal{A} \cap \mathbb{R}_+\mathcal{A},$$

where $\mathbb{N}\mathcal{A}$ is the subsemigroup of $\mathbb{N}^n$ generated by $\mathcal{A}$. The description of the integral closure given in Eq. (1.1) can of course be applied to any of the monomial algebras considered here. In particular if $\mathcal{A}'$ is the set

$$\mathcal{A}' = \{e_1, \ldots, e_n, (v_1, 1), \ldots, (v_q, 1)\},$$

where $e_i$ is the *ith* unit vector, then $\mathbb{Z}\mathcal{A}' = \mathbb{Z}^{n+1}$ and $R[It]$ is normal if and only if $\mathbb{N}\mathcal{A}' = \mathbb{Z}^{n+1} \cap \mathbb{R}_+\mathcal{A}'$. A dual characterization of the normality of $R[It]$ will be given in Proposition 1.2.8. Recall that the *integral closure* of $I^i$, denoted by $\overline{I^i}$, is the ideal of $R$ given by

$$\overline{I^i} = (\{x^a \in R \mid \exists\, p \in \mathbb{N} \setminus \{0\}; (x^a)^p \in I^{pi}\}),$$

see for instance [87, Proposition 7.3.3]. The ideal $I$ is called *normal* if $I^i = \overline{I^i}$ for all $i$. If $\overline{I} = I$, the ideal $I$ is called *integrally closed*. Note that the ideal $I$ is normal if and only if $R[It]$ is normal [87].

Recall that the Ehrhart ring $A(P)$ is always normal [10]. A finite set $\mathcal{A} \subset \mathbb{R}^n$ is called an *Hilbert basis* if $\mathbb{N}\mathcal{A} = \mathbb{R}_+\mathcal{A} \cap \mathbb{Z}^n$. An *integral Hilbert basis* is a Hilbert basis consisting of integral vectors only. Note that if $\mathcal{A} = \{v_1, ..., v_q\}$ is an integral Hilbert basis, then the ring $K[x^{v_1}, ..., x^{v_q}]$ is normal.

The contents of this chapter are as follows. In Section 1.1 we study the normality of general monomial ideals. We are able to characterize this property in terms of blocking polyhedra and the integer decomposition property (see Theorem 1.1.2). For integrally closed ideals and for edge ideals this property completely characterizes their normality. As an interesting consequence using a result of Baum and Trotter [4] we describe the normality of a monomial ideal in terms of the integer rounding property:

**Corollary 1.1.8** *The ideal $I$ is normal if and only if the system $xA \geq \mathbf{1}; x \geq 0$ has the integer rounding property.*

There are already some characterizations of the normality of monomial ideals [29, Theorem 4.4]. There are also useful membership tests based on linear optimization to decide whether a given monomial lies in the integral closure of a monomial ideal [19, Proposition 3.5]. A combinatorial description of the integral closure of a monomial ideal is given in [83, Section 6.6] (cf. [26, p. 141]). This description has been used in [19] to find multiplicities of ideals and volumes of lattice polytopes based on a probabilistic approach using a Monte Carlo method.

We use the theory of blocking and antiblocking polyhedra [4, 37, 70] to describe when the systems

$$x \geq 0;\ xA \leq \mathbf{1}, \quad x \geq 0; xA \geq \mathbf{1}, \quad xA \leq \mathbf{1},$$

have the *integer rounding property* (see Definitions 1.2.2, 1.1.6, 1.2.18) in terms of the normality of the monomial algebras considered here. As usual $\mathbf{1}$ denotes the vector whose entries are equal to 1.

One of the main results of Section 1.2 is:

**Theorem 1.2.5** *The system $x \geq 0$; $xA \leq \mathbf{1}$ has the integer rounding property if and only if the subring $S = K[x^{w_1}t, \ldots, x^{w_r}t]$ is normal.*

This result was shown in [23] when $A$ is the incidence matrix of a clutter, i.e., when the entries of $A$ are in $\{0, 1\}$. Recall that a *clutter* $\mathcal{C}$ with finite vertex set $X = \{x_1, \ldots, x_n\}$ is a family of subsets of $X$, called edges, none of which is included in another. The *incidence matrix* of a clutter $\mathcal{C}$ is the vertex-edge matrix whose columns are the characteristic vectors of the edges of $\mathcal{C}$. The *edge ideal* of a clutter $\mathcal{C}$, denoted by $I(\mathcal{C})$, is the ideal of $R$ generated by all monomials $x_e = \prod_{x_i \in e} x_i$ such that $e$ is an edge of $\mathcal{C}$. The *Alexander dual* of $I(\mathcal{C})$ is the ideal of $R$ given by $I(\mathcal{C})^\vee = \cap_{e \in E}(e)$, where $E = E(\mathcal{C})$ is the edge set of $\mathcal{C}$.

As a consequence we show that if $A$ is the incidence matrix of a clutter all of whose edges have the same number of elements and either system $x \geq 0; xA \leq \mathbf{1}$ or $x \geq 0; xA \geq \mathbf{1}$ has the integer rounding property, then the subring $K[x^{v_1}t, \ldots, x^{v_q}t]$ is normal (see Corollary 1.2.6).

We present a duality between the integer rounding property of the systems $x \geq 0; xA \geq \mathbf{1}$ and $x \geq 0; xA^* \leq \mathbf{1}$ valid for matrices with entries in $\{0, 1\}$, where $a_{ij}^* = 1 - a_{ij}$ is the $ij$-entry of $A^*$. This duality is extended to a duality between monomial subrings.

Altogether another main result of Section 1.2 is:

**Theorem 1.2.11** *Let $A$ be the incidence matrix of a clutter. If $v_i^* = \mathbf{1} - v_i$*

*and* $A^*$ *is the matrix with column vectors* $v_1^*,\ldots,v_q^*$*, then the following are equivalent:*

(a) $R[I^*t]$ *is normal, where* $I^*=(x^{v_1^*},\ldots,x^{v_q^*})$.

(b) $S=K[x^{w_1}t,\ldots,x^{w_r}t]$ *is normal.*

(c) $\{-e_1,\ldots,-e_n,(v_1,1),\ldots,(v_q,1)\}$ *is a Hilbert basis.*

(d) $x\geq 0; xA^*\geq \mathbf{1}$ *has the integer rounding property.*

(e) $x\geq 0; xA\leq \mathbf{1}$ *has the integer rounding property.*

Then we present some interesting consequences of this duality. First of all we recover one of the main results of [91] showing that if

$$P=\{x|x\geq 0; xA\leq \mathbf{1}\}$$

is an integral polytope, i.e., $P$ has only integer vertices, and $A$ is a $\{0,1\}$-matrix, then the Rees algebra $R[I^*t]$ is normal (see Corollary 1.2.14). This result is related to perfect graphs. Indeed if $P$ is integral, then $v_1,\ldots,v_q$ correspond to the maximal cliques (maximal complete subgraphs) of a perfect graph $H$ [14, 62], and $v_1^*,\ldots,v_q^*$ correspond to the minimal vertex covers of the complement of $H$. Second we show that if $A$ is the incidence matrix of the collection of basis of a matroid, then all systems

$$x\geq 0; xA\geq \mathbf{1},\quad x\geq 0; xA^*\geq \mathbf{1},\quad x\geq 0; xA\leq \mathbf{1},\quad x\geq 0; xA^*\leq \mathbf{1}$$

have the integer rounding property (see Corollary 1.2.15). Third we show that if $A$ is the incidence matrix of a graph, then $R[It]$ is normal if and only if $R[I^*t]$ is normal (see Corollary 1.3.4). We give an example to show that this result does not extends to arbitrary uniform clutters (see Example 1.3.5). If $A$ is the incidence matrix of a graph $G$, we characterize when $I^*$ is the Alexander dual of the edge ideal of the complement of $G$ (see Proposition 1.3.6). If $G$ is a triangle-free graph, we show a duality between the normality of $I=I(G)$ and that of the Alexander dual of the edge ideal of the complement of $G$ (see Corollary 1.3.7). We show an example of an edge ideal of a graph whose Alexander dual is not normal (see Example 1.3.8). In [91] it is shown that this is never the case if the graph is perfect, i.e., the Alexander dual of the edge ideal of a perfect graph is always normal. Finally we recover one of the main results of [45] showing that if $A$ is the incidence matrix of a clutter $\mathcal{C}$, then $\mathcal{C}$ satisfies the max-flow min-cut property if and only if the set covering polyhedron

$$Q(A) = \{x |\, x \geq 0; xA \geq \mathbf{1}\}$$

is integral and $R[It]$ is normal (see Corollary 1.2.17).

Another main result of Section 1.2 is:

**Theorem 1.2.20** *If the system $xA \leq \mathbf{1}$ has the integer rounding property, then $K[F]$ is normal and $\mathbb{Z}^n/\mathbb{Z}\mathcal{A}$ is a torsion-free group. The converse holds if $|v_i| = d$ for all $i$. Here $|v_i| = \langle v_i, \mathbf{1}\rangle$.*

As a consequence of this result we prove: (i) If $A$ is the incidence matrix of a connected graph $G$, then the system $xA \leq \mathbf{1}$ has the integer rounding property if and only if $G$ is a bipartite graph (see Corollary 1.2.21), and (ii) Let $A$ be the incidence matrix of a clutter $\mathcal{C}$. If $\mathcal{C}$ is uniform, i.e., all its edges have the same number of elements, and $\mathcal{C}$ has the max-flow min-cut property (see Definition 2.2.3), then the system $xA \leq \mathbf{1}$ has the integer rounding property (see Corollary 1.2.22).

If $A$ is the incidence matrix of a bipartite graph, a remarkable result of Section 1.3 shows that the system $x \geq 0; xA \leq \mathbf{1}$ has the integer rounding property if and only if the extended Rees algebra $R[It, t^{-1}]$ is normal.

Let $G$ be a graph with incidence matrix $A$. If $G$ is a connected graph, in Section 1.3 we are able to prove that the system $x \geq 0; xA \leq \mathbf{1}$ has the integer rounding property if and only if the induced subgraph of the vertices of any two vertex disjoint odd cycles of $G$ is connected. Other equivalent descriptions of this property are also presented.

Before stating the main results of Section 1.5, we need to introduce the canonical module and the $a$-invariant (see Section 1.4 for additional details). Below we briefly explain the important role that these two objects play in the general theory. The subring $S$ is a standard $K$-algebra because $\langle (w_i, 1), e_{n+1}\rangle = 1$ for all $i$. Here $\langle\, ,\rangle$ is the standard inner product and $e_i$ is the *ith* unit vector. If $S$ is normal, then according to a formula of Danilov and Stanley [18] the canonical module of $S$ is the ideal of $S$ given by

$$\omega_S = (\{x^a t^b |\, (a, b) \in \mathbb{N}\mathcal{B} \cap (\mathbb{R}_+\mathcal{B})^{\mathrm{o}}\}), \tag{1.2}$$

where $\mathcal{B} = \{(w_1, 1), \ldots, (w_r, 1)\}$ and $(\mathbb{R}_+\mathcal{B})^{\mathrm{o}}$ is the relative interior of $\mathbb{R}_+\mathcal{B}$. This expression for the canonical module of $S$ is central for our purposes. Recall that the *a-invariant* of $S$, denoted by $a(S)$, is the degree as a rational function of the Hilbert series of $S$ [87, p. 99]. Thus we may compute $a$-invariants using the program *Normaliz* [11]. Let $H_S$ and $\varphi_S$ be the Hilbert function and the Hilbert polynomial of $S$ respectively. The *index of regularity* of $S$, denoted by $\text{reg}(S)$, is the least positive integer such that $H_S(i) = \varphi_S(i)$ for $i \geq \text{reg}(S)$. The $a$-invariant plays a fundamental role in algebra and geometry because one has:

$$\text{reg}(S) = \begin{cases} 0 & \text{if } a(S) < 0, \\ a(S) + 1 & \text{otherwise,} \end{cases}$$

see [87, Corollary 4.1.12]. If $S$ is normal, then $S$ is Cohen-Macaulay [56] and its $a$-invariant is given by

$$a(S) = -\min\{\, i \,|\, (\omega_S)_i \neq 0\}, \tag{1.3}$$

see [10, p. 141] and [87, Proposition 4.2.3].

In Section 1.4 we give a general technique to compute the canonical module and the $a$-invariant of a wide class of monomial subrings (see Theorem 1.4.1).

Then in Section 1.5 we study the canonical module and the $a$-invariant of $S$, when $S$ arises from an integer rounding problem. This invariant plays a key role in the theory of Hilbert functions [10]. There are some methods, based on combinatorial optimization, that have been used to study canonical modules of edge subrings of bipartite graphs [79]. Our approach to study canonical modules is inspired by these methods. If $S$ is a normal domain, we express the canonical module of $S$ and its $a$-invariant in terms of the vertices of the polytope $\{x|x \geq 0; xA \leq \mathbf{1}\}$ (see Theorem 1.5.2). We give necessary and sufficient conditions for $S$ to be Gorenstein and give a formula for the $a$-invariant of $S$ in terms of the vertices of the polytope $P = \{x|\, x \geq 0; xA \leq \mathbf{1}\}$. We are able to give an explicit description of the canonical module of $S$ and its $a$-invariant when $\mathcal{C}$ is the clutter of maximal cliques of a perfect graph (Theorem 1.5.7). The $a$-invariant of general subrings arising from clutters seems to be closely related to the combinatorial structure of the clutter (see [79, Proposition 4.2] and Theorem 1.5.7).

For use below let $\text{vert}(P)$ be the set of vertices of $P$ and let $\ell_1, \ldots, \ell_p$ be the set of all maximal elements of $\text{vert}(P)$ (maximal with respect to $<$). For each $1 \leq i \leq p$ there is a unique positive integer $d_i$ such that the non-zero entries of $(-d_i\ell_i, d_i)$ are relatively prime.

The main results of Section 1.5 are as follows.

**Theorem 1.5.2** *If the system $x \geq 0; xA \leq \mathbf{1}$ has the integer rounding property, then the canonical module of $S = K[x^{w_1}t, \ldots, x^{w_r}t]$ is given by*

$$\omega_S = \left(\left\{ x^a t^b \middle| \, (a,b) \begin{pmatrix} -d_1\ell_1 & \cdots & -d_p\ell_p & e_1 & \cdots & e_n \\ d_1 & \cdots & d_p & 0 & \cdots & 0 \end{pmatrix} \geq \mathbf{1} \right\}\right), \qquad (1.4)$$

*and the a-invariant of $S$ is equal to $-\max_i\{\lceil 1/d_i + |\ell_i| \rceil\}$. Here $|\ell_i| = \langle \ell_i, \mathbf{1}\rangle$.*

This result complements and generalizes a result of [23] valid only for incidence matrices of clutters. If $S$ is normal, the last Betti number in the homogeneous free resolution of the toric ideal $P_S$ of $S$ is equal to $\nu(\omega_S)$, the minimum number of generators of $\omega_S$. This number is called the *type* of $P_S$. Thus by describing the canonical module of $S$ we are in fact providing a device to compute the type of $P_S$. According to [79] the number of integral vertices of the polyhedron that defines $\omega_S$ (see Eq. (1.4)) is a lower bound for $\nu(\omega_S)$.

Using the description above for $\omega_S$ we then prove:

**Theorem 1.5.3** *Assume that the system $x \geq 0$; $xA \leq \mathbf{1}$ has the integer rounding property. If $S$ is Gorenstein and $c_0 = \max\{|\ell_i| \colon 1 \leq i \leq p\}$ is an integer, then $|\ell_k| = c_0$ for each $1 \leq k \leq p$ such that $\ell_k$ has integer entries.*

**Theorem 1.5.4** *Assume that the system $x \geq 0; xA \leq \mathbf{1}$ has the integer rounding property. If $-a(S) = 1/d_i + |\ell_i|$ for $i = 1, \ldots, p$, then $S$ is Gorenstein.*

As a consequence of Theorems 1.5.3 and 1.5.4 we obtain that if $P$ is an integral polytope, i.e., it has only integral vertices, then $S$ is Gorenstein if and only if $a(S) = -(|\ell_i| + 1)$ for $i = 1, \ldots, p$ (see Corollary 1.5.5).

We also examine the Gorenstein and complete intersection properties of subrings arising from systems with the integer rounding property of incidence matrices of graphs. Based on a computer analysis, using the program *Normaliz* [11], we conjecture a possible description of all Gorenstein subrings $S$ in terms of the vertices of $P$ (see Problem 1.5.6). Let $G$ be a connected graph with $n$ vertices and $q$ edges and let $A$ be its incidence matrix. If the system $xA \leq \mathbf{1}$ has the integer rounding property, then we show that $K[Ft \cup \{t\}]$ is a complete intersection if and only if $G$ is bipartite and the number of primitive cycles of $G$ is equal to $q - n + 1$ (see Proposition 1.5.11).

Let $G$ be a bipartite graph and let $A$ be its incidence matrix. A constructive description of all bipartite graphs such that $K[G] = K[x^{v_1}, \ldots, x^{v_q}]$ is a complete intersection is given in [41]. The Gorenstein property of $K[G]$ has been studied in [44, 55]. Thus by Lemma 1.5.10 and [10, Proposition 3.1.19] the Gorenstein property and the complete intersection property of $K[x^{v_1}t,\ldots,x^{v_q}t,t]$ are well understood in this particular case. The $a$-invariant of $K[G]$ has a combinatorial expression in terms of directed cuts and can be computed using linear programming [79]. Some other expressions for $a(K[G])$ can be found in [17, 44, 86].

# 1.1 Normality of monomial ideals

Let $R=K[x_1,\ldots,x_n]$ be a polynomial ring over a field $K$, let $I$ be a monomial ideal of $R$ generated by $x^{v_1},\ldots,x^{v_q}$, and let $A$ be the $n\times q$ matrix with column vectors $v_1,\ldots,v_q$. Recall that the *integral closure* of $I^i$, denoted by $\overline{I^i}$, is the ideal of $R$ given by

$$\overline{I^i}=(\{x^a\in R|\,\exists\, p\in\mathbb{N}\setminus\{0\};(x^a)^p\in I^{pi}\}),$$

see for instance [87, Proposition 7.3.3]. The ideal $I$ is called *normal* if $I^i=\overline{I^i}$ for all $i$. If $\overline{I}=I$, the ideal $I$ is called *integrally closed*. Note that the ideal $I$ is normal if and only if $R[It]$ is normal [87].

For a rational polyhedron $Q$ in $\mathbb{R}^n$, i.e., with rational vertices, define its *blocking polyhedron* $B(Q)$ by:

$$B(Q)=\{z\in\mathbb{R}^n|\, z\geq 0;\ \langle z,x\rangle\geq 1 \text{ for all } x \text{ in } Q\}.$$

For a matrix $A$ with entries in $\mathbb{N}$, its *covering polyhedron* $Q(A)$ is defined by:

$$Q(A)=\{x|\, x\geq 0; xA\geq \mathbf{1}\}.$$

If $A$ is the incidence matrix of a clutter $\mathcal{C}$, then the integral vectors of $Q(A)$ are in one to one correspondence with the minimal vertex covers of $\mathcal{C}$ [71].

The blocking polyhedron of $Q(A)$ can be expressed as follows.

**Lemma 1.1.1** *If $Q=Q(A)$, then $B(Q)=\mathbb{R}^n_+ +\mathrm{conv}(v_1,\ldots,v_q)$.*

**Proof.** The right hand side is clearly contained in the left hand side. Conversely take $z$ in $B(Q)$, then $\langle z,x\rangle\geq 1$ for all $x\in Q$ and $z\geq 0$. Let $\ell_1,\ldots,\ell_r$ be the vertex set of $Q$. In particular $\langle z,\ell_i\rangle\geq 1$ for all $i$. Then $\langle (z,1),(\ell_i,-1)\rangle\geq 0$ for all $i$. From [45, Theorem 3.2] we get that $(z,1)$ belongs to the cone generated by

$$\mathcal{A}'=\{e_1,\ldots,e_n,(v_1,1),\ldots,(v_q,1)\}.$$

Thus $z$ is in $\mathbb{R}_+^n + \mathrm{conv}(v_1, \ldots, v_q)$. This completes the proof of the asserted equality. $\Box$

A rational polyhedron $Q$ is said to have the *integer decomposition property* if for each natural number $k$ and for each integer vector $a$ in $kQ$, $a$ is the sum of $k$ integer vectors in $Q$; see [71, pp. 66–82]. Recall that $kQ$ is equal to $\{ka \mid a \in Q\}$.

**Theorem 1.1.2** *The ideal $I$ is normal if and only if the blocking polyhedron $B(Q)$ of $Q = Q(A)$ has the integer decomposition property and all minimal integer vectors of $B(Q)$ are columns of $A$ (minimal with respect to $\leq$).*

**Proof.** By Lemma 1.1.1, we have the equality $B(Q) = \mathbb{R}_+^n + \mathrm{conv}(v_1, \ldots, v_q)$. Hence $B(Q) \cap \mathbb{Q}^n = \mathbb{Q}_+^n + \mathrm{conv}_{\mathbb{Q}}(v_1, \ldots, v_q)$ because the polyhedron $B(Q)$ is rational. From this equality we readily obtain

$$\overline{I^k} = (\{x^a \mid a \in kB(Q) \cap \mathbb{Z}^n\}) \tag{1.5}$$

for $0 \neq k \in \mathbb{N}$. Assume that $I$ is normal, i.e., $\overline{I^k} = I^k$ for $k \geq 1$. Let $a$ be an integral vector in $kB(Q)$. Then $x^a \in I^k$ and consequently $a$ is the sum of $k$ integral vectors in $B(Q)$, that is, $B(Q)$ has the integer decomposition property. Take a minimal integer vector $a$ in $B(Q)$. Then $x^a \in \overline{I} = I$ and we can write $a = \delta + v_i$ for some $v_i$ and for some $\delta \in \mathbb{N}^n$. Thus $a = v_i$ by the minimality of $a$. Conversely assume that $B(Q)$ has the integer decomposition property and all minimal integer vectors of $B(Q)$ are columns of $A$. Take $x^a \in \overline{I^k}$, i.e., $a$ is an integral vector of $kB(Q)$. Hence $a$ is the sum of $k$ integral vectors $\alpha_1, \ldots, \alpha_k$ in $B(Q)$. Since any minimal vector of $B(Q)$ is a column of $A$ we may assume that $\alpha_i = c_i + v_i$ for $i = 1, \ldots, k$. Hence $x^a \in I^k$, as required. $\Box$

**Theorem 1.1.3** *If $I = \overline{I}$, then $I$ is normal if and only if the blocking polyhedron $B(Q)$ has the integer decomposition property.*

**Proof.** $\Rightarrow$) If $I$ is normal, by Theorem 1.1.2 the blocking polyhedron $B(Q)$ has the integer decomposition property.

$\Leftarrow$) Take $x^a \in \overline{I^k}$. From Eq. (1.5) we get that $a$ is an integral vector of $kB(Q)$. Hence $a$ is the sum of $k$ integral vectors $\alpha_1, \ldots, \alpha_k$ in $B(Q)$. Using Eq. (1.5) with $k = 1$, we get that $\alpha_1, \ldots, \alpha_k$ are in $\overline{I} = I$. Hence $x^a \in I^k$, as required. $\Box$

**Corollary 1.1.4** *If $I = I(\mathcal{C})$ is the edge ideal of a clutter $\mathcal{C}$, then $I$ is normal if and only if the blocking polyhedron $B(Q)$ has the integer decomposition property.*

**Proof.** Recall that $I$ is an intersection of prime ideals (see Eq. (4.2)). Thus it is seen that $\overline{I} = I$ and the result follows from Theorem 1.1.3. □

The general definition of integer rounding property is the following:

**Definition 1.1.5** A rational system $Ax \leq b$ of linear inequalities is said to have the *integer rounding property* if

$$\min\{\langle y, b\rangle \,|\, y \geq 0; yA = c; y \text{ integral}\} = \lceil \min\{\langle y, b\rangle \,|\, y \geq 0; yA = c\}\rceil \quad (1.6)$$

for each integral vector $c$ for which $\min\{\langle y, b\rangle \,|\, y \geq 0; yA = c\}$ is finite.

For systems of the form $x \geq 0; xA \geq \mathbf{1}$ this definition can be stated as:

**Definition 1.1.6** Let $A$ be a matrix with entries in $\mathbb{N}$. The linear system $x \geq 0; xA \geq \mathbf{1}$ has the *integer rounding property* if

$$\max\{\langle y, \mathbf{1}\rangle \,|\, y \geq 0; Ay \leq w; y \in \mathbb{N}^q\} = \lfloor \max\{\langle y, \mathbf{1}\rangle \,|\, y \geq 0; Ay \leq w\}\rfloor \quad (1.7)$$

for each integral vector $w$ for which the right hand side is finite.

Systems with the integer rounding property have been widely studied; see for instance [70, Chapter 22], [71, Chapter 5], and the references there.

**Theorem 1.1.7** ([4], [71, p. 82]) *The system $x \geq 0; xA \geq \mathbf{1}$ has the integer rounding property if and only if the blocking polyhedron $B(Q)$ of $Q = Q(A)$ has the integer decomposition property and all minimal integer vectors of $B(Q)$ are columns of $A$ (minimal with respect to $\leq$).*

The next result is interesting because it allows us to determine whether a given system $xA \geq \mathbf{1}; x \geq 0$, with $A$ an integral non-negative matrix, has the rounding property using the program [11].

**Corollary 1.1.8** *Let $I = (x^{v_1}, \ldots, x^{v_q})$ be a monomial ideal and let $A$ be the matrix with column vectors $v_1, \ldots, v_q$. Then $I$ is a normal ideal if and only if the system $xA \geq \mathbf{1}; x \geq 0$ has the integer rounding property.*

**Proof.** According to Theorem 1.1.7, the system $xA \geq \mathbf{1}; x \geq 0$ has the integer rounding property if and only if the blocking polyhedron $B(Q)$ of $Q = Q(A)$ has the integer decomposition property and all minimal integer vectors of $B(Q)$ are columns of $A$ (minimal with respect to $\leq$) (cf. [71, p. 82]). Thus the result follows at once from Theorem 1.1.2. $\Box$

We have been informed that Corollary 1.1.8 was observed by Trung when $I$ is the edge ideal of a hypergraph.

## 1.2 Integer rounding properties

We continue to use the notation and definitions used in the introduction. In this section we introduce and study integer rounding properties, describe some of their properties, present a duality theorem and show several applications.

Let $P$ be a rational polyhedron in $\mathbb{R}^n$, i.e., with rational vertices. Recall that the *antiblocking polyhedron* of $P$ is defined as:

$$T(P) := \{z |\, z \geq 0; \langle z, x\rangle \leq 1 \text{ for all } x \in P\}.$$

**Lemma 1.2.1** *Let $A$ be a matrix of size $n \times q$ with entries in $\mathbb{N}$, let $v_1, \ldots, v_q$ be the column vectors of $A$ and let $\{w_1, \ldots, w_r\}$ be the set of all $\alpha$ in $\mathbb{N}^n$ such that $\alpha \leq v_i$ for some $i$. If $P = \{x|\, x \geq 0; xA \leq \mathbf{1}\}$, then*

$$T(P) = \text{conv}(w_1, \ldots, w_r).$$

**Proof.** First we show the following equality which is interesting in its own right:

$$\text{conv}(w_1, \ldots, w_r) = \mathbb{R}^n_+ \cap (\text{conv}(w_1, \ldots, w_r) + \mathbb{R}_+\{-e_1, \ldots, -e_n\}). \quad (1.8)$$

Clearly the left hand side is contained in the right hand side. Conversely let $z$ be a vector in the right hand side. Then $z \geq 0$ and we can write

$$z = \lambda_1 w_1 + \cdots + \lambda_r w_r - \delta_1 e_1 - \cdots - \delta_n e_n, \quad (\lambda_i \geq 0; \textstyle\sum_i \lambda_i = 1; \delta_i \geq 0). \quad (1.9)$$

Consider the vector $z' = \lambda_1 w_1 + \cdots + \lambda_r w_r - \delta_1 e_1$. We set $T' = \text{conv}(w_1, \ldots, w_r)$ and $w_i = (w_{i1}, \ldots, w_{in})$. We claim that $z'$ is in $T'$. We may assume that $\delta_1 > 0$, $\lambda_i > 0$ for all $i$, and that the first entry $w_{i1}$ of $w_i$ is positive for $1 \leq i \leq s$ and is equal to zero for $i > s$. From Eq. (1.9) we get $\lambda_1 w_{11} + \cdots + \lambda_s w_{s1} \geq \delta_1$.

Case (I): $\lambda_1 w_{11} \geq \delta_1$. Then we can write

$$z' = \frac{\delta_1}{w_{11}}(w_1 - w_{11}e_1) + \left(\lambda_1 - \frac{\delta_1}{w_{11}}\right) w_1 + \lambda_2 w_2 + \cdots + \lambda_r w_r.$$

Notice that $w_1 - w_{11}e_1$ is again in $\{w_1, \ldots, w_r\}$. Thus $z'$ is a convex combination of $w_1, \ldots, w_r$, i.e., $z' \in T'$.

Case (II): $\lambda_1 w_{11} < \delta_1$. Let $m$ be the largest integer less than or equal to $s$ such that $\lambda_1 w_{11} + \cdots + \lambda_{m-1} w_{(m-1)1} < \delta_1 \leq \lambda_1 w_{11} + \cdots + \lambda_m w_{m1}$. Then

$$\begin{aligned} z' \;\; = \;\; & \sum_{i=1}^{m-1} \lambda_i (w_i - w_{i1}e_1) + \left[\frac{\delta_1}{w_{m1}} - \left(\sum_{i=1}^{m-1} \frac{\lambda_i w_{i1}}{w_{m1}}\right)\right](w_m - w_{m1}e_1) + \\ & \left[\lambda_m - \frac{\delta_1}{w_{m1}} + \left(\sum_{i=1}^{m-1} \frac{\lambda_i w_{i1}}{w_{m1}}\right)\right] w_m + \sum_{i=m+1}^{r} \lambda_i w_i. \end{aligned}$$

Notice that $w_i - w_{i1}e_1$ is again in $\{w_1, \ldots, w_r\}$ for $i = 1, \ldots, m$. Thus $z'$ is a convex combination of $w_1, \ldots, w_r$, i.e., $z' \in T'$. This completes the proof of the claim. Note that we can apply the argument above to any entry of $z$ or $z'$ thus we obtain that $z' - \delta_2 e_2 \in T'$. Thus by induction we obtain that $z \in T'$, as required. This completes the proof of Eq. (1.8).

Clearly one has the equality $P = \{z \mid z \geq 0; \langle z, w_i \rangle \leq 1 \, \forall i\}$ because for each $w_i$ there is $v_j$ such that $w_i \leq v_j$. Hence by the finite basis theorem [70] we can write

$$P = \{z \mid z \geq 0; \langle z, w_i \rangle \leq 1 \, \forall i\} = \text{conv}(\ell_0, \ell_1, \ldots, \ell_m) \tag{1.10}$$

for some $\ell_1, \ldots, \ell_m$ in $\mathbb{Q}_+^n$ and $\ell_0 = 0$. From Eq. (1.10) we readily get the equality

$$\{z \mid z \geq 0; \langle z, \ell_i \rangle \leq 1 \, \forall i\} = T(P). \tag{1.11}$$

Using Eq. (1.10) and noticing that $\langle \ell_i, w_j \rangle \leq 1$ for all $i, j$, we get

$$\mathbb{R}_+^n \cap (\text{conv}(\ell_0, \ldots, \ell_m) + \mathbb{R}_+\{-e_1, \ldots, -e_n\}) = \{z \mid z \geq 0; \langle z, w_i \rangle \leq 1 \, \forall i\}.$$

Hence using this equality and [70, Theorem 9.4] we obtain

$$\mathbb{R}_+^n \cap (\text{conv}(w_1, \ldots, w_r) + \mathbb{R}_+\{-e_1, \ldots, -e_n\}) = \{z \mid z \geq 0; \langle z, \ell_i \rangle \leq 1 \, \forall i\}. \tag{1.12}$$

Therefore by Eq. (1.8) together with Eqs. (1.11) and (1.12) we conclude that $T(P)$ is equal to $\text{conv}(w_1, \ldots, w_r)$, as required. $\square$

If $v_1, \ldots, v_q$ are $\{0, 1\}$-vectors, then the equality of Lemma 1.2.1 follows directly from [37, Theorem 8].

**Definition 1.2.2** Let $A$ be a matrix with entries in $\mathbb{N}$. The linear system $x \geq 0; xA \leq \mathbf{1}$ has the *integer rounding property* if

$$\lceil \min\{\langle y, \mathbf{1}\rangle |\, y \geq 0;\ Ay \geq a\}\rceil = \min\{\langle y, \mathbf{1}\rangle |\, Ay \geq a;\ y \in \mathbb{N}^q\}$$

for each integral vector $a$ for which $\min\{\langle y, \mathbf{1}\rangle |\, y \geq 0;\ Ay \geq a\}$ is finite.

If $a \in \mathbb{R}^n$, its *support* is given by $\mathrm{supp}(a) = \{i|\, a_i \neq 0\}$. Note that $a = a_+ - a_-$, where $a_+$ and $a_-$ are two non negative vectors with disjoint support called the *positive* and *negative* part of $a$ respectively.

**Remark 1.2.3** Let $A$ be a matrix with entries in $\mathbb{N}$. Then the linear system $x \geq 0; xA \leq \mathbf{1}$ has the *integer rounding property* if and only if

$$\lceil \min\{\langle y, \mathbf{1}\rangle |\, y \geq 0;\ Ay \geq a\}\rceil = \min\{\langle y, \mathbf{1}\rangle |\, Ay \geq a;\ y \in \mathbb{N}^q\}$$

for each vector $a \in \mathbb{N}^n$ for which $\min\{\langle y, \mathbf{1}\rangle |\, y \geq 0;\ Ay \geq a\}$ is finite. This follows from decomposing an integral vector $a$ as $a = a_+ - a_-$ and noticing that for $y \geq 0$ we have that $Ay \geq a$ if and only if $Ay \geq a_+$

**Theorem 1.2.4** ([4], [71, p. 82]) *Let $A$ be a non-negative integer matrix and let $P = \{x|\, x \geq 0;\ xA \leq \mathbf{1}\}$. The system $xA \leq \mathbf{1}; x \geq 0$ has the integer rounding property if and only if $T(P)$ has the integer decomposition property and all maximal integer vectors of $T(P)$ are columns of $A$ (maximal with respect to $\leq$).*

**Theorem 1.2.5** *Let $A$ be a matrix with entries in $\mathbb{N}$ and let $v_1, \ldots, v_q$ be the columns of $A$. If $w_1, \ldots, w_r$ is the set of all $\alpha \in \mathbb{N}^n$ such that $\alpha \leq v_i$ for some $i$, then the system $x \geq 0;\ xA \leq \mathbf{1}$ has the integer rounding property if and only if the subring $K[x^{w_1}t, \ldots, x^{w_r}t]$ is normal.*

**Proof.** Let $P = \{x|\, x \geq 0;\ xA \leq \mathbf{1}\}$ and let $T(P)$ be its antiblocking polyhedron. By Lemma 1.2.1 one has

$$T(P) = \mathrm{conv}(w_1, \ldots, w_r). \tag{1.13}$$

Let $\overline{S}$ be the integral closure of $S = K[x^{w_1}t, \ldots, x^{w_r}t]$ in its field of fractions. By the description of $\overline{S}$ given in Eq. (1.1) one has

$$\overline{S} = K[\{x^a t^b \mid (a, b) \in \mathbb{Z}\mathcal{B} \cap \mathbb{R}_+\mathcal{B}\}],$$

where $\mathcal{B} = \{(w_1, 1), \ldots, (w_r, 1)\}$. By Theorem 1.2.4 it suffices to prove that $S$ is normal if and only if $T(P)$ has the integer decomposition property and all maximal integer vectors of $T(P)$ are columns of $A$ (maximal with respect to $\leq$).

Assume that $S$ is normal, i.e., $S = \overline{S}$. Let $b$ be a natural number and let $a$ be an integer vector in $bT(P)$. Then using Eq. (1.13) it is seen that $(a, b)$ is in $\mathbb{R}_+\mathcal{B}$. Since $S$ is normal we have $\mathbb{R}_+\mathcal{B} \cap \mathbb{Z}\mathcal{B} = \mathbb{N}\mathcal{B}$. In our situation one has $\mathbb{Z}\mathcal{B} = \mathbb{Z}^{n+1}$. Hence $(a, b) \in \mathbb{N}\mathcal{B}$ and $a$ is the sum of $b$ integer vectors in $T(P)$. Thus $T(P)$ has the integer decomposition property. Assume that $a$ is a maximal integer vector of $T(P)$. It is not hard to see that $(a, 1)$ is in $\mathbb{R}_+\mathcal{B}$, i.e., $x^a t \in \overline{S} = S$. Thus $(a, 1)$ is a linear combination of vectors in $\mathcal{B}$ with coefficients in $\mathbb{N}$. Hence $(a, 1)$ is equal to $(w_j, 1)$ for some $j$. There exists $v_i$ such that $a = w_j \leq v_i$. Therefore by the maximality of $a$, we get $a = v_i$ for some $i$. Thus $a$ is a column of $A$ as required.

Conversely assume that $T(P)$ has the integer decomposition property and that all maximal integer vectors of $T(P)$ are columns of $A$. Let $x^a t^b \in \overline{S}$. Then $(a, b)$ is in the cone $\mathbb{R}_+\mathcal{B}$. Hence, using Eq. (1.13), we get $a \in bT(P)$. Thus $a = \alpha_1 + \cdots + \alpha_b$, where $\alpha_i$ is an integral vector of $T(P)$ for all $i$. Since each $\alpha_i$ is less than or equal to a maximal integer vector of $T(P)$, we get that $\alpha_i \in \{w_1, \ldots, w_r\}$. Then $x^a t^b \in S$. This proves that $S = \overline{S}$. □

Recall that $A(P) = K[\{x^a t^i \,|\, a \in \mathbb{Z}^n \cap iP; i \in \mathbb{N}\}] \subset R[t]$, is the Ehrhart ring of a lattice polytope $P$.

**Corollary 1.2.6** *Let $\mathcal{C}$ be a uniform clutter, let $A$ be its incidence matrix, and let $v_1, \ldots, v_q$ be the columns of $A$. If either system $x \geq 0; xA \leq \mathbf{1}$ or $x \geq 0; xA \geq \mathbf{1}$ has the integer rounding property and $P = \mathrm{conv}(v_1, \ldots, v_q)$, then*

$$K[x^{v_1}t, \ldots, x^{v_q}t] = A(P).$$

**Proof.** In general the subring $K[x^{v_1}t, \ldots, x^{v_q}t]$ is contained in $A(P)$. Assume that $x \geq 0; xA \leq \mathbf{1}$ has the integer rounding property and that every edge of $\mathcal{C}$ has $d$ elements. Let $w_1, \ldots, w_r$ be the set of all $\alpha \in \mathbb{N}^n$ such that $\alpha \leq v_i$ for some $i$. Then by Theorem 1.2.5 the subring $K[x^{w_1}t, \ldots, x^{w_r}t]$ is normal. Using that $v_1, \ldots, v_q$ is the set of $w_i$ with $|w_i| = d$, it is not hard to see that $A(P)$ is contained in $K[x^{v_1}t, \ldots, x^{v_q}t]$.

Assume that $x \geq 0; xA \geq \mathbf{1}$ has the integer rounding property. Let $I = I(\mathcal{C})$ be the edge ideal of $\mathcal{C}$ and let $R[It]$ be its Rees algebra. By Corollary 1.1.8, $R[It]$ is a normal domain. Since the clutter $\mathcal{C}$ is uniform the required equality follows at once from [28, Theorem 3.15]. $\Box$

The converse of Corollary 1.2.6 fails as the following example shows.

**Example 1.2.7** Let $\mathcal{C}$ be the uniform clutter with vertex set $X = \{x_1, \ldots, x_8\}$ and edge set

$$E(\mathcal{C}) = \{\{x_3, x_4, x_6, x_8\},\ \{x_2, x_5, x_6, x_7\},\ \{x_1, x_4, x_5, x_8\},\ \{x_1, x_2, x_3, x_8\}\}.$$

The characteristic vectors of the edges of $\mathcal{C}$ are

$$v_1 = (0,0,1,1,0,1,0,1),\quad v_2 = (0,1,0,0,1,1,1,0),$$
$$v_3 = (1,0,0,1,1,0,0,1),\quad v_4 = (1,1,1,0,0,0,0,1).$$

Let $A$ be the incidence matrix of $\mathcal{C}$ with column vectors $v_1, \ldots, v_4$ and let $P$ be the convex hull of $\{v_1, \ldots, v_4\}$. It is not hard to verify that the set

$$\{(v_1, 1), (v_2, 1), (v_3, 1), (v_4, 1)\}$$

is a Hilbert basis in the sense of [70]. Therefore we have the equality

$$K[x^{v_1}t, x^{v_2}t, x^{v_3}t, x^{v_4}t] = A(P).$$

Using Corollary 1.1.8 and [9, Theorem 2.12] it is seen that none of the two systems $x \geq 0; xA \leq \mathbf{1}$ and $x \geq 0; xA \geq \mathbf{1}$ have the integer rounding property.

A set $\mathcal{A} \subset \mathbb{Z}^n$ is called an *integral Hilbert basis* if $\mathbb{N}\mathcal{A} = \mathbb{R}_+\mathcal{A} \cap \mathbb{Z}^n$. Note that if $\mathcal{A}$ is an integral Hilbert basis, then the semigroup ring $K[\mathbb{N}\mathcal{A}]$ is normal.

**Proposition 1.2.8** *Let $I = (x^{v_1}, \ldots, x^{v_q})$ be a monomial ideal and let $v_i^* = \mathbf{1} - v_i$. Then $R[It]$ is normal if and only if the set*

$$\Gamma = \{-e_1, \ldots, -e_n, (v_1^*, 1), \ldots, (v_q^*, 1)\}$$

*is a Hilbert basis.*

**Proof.** Let $\mathcal{A}' = \{e_1, \ldots, e_n, (v_1, 1), \ldots, (v_q, 1)\}$. Assume that $R[It]$ is normal. Then $\mathcal{A}'$ is an integral Hilbert basis. Let $(a, b)$ be an integral vector in $\mathbb{R}_+\Gamma$, with $a \in \mathbb{Z}^n$ and $b \in \mathbb{Z}$. Then we can write

$$(a, b) = \mu_1(-e_1) + \cdots + \mu_n(-e_n) + \lambda_1(v_1^*, 1) + \cdots + \lambda_q(v_q^*, 1),$$

where $\mu_i \geq 0$ and $\lambda_j \geq 0$ for all $i, j$. Therefore

$$-(a, b) + b\mathbf{1}^* = \mu_1 e_1 + \cdots + \mu_n e_n + \lambda_1(v_1, -1) + \cdots + \lambda_q(v_q, -1),$$

where $\mathbf{1}^* = e_1 + \cdots + e_n$. This equality is equivalent to

$$-(a, -b) + b\mathbf{1}^* = \mu_1 e_1 + \cdots + \mu_n e_n + \lambda_1(v_1, 1) + \cdots + \lambda_q(v_q, 1).$$

As $\mathcal{A}'$ is an integral Hilbert basis we can write

$$-(a, -b) + b\mathbf{1}^* = \mu_1' e_1 + \cdots + \mu_n' e_n + \lambda_1'(v_1, 1) + \cdots + \lambda_q'(v_q, 1),$$

where $\mu_i' \in \mathbb{N}$ and $\lambda_j' \in \mathbb{N}$ for all $i, j$. Thus $(a, b) \in \mathbb{N}\Gamma$. This proves that $\Gamma$ is an integral Hilbert basis. The converse can be shown using similar arguments. $\Box$

A *clutter* $\mathcal{C}$ with finite vertex set $X = \{x_1, \ldots, x_n\}$ is a family of subsets of $X$, called edges, none of which is included in another. Let $f_1, \ldots, f_q$ be the edges of $\mathcal{C}$ and let $v_k = \sum_{x_i \in f_k} e_i$ be the *characteristic vector* of $f_k$. The *incidence matrix* of $\mathcal{C}$ is the $n \times q$ matrix with column vectors $v_1, \ldots, v_q$.

**Definition 1.2.9** *Let $A = (a_{ij})$ be a matrix with entries in $\{0, 1\}$. Its dual is the matrix $A^* = (a_{ij}^*)$, where $a_{ij}^* = 1 - a_{ij}$.*

The following duality is valid for incidence matrices of clutters. It will be used later to establish a duality theorem for monomial subrings.

**Theorem 1.2.10** *Let $A$ be the incidence matrix of a clutter and let $v_1, \ldots, v_q$ be its column vectors. If $v_i^* = \mathbf{1} - v_i$ and $A^*$ is the matrix with column vectors $v_1^*, \ldots, v_q^*$, then the system $x \geq 0; xA \geq \mathbf{1}$ has the integer rounding property if and only if the system $x \geq 0; xA^* \leq \mathbf{1}$ has the integer rounding property.*

**Proof.** Consider $Q = \{x | x \geq 0; xA \geq \mathbf{1}\}$ and $P^* = \{x | x \geq 0; xA^* \leq \mathbf{1}\}$. Let $w_1^*, \ldots, w_s^*$ be the set of all $\alpha \in \mathbb{N}^n$ such that $\alpha \leq v_i^*$ for some $i$. Then, using Lemmas 1.1.1 and 1.2.1, we obtain that the blocking polyhedron of $Q$ and the antiblocking polyhedron of $P^*$ are given by

$$B(Q) = \mathbb{R}_+^n + \mathrm{conv}(v_1, \ldots, v_q) \quad \text{and} \quad T(P^*) = \mathrm{conv}(w_1^*, \ldots, w_s^*)$$

respectively.

$\Rightarrow$) By Theorem 1.2.4 it suffices to show that $T(P^*)$ has the integer decomposition property and all maximal integer vectors of $T(P^*)$ are columns of $A^*$. Let $b$ be an integer and let $a$ be an integer vector in $bT(P^*)$. Then we can write

$$a = b(\lambda_1 w_1^* + \cdots + \lambda_s w_s^*), \quad (\textstyle\sum_i \lambda_i = 1; \lambda_i \geq 0).$$

For each $1 \leq i \leq s$ there is $v_{j_i}^*$ in $\{v_1^*, \ldots, v_q^*\}$ such that $w_i^* \leq v_{j_i}^*$. Thus for each $i$ we can write $\mathbf{1} - w_i^* = v_{j_i} + \delta_i$, where $\delta_i \in \mathbb{N}^n$. Therefore

$$\mathbf{1} - a/b = \lambda_1(v_{j_1} + \delta_1) + \cdots + \lambda_s(v_{j_s} + \delta_s).$$

This means that $\mathbf{1} - a/b \in B(Q)$, i.e., $b\mathbf{1} - a$ is an integer vector in $bB(Q)$. Hence by Theorem 1.1.7 we can write $b\mathbf{1} - a = \alpha_1 + \cdots + \alpha_b$ for some $\alpha_1, \ldots, \alpha_b$ integer vectors in $B(Q)$, and for each $\alpha_i$ there is $v_{k_i}$ in $\{v_1, \ldots, v_q\}$ such that $v_{k_i} \leq \alpha_i$. Thus $\alpha_i = v_{k_i} + \epsilon_i$ for some $\epsilon_i \in \mathbb{N}^n$ and consequently:

$$a = (\mathbf{1} - v_{k_1}) + \cdots + (\mathbf{1} - v_{k_b}) - c = v_{k_1}^* + \cdots + v_{k_b}^* - c,$$

where $c = (c_1, \ldots, c_n) \in \mathbb{N}^n$. Notice that $v_{k_1}^* + \cdots + v_{k_b}^* \geq c$ because $a \geq 0$. If $c_1 \geq 1$, then the first entry of $v_{k_i}^*$ is non-zero for some $i$ and we can write

$$a = v_{k_1}^* + \cdots + v_{k_{i-1}}^* + (v_{k_i}^* - e_1) + v_{k_{i+1}}^* + \cdots + v_{k_b}^* - (c - e_1)$$

Since $v_{k_i}^* - e_1$ is again in $\{w_1^*, \ldots, w_s^*\}$, we can apply this argument recursively to obtain that $a$ is the sum of $b$ integer vectors in $\{w_1^*, \ldots, w_s^*\}$. This proves that $T(P^*)$ has the integer decomposition property. Let $a$ be a maximal integer vector of $T(P^*)$. Since the vectors $w_1^*, \ldots, w_s^*$ have entries in $\{0, 1\}$, we get $T(P^*) \cap \mathbb{Z}^n = \{w_1^*, \ldots, w_s^*\}$. Then $a = w_i^*$ for some $i$. As $w_i^* \leq v_j^*$ for some $j$, we conclude that $a = v_j^*$, i.e., $a$ is a column of $A^*$, as required.

$\Leftarrow$) According Corollary 1.1.8, the system $x \geq 0; xA \geq \mathbf{1}$ has the integer rounding property if and only if $R[It]$ is normal. Thus by Proposition 1.2.8 we need only show that the set $\Gamma = \{-e_1, \ldots, -e_n, (v_1^*, 1), \ldots, (v_q^*, 1)\}$ is an integral Hilbert basis. Let $(a, b)$ be an integral vector in $\mathbb{R}_+\Gamma$, with $a \in \mathbb{Z}^n$ and $b \in \mathbb{Z}$. Then we can write

$$(a, b) = \mu_1(-e_1) + \cdots + \mu_n(-e_n) + \lambda_1(v_1^*, 1) + \cdots + \lambda_q(v_q^*, 1),$$

where $\mu_i \geq 0$, $\lambda_j \geq 0$ for all $i, j$. Hence $A^*\lambda \geq a$, where $\lambda = (\lambda_i)$. By hypothesis the system $x \geq 0; xA^* \leq \mathbf{1}$ has the integer rounding property. Then one has

$$b \geq \lceil \min\{\langle y, \mathbf{1}\rangle \,|\, y \geq 0;\ A^*y \geq a\} \rceil = \min\{\langle y, \mathbf{1}\rangle \,|\, A^*y \geq a;\ y \in \mathbb{N}^q\} = \langle y_0, \mathbf{1}\rangle$$

for some $y_0 = (y_i) \in \mathbb{N}^q$ such that $|y_0| = \langle y_0, \mathbf{1}\rangle \leq b$ and $a \leq A^* y_0$. Then

$$a = y_1 v_1^* + \cdots + y_q v_q^* - \delta_1 e_1 - \cdots - \delta_n e_n,$$

where $\delta_1, \ldots, \delta_n$ are in $\mathbb{N}$. Hence we can write

$$(a, b) = y_1(v_1^*, 1) + \cdots + y_{q-1}(v_{q-1}^*, 1) + (y_q + b - |y_0|)(v_q^*, 1) - (b - |y_0|)v_q^* - \delta,$$

where $\delta = (\delta_i)$. As the entries of $A^*$ are in $\mathbb{N}$, the vector $-v_q^*$ can be written as a non-negative integer combination of $-e_1, \ldots, -e_n$. Thus $(a, b) \in \mathbb{N}\Gamma$. This proves that $\Gamma$ is an integral Hilbert basis. $\Box$

We come to one of the main result of this section. It establishes a duality for monomial subrings.

**Theorem 1.2.11** *Let $A$ be the incidence matrix of a clutter, let $v_1, \ldots, v_q$ be its column vectors and let $v_i^* = \mathbf{1} - v_i$. If $w_1^*, \ldots, w_s^*$ is the set of all $\alpha \in \mathbb{N}^n$ such that $\alpha \leq v_i^*$ for some $i$, then the following conditions are equivalent:*

(a) *$R[It]$ is normal, where $I = (x^{v_1}, \ldots, x^{v_q})$.*

(b) *$S^* = K[x^{w_1^*}t, \ldots, x^{w_s^*}t]$ is normal.*

(c) *$\{-e_1, \ldots, -e_n, (v_1^*, 1), \ldots, (v_q^*, 1)\}$ is a Hilbert basis.*

(d) *$x \geq 0; xA \geq \mathbf{1}$ has the integer rounding property.*

(e) *$x \geq 0; xA^* \leq \mathbf{1}$ has the integer rounding property.*

**Proof.** (a) $\Leftrightarrow$ (c): This was shown in Proposition 1.2.8. (a) $\Leftrightarrow$ (d): This was shown in Corollary 1.1.8. (b) $\Leftrightarrow$ (e): This was shown in Theorem 1.2.5. (d) $\Leftrightarrow$ (e): This follows from Theorem 1.2.10. $\Box$

To illustrate the usefulness of this duality, below we show various results that follow from there.

**Definition 1.2.12** *Let $\mathcal{C}$ be a clutter on the vertex set $X = \{x_1, \ldots, x_n\}$. The edge ideal of $\mathcal{C}$, denoted by $I(\mathcal{C})$, is the ideal of $R$ generated by all monomials $x_e = \prod_{x_i \in e} x_i$ such that $e$ is an edge of $\mathcal{C}$. The dual $I^*$ of an edge ideal $I$ is the ideal of $R$ generated by all monomials $x_1 \cdots x_n / x_e$ such that $e$ is an edge of $\mathcal{C}$.*

**Definition 1.2.13** Let $A$ be an integral matrix and let $b$ an integral vector. The system $Ax \leq b$ is *totally dual integral* (TDI) if the minimum in the LP-duality equation

$$\max\{\langle c, x\rangle \,|\, Ax \leq b\} = \min\{\langle y, b\rangle \,|\, y \geq 0; yA = c\} \tag{1.14}$$

has an integral optimum solution $y$ for each integral vector $c$ for which the minimum is finite.

**Corollary 1.2.14** ([91, Theorem 2.10]) *Let $\mathcal{C}$ be a clutter and let $A$ be its incidence matrix. If $P = \{x | x \geq 0; xA \leq \mathbf{1}\}$ is an integral polytope and $I = I(\mathcal{C})$, then*

(i) *$R[I^*t]$ is normal.*

(ii) *$S = K[x^{w_1}t, \ldots, x^{w_r}t]$ is normal.*

**Proof.** Since $P$ has only integral vertices, by a result of Lovász [62] the system $x \geq 0; xA \leq \mathbf{1}$ is totally dual integral, i.e., the minimum in the LP-duality equation

$$\max\{\langle a, x\rangle \,|\, x \geq 0; xA \leq \mathbf{1}\} = \min\{\langle y, \mathbf{1}\rangle \,|\, y \geq 0; Ay \geq a\} \tag{1.15}$$

has an integral optimum solution $y$ for each integral vector $a$ with finite minimum. In particular the system $x \geq 0; xA \leq \mathbf{1}$ satisfies the integer rounding property. Therefore $R[I^*t]$ and $K[x^{w_1}t, \ldots, x^{w_r}t]$ are normal by Theorem 1.2.11. □

The normality assertion of part (ii) is well known and it can also be shown directly using the fact that the system $x \geq 0; xA \leq \mathbf{1}$ is TDI if $P$ is integral. Part (ii) is related to perfect graphs. Indeed if $P$ is integral, the $w_i's$ correspond to the cliques (complete subgraphs) of a perfect graph [14, 62].

**Corollary 1.2.15** *Let $B_1, \ldots, B_q$ be the collection of basis of a matroid $M$ with vertex set $X$ and let $v_1, \ldots, v_q$ be their characteristic vectors. If $A$ is the matrix with column vectors $v_1, \ldots, v_q$, then all systems*

$$x \geq 0; xA \geq \mathbf{1}, \quad x \geq 0; xA^* \geq \mathbf{1}, \quad x \geq 0; xA \leq \mathbf{1}, \quad x \geq 0; xA^* \leq \mathbf{1}$$

*have the integer rounding property.*

**Proof.** Consider the basis monomial ideal $I = (x^{v_1}, \ldots, x^{v_q})$ of the matroid $M$. By [67, Theorem 2.1.1], the collection of basis of the dual matroid $M^*$ of $M$ is given by $X \setminus B_1, \ldots, X \setminus B_q$. Now, the basis monomial ideal of a matroid is normal [92, Corollary 3.8], thus the result follows at once from the duality of Theorem 1.2.11. $\Box$

**Definition 1.2.16** A clutter $\mathcal{C}$ satisfies the *max-flow min-cut* (MFMC) property if both sides of the LP-duality equation

$$\min\{\langle a, x\rangle \,|\, x \geq 0; xA \geq \mathbf{1}\} = \max\{\langle y, \mathbf{1}\rangle \,|\, y \geq 0; Ay \leq a\} \tag{1.16}$$

have integral optimum solutions $x$ and $y$ for each non-negative integral vector $a$. The system $xA \geq \mathbf{1}$; $x \geq 0$ is called *totally dual integral* (TDI) if the maximum has an integral optimum solution $y$ for each integral vector $a$ with finite maximum.

**Corollary 1.2.17** ([45, Theorem 3.4]) *Let $A$ be the incidence matrix of a clutter $\mathcal{C}$ and let $I = I(\mathcal{C})$ be its edge ideal. Then $\mathcal{C}$ satisfies the max-flow min-cut property if and only if $Q(A)$ is integral and $R[It]$ is normal.*

**Proof.** Notice that if $\mathcal{C}$ has the max-flow min-cut property, then $Q(A)$ is integral [70, Corollary 22.1c]. Therefore the result follows directly from Eqs. (1.7), (1.16), and Theorem 1.2.11. $\Box$

We now turn our attention to the integer rounding property of systems of the form $xA \leq \mathbf{1}$.

**Definition 1.2.18** Let $A$ be a matrix with entries in $\mathbb{N}$. The system $xA \leq \mathbf{1}$ is said to have the *integer rounding property* if

$$\lceil \min\{\langle y, \mathbf{1}\rangle \,|\, y \geq 0;\ Ay = a\} \rceil = \min\{\langle y, \mathbf{1}\rangle \,|\, Ay = a;\ y \in \mathbb{N}^q\}$$

for each integral vector $a$ for which $\min\{\langle y, \mathbf{1}\rangle \,|\, y \geq 0;\ Ay = a\}$ is finite.

The next result is just a reinterpretation of an unpublished result of Giles and Orlin [70, Theorem 22.18] that characterizes the integer rounding property in terms of Hilbert bases.

**Proposition 1.2.19** *Let $v_1, \ldots, v_q$ be the column vectors of a non-negative integer matrix $A$ and let $A(P)$ be the Ehrhart ring of $P = \mathrm{conv}(0, v_1, \ldots, v_q)$. Then the system $xA \leq \mathbf{1}$ has the integer rounding property if and only if*

$$K[x^{v_1}t, \ldots, x^{v_q}t, t] = A(P).$$

**Proof.** By [70, Theorem 22.18], we have that the system $xA \leq \mathbf{1}$ has the integer rounding property if and only if the set $\mathcal{B} = \{(v_1, 1), \ldots, (v_q, 1), (0, 1)\}$ is an integral Hilbert basis. Thus the proposition follows readily by noticing the equality

$$A(P) = K[\{x^a t^b | (a, b) \in \mathbb{R}_+\mathcal{B} \cap \mathbb{Z}^{n+1}\}]$$

and the inclusion $K[x^{v_1}t, \ldots, x^{v_q}t, t] \subset A(P)$. $\Box$

**Theorem 1.2.20** *Let $\mathcal{A} = \{v_1, \ldots, v_q\}$ be the set of column vectors of a matrix $A$ with entries in $\mathbb{N}$. If the system $xA \leq \mathbf{1}$ has the integer rounding property, then*

(a) *$K[F]$ is normal, where $F = \{x^{v_1}, \ldots, x^{v_q}\}$, and*

(b) *$\mathbb{Z}^n/\mathbb{Z}\mathcal{A}$ is a torsion-free group.*

*The converse holds if $|v_i| = d$ for all $i$.*

**Proof.** For use below we set $\mathcal{B} = \{(v_1, 1), \ldots, (v_q, 1), (0, 1)\}$. First we prove (a). Let $x^a \in \overline{K[F]}$. Then $a \in \mathbb{Z}\mathcal{A}$ and we can write

$$a = \lambda_1 v_1 + \cdots + \lambda_q v_q,$$

for some $\lambda_1, \ldots, \lambda_q$ in $\mathbb{R}_+$. Hence

$$(a, \lceil \textstyle\sum_i \lambda_i \rceil) = \lambda_1(v_1, 1) + \cdots + \lambda_q(v_q, 1) + \delta(0, 1),$$

where $\delta \geq 0$. Therefore by Proposition 1.2.19, there are $\lambda'_1, \ldots \lambda'_q \in \mathbb{N}$ and $\delta' \in \mathbb{N}$ such that

$$(a, \lceil \textstyle\sum_i \lambda_i \rceil) = \lambda'_1(v_1, 1) + \cdots + \lambda'_q(v_q, 1) + \delta'(0, 1),$$

Thus $x^a \in K[F]$, as required. Next we show (b). From Proposition 1.2.19, we get

$$K[x^{v_1}t, \ldots, x^{v_q}t, t] = A(P).$$

Hence using [29, Theorem 3.9] we obtain that the group $M = \mathbb{Z}^{n+1}/\mathbb{Z}\mathcal{B}$ is torsion free. Let $\overline{a}$ be an element of $T(\mathbb{Z}^n/\mathbb{Z}\mathcal{A})$, the torsion subgroup of $\mathbb{Z}^n/\mathbb{Z}\mathcal{A}$. Thus there is a positive integer $s$ so that

$$sa = \lambda_1 v_1 + \cdots + \lambda_q v_q$$

for some $\lambda_1,\dots,\lambda_q$ in $\mathbb{Z}$. From the equality

$$s(a,|a|)=\lambda_1(v_1,1)+\cdots+\lambda_q(v_q,1)+(s|a|-\lambda_1-\cdots-\lambda_q)(0,1)$$

we obtain that the image of $(a,|a|)$ in $M$, denoted by $\overline{(a,|a|)}$, is a torsion element, i.e., $\overline{(a,|a|)}\in T(M)=(\overline{0})$. Hence it is readily seen that $a\in\mathbb{Z}\mathcal{A}$, i.e., $\overline{a}=\overline{0}$. Altogether we have $T(\mathbb{Z}^n/\mathbb{Z}\mathcal{A})=(0)$.

Conversely assume that $|v_i|=d$ for all $i$ and that (a) and (b) hold. We need only show that $\mathcal{B}$ is an integral Hilbert basis. Let $(a,b)$ be an integral vector in $\mathbb{R}_+\mathcal{B}$, where $a\in\mathbb{N}^n$ and $b\in\mathbb{N}$. Then we can write

$$(a,b)=\lambda_1(v_1,1)+\cdots+\lambda_q(v_q,1)+\mu(0,1), \tag{1.17}$$

for some $\lambda_1,\dots,\lambda_q,\mu$ in $\mathbb{Q}_+$. Hence using this equality together with (b) gives that $a$ is in $\mathbb{R}_+\mathcal{A}\cap\mathbb{Z}\mathcal{A}$. Hence $x^a\in\overline{K[F]}=K[F]$, i.e., $a\in\mathbb{N}\mathcal{A}$. Then we can write

$$a=\eta_1v_1+\cdots+\eta_qv_q$$

for some $\eta_1,\dots,\eta_q$ in $\mathbb{N}$. Since $|v_i|=d$ for all $i$, one has $\sum_i\lambda_i=\sum_i\eta_i$. Therefore using Eq. (1.17), we get $\mu\in\mathbb{N}$. Consequently from the equality

$$(a,b)=\eta_1(v_1,1)+\cdots+\eta_q(v_q,1)+\mu(0,1),$$

we conclude that $(a,b)\in\mathbb{N}\mathcal{B}$. This proves that $\mathcal{B}$ is an integral Hilbert basis. $\Box$

**Corollary 1.2.21** *Let $A$ be the incidence matrix of a connected graph $G$. Then the system $xA\leq\mathbf{1}$ has the integer rounding property if and only if $G$ is a bipartite graph.*

**Proof.** $\Rightarrow$) Let $\mathcal{A}=\{v_1,\dots,v_q\}$ be the set of columns of $A$. If $G$ is not bipartite, then according to [89, Corollary 3.4] one has $\mathbb{Z}^n/\mathbb{Z}\mathcal{A}\simeq\mathbb{Z}_2$, a contradiction to Theorem 1.2.20(b).

$\Leftarrow$) By [89, Theorem 2.15, Corollary 3.4] we get that the ring $K[x^{v_1},\dots,x^{v_q}]$ is normal and that $\mathbb{Z}^n/\mathbb{Z}\mathcal{A}\simeq\mathbb{Z}$. Thus by Theorem 1.2.20 the system $xA\leq\mathbf{1}$ has the integer rounding property, as required. This part of the proof also follows directly from the fact that the incidence matrix of a bipartite graph is totally unimodular. Indeed, since $A$ is totally unimodular, both problems of the LP-duality equation

$$\max\{\langle a,x\rangle\,|\,xA\leq\mathbf{1}\}=\min\{\langle y,\mathbf{1}\rangle\,|\,y\geq 0; Ay=a\}$$

have integral optimum solutions for each integral vector $a$ for which the minimum is finite, see [70, Corollary 19.1a]. Thus the system $xA \leq \mathbf{1}$ has the integer rounding property. □

**Corollary 1.2.22** *Let $A$ be the incidence matrix of a clutter $\mathcal{C}$. If $\mathcal{C}$ is uniform and has the max-flow min-cut property, then the system $xA \leq \mathbf{1}$ has the integer rounding property.*

**Proof.** Since all edges of $\mathcal{C}$ have the same number of elements, it suffices to observe that conditions (a) and (b) of Theorem 1.2.20 are satisfied because of [22, Theorem 3.6]. □

## 1.3 Rounding properties and graphs

Let $G$ be a graph with vertex set $X = \{x_1, \ldots, x_n\}$ and let $v_1, \ldots, v_q$ be the column vectors of the incidence matrix $A$ of $G$. If $G$ is a connected graph, we are able to prove that the system $x \geq 0; xA \leq \mathbf{1}$ has the integer rounding property if and only if the induced subgraph of the vertices of any two vertex disjoint odd cycles of $G$ is connected. Other equivalent descriptions of this property are also presented.

Let $R = K[x_1, \ldots, x_n]$ be a polynomial ring over a field $K$. We will examine the integer rounding property of the system $x \geq 0; xA \leq \mathbf{1}$ using the monomial subring:

$$S = K[x_1t, \ldots, x_nt, x^{v_1}t, \ldots, x^{v_q}t, t] \subset R[t],$$

where $t$ is a new variable. The main results of this section show that the system $x \geq 0; xA \leq \mathbf{1}$ has the integer rounding property if and only if any of the following equivalent conditions hold (see Theorems 1.3.2 and 1.3.3):

(a) $x \geq 0; xA \geq \mathbf{1}$ is a system with the integer rounding property.

(b) $R[x^{v_1}t, \ldots, x^{v_q}t]$ is normal, where $v_1, \ldots, v_q$ are the column vectors of $A$.

(c) $K[x^{v_1}t, \ldots, x^{v_q}t]$ is normal.

(d) $K[t, x_1t, \ldots, x_nt, x^{v_1}t, \ldots, x^{v_q}t]$ is normal.

(e) The induced subgraph of the vertices of any two vertex disjoint odd cycles of $G$ is connected.

The equivalences between (b), (c), and (e) follow from [74]. That (a) is equivalent to (b) is shown in Corollary 1.1.8. We are able to prove that the ring in (d) is isomorphic to the extended Rees algebra of the edge ideal of $G$ (see Lemma 1.3.1).

Let $I = I(G)$ be the edge ideal of $G$. Recall that the *extended Rees algebra* of $I$ is the subring

$$R[It, t^{-1}] := R[It][t^{-1}] \subset R[t, t^{-1}],$$

where $R[It]$ is the Rees algebra of $I$.

**Lemma 1.3.1** $R[It, t^{-1}] \simeq K[t, x_1t, \ldots, x_nt, x^{v_1}t, \ldots, x^{v_q}t]$.

**Proof.** We set $S = K[t, x_1t, \ldots, x_nt, x^{v_1}t, \ldots, x^{v_q}t]$. Note that $S$ and $R[It, t^{-1}]$ are both integral domains of the same Krull dimension, this follows from the dimension formula given in [77, Lemma 4.2]. Thus it suffices to prove that there is an epimorphism $\overline{\psi}\colon S \to R[It, I^{-1}]$ of $K$-algebras.

Let $u_0, u_1, \ldots, u_n, t_1, \ldots, t_q$ be a new set of variables and let $\varphi$, $\psi$ be the maps of $K$-algebras defined by the diagram

$$\begin{array}{ccc} K[u_0, u_1, \ldots, u_n, t_1, \ldots, t_q] & \xrightarrow{\psi} & R[It, t^{-1}] \\ \downarrow \varphi & \overset{\overline{\psi}}{\nearrow} & \\ S & & \end{array} \qquad \begin{array}{ll} u_0 \overset{\varphi}{\longmapsto} t, & u_0 \overset{\psi}{\longmapsto} t^{-1}, \\ u_i \overset{\varphi}{\longmapsto} x_it, & u_i \overset{\psi}{\longmapsto} x_i, \\ t_i \overset{\varphi}{\longmapsto} x^{v_i}t, & t_i \overset{\psi}{\longmapsto} x^{v_i}t. \end{array}$$

To complete the proof we will show that there is an epimorphism $\overline{\psi}$ of $K$-algebras that makes this diagram commutative, i.e., $\psi = \overline{\psi}\varphi$. To show the existence of $\overline{\psi}$ we need only show the inclusion $\ker(\varphi) \subset \ker(\psi)$. As $\ker(\varphi)$, being a toric ideal, is generated by binomials [77], it suffices to prove that any binomial of $\ker(\varphi)$ belongs to $\ker(\psi)$. Let

$$f = u_0^{a_0}u_1^{a_1}\cdots u_n^{a_n}t_1^{b_1}\cdots t_q^{b_q} - u_0^{c_0}u_1^{c_1}\cdots u_n^{c_n}t_1^{d_1}\cdots t_q^{d_q}$$

be a binomial in $\ker(\varphi)$. Then

$$t^{a_0}(x_1t)^{a_1}\cdots(x_nt)^{a_n}(x^{v_1}t)^{b_1}\cdots(x^{v_q}t)^{b_q} = t^{c_0}(x_1t)^{c_1}\cdots(x_nt)^{c_n}(x^{v_1}t)^{d_1}\cdots(x^{v_q}t)^{d_q}$$

Taking degrees in $t$ and $\underline{x} = \{x_1, \ldots, x_n\}$ we obtain

$$\begin{aligned} a_0 + (a_1 + \cdots + a_n) + (b_1 + \cdots + b_q) &= c_0 + (c_1 + \cdots + c_n) + (d_1 + \cdots + d_q), \\ a_1 + \cdots + a_n + 2(b_1 + \cdots + b_q) &= c_1 + \cdots + c_n + 2(d_1 + \cdots + d_q). \end{aligned}$$

Thus $-a_0 + b_1 + \cdots + b_q = -c_0 + d_1 + \cdots + d_q$, and we obtain the equality

$$t^{-a_0} x_1^{a_1} \cdots x_n^{a_n} (x^{v_1} t)^{b_1} \cdots (x^{v_q} t)^{b_q} = t^{-c_0} x_1^{c_1} \cdots x_n^{c_n} (x^{v_1} t)^{d_1} \cdots (x^{v_q} t)^{d_q},$$

i.e., $f \in \ker(\psi)$, as required. □

We come to the main result of this section.

**Theorem 1.3.2** *Let $G$ be a connected graph and let $A$ be its incidence matrix. Then the system*

$$x \geq 0; xA \leq \mathbf{1}$$

*has the integer rounding property if and only if the induced subgraph of the vertices of any two vertex disjoint odd cycles of $G$ is connected.*

**Proof.** Let $v_1, \ldots, v_q$ be the column vectors of $A$. According to [74, Theorem 1.1] (cf. [89, Corollary 3.10]), the subring $K[Gt] := K[x^{v_1}t, \ldots, x^{v_q}t]$ is normal if and only if the induced subgraph of the vertices of any two vertex disjoint odd cycles of $G$ is connected. Thus we need only show that $K[Gt]$ is normal if and only the system $x \geq 0; xA \leq \mathbf{1}$ has the integer rounding property. Let $I = I(G)$ be the edge ideal of $G$. Since $G$ is connected, by [74, Corollary 2.8] the subring $K[Gt]$ is normal if and only if the Rees algebra $R[It]$ of $I$ is normal. By a result of [54], $R[It]$ is normal if and only if $R[It, t^{-1}]$ is normal. By Lemma 1.3.1, $R[It, t^{-1}]$ is normal if and only if the subring

$$S = K[t, x_1 t, \ldots, x_n t, x^{v_1} t, \ldots, x^{v_q} t]$$

is normal. Thus we can apply Theorem 1.2.5 to conclude that $S$ is normal if and only if the system $x \geq 0; xA \leq \mathbf{1}$ has the integer rounding property. □

**Theorem 1.3.3** *Let $G$ be a connected graph and let $A$ be its incidence matrix. Then the system $x \geq 0; xA \leq \mathbf{1}$ has the integer rounding property if and only if any of the following equivalent conditions hold*

(a) $x \geq 0; xA \geq \mathbf{1}$ *is a system with the integer rounding property.*

(b) $R[It]$ *is a normal domain, where $I = I(G)$ is the edge ideal of $G$.*

(c) $K[x^{v_1}t, \ldots, x^{v_q}t]$ *is normal, where $v_1, \ldots, v_q$ are the column vectors of $A$*

(d) $K[t, x_1 t, \ldots, x_n t, x^{v_1} t, \ldots, x^{v_q} t]$ *is normal.*

**Proof.** According Corollary 1.1.8, the system $x \geq 0; xA \geq \mathbf{1}$ has the integer rounding property if and only if the Rees algebra $R[It]$ is normal. Thus the result follows from the proof of Theorem 1.3.2. $\Box$

**Corollary 1.3.4** *Let $G$ be a connected graph and let $I = I(G)$ be its edge ideal. Then $R[It]$ is normal if and only if $R[I^*t]$ is normal.*

**Proof.** By Theorem 1.3.3, the system $x \geq 0; xA \geq \mathbf{1}$ has the integer rounding property if and only if the system $x \geq 0; xA \leq \mathbf{1}$ does. Therefore the result follows at once using Theorem 1.2.11. $\Box$

This result is valid even if the graph is not connected. To prove it recall that (i) the Rees algebra $R[It]$ is normal if and only if the extended Rees algebra $R[It, t^{-1}]$ is normal, and (ii) $R[It, t^{-1}]$ is isomorphic to

$$K[t, x_1t, \ldots, x_nt, x^{v_1}t, \ldots, x^{v_q}t],$$

when $I$ is the edge ideal of a graph (see [23]). Then apply Theorem 1.2.11. The next example shows that Corollary 1.3.4 does not extend to arbitrary uniform clutters.

**Example 1.3.5** Consider the clutter $\mathcal{C}$ whose incidence matrix $A$ is the transpose of the matrix:

$$\begin{bmatrix} 0&0&1&1&0&1&1&1&1&1\\ 0&0&1&0&1&1&1&1&1&1\\ 0&1&1&0&0&1&1&1&1&1\\ 1&1&0&0&0&1&1&1&1&1\\ 0&1&1&0&1&0&1&1&1&1\\ 1&1&1&1&1&0&0&1&1&0\\ 1&1&1&1&1&0&0&1&0&1\\ 1&1&1&1&1&0&1&1&0&0\\ 1&1&1&1&1&1&1&0&0&0\\ 1&1&1&1&0&0&1&1&0&1 \end{bmatrix}$$

Let $I = I(\mathcal{C})$ be the edge ideal of $\mathcal{C}$. Note that all edges of $\mathcal{C}$ have 7 vertices. Using *Normaliz* [11] it is seen that $R[It]$ is normal and that $R[I^*t]$ is not normal.

Let $\mathcal{C}$ be a clutter with vertex set $X$. A vertex $x$ of $\mathcal{C}$ is called *isolated* if $x$ does not occur in any edge of $\mathcal{C}$. A subset $C \subset X$ is a *minimal vertex cover* of $\mathcal{C}$ if: ($c_1$) every edge of $\mathcal{C}$ contains at least one vertex of $C$, and ($c_2$) there is no proper subset of $C$ with the first property. If $C$ only satisfies condition ($c_1$), then $C$ is called a *vertex cover* of $\mathcal{C}$. The *Alexander dual* of $\mathcal{C}$, denoted by $\mathcal{C}^\vee$, is the clutter whose edges are the minimal vertex covers of $\mathcal{C}$. The edge ideal of $\mathcal{C}^\vee$, denoted by $I(\mathcal{C})^\vee$, is called the *Alexander dual* of $I(\mathcal{C})$. In combinatorial optimization the Alexander dual of a clutter is referred to as the *blocker* of the clutter [71].

**Proposition 1.3.6** *Let $G$ be a graph without isolated vertices and let $G'$ be its complement. Then $I(G')^\vee = I(G)^*$ if and only if $G$ is triangle free.*

**Proof.** $\Rightarrow$) Let $X = \{x_1, \ldots, x_n\}$ be the vertex set of $G$. Assume that $G$ has a triange $\mathcal{C}_3 = \{x_1, x_2, x_3\}$, i.e., $\{x_i, x_j\}$ are edges of $G$ for $1 \leq i < j \leq 3$. Clearly we may assume $n \geq 4$. Notice that $C' = \{x_4, \ldots, x_n\}$ is a vertex cover of $G'$, i.e., $x_4 \cdots x_n$ belongs to $I(G')^\vee$ and consequently it belongs to $I(G)^*$, a contradiction because $I(G)^*$ is generated by monomials of degree $n-2$.

$\Leftarrow$) Let $x^a = x_1 \cdots x_r$ be a minimal generator of $I(G')^\vee$. Then $C = \{x_1, \ldots, x_r\}$ is a minimal vertex cover of $G'$. Hence $X \setminus C$ is a maximal complete subgraph of $G$. Thus by hypothesis $X \setminus C$ is an edge of $G$, i.e., $x^a \in I(G)^*$. This proves the inclusion $I(G')^\vee \subset I(G)^*$. Conversely, let $x^a$ be a minimal generator of $I(G)^*$. There is an edge $\{x_1, x_2\}$ of $G$ such that $x^a = x_3 \cdots x_n$. Every edge of $G'$ must intersect $C = \{x_3, \ldots, x_n\}$, i.e., $x^a \in I(G')^\vee$. □

This formula applies for instance if $G$ is a bipartite graph.

**Corollary 1.3.7** *Let $G$ be a free triangle graph without isolated vertices. Then $R[I(G)t]$ is normal if and only if $R[I(G')^\vee t]$ is normal.*

**Proof.** It follows directly from Corollary 1.3.4 and Proposition 1.3.6. □

To the best of our knowledge the following is the first example of an edge ideal of a graph whose Alexander dual is not normal.

**Example 1.3.8** Let $G$ be the graph consisting of two vertex disjoint odd cycles of length 5 and let $G'$ be its complement. According to [74] the Rees algebra of $I(G)$ is not normal. Thus $R[I(G')^\vee t]$ is not normal by Corollary 1.3.7.

## 1.4 An expression for the canonical module

Let $R = K[x_1, \ldots, x_n]$ be a polynomial ring over an arbitrary field $K$ and let $K[F] = K[x^{v_1}, \ldots, x^{v_q}]$ be a homogeneous monomial subring, i.e., there exists $0 \neq x_0 \in \mathbb{Q}^n$ satisfying $\langle x_0, v_i \rangle = 1$ for all $i$. Then $K[F]$ is a standard graded $K$-algebra with the grading induced by declaring that a monomial $x^a \in K[F]$ has degree $i$ if and only if $\langle a, x_0 \rangle = i$. Recall that the *a-invariant* of $K[F]$, denoted by $a(K[F])$, is the degree as a rational function of the Hilbert series of $K[F]$, see for instance [87, p. 99]. Let $H$ and $\varphi$ be the Hilbert function and the Hilbert polynomial of $K[F]$ respectively. The index of regularity of $K[F]$, denoted by $\text{reg}(K[F])$, is the least positive integer such that $H(i) = \varphi(i)$ for $i \geq \text{reg}(K[F])$. The $a$-invariant plays a fundamental role in algebra and geometry because one has: $\text{reg}(K[F]) = 0$ if $a(K[F]) < 0$ and $\text{reg}(K[F]) = a(K[F]) + 1$ otherwise [87, Corollary 4.1.12].

If $K[F]$ is Cohen-Macaulay and $\omega_{K[F]}$ is the canonical module of $K[F]$, then

$$a(K[F]) = -\min\{\, i \mid (\omega_{K[F]})_i \neq 0\}, \tag{1.18}$$

see [10, p. 141] and [87, Proposition 4.2.3]. This formula applies if $K[F]$ is normal because normal monomial subrings are Cohen-Macaulay [56]. If $K[F]$ is normal, then by a formula of Danilov and Stanley (see [10, Theorem 6.3.5] and [18]) the canonical module of $K[F]$ is the ideal given by

$$\omega_{K[F]} = (\{x^a \mid a \in \mathbb{N}\mathcal{A} \cap (\mathbb{R}_+\mathcal{A})^{\circ}\}), \tag{1.19}$$

where $\mathcal{A} = \{v_1, \ldots, v_q\}$ and $(\mathbb{R}_+\mathcal{A})^{\circ}$ is the interior of $\mathbb{R}_+\mathcal{A}$ relative to $\text{aff}(\mathbb{R}_+\mathcal{A})$, the affine hull of $\mathbb{R}_+\mathcal{A}$.

The *dual cone* of $\mathbb{R}_+\mathcal{A}$ is the polyhedral cone given by

$$(\mathbb{R}_+\mathcal{A})^* = \{x \mid \langle x, y \rangle \geq 0;\ \forall\, y \in \mathbb{R}_+\mathcal{A}\}.$$

A set $\mathcal{H} \subset \mathbb{R}^n \setminus \{0\}$ is called an *integral basis* of $(\mathbb{R}_+\mathcal{A})^*$ if $(\mathbb{R}_+\mathcal{A})^* = \mathbb{R}_+\mathcal{H}$ and $\mathcal{H} \subset \mathbb{Z}^n$. Let $0 \neq a \in \mathbb{R}^n$. In what follows $H_a^+$ denotes the closed halfspace $H_a^+ = \{x \mid \langle x, a \rangle \geq 0\}$ and $H_a$ stands for the hyperplane through the origin with normal vector $a$.

The next result gives a general technique to compute the canonical module and the $a$-invariant of a wide class of monomial subrings. Another technique is given in [79]. In Section 1.5 we give some more precise expressions for the canonical module and the $a$-invariant of special families of monomial subrings arising from integer rounding properties.

**Theorem 1.4.1** *Let $c_1,\ldots,c_r$ be an integral basis of $(\mathbb{R}_+\mathcal{A})^*$ and let $b=(b_i)$ be the $\{0,-1\}$-vector given by $b_i=0$ if $\mathbb{R}_+\mathcal{A}\subset H_{c_i}$ and $b_i=-1$ if $\mathbb{R}_+\mathcal{A}\not\subset H_{c_i}$. If $\mathbb{N}\mathcal{A}=\mathbb{Z}^n\cap\mathbb{R}_+\mathcal{A}$ and $B$ is the matrix with column vectors $-c_1,\ldots,-c_r$, then*

(a) $\omega_{K[F]}=(\{x^a|\, a\in\mathbb{Z}^n\cap\{x|xB\leq b\}).$

(b) $a(K[F])=-\min\{\langle x_0,x\rangle|\; x\in\mathbb{Z}^n\cap\{x|\,xB\leq b\}\}.$

**Proof.** Let $\mathcal{H}=\{c_1,\ldots,c_r\}$. By duality [70, Corollary 7.1a], we have the equality

$$\mathbb{R}_+\mathcal{A}=H^+_{c_1}\cap\cdots\cap H^+_{c_r}. \tag{1.20}$$

Observe that $\mathbb{R}_+\mathcal{A}\cap H_{c_i}$ is a proper face if $b_i=-1$ and it is an improper face otherwise. From Eq. (1.20) we get that each facet of $\mathbb{R}_+\mathcal{A}$ has the form $\mathbb{R}_+\mathcal{A}\cap H_{c_i}$ for some $i$. The relative interior of the cone $\mathbb{R}_+\mathcal{A}$ is the union of its facets. Hence, using that $\mathcal{H}$ is an integral basis, we obtain the equality

$$\mathbb{Z}^n\cap(\mathbb{R}_+\mathcal{A})^{\mathrm{o}}=\mathbb{Z}^n\cap\{x|\,xB\leq b\}. \tag{1.21}$$

Now, part (a) follows readily from Eqs. (1.19) and (1.21). Part (b) follows from Eq. (1.18) and part (a). □

## 1.5 The canonical module and the $a$-invariant

In this section we give a description of the canonical module and the $a$-invariant for subrings arising from systems with the integer rounding property.

Let $A$ be a matrix of size $n\times q$ with entries in $\mathbb{N}$ such that $A$ has non-zero rows and non-zero columns. Let $v_1,\ldots,v_q$ be the columns of $A$. For use below consider the set $w_1,\ldots,w_r$ of all $\alpha\in\mathbb{N}^n$ such that $\alpha\leq v_i$ for some $i$. Let $R=K[x_1,\ldots,x_n]$ be a polynomial ring over a field $K$ and let

$$S=K[x^{w_1}t,\ldots,x^{w_r}t]\subset R[t]$$

be the subring of $R[t]$ generated by $x^{w_1}t,\ldots,x^{w_r}t$, where $t$ is a new variable. As $(w_i,1)$ lies in the hyperplane $x_{n+1}=1$ for all $i$, $S$ is a standard $K$-algebra. Thus a monomial $x^at^b$ in $S$ has degree $b$. In what follows we assume that $S$ has this grading. If $S$ is normal, then according to Eq. (1.19) the canonical module of $S$ is the ideal given by

$$\omega_S=(\{x^at^b|\,(a,b)\in\mathbb{N}\mathcal{B}\cap(\mathbb{R}_+\mathcal{B})^{\mathrm{o}}\}), \tag{1.22}$$

where $\mathcal{B} = \{(w_1, 1), \ldots, (w_r, 1)\}$ and $(\mathbb{R}_+\mathcal{B})^{\mathrm{o}}$ is the interior of $\mathbb{R}_+\mathcal{B}$ relative to $\mathrm{aff}(\mathbb{R}_+\mathcal{B})$, the affine hull of $\mathbb{R}_+\mathcal{B}$. In our case $\mathrm{aff}(\mathbb{R}_+\mathcal{B}) = \mathbb{R}^{n+1}$.

Let $\ell_0, \ell_1, \ldots, \ell_m$ be the vertices of $P = \{x \mid x \geq 0; xA \leq \mathbf{1}\}$, where $\ell_0 = 0$, and let $\ell_1, \ldots, \ell_p$ be the set of all maximal elements of $\ell_0, \ell_1, \ldots, \ell_m$ (maximal with respect to $<$).

The following lemma is not hard to prove.

**Lemma 1.5.1** *For each $1 \leq i \leq p$ there is a unique positive integer $d_i$ such that the non-zero entries of $(-d_i\ell_i, d_i)$ are relatively prime.*

*Notation* In what follows $\{\ell_1, \ldots, \ell_p\}$ is the set of maximal elements of $\{\ell_0, \ldots, \ell_m\}$ and $d_1, \ldots, d_p$ are the unique positive integers in Lemma 1.5.1.

The next result complements and generalizes a result of [23].

**Theorem 1.5.2** *If the system $x \geq 0; xA \leq \mathbf{1}$ has the integer rounding property, then the subring $S = K[x^{w_1}t, \ldots, x^{w_r}t]$ is normal, the canonical module of $S$ is given by*

$$\omega_S = \left(\left\{ x^a t^b \middle| \, (a, b) \begin{pmatrix} -d_1\ell_1 & \cdots & -d_p\ell_p & e_1 & \cdots & e_n \\ d_1 & \cdots & d_p & 0 & \cdots & 0 \end{pmatrix} \geq \mathbf{1} \right\}\right), \quad (1.23)$$

*and the $a$-invariant of $S$ is equal to $-\max_i\{\lceil 1/d_i + |\ell_i| \rceil\}$. Here $|\ell_i| = \langle \ell_i, \mathbf{1} \rangle$.*

**Proof.** Note that in Eq. (1.23) we regard $(-d_i\ell_i, d_i)$ and $e_j$ as column vectors for all $i, j$. The normality of $S$ follows from Theorem 1.2.5. Recall that we have the following duality (see Section 1.2):

$$\begin{aligned} P = \{x \mid x \geq 0; \langle x, w_i \rangle \leq 1\, \forall i\} &= \mathrm{conv}(\ell_0, \ell_1, \ldots, \ell_m), \\ \mathrm{conv}(w_1, \ldots, w_r) &= \{x \mid x \geq 0; \langle x, \ell_i \rangle \leq 1 \forall i\} = T(P), \end{aligned} \quad (1.24)$$

where $\{\ell_0, \ell_1, \ldots, \ell_m\} \subset \mathbb{Q}_+^n$ is the set of vertices of $P$ and $\ell_0 = 0$. Therefore using Eq. (1.24) and the maximality of $\ell_1, \ldots, \ell_p$ we obtain

$$\mathrm{conv}(w_1, \ldots, w_r) = \{x \mid x \geq 0; \langle x, \ell_i \rangle \leq 1, \, \forall \, i = 1, \ldots, p\}. \quad (1.25)$$

We set $\mathcal{B} = \{(w_1, 1), \ldots, (w_r, 1)\}$. Note that $\mathbb{Z}\mathcal{B} = \mathbb{Z}^{n+1}$. From Eq. (1.25) it is seen that

$$\mathbb{R}_+\mathcal{B} = H_{e_1}^+ \cap \cdots \cap H_{e_n}^+ \cap H_{(-d_1\ell_1, d_1)}^+ \cap \cdots \cap H_{(-d_p\ell_p, d_p)}^+. \quad (1.26)$$

Here $H_a^+$ denotes the closed halfspace $H_a^+ = \{x|\,\langle x,a\rangle \geq 0\}$ and $H_a$ stands for the hyperplane through the origin with normal vector $a$. Notice that

$$H_{e_1}\cap\mathbb{R}_+\mathcal{B},\ldots,H_{e_n}\cap\mathbb{R}_+\mathcal{B},H_{(-d_1\ell_1,d_1)}\cap\mathbb{R}_+\mathcal{B},\ldots,H_{(-d_p\ell_p,d_p)}\cap\mathbb{R}_+\mathcal{B}$$

are proper faces of $\mathbb{R}_+\mathcal{B}$. Hence from Eq. (1.26) we get that a vector $(a,b)$, with $a\in\mathbb{Z}^n$, $b\in\mathbb{Z}$, is in the relative interior of $\mathbb{R}_+\mathcal{B}$ if and only if the entries of $a$ are positive and $\langle(a,b),(-d_i\ell_i,d_i)\rangle\geq 1$ for all $i$. Thus the required expression for $\omega_S$, i.e., Eq. (1.23), follows using the normality of $S$ and the Danilov-Stanley formula given in Eq. (1.22).

It remains to prove the formula for $a(S)$, the $a$-invariant of $S$. Consider the vector $(\mathbf{1},b_0)$, where $b_0=\max_i\{\lceil 1/d_i+|\ell_i|\rceil\}$. Using Eq. (1.23), it is not hard to see (by direct substitution of $(\mathbf{1},b_0)$), that the monomial $x^{\mathbf{1}}t^{b_0}$ is in $\omega_S$. Thus from Eq. (1.18) we get $a(S)\geq -b_0$. Conversely if the monomial $x^at^b$ is in $\omega_S$, then again from Eq. (1.23) we get $\langle(-d_i\ell_i,d_i),(a,b)\rangle\geq 1$ for all $i$ and $a_i\geq 1$ for all $i$, where $a=(a_i)$. Hence

$$bd_i\geq 1+d_i\langle a,\ell_i\rangle\geq 1+d_i\langle\mathbf{1},\ell_i\rangle = 1+d_i|\ell_i|.$$

Since $b$ is an integer we obtain $b\geq\lceil 1/d_i+|\ell_i|\rceil$ for all $i$. Therefore $b\geq b_0$, i.e., $\deg(x^at^b)=b\geq b_0$. As $x^at^b$ was an arbitrary monomial in $\omega_S$, by the formula for the $a$-invariant of $S$ given in Eq. (1.18) we obtain that $a(S)\leq -b_0$. Altogether one has $a(S)=-b_0$, as required. □

**Theorem 1.5.3** *Assume that the system $x\geq 0$; $xA\leq\mathbf{1}$ has the integer rounding property. If $S=K[x^{w_1}t,\ldots,x^{w_r}t]$ is Gorenstein and $c_0=\max\{|\ell_i|\colon 1\leq i\leq p\}$ is an integer, then $|\ell_k|=c_0$ for each $1\leq k\leq p$ such that $\ell_k$ has integer entries.*

**Proof.** We proceed by contradiction. Assume that $|\ell_k|<c_0$ for some integer $1\leq k\leq p$ such that $\ell_k$ is integral. We may assume that $\ell_k=(1,\ldots,1,0,\ldots,0)$ and $|\ell_k|=s$. From Eq. (1.26) it follows that the monomial $x^{\ell_k}t^{s-1}$ cannot be in $S$ because $(\ell_k,s-1)$ does not belong to $H^+_{(-d_k\ell_k,d_k)}$. Consider the monomial $x^at^b$, where $a=\ell_k+\mathbf{1}$, $b=b_0+s-1$ and $b_0=-a(S)$. We claim that the monomial $x^at^b$ is in $\omega_S$. By Theorem 1.5.2 it suffices to show that $\langle(a,b),(-d_j\ell_j,d_j)\rangle\geq 1$ for $1\leq j\leq p$. Thus we need only show that $\langle(a,b),(-\ell_j,1)\rangle>0$ for $1\leq j\leq p$. From the proof of Theorem 1.5.2, it is seen that $-a(S)=\max_i\{\lfloor|\ell_i|\rfloor\}+1$. Hence we get $b_0=c_0+1$. One has the following equalities

$$\langle(a,b),(-\ell_j,1)\rangle = -|\ell_j|-\langle\ell_k,\ell_j\rangle+b_0+s-1 = -|\ell_j|-\langle\ell_k,\ell_j\rangle+c_0+s.$$

Set $\ell_j = (\ell_{j_1}, \ldots, \ell_{jn})$. From Eq. (1.26) we get that the entries of each $\ell_j$ are less than or equal to 1. Case (I): If $\ell_{ji} < 1$ for some $1 \leq i \leq s$, then $s - \langle \ell_k, \ell_j \rangle > 0$ and $c_0 \geq |\ell_j|$. Case (II): $\ell_{ji} = 1$ for $1 \leq i \leq s$. Then $\ell_j \geq \ell_k$. Thus by the maximality of $\ell_k$ we obtain $\ell_j = \ell_k$. In both cases we obtain $\langle (a, b), (-\ell_j, 1) \rangle > 0$, as required. Hence the monomial $x^a t^b$ is in $\omega_S$. Since $S$ is Gorenstein and $\omega_S$ is generated by $x^{\mathbf{1}} t^{b_0}$, we obtain that $x^a t^b$ is a multiple of $x^{\mathbf{1}} t^{b_0}$, i.e., $x^{\ell_k} t^{s-1}$ must be in $S$, a contradiction. $\Box$

**Theorem 1.5.4** *Assume that the system $x \geq 0; xA \leq \mathbf{1}$ has the integer rounding property. If $S = K[x^{w_1}t, \ldots, x^{w_r}t]$ and $-a(S) = 1/d_i + |\ell_i|$ for $i = 1, \ldots, p$, then $S$ is Gorenstein.*

**Proof.** We set $b_0 = -a(S)$ and $\mathcal{B} = \{(w_1, 1), \ldots, (w_r, 1)\}$. The ring $S$ is normal by Theorem 1.2.5. Since the monomial $x^1 t^{b_0} = x_1 \cdots x_n t^{b_0}$ is in $\omega_S$, we need only show that $\omega_S = (x^1 t^{b_0})$. Take $x^a t^b \in \omega_S$. It suffices to prove that $x^{a-\mathbf{1}} t^{b-b_0}$ is in $S$. Using Theorem 1.2.5, one has $\mathbb{R}_+\mathcal{B} \cap \mathbb{Z}^{n+1} = \mathbb{N}\mathcal{B}$. Thus we need only show that the vector $(a - \mathbf{1}, b - b_0)$ is in $\mathbb{R}_+\mathcal{B}$. From Eq. (1.26), the proof reduces to showing that the vector $(a - \mathbf{1}, b - b_0)$ is in $H^+_{(-\ell_i,1)}$ for $i = 1, \ldots, p$.

As $x^a t^b \in \omega_S$, from the description of $\omega_S$ given in Theorem 1.5.2 we get

$$\langle (a, b), (-d_i\ell_i, d_i) \rangle = -\langle a, d_i\ell_i \rangle + bd_i \geq 1 \implies -\langle a, \ell_i \rangle \geq -b + 1/d_i$$

for $i = 1, \ldots, p$. Therefore

$$\langle (a - \mathbf{1}, b - b_0), (-\ell_i, 1) \rangle = -\langle a, \ell_i \rangle + |\ell_i| + b - b_0 \geq -b + 1/d_i + |\ell_i| + b - b_0 = 0$$

for all $i$, as required. $\Box$

**Corollary 1.5.5** *If $P = \{x \mid x \geq 0; xA \leq \mathbf{1}\}$ is an integral polytope, then the monomial subring $S = K[x^{w_1}t, \ldots, x^{w_r}t]$ is Gorenstein if and only if $a(S) = -(|\ell_i| + 1)$ for $i = 1, \ldots, p$.*

**Proof.** Notice that if $P$ is integral, then $\ell_i$ has entries in $\{0, 1\}$ for $1 \leq i \leq p$ and consequently $d_i = 1$ for $1 \leq i \leq p$. Thus the result follows from Theorems 1.5.3 and 1.5.4. $\Box$

**Problem 1.5.6** If $A$ is the incidence matrix of a connected graph and the system $x \geq 0$; $xA \leq \mathbf{1}$ has the integer rounding property, then the subring $S = K[x^{w_1}t, \ldots, x^{w_r}t]$ is Gorenstein if and only if $-a(S) = 1/d_i + |\ell_i|$ for $i = 1, \ldots, p$.

Note that the answer to this problem is positive if $A$ is the incidence matrix of a bipartite graph because in this case $P$ is an integral polytope and we may apply Corollary 1.5.5. If $A$ is the incidence matrix of a graph, then it is seen that $d_i = 1$ or $d_i = 1/2$ for each $i$.

**Monomial subrings of cliques of perfect graphs** Let $S$ be a set of vertices of a graph $G$, the *induced subgraph* $\langle S\rangle$ is the maximal subgraph of $G$ with vertex set $S$. A *clique* of a graph $G$ is a subset of the set of vertices that induces a complete subgraph. Let $G$ be a graph with vertex set $X = \{x_1, \ldots, x_n\}$. A *colouring* of the vertices of $G$ is an assignment of colours to the vertices of $G$ in such a way that adjacent vertices have distinct colours. The *chromatic number* of $G$ is the minimal number of colours in a colouring of $G$. A graph is *perfect* if for every induced subgraph $H$, the chromatic number of $H$ equals the size of the largest complete subgraph of $H$. Let $S$ be a subset of the vertices of $G$. The set $S$ is called *independent* if no two vertices of $S$ are adjacent.

For use below we consider the empty set as a clique whose vertex set is empty. The *support* of a monomial $x^a$ is given by $\mathrm{supp}(x^a) = \{x_i \mid a_i > 0\}$. Note that $\mathrm{supp}(x^a) = \emptyset$ if and only if $a = 0$.

**Theorem 1.5.7** *Let $G$ be a perfect graph and let $S = K[x^{\omega_1}t, \ldots, x^{\omega_r}t]$ be the subring generated by all square-free monomials $x^a t$ such that $\mathrm{supp}(x^a)$ is a clique of $G$. Then the canonical module of $S$ is given by*

$$\omega_S = \left(\left\{ x^a t^b \middle| \, (a,b) \begin{pmatrix} -a_1 & \cdots & -a_m & e_1 & \cdots & e_n \\ 1 & \cdots & 1 & 0 & \cdots & 0 \end{pmatrix} \geq \mathbf{1} \right\}\right),$$

*where $a_1, \ldots, a_m$ are the characteristic vectors of the maximal independent sets of $G$, and the a-invariant of $S$ is equal to $-(\max_i\{|a_i|\} + 1)$.*

**Proof.** Let $v_1, \ldots, v_q$ be the set of characteristic vectors of the maximal cliques of $G$. Note that $w_1, \ldots, w_r$ is the set of all $\alpha \in \mathbb{N}^n$ such that $\alpha \leq v_i$ for some $i$. Since $G$ is a perfect graph, by [60, Theorem 16.14] we have the equality

$$P = \{x | x \geq 0; xA \leq \mathbf{1}\} = \mathrm{conv}(a_0, a_1, \ldots, a_p),$$

where $a_0 = 0$ and $a_1, \ldots, a_p$ are the characteristic vectors of the independent sets of $G$. We may assume that $a_1, \ldots, a_m$ correspond to the maximal independent sets of $G$. Furthermore, since $P$ has only integral vertices, by a result of Lovász [62] the system $x \geq 0; xA \leq \mathbf{1}$ is totally dual integral, i.e., the minimum in the LP-duality equation

$$\max\{\langle \alpha, x\rangle \,|\, x \geq 0; xA \leq \mathbf{1}\} = \min\{\langle y, \mathbf{1}\rangle \,|\, y \geq 0; Ay \geq \alpha\} \tag{1.27}$$

has an integral optimum solution $y$ for each integral vector $\alpha$ with finite minimum. In particular the system $x \geq 0; xA \leq \mathbf{1}$ satisfies the integer rounding property. Therefore the result follows readily from Theorem 1.5.2. □

**Proposition 1.5.8** Let $G$ be a perfect graph, let $w_1, \ldots, w_r$ be the incidence vectors of the cliques of $G$, and let $a_1, \ldots, a_m$ be the incidence vectors of the maximal independent sets of $G$. Then the set

$$\Gamma = \{(-a_1, 1), \ldots, (-a_m, 1), e_1, \ldots, e_n\}$$

is a Hilbert basis of $(\mathbb{R}_+\mathcal{B})^*$, where $\mathcal{B} = \{(w_1, 1), \ldots, (w_r, 1)\}$.

**Proof.** The incidence vector of $\emptyset$ is set to be equal to zero. As $G$ is perfect we have

$$\mathrm{conv}(w_1, \ldots, w_r) = \{x \,|\, x \geq 0; x(a_1 \cdots a_m) \leq \mathbf{1}\}.$$

Therefore it is seen that

$$\mathbb{R}_+(w_1, 1) + \cdots + \mathbb{R}_+(w_r, 1) = H_{e_1}^+ \cap \cdots \cap H_{e_n}^+ \cap \mathbb{R}_+(-a_1, 1) \cap \cdots \mathbb{R}_+(-a_m, 1).$$

Thus by duality we obtain that $(\mathbb{R}_+\mathcal{B})^* = \mathbb{R}_+\Gamma$. Using that the system $x \geq 0$; $x(a_1, \ldots, a_m) \leq \mathbf{1}$ is TDI it follows that $\mathbb{Z}^n \cap \mathbb{R}_+\Gamma = \mathbb{N}\Gamma$, i.e., $\Gamma$ is a Hilbert basis, as required. □

For use below recall that a graph $G$ is called *unmixed* if all maximal independent sets of $G$ have the same cardinality. Unmixed bipartite graphs have been nicely characterized in [68, 90].

**Corollary 1.5.9** *Let $G$ be a connected bipartite graph and let $I = I(G)$ be its edge ideal. Then the extended Rees algebra $R[It, t^{-1}]$ is a Gorenstein standard $K$-algebra if and only if $G$ is unmixed.*

**Proof.** Let $\omega_S$ be the canonical module of $S = R[It, t^{-1}]$. Recall that $S$ is Gorenstein if and only if $\omega_S$ is a principal ideal [10]. Since $G$ is a perfect graph, the result follows using Lemma 1.3.1 together with the description of the canonical module given in Theorem 1.5.7. □

**Subrings associated to the system $xA \leq \mathbf{1}$**

Let $A$ be a matrix with entries in $\mathbb{N}$ such that the system $xA \leq \mathbf{1}$ has integer rounding property. As before we assume that the rows and columns of $A$ are different from zero and that $v_1, \ldots, v_q$ are the columns of $A$. In what follows we assume that $|v_i| = d$ for all $i$.

The following lemma is not hard to show.

**Lemma 1.5.10** *If $|v_i| = d$ for all $i$. Then there are isomorphisms*

$$K[x^{v_1}t, \ldots, x^{v_q}t, t] \simeq K[x^{v_1}t, \ldots, x^{v_q}t][T] \text{ and } K[x^{v_1}t, \ldots, x^{v_q}t] \simeq K[x^{v_1}, \ldots, x^{v_q}]$$

*induced by $x^{v_i}t \mapsto x^{v_i}t$, $t \mapsto T$ and $x^{v_i}t \mapsto x^{v_i}$ respectively, where $T$ is a new variable.*

Let $S$ be a homogeneous monomial subring and let $P_S$ be its toric ideal. Recall that $S$ is called a *complete intersection* if $P_S$ is a complete intersection, i.e., $P_S$ can be generated by $\mathrm{ht}(P_S)$ binomials, where $\mathrm{ht}(P_S)$ is the height of $P_S$. Let $c$ be a cycle of a graph $G$. A *chord* of $c$ is any edge of $G$ joining two non adjacent vertices of $c$. A cycle without chords is called *primitive*.

**Proposition 1.5.11** *Let $G$ be a connected graph with $n$ vertices and $q$ edges and let $A$ be its incidence matrix. If the system $xA \leq \mathbf{1}$ has the integer rounding property, then $K[x^{v_1}t, \ldots, x^{v_q}t, t]$ is a complete intersection if and only if $G$ is bipartite and the number of primitive cycles of $G$ is equal to $q - n + 1$.*

**Proof.** $\Rightarrow$) By Corollary 1.2.21 the graph $G$ is bipartite. From Lemma 1.5.10 it follows that $K[x^{v_1}t, \ldots, x^{v_q}t, t]$ is a complete intersection if and only if $K[G] = K[x^{v_1}, \ldots, x^{v_q}]$ is a complete intersection. Therefore by [72] we get that $K[G]$ is a complete intersection if and only if the number of primitive cycles of $G$ is equal to $q - n + 1$.

$\Leftarrow$) By [72] the ring $K[G]$ is a complete intersection. Hence $K[x^{v_1}t, \ldots, x^{v_q}t, t]$ is a complete intersection by Lemma 1.5.10. $\square$

# Chapter 2

# Algebraic and Combinatorial Properties of Ideals and Algebras of Uniform Clutters of TDI Systems

Let $\mathcal{C}$ be a uniform clutter, i.e., all the edges of $\mathcal{C}$ have the same number of elements, and let $A$ be the incidence matrix of $\mathcal{C}$. We denote the column vectors of $A$ by $v_1, \ldots, v_q$. The *vertex covering number* of $\mathcal{C}$, denoted by $\alpha_0(\mathcal{C})$, is the smallest number of vertices in any minimal vertex cover of $\mathcal{C}$. The clutter obtained from $\mathcal{C}$ by deleting a vertex $x_i$ and removing all edges containing $x_i$ is denoted by $\mathcal{C} \setminus \{x_i\}$. A clutter $\mathcal{C}$ is called *vertex critical* if $\alpha_0(\mathcal{C} \setminus \{x_i\}) < \alpha_0(\mathcal{C})$ for all $i$. Under certain conditions we prove that $\mathcal{C}$ is vertex critical. If $\mathcal{C}$ satisfies the max-flow min-cut property, we prove that $A$ diagonalizes over $\mathbb{Z}$ to an identity matrix and that $v_1, \ldots, v_q$ is an integral Hilbert basis. We also show that if $\mathcal{C}$ has a perfect matching such that $\mathcal{C}$ has the packing property and $\alpha_0(\mathcal{C}) = 2$, then $A$ diagonalizes over $\mathbb{Z}$ to an identity matrix. If $A$ is a balanced matrix we prove that any regular triangulation of the cone generated by $v_1, \ldots, v_q$ is unimodular. Some examples are presented to show that our results only hold for uniform clutters. These results are closely related to certain algebraic properties, such as the normality or torsion freeness, of blowup algebras of edge ideals and to finitely generated abelian groups. They are also related to the theory of Gröbner bases of toric ideals and to Ehrhart rings.

Let $R = K[x_1, \ldots, x_n]$ be a polynomial ring over a field $K$ and let $I$ be an ideal of $R$ of height $g$ minimally generated by a finite set

$$F = \{x^{v_1}, \ldots, x^{v_q}\}$$

of square-free monomials. As usual we use the notation $x^a := x_1^{a_1} \cdots x_n^{a_n}$, where $a = (a_1, \ldots, a_n) \in \mathbb{N}^n$. The *support* of $x^a$ is given by $\mathrm{supp}(x^a) = \{x_i \,|\, a_i > 0\}$. For technical reasons we shall assume that each variable $x_i$ occurs in at least one monomial of $F$.

A *clutter* ($\mathcal{C}$) with finite vertex set $X$ is a family of subsets of $X$, called edges, none of which is included in another. The set of vertices and edges of $\mathcal{C}$ are denoted by $X = V(\mathcal{C})$ and $E(\mathcal{C})$ respectively. Clutters are special types of hypergraphs and are sometimes called *Sperner families* in the literature. One example of a clutter is a graph with the vertices and edges defined in the usual way for graphs. For a thorough study of clutters and hypergraphs from the point of view of combinatorial optimization see [16, 71].

We associate to the ideal $I$ a *clutter* $\mathcal{C}$ by taking the set of indeterminates $X = \{x_1, \ldots, x_n\}$ as vertex set and $E = \{f_1, \ldots, f_q\}$ as edge set, where $f_k$ is the support of $x^{v_k}$. The assignment $I \mapsto \mathcal{C}$ gives a natural one to one correspondence between the family of square-free monomial ideals and the family of clutters. The ideal $I$ is called the *edge ideal* of $\mathcal{C}$. To stress the relationship between $I$ and $\mathcal{C}$ we will use the notation $I = I(\mathcal{C})$. The $\{0,1\}$-vector $v_k$ is the so called *characteristic vector* or *incidence vector* of $f_k$, i.e., $v_k = \sum_{x_i \in f_k} e_i$, where $e_i$ is the *ith* unit vector. We shall always assume that $\mathcal{C}$ has no isolated vertices, i.e., each vertex $x_i$ occurs in at least one edge of $\mathcal{C}$.

Let $A$ be the *incidence matrix* of $\mathcal{C}$ whose column vectors are $v_1, \ldots, v_q$. The *set covering polyhedron* of $\mathcal{C}$ is given by:

$$Q(A) = \{x \in \mathbb{R}^n |\, x \geq 0;\, xA \geq \mathbf{1}\},$$

where $\mathbf{1} = (1, \ldots, 1)$. A subset $C \subset X$ is called a *minimal vertex cover* of $\mathcal{C}$ if: (i) every edge of $\mathcal{C}$ contains at least one vertex of $C$, and (ii) there is no proper subset of $C$ with the first property. The map $C \mapsto \sum_{x_i \in C} e_i$ gives a bijection between the minimal vertex covers of $\mathcal{C}$ and the integral vectors of $Q(A)$. A polyhedron is called an *integral polyhedron* if it has only integral vertices. A clutter is called *d-uniform* or *uniform* if all its edges have exactly $d$ vertices.

We begin in Section 2.1 by introducing various combinatorial properties of clutters. We then give a simple combinatorial proof of the following result of [42]:

**Proposition 2.1.2** *If $\mathcal{C}$ is a $d$-uniform clutter whose set covering polyhedron $Q(A)$ is integral, then there are $X_1,\ldots,X_d$ mutually disjoint minimal vertex covers of $\mathcal{C}$ such that $X=\cup_{i=1}^{d}X_i$. In particular $|\mathrm{supp}(x^{v_i})\cap X_k|=1$ for all $i,k$.*

The original proof of this result was algebraic. It was based on the fact that the radical of the ideal $IR[It]$ can be expressed in terms of the minimal primes of $I$, where $R[It]$ is the Rees algebra of $I$ (see Section 2.2). Example 2.1.3 shows that this result fails if we drop the uniformity hypothesis. For use below we denote the smallest number of vertices in any minimal vertex cover of $\mathcal{C}$ by $\alpha_0(\mathcal{C})$. A set of pairwise disjoint edges of $\mathcal{C}$ is called *independent* or a *matching* and a set of independent edges of $\mathcal{C}$ whose union is $X$ is called a *perfect matching*. We then prove:

**Proposition 2.1.6** *Let $\mathcal{C}$ be a $d$-uniform clutter with a perfect matching such that $Q(A)$ is integral. Then $\mathcal{C}$ is vertex critical.*

A simple example is shown to see that this result fails for non uniform clutters with integral set covering polyhedron (Remark 2.1.7).

In Section 2.2 we introduce Rees algebras, Ehrhart rings, and edge subrings. Certain algebraic properties of these graded algebras such as the normality and torsion–freeness are related to combinatorial optimization properties of clutters such as the max-flow min-cut property (see Definition 2.2.3) and the integrality of $Q(A)$ [42, 45]. This relation between algebra and combinatorics will be quite useful here. In Theorems 2.2.2, 2.2.4, and Proposition 2.2.5 we summarize the algebro-combinatorial facts needed to show the main result of Section 2.2:

**Theorem 2.2.6** *If $\mathcal{C}$ is a uniform clutter with the max-flow min-cut property, then*

(a) $\Delta_r(A)=1$, *where* $r=\mathrm{rank}(A)$.

(b) $\mathbb{N}\mathcal{A}=\mathbb{R}_+\mathcal{A}\cap\mathbb{Z}^n$, *where* $\mathcal{A}=\{v_1,\ldots,v_q\}$.

Here $\Delta_r(A)$ denotes the greatest common divisor of all the nonzero $r\times r$ sub-determinants of $A$, $\mathbb{N}\mathcal{A}$ denotes the semigroup generated by $\mathcal{A}$, and $\mathbb{R}_+\mathcal{A}$ denotes the cone generated by $\mathcal{A}$. Condition (b) means that $\mathcal{A}$ is a *Hilbert basis* for $\mathbb{R}_+\mathcal{A}$. As interesting consequences we obtain that if $\mathcal{C}$ is a $d$-uniform clutter with the max-flow min-cut property, then $A$ diagonalizes over $\mathbb{Z}$—using row and column operations—to an identity matrix (see Corollary 2.2.8) and $\mathcal{C}$ has a perfect matching if and only if $n = d\alpha_0(\mathcal{C})$ (see Corollary 2.2.10). In Example 2.2.7 we show that the uniformity hypothesis is essential in the two statements of Theorem 2.2.6.

Section 2.3 deals with the diagonalization problem (see Conjecture 2.2.16) for clutters with the packing property (see Definition 2.2.11). The following is one of the main results of this section. It gives some support to Conjecture 2.2.16.

**Theorem 2.3.1** *Let $\mathcal{C}$ be a $d$-uniform clutter with a perfect matching such that $\mathcal{C}$ has the packing property and $\alpha_0(\mathcal{C}) = 2$. If $A$ has rank $r$, then*

$$\Delta_r \begin{pmatrix} A \\ \mathbf{1} \end{pmatrix} = 1.$$

As an application we obtain the next result which gives some support to a Conjecture of Conforti and Cornuéjols [15] (see Conjecture 2.2.14).

**Corollary 2.3.3** *Let $\mathcal{C}$ be a $d$-uniform clutter with a perfect matching such that $\mathcal{C}$ has the packing property and $\alpha_0(\mathcal{C}) = 2$. If $v_1, \ldots, v_q$ are linearly independent, then $\mathcal{C}$ has the max-flow min-cut property.*

All clutters of Section 2.3 satisfy the hypotheses of Theorem 2.3.4. We call this type of clutter 2-*partitionable* (see Example 2.3.5). They occur naturally in the theory of blockers of unmixed bipartite graphs (see Corollary 2.3.6). The other main result of Section 2.3 is about some of the properties of this family of clutters:

**Theorem 2.3.4.** *Let $\mathcal{C}$ be a $d$-uniform clutter with a partition $X_1, \ldots, X_d$ of $X$ such that $X_i$ is a minimal vertex cover of $\mathcal{C}$ and $|X_i| = 2$ for all $i$. Then* (a) $\mathrm{rank}(A) \leq d + 1$. (b) *If $C$ is a minimal vertex cover of $\mathcal{C}$, then $2 \leq |C| \leq d$.* (c) *If $\mathcal{C}$ satisfies the König property and there is a minimal vertex cover $C$ of $\mathcal{C}$ with $|C| = d \geq 3$, then* $\mathrm{rank}(A) = d + 1$. (d) *If $I = I(\mathcal{C})$ is minimally non-normal and $\mathcal{C}$ satisfies the packing property, then* $\mathrm{rank}(A) = d + 1$.

Regular and unimodular triangulations are introduced in Section 2.4. There is a relationship between the Gröbner bases of the toric ideal of $K[F]$ and the triangulations of $\mathcal{A}$, which has many interesting applications. We make use of the theory of Gröbner bases and convex polytopes, which was created and developed by Sturmfels [77], to prove the following main result of Section 2.4:

**Theorem 2.4.6** *Let $A$ be a balanced matrix with distinct column vectors $v_1, \ldots, v_q$. If $|v_i| = d$ for all $i$, then any regular triangulation of the cone $\mathbb{R}_+\{v_1, \ldots, v_q\}$ is unimodular.*

Here $|v_i|$ denotes the sum of the entries of the vector $v_i$. Recall that a matrix $A$ with entries in $\{0,1\}$ is called *balanced* if $A$ has no square submatrix of odd size with exactly two 1's in each row and column. If we do not require the uniformity condition $|v_i| = d$ for all $i$ this result is false, as is seen in Example 2.4.7. What makes this result surprising is the fact that not all balanced matrices are unimodular (see Example 2.4.7). This result gives some support to Conjecture 2.4.3: If $\mathcal{C}$ is a uniform clutter that satisfies the max-flow min-cut property, then the rational polyhedral cone $\mathbb{R}_+\{v_1,\ldots,v_q\}$ has a unimodular regular triangulation.

Throughout the chapter we introduce most of the notions that are relevant for our purposes. For unexplained terminology and notation we refer to [71] (for the theory of combinatorial optimization) and [10, 84] (for the theory of blowup algebras and integral closures). See [65] for additional information about commutative rings and ideals.

## 2.1 Structure of ideals of uniform clutters

We continue to use the notation and definitions used in the introduction. In what follows $\mathcal{C}$ denotes a $d$-uniform clutter with vertex set $X = \{x_1,\ldots,x_n\}$, edge set $E(\mathcal{C})$, edge ideal $I = I(\mathcal{C})$, and incidence matrix $A$. The column vectors of $A$ are denoted by $v_1,\ldots,v_q$ and the edge ideal of $\mathcal{C}$ is given by $I = (x^{v_1},\ldots,x^{v_q})$.

In this section we study the structure of uniform clutters whose set covering polyhedron is integral. Examples of this type of clutter include bipartite graphs and uniform clutters with the max-flow min-cut property. A clutter whose set covering polyhedron is integral is called *ideal* in the literature [16]. We denote the smallest number of vertices in any minimal vertex cover of $\mathcal{C}$ by $\alpha_0(\mathcal{C})$ and the maximum number of independent edges of $\mathcal{C}$ by $\beta_1(\mathcal{C})$. These two numbers are called the *vertex covering number* and the *edge independence number* respectively. Notice that in general $\beta_1(\mathcal{C}) \leq \alpha_0(\mathcal{C})$. If equality occurs we say that $\mathcal{C}$ has the *König property.*

Recall that $\mathfrak{p}$ is a minimal prime of $I = I(\mathcal{C})$ if and only if $\mathfrak{p} = (C)$ for some minimal vertex cover $C$ of $\mathcal{C}$ [87, Proposition 6.1.16]. Thus the primary decomposition of the edge ideal of $\mathcal{C}$ is given by

$$I(\mathcal{C}) = (C_1) \cap (C_2) \cap \cdots \cap (C_s),$$

where $C_1,\ldots,C_s$ are the minimal vertex covers of $\mathcal{C}$ and $(C_i)$ denotes the prime ideal of $R$ generated by $C_i$. In particular observe that $\operatorname{ht} I(\mathcal{C})$, the height of $I(\mathcal{C})$, equals the number of vertices in a minimum vertex cover of $\mathcal{C}$, i.e., $\operatorname{ht} I(\mathcal{C}) = \alpha_0(\mathcal{C})$. This is a hint of the rich interaction between the combinatorics of $\mathcal{C}$ and the algebra of $I(\mathcal{C})$.

The next result was shown in [42] using commutative algebra methods. Here we give a simple combinatorial proof.

**Lemma 2.1.1** ([42]) *If $\mathcal{C}$ is a $d$-uniform clutter such that $Q(A)$ is integral, then there exists a minimal vertex cover of $\mathcal{C}$ intersecting every edge of $\mathcal{C}$ in exactly one vertex.*

**Proof.** Let $B$ be the integral matrix whose columns are the vertices of $Q(A)$. It is not hard to show that *a vector $\alpha \in \mathbb{R}^n$ is an integral vertex of $Q(A)$ if and only if $\alpha = \sum_{x_i \in C} e_i$ for some minimal vertex cover $C$ of $\mathcal{C}$*. Thus the columns of $B$ are the characteristic vectors of the minimal vertex covers of $\mathcal{C}$. Using [16, Theorem 1.17] we get that

$$Q(B) = \{x \mid x \geq 0; xB \geq \mathbf{1}\}$$

is an integral polyhedron whose vertices are the columns of $A$. Therefore we have the equality

$$Q(B) = \mathbb{R}^n_+ + \operatorname{conv}(v_1, \ldots, v_q). \tag{2.1}$$

We proceed by contradiction. Assume that for each column $u_k$ of $B$ there exists a vector $v_{i_k}$ in $\{v_1, \ldots, v_q\}$ such that $\langle v_{i_k}, u_k \rangle \geq 2$. Then

$$v_{i_k} B \geq \mathbf{1} + e_k.$$

Consider the vector $\alpha = v_{i_1} + \cdots + v_{i_s}$, where $s$ is the number of columns of $B$. From the inequality

$$\alpha B \geq (\mathbf{1} + e_1) + \cdots + (\mathbf{1} + e_s) = (s+1, \ldots, s+1)$$

we obtain that $\alpha/(s+1) \in Q(B)$. Thus, using Eq. (2.1), we can write

$$\alpha/(s+1) = \mu_1 e_1 + \cdots + \mu_n e_n + \lambda_1 v_1 + \cdots + \lambda_q v_q \quad (\mu_i, \lambda_j \geq 0;\ \textstyle\sum \lambda_i = 1). \tag{2.2}$$

Therefore taking inner products with $\mathbf{1}$ in Eq. (2.2) and using the fact that $\mathcal{C}$ is $d$-uniform we get that $|\alpha| \geq (s+1)d$. Then using the equality $\alpha = v_{i_1} + \cdots + v_{i_s}$ we conclude

$$sd = |v_{i_1}| + \cdots + |v_{i_s}| = |\alpha| \geq (s+1)d,$$

a contradiction because $d \geq 1$. $\Box$

A graph $G$ is called *strongly perfect* if every induced subgraph $H$ of $G$ has a maximal independent set of vertices $F$ such that $|F \cap K| = 1$ for any maximal clique $K$ of $H$. Bipartite and chordal graphs are strongly perfect. If $A$ is the vertex-clique matrix of $G$, then $G$ being strongly perfect implies that the clique polytope of $G$, $\{x | \, x \geq 0; \, xA \leq \mathbf{1}\}$, has a vertex that intersects every maximal clique. In this sense, uniform clutters such that $Q(A)$ is integral can be thought of as being analogous to strongly perfect graphs.

The notion of a minor plays a prominent role in combinatorial optimization [16]. Recall that a proper ideal $I'$ of $R$ is called a *minor* of $I = I(\mathcal{C})$ if there is a subset

$$X' = \{x_{i_1}, \ldots, x_{i_r}, x_{j_1}, \ldots, x_{j_s}\}$$

of the set of variables $X$ such that $I'$ is obtained from $I$ by making $x_{i_k} = 0$ and $x_{j_\ell} = 1$ for all $k, \ell$. Notice that a set of generators $x^{v_1'}, \ldots, x^{v_q'}$ of $I'$ is obtained from a set of generators $x^{v_1}, \ldots, x^{v_q}$ of $I$ by making $x_{i_k} = 0$ and $x_{j_\ell} = 1$ for all $k, \ell$. The ideal $I$ is considered itself a minor. A clutter $\mathcal{C}'$ is called a *minor* of $\mathcal{C}$ if $\mathcal{C}'$ corresponds to a minor $I'$ of $I$ under the correspondence between square-free monomial ideals and clutters. This terminology is consistent with that of [16, p. 23]. Also notice that $\mathcal{C}'$ is obtained from $I'$ by considering the unique set of square-free monomials of $R$ that minimally generate $I'$. The clutter $\mathcal{C} \setminus \{x_i\}$ corresponds to the ideal $I'$ obtained from $I$ by making $x_i = 0$, i.e., $\mathcal{C} \setminus \{x_i\}$ is a special type of a minor which is called a *deletion*.

The notion of a minor of a clutter is not a generalization of the notion of a minor of a graph in the sense of graph theory [71, p. 25]. For instance if $G$ is a cycle of length four and we contract an edge we obtain that a triangle is a minor of $G$, but a triangle cannot be a minor of $G$ in our sense.

**Proposition 2.1.2** *If $\mathcal{C}$ is a $d$-uniform clutter whose set covering polyhedron $Q(A)$ is integral, then there are $X_1, \ldots, X_d$ mutually disjoint minimal vertex covers of $\mathcal{C}$ such that $X = \cup_{i=1}^{d} X_i$.*

**Proof.** By induction on $d$. If $d = 1$, then $E(\mathcal{C}) = \{\{x_1\}, \ldots, \{x_n\}\}$ and $X$ is a minimal vertex cover of $\mathcal{C}$. In this case we set $X_1 = X$. Assume $d \geq 2$. By Lemma 2.1.1 there is a minimal vertex cover $X_1$ of $\mathcal{C}$ such that $|\mathrm{supp}(x^{v_i}) \cap X_1| = 1$ for all $i$. Consider the ideal $I'$ obtained from $I$ by making $x_i = 1$ for $x_i \in X_1$. Let $\mathcal{C}'$ be the clutter corresponding to $I'$ and let $A'$ be the incidence matrix of $\mathcal{C}'$. The ideal $I'$ (resp. the clutter $\mathcal{C}'$) is a minor of $I$ (resp.

$\mathcal{C}$). Recall that the integrality of $Q(A)$ is preserved under taking minors [71, Theorem 78.2], so $Q(A')$ is integral. Then $\mathcal{C}'$ is a $(d-1)$-uniform clutter whose set covering polyhedron $Q(A')$ is integral. Note that $V(\mathcal{C}') = X\setminus X_1$. Therefore by induction hypothesis there are $X_2, \ldots, X_d$ pairwise disjoint minimal vertex covers of $\mathcal{C}'$ such that $X \setminus X_1 = X_2 \cup \cdots \cup X_d$. To complete the proof observe that $X_2, \ldots, X_d$ are minimal vertex covers of $\mathcal{C}$. Indeed if $e$ is an edge of $\mathcal{C}$ and $2 \leq k \leq d$, then $e \cap X_1 = \{x_i\}$ for some $i$. Since $e \setminus \{x_i\}$ is an edge of $\mathcal{C}'$, we get $(e \setminus \{x_i\}) \cap X_k \neq \emptyset$. Hence $X_k$ is a vertex cover of $\mathcal{C}$. Furthermore if $x \in X_k$, then by the minimality of $X_k$ relative to $\mathcal{C}'$ there is an edge $e'$ of $\mathcal{C}'$ disjoint from $X_k \setminus \{x\}$. Since $e = e' \cup \{y\}$ is an edge of $\mathcal{C}$ for some $y \in X_1$, we obtain that $e$ is an edge of $\mathcal{C}$ disjoint from $X_k \setminus \{x\}$. Therefore $X_k$ is a minimal vertex cover of $\mathcal{C}$, as required. □

**Example 2.1.3** Consider the clutter $\mathcal{C}$ with vertex set $X = \{x_1, \ldots, x_9\}$ whose edges are

$$\begin{array}{llll} f_1 = \{x_1, x_2\}, & f_2 = \{x_3, x_4, x_5, x_6\}, & f_3 = \{x_7, x_8, x_9\}, & \\ f_4 = \{x_1, x_3\}, & f_5 = \{x_2, x_4\}, & f_6 = \{x_5, x_7\}, & f_7 = \{x_6, x_8\}. \end{array}$$

In this example $Q(A)$ is integral because the incidence matrix of $\mathcal{C}$ is a balanced matrix. However $|C \cap f_i| \geq 2$ for any minimal vertex cover $C$ and for some $i$. Thus the uniformity hypothesis is essential in Proposition 2.1.2. This clutter can be represented as:

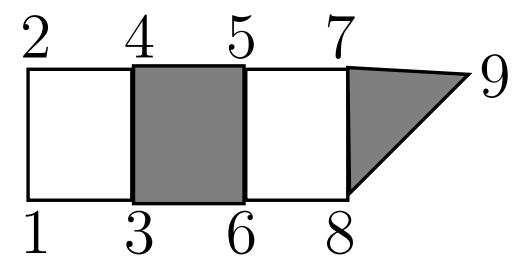

**Definition 2.1.4** Let $X_1, \ldots, X_d$ be a partition of $X$. The matroid whose collection of bases is

$$\mathcal{B} = \{\{y_1, \ldots, y_d\} | \, y_i \in X_i \text{ for } i = 1, \ldots, d\}.$$

is called the *transversal matroid* defined by $X_1, \ldots, X_d$ and is denoted by $\mathcal{M}$..

Recall that the set covering polyhedron of the clutter of bases of any transversal matroid is integral [42, p. 92]. The next result generalizes the fact that any bipartite graph is a subgraph of a complete bipartite graph.

**Corollary 2.1.5** *If $\mathcal{C}$ is a d-uniform clutter and $Q(A)$ is integral, then there is a partition $X_1,\ldots,X_d$ of $X$ such that $\mathcal{C}$ is a subclutter of the clutter of bases of the transversal matroid $\mathcal{M}$ defined by $X_1,\ldots,X_d$.*

The notion of a vertex critical clutter is the natural generalization of the corresponding notion for graph [50]. The family of vertex critical graphs has many nice properties, for instance if $G$ is a vertex critical graph with $n$ vertices, then $\alpha_0(G) \geq n/2$ (see [43, Theorem 3.11]). If a clutter $\mathcal{C}$ is vertex critical, then $\alpha_0(\mathcal{C}) = \alpha_0(\mathcal{C}\setminus\{x_i\})+1$ for all $i$.

**Proposition 2.1.6** *Let $\mathcal{C}$ be a d-uniform clutter with a perfect matching such that $Q(A)$ is integral. Then $\mathcal{C}$ is vertex critical.*

**Proof.** First, we claim that $n = gd$, where $g = \alpha_0(\mathcal{C}) = \operatorname{ht} I(\mathcal{C})$. By Proposition 2.1.2 there are $X_1,\ldots,X_d$ mutually disjoint minimal vertex covers of $\mathcal{C}$ such that $X = \cup_{i=1}^d X_i$. First we show that $n \geq gd$. Notice that

$$n = |X| = \sum_{i=1}^{d} |X_i| \geq \sum_{i=1}^{d} g = dg.$$

Thus $n \geq gd$. By hypothesis there are mutually disjoint edges $f_1,\ldots,f_r$ such that $X$ is equal to $f_1\cup\cdots\cup f_r$. Consequently $n = rd$. So $n = rd \geq gd$, i.e., $r \geq g$. On the other hand $g = \operatorname{ht} I(\mathcal{C}) \geq \beta_1(\mathcal{C}) \geq r$. Thus $r = g$ and $n = gd$ as claimed. In particular $\mathcal{C}$ has the König property. We now prove that $\mathcal{C}$ is vertex critical. Notice that

$$n = gd = |X_1| + \cdots + |X_d|.$$

As $|X_i| \geq g$ for all $i$, we get $|X_i| = g$ for all $i$. It follows rapidly that $\mathcal{C}$ is vertex critical. Indeed notice that each vertex $x_i$ belongs to a minimal vertex cover $C_i$ of $\mathcal{C}$ with $g$ vertices. The set $C_i \setminus \{x_i\}$ is a vertex cover of $\mathcal{C}\setminus\{x_i\}$ of size $g-1$. Hence $\alpha_0(\mathcal{C}\setminus\{x_i\}) < \alpha_0(\mathcal{C})$. □

**Remark 2.1.7** Consider the clutter $\mathcal{C}$ of Example 2.1.3. This clutter has a perfect matching and $Q(A)$ is integral, but it is not vertex critical because $\alpha_0(\mathcal{C}\setminus\{x_9\}) = \alpha_0(\mathcal{C}) = 4$. Thus the uniformity condition is essential in Proposition 2.1.6.

From the proof of Proposition 2.1.6 we get:

**Proposition 2.1.8** *Let $\mathcal{C}$ be a d-uniform clutter with a perfect matching $f_1,\ldots,f_r$. If $Q(A)$ is integral, then $r = \alpha_0(\mathcal{C})$ and there are $X_1,\ldots,X_d$ mutually disjoint minimal vertex covers of $\mathcal{C}$ of size $\alpha_0(\mathcal{C})$ such that $X = \cup_{i=1}^d X_i$.*

## 2.2 Algebras and TDI systems of clutters

As before let $R = K[x_1, \ldots, x_n]$ be a polynomial ring over a field $K$ and let $\mathcal{C}$ be a clutter with vertex set $X = \{x_1, \ldots, x_n\}$, edge set $E(\mathcal{C})$, edge ideal $I = I(\mathcal{C})$, and incidence matrix $A$. The column vectors of $A$ are denoted by $v_1, \ldots, v_q$. Thus the edge ideal of $\mathcal{C}$ is the ideal of $R$ generated by the set $F = \{x^{v_1}, \ldots, x^{v_q}\}$.

First we examine the interaction between combinatorial optimization properties of clutters and algebraic properties of monomial algebras. The *monomial algebras* considered here are: (a) the *Rees algebra*

$$R[It] := R \oplus It \oplus \cdots \oplus I^i t^i \oplus \cdots \subset R[t],$$

where $t$ is a new variable, (b) the *homogeneous monomial subring*

$$K[Ft] = K[x^{v_1}t, \ldots, x^{v_q}t] \subset R[t]$$

spanned by $Ft = \{x^{v_1}t, \ldots, x^{v_q}t\}$, (c) the *edge subring*

$$K[F] = K[x^{v_1}, \ldots, x^{v_q}] \subset R$$

spanned by $F$, and (d) the *Ehrhart ring*

$$A(P) = K[\{x^a t^i \,|\, a \in \mathbb{Z}^n \cap iP; i \in \mathbb{N}\}] \subset R[t]$$

of the lattice polytope $P = \text{conv}(v_1, \ldots, v_q)$.

The Rees algebra of the edge ideal $I$ can be written as

$$R[It] \;\; = \;\; K[\{x^a t^b \,|\, (a,b) \in \mathbb{N}\mathcal{A}'\}]$$

where $\mathcal{A}' = \{(v_1, 1), \ldots, (v_q, 1), e_1, \ldots, e_n\}$, $e_i$ is the $ith$ unit vector in $\mathbb{R}^{n+1}$, and $\mathbb{N}\mathcal{A}'$ is the subsemigroup of $\mathbb{N}^{n+1}$ spanned by $\mathcal{A}'$. According to [87, Theorem 7.2.28] the integral closure of $R[It]$ in its field of fractions can be expressed as

$$\overline{R[It]} \;\; = \;\; K[\{x^a t^b \,|\, (a,b) \in \mathbb{Z}\mathcal{A}' \cap \mathbb{R}_+\mathcal{A}'\}]$$

where $\mathbb{R}_+\mathcal{A}'$ is the cone spanned by $\mathcal{A}'$ and $\mathbb{Z}\mathcal{A}'$ is the subgroup spanned by $\mathcal{A}'$. The cone $\mathbb{R}_+\mathcal{A}'$ is called the *Rees cone* of $I$. The Rees algebra of $I$ is called *normal* if $R[It] = \overline{R[It]}$. Notice that $\mathbb{Z}\mathcal{A}' = \mathbb{Z}^{n+1}$. Hence we obtain the following well known fact:

**Lemma 2.2.1** *$R[It]$ is normal if and only if $\mathbb{N}\mathcal{A}' = \mathbb{Z}^{n+1} \cap \mathbb{R}_+\mathcal{A}'$.*

The Rees cone of $I$ has dimension $n+1$ because $\mathbb{Z}\mathcal{A}' = \mathbb{Z}^{n+1}$. According to [93, Theorem 4.1.1] there is a unique irreducible representation

$$\mathbb{R}_+\mathcal{A}' = H^+_{e_1} \cap H^+_{e_2} \cap \cdots \cap H^+_{e_{n+1}} \cap H^+_{\ell_1} \cap H^+_{\ell_2} \cap \cdots \cap H^+_{\ell_r}$$

such that each $\ell_k$ is in $\mathbb{Z}^{n+1}$, the non-zero entries of each $\ell_k$ are relatively prime, and none of the closed halfspaces $H^+_{e_1}, \ldots, H^+_{e_{n+1}}, H^+_{\ell_1}, \ldots, H^+_{\ell_r}$ can be omitted from the intersection. Here $H^+_a$ denotes the closed halfspace

$$H^+_a = \{x \in \mathbb{R}^{n+1} | \, \langle x, a\rangle \geq 0\},$$

$H_a$ stands for the hyperplane through the origin with normal vector $a$, and $\langle\ ,\ \rangle$ denotes the standard inner product. Irreducible representations of Rees cones were first introduced and studied in [29]. There are some interesting links between these representations, edge ideals [42], perfect graphs [91], and bases monomial ideals of matroids or polymatroids [92].

The Rees cone of $I$ and the set covering polyhedron of $\mathcal{C}$ are closely related:

**Theorem 2.2.2** [45, Corollary 3.13] *Let $C_1, \ldots, C_s$ be the minimal vertex covers of a clutter $\mathcal{C}$ and let $u_k = \sum_{x_i \in C_k} e_i$ for $1 \leq k \leq s$. Then $Q(A)$ is integral if and only if the irreducible representation of the Rees cone is*:

$$\mathbb{R}_+\mathcal{A}' = H^+_{e_1} \cap H^+_{e_2} \cap \cdots \cap H^+_{e_{n+1}} \cap H^+_{\ell_1} \cap H^+_{\ell_2} \cap \cdots \cap H^+_{\ell_s}, \tag{2.3}$$

*where $\ell_k = (u_k, -1)$ for $1 \leq k \leq s$.*

**Definition 2.2.3** The clutter $\mathcal{C}$ satisfies the *max-flow min-cut* (MFMC) property if both sides of the LP-duality equation

$$\min\{\langle \alpha, x\rangle | \, x \geq 0; xA \geq \mathbf{1}\} = \max\{\langle y, \mathbf{1}\rangle | \, y \geq 0; Ay \leq \alpha\} \tag{2.4}$$

have integral optimum solutions $x$ and $y$ for each non-negative integral vector $\alpha$. The system $xA \geq \mathbf{1}$; $x \geq 0$ is called *totally dual integral* (TDI) if the maximum has an integral optimum solution $y$ for each integral vector $\alpha$ with finite maximum.

Recall that the *ith symbolic power* of $I$ is given by

$$I^{(i)} = S^{-1}I^i \cap R \ \text{ for } i \geq 1,$$

where $S = R \setminus \cup_{k=1}^{s}\mathfrak{p}_i$ and $S^{-1}I^i$ is the localization of $I^i$ at $S$. In our situation the $ith$ symbolic power of $I$ has a simple expression: $I^{(i)} = \mathfrak{p}_1^i \cap \cdots \cap \mathfrak{p}_s^i$, see [87].

**Theorem 2.2.4** ([29, 42, 45, 57, 71]) *The following statements are equivalent*:

(i) *The associated graded ring* $\mathrm{gr}_I(R) = R[It]/IR[It]$ *is reduced.*

(ii) $R[It]$ *is normal and* $Q(A)$ *is an integral polyhedron.*

(iii) $I^i = I^{(i)}$ *for* $i \geq 1$, *where* $I^{(i)}$ *is the ith symbolic power of* $I$.

(iv) $\mathcal{C}$ *has the max-flow min-cut property.*

(v) $x \geq 0$; $xA \geq 1$ *is a* TDI system.

*Notation* For an integral matrix $B \neq (0)$, the greatest common divisor of all the nonzero $r \times r$ sub-determinants of $B$ will be denoted by $\Delta_r(B)$.

**Proposition 2.2.5** *Let* $\mathcal{C}$ *be a clutter and let* $B$ *be the matrix with column vectors* $(v_1, 1), \ldots, (v_q, 1)$. *The following statements hold*:

(i) [42, Proposition 4.4] *If* $\mathcal{C}$ *is uniform, then* $\mathcal{C}$ *has the max-flow min-cut property if and only if* $Q(A)$ *is integral and* $K[Ft] = A(P)$.

(ii) [28, Theorem 3.9] $\Delta_r(B) = 1$ *if and only if* $\overline{K[Ft]} = A(P)$, *where* $r$ *is the rank of* $B$.

We come to the main result of this section.

**Theorem 2.2.6** *If* $\mathcal{C}$ *is a uniform clutter with the max-flow min-cut property, then*

(a) $\Delta_r(A) = 1$, *where* $r = \mathrm{rank}(A)$.

(b) $\mathbb{N}\mathcal{A} = \mathbb{R}_+\mathcal{A} \cap \mathbb{Z}^n$, *where* $\mathcal{A} = \{v_1, \ldots, v_q\}$.

**Proof.** (a) Let $\widetilde{A}$ be the matrix with column vectors $(v_1,0),\ldots,(v_q,0)$. We need only show that $\Delta_r(\widetilde{A})=1$ because $A$ and $\widetilde{A}$ have the same rank and $\Delta_r(A)=\Delta_r(\widetilde{A})$. Let $B$ be the matrix with column vectors $(v_1,1),\ldots,(v_q,1)$. Since the clutter $\mathcal{C}$ is uniform, the last row vector of $B$, i.e., the vector $\mathbf{1}=(1,\ldots,1)$, is a $\mathbb{Q}$-linear combination of the first $n$ rows of $B$. Thus $\widetilde{A}$ and $B$ have the same rank. By Proposition 2.2.5(i) we obtain $K[Ft]=A(P)$. In particular, taking integral closures, one has $\overline{K[Ft]}=A(P)$ because $A(P)$ is always a normal domain. Hence by Proposition 2.2.5(ii) we have $\Delta_r(B)=1$. Recall that $\Delta_r(\widetilde{A})=1$ if and only if $\widetilde{A}$ is equivalent over $\mathbb{Z}$ to an identity matrix. In other words $\Delta_r(\widetilde{A})=1$ if and only if all the invariant factors of $\widetilde{A}$ are equal to 1. Thus it suffices to prove that $B$ is equivalent to $\widetilde{A}$ over $\mathbb{Z}$. Notice that in general $B$ and $\widetilde{A}$ are not equivalent over $\mathbb{Z}$ (for instance if $\mathcal{C}$ is a cycle of length three, then $\widetilde{A}$ and $B$ have rank 3, $\Delta_3(\widetilde{A})=2$ and $\Delta_3(B)=1$). By Proposition 2.1.2, there are $X_1,\ldots,X_d$ mutually disjoint minimal vertex covers of $\mathcal{C}$ such that $X=\cup_{i=1}^{d}X_i$ and

$$|\mathrm{supp}(x^{v_i})\cap X_k|=1 \quad \forall\ i,k. \tag{2.5}$$

By permuting the variables we may assume that $X_1$ is equal to $\{x_1,\ldots,x_r\}$. Hence the last row of $B$, which is the vector $\mathbf{1}$, is the sum of the first $|X_1|$ rows of $B$, i.e., the matrix $B$ is equivalent to $\widetilde{A}$ over $\mathbb{Z}$.

(b) It suffices to prove the inclusion $\mathbb{R}_+\mathcal{A}\cap\mathbb{Z}^n\subset\mathbb{N}\mathcal{A}$. Let $a$ be an integral vector in $\mathbb{R}_+\mathcal{A}$. Then $a=\lambda_1v_1+\cdots+\lambda_qv_q$, $\lambda_i\geq 0$ for all $i$. Set $b=\sum_i\lambda_i$ and denote the *ceiling* of $b$ by $\lceil b\rceil$. Recall that $\lceil b\rceil=b$ if $b\in\mathbb{N}$ and $\lceil b\rceil=\lfloor b\rfloor+1$ if $b\notin\mathbb{N}$, where $\lfloor b\rfloor$ is the integer part of $b$. Then $|a|=bd$. We claim that $(a,\lceil b\rceil)$ belongs to $\mathbb{R}_+\mathcal{A}'$, where $\mathcal{A}'$ is the set $\{e_1,\ldots,e_n,(v_1,1),\ldots,(v_q,1)\}$. Let $C_1,\ldots,C_s$ be the minimal vertex covers of $\mathcal{C}$ and let $u_i$ be the incidence vector of $C_i$ for $1\leq i\leq s$. Since $Q(A)$ is integral, by Theorem 2.2.2, we can write

$$\mathbb{R}_+\mathcal{A}'=H_{e_1}^+\cap H_{e_2}^+\cap\cdots\cap H_{e_{n+1}}^+\cap H_{\ell_1}^+\cap H_{\ell_2}^+\cap\cdots\cap H_{\ell_s}^+, \tag{2.6}$$

where $\ell_i=(u_i,-1)$ for $1\leq i\leq s$. Notice that $(a,b)\in\mathbb{R}_+\mathcal{A}'$, thus using Eq. (2.6) we get that $\langle a,u_i\rangle\geq b$ for all $i$. Hence $\langle a,u_i\rangle\geq\lceil b\rceil$ for all $i$ because $\langle a,u_i\rangle$ is an integer for all $i$. Using Eq. (2.6) again we get that $(a,\lceil b\rceil)\in\mathbb{R}_+\mathcal{A}'$, as claimed. By Theorem 2.2.4 the Rees ring $R[It]$ is normal. Consequently applying Lemma 2.2.1, we obtain that $(a,\lceil b\rceil)\in\mathbb{N}\mathcal{A}'$. There are non-negative integers $\eta_1,\ldots,\eta_q$ and $\rho_1,\ldots,\rho_n$ such that

$$(a,\lceil b\rceil)=\eta_1(v_1,1)+\cdots+\eta_q(v_q,1)+\rho_1e_1+\cdots+\rho_ne_n.$$

Hence it is seen that $|a| = \lceil b \rceil d + \sum_i \rho_i = bd$. Consequently $\rho_i = 0$ for all $i$ and $b = \lceil b \rceil$. It follows at once that $a \in \mathbb{N}\mathcal{A}$ as required. □

The next example shows that the uniformity hypothesis is essential in the two statements of Theorem 2.2.6.

**Example 2.2.7** Consider the clutter $\mathcal{C}$ whose incidence matrix is

$$\begin{bmatrix} 1 & 0 & 0 & 1 \\ 0 & 1 & 0 & 1 \\ 0 & 0 & 1 & 1 \\ 0 & 1 & 1 & 0 \\ 1 & 0 & 1 & 0 \end{bmatrix}.$$

Let $v_1, v_2, v_3, v_4$ be the columns of $A$. This clutter is not uniform, satisfies max-flow min-cut, $A$ is not equivalent over $\mathbb{Z}$ to an identity matrix, and $\{v_1, \ldots, v_4\}$ is not a Hilbert basis for the cone it generates. This clutter can be represented as:

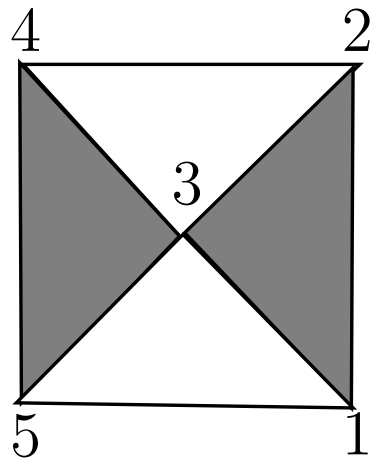

**Corollary 2.2.8** *If $\mathcal{C}$ is a uniform clutter with the max-flow min-cut property, then its incidence matrix diagonalizes over $\mathbb{Z}$ to an identity matrix.*

**Proof.** By Theorem 2.2.6 one has $\Delta_r(A) = 1$, where $r$ is the rank of $A$. Thus the invariant factors of $A$ are all equal to 1 (see [59, Theorem 3.9]), i.e., the Smith normal form of $A$ is an identity matrix. □

**Corollary 2.2.9** *Let $\mathcal{C}$ be a uniform clutter. Then the following are equivalent*:

(i) *$\mathcal{C}$ has the max-flow min-cut property.*

(ii) *$Q(A)$ is an integral polyhedron and $\mathbb{N}\mathcal{A} = \mathbb{R}_+\mathcal{A} \cap \mathbb{Z}^n$, where $\mathcal{A} = \{v_1, \ldots, v_q\}$.*

**Proof.** By Theorems 2.2.4 and 2.2.6 we obtain that (i) implies (ii). Next we prove that (ii) implies (i). By Proposition 2.2.5(i) it suffices to prove that $K[Ft] = A(P)$. Clearly $K[Ft] \subset A(P)$. To show the other inclusion take $x^a t^b \in A(P)$, i.e., $a \in bP \cap \mathbb{Z}^n$. Then from the equality $\mathbb{N}\mathcal{A} = \mathbb{R}_+\mathcal{A} \cap \mathbb{Z}^n$ it is seen that $a = \eta_1 v_1 + \cdots + \eta_q v_q$ for some $\eta_i$'s in $\mathbb{N}$ such that $\sum_i \eta_i = b$. Thus $x^a t^b \in K[Ft]$, as required. □

**Corollary 2.2.10** *Let $\mathcal{C}$ be a $d$-uniform clutter with $n$ vertices. If $\mathcal{C}$ has the max-flow min-cut property, then $\mathcal{C}$ has a perfect matching if and only if $n = d\alpha_0(\mathcal{C})$.*

**Proof.** $\Rightarrow$) As $Q(A)$ is integral and $\mathcal{C}$ has a perfect matching, from the proof of Proposition 2.1.6 we obtain the equality $n = d\alpha_0(\mathcal{C})$.

$\Leftarrow$) We set $g = \alpha_0(\mathcal{C})$. Let $u_1, \ldots, u_s$ be the characteristic vectors of the minimal vertex covers of $\mathcal{C}$ and let $B$ be the matrix with column vectors $u_1, \ldots, u_s$. Then $\langle \mathbf{1}, u_i \rangle \geq g$ for all $i$ because any minimal vertex cover of $\mathcal{C}$ has at least $g$ vertices. Thus the vector $\mathbf{1}/g$ belongs to the polyhedron $Q(B) = \{x | x \geq 0; xB \geq \mathbf{1}\}$. As $Q(A)$ is integral, by [16, Theorem 1.17] we get that $Q(B)$ is an integral polyhedron. Consequently

$$Q(B) = \mathbb{R}_+^n + \mathrm{conv}(v_1, \ldots, v_q),$$

where $\mathcal{A} = \{v_1, \ldots, v_q\}$ is the set of column vectors of the incidence matrix of $\mathcal{C}$. Therefore since the rational vector $\mathbf{1}/g$ is in $Q(B)$ we can write

$$\mathbf{1}/g = \delta + \mu_1 v_1 + \cdots + \mu_q v_q \quad (\delta \in \mathbb{R}_+^n;\ \mu_i \geq 0;\ \mu_1 + \cdots + \mu_q = 1).$$

Hence $n/g = |\delta| + (\sum_{i=1}^q \mu_i)d = |\delta| + d$. Since $n = dg$, we obtain that $\delta = 0$. Thus the vector $\mathbf{1}$ is in $\mathbb{R}_+\mathcal{A} \cap \mathbb{Z}^n$. By Theorem 2.2.6(b) this intersection is equal to $\mathbb{N}\mathcal{A}$. Then we can write $\mathbf{1} = \eta_1 v_1 + \cdots + \eta_q v_q$, for some $\eta_1, \ldots, \eta_q$ in $\mathbb{N}$. Hence it is readily seen that $\mathcal{C}$ has a perfect matching. □

Packing problems of clutters occur in many contexts of combinatorial optimization [16, 71], especially where the question of whether a linear program and its dual have integral optimum solutions is fundamental.

**Definition 2.2.11** A clutter $\mathcal{C}$ satisfies the *packing property* if all its minors satisfy the König property, i.e., $\alpha_0(\mathcal{C}') = \beta_1(\mathcal{C}')$ for any minor $\mathcal{C}'$ of $\mathcal{C}$.

To study linear algebra properties and ring theoretical properties of uniform clutters with the packing property we need the following result. This interesting result of Lehman is essential in the proofs of Theorem 2.3.1 and Corollary 2.3.2 because it allows the use of the structure theorems presented in Section 2.1, it also allows to state some conjectures about the packing property.

**Theorem 2.2.12** (A. Lehman [61]; see [16, Theorem 1.8]) *If a clutter $\mathcal{C}$ has the packing property, then $Q(A)$ is integral.*

**Proposition 2.2.13** ([16]) *If a clutter $\mathcal{C}$ has the max-flow min-cut property, then $\mathcal{C}$ has the packing property.*

**Proof.** It suffices to prove that $\mathcal{C}$ has the König property because the max-flow min-cut property is closed under taking minors. Making $\alpha = \mathbf{1}$ in Eq. (2.4), we get that the LP-duality equation:

$$\min\{\langle \mathbf{1}, x\rangle \,|\, x \geq 0; xA \geq \mathbf{1}\} = \max\{\langle y, \mathbf{1}\rangle \,|\, y \geq 0; Ay \leq \mathbf{1}\}$$

has optimum integral solutions $x$, $y$. To complete the proof notice that the left hand side of this equality is $\alpha_0(\mathcal{C})$ and the right hand side is $\beta_1(\mathcal{C})$. □

Conforti and Cornuéjols conjecture that the converse is also true:

**Conjecture 2.2.14** ([15], [16, Conjecture 1.6]) If a clutter $\mathcal{C}$ satisfies the packing property, then $\mathcal{C}$ has the max-flow min-cut property.

To the best of our knowledge this conjecture is open. For uniform clutters, using Proposition 2.2.5 and Theorem 2.2.12, we obtain the following algebraic version of this conjecture:

**Conjecture 2.2.15** If $\mathcal{C}$ is a uniform clutter with the packing property, then the homogeneous monomial ring $K[Ft]$ equals the Erhart ring $A(P)$.

Conjecture 2.2.14 together with Proposition 2.2.5 suggest the following:

**Conjecture 2.2.16** ([42]) If $\mathcal{C}$ is a uniform clutter with the packing property, then any of the following equivalent conditions hold

(a) $\mathbb{Z}^{n+1}/((v_1, 1), \ldots, (v_q, 1))$ is a free group.

(b) $\Delta_r(B) = 1$, where $B$ is the matrix with column vectors $(v_1, 1), \ldots, (v_q, 1)$ and $r$ is the rank of $B$.

(c) $B$ diagonalizes over $\mathbb{Z}$ to an identity matrix.

(d) $\overline{K[Ft]} = A(P)$, where $A(P)$ is the Ehrhart ring of $P = \text{conv}(v_1, \ldots, v_q)$.

This conjecture will be proved for $d$-uniform clutters with a perfect matching and vertex covering number equal to 2 (see Theorem 2.3.1).

## 2.3 Minors and the packing property

Let $\mathcal{C}$ be a $d$-uniform clutter with vertex set $X = \{x_1, \ldots, x_n\}$, let $x^{v_1}, \ldots, x^{v_q}$ be the minimal set of generators of the edge ideal $I = I(\mathcal{C})$, and let $A$ be the incidence matrix of $\mathcal{C}$ with column vectors $v_1, \ldots, v_q$. We denote the transpose of $A$ by $A^t$. All clutters considered in this section have vertex covering number equal to 2. We shall be interested in studying the relationships between the combinatorics of $\mathcal{C}$, the algebraic properties of the Ehrhart ring $A(P)$, the linear algebra of the incidence matrix $A$, and the algebra of the edge ideal $I(\mathcal{C})$.

We begin by showing that nice combinatorial properties of $\mathcal{C}$ are reflected in nice linear algebra properties of $A$. The following is one of the main results of this section. It gives some support to Conjecture 2.2.16.

**Theorem 2.3.1** *Let $\mathcal{C}$ be a $d$-uniform clutter with a perfect matching such that $\mathcal{C}$ has the packing property and $\alpha_0(\mathcal{C}) = 2$. If $A$ has rank $r$, then*

$$\Delta_r\begin{pmatrix} A \\ \mathbf{1} \end{pmatrix} = 1.$$

**Proof.** By the Lehman theorem the polyhedron $Q(A)$ is integral (see Theorem 2.2.12). Thus by Proposition 2.1.8 there is a perfect matching $f_1, f_2$ of $\mathcal{C}$ with $X = f_1 \cup f_2$ and there is a partition $X_1, \ldots, X_d$ of $X$ such that $X_i$ is a minimal vertex cover of $\mathcal{C}$ for all $i$, $|X_i| = 2$ for all $i$, and

$$|\text{supp}(x^{v_i}) \cap X_k| = 1 \quad \forall\ i, k. \tag{2.7}$$

Thus we may assume that $X_i = \{x_{2i-1}, x_{2i}\}$ for $i = 1, \ldots, d$. Notice that $n = 2d$ and $\text{rank}(A) \geq 2$ because $X = f_1 \cup f_2$.

We proceed by induction on $r = \mathrm{rank}(A)$. Since the sum of the first two rows of $A$ is equal to $\mathbf{1}$, it suffices to prove that 1 is the only invariant factor of $A^t$ or equivalently that the Smith normal form of $A^t$ is the identity.

First we show the case $r = 2$ which is the base case for the induction process. We may assume that $v_1$ and $v_2$ are the characteristic vectors of $f_1$ and $f_2$ respectively. The matrix $A_0$ with rows $v_1, v_2$ is equivalent over $\mathbb{Z}$ to

$$\begin{pmatrix} 1 & 0 & 0 & \cdots & 0 \\ 0 & 1 & 0 & \cdots & 0 \end{pmatrix}$$

and this is the Smith normal form of $A_0$. Thus the quotient group $M = \mathbb{Z}^n/\mathbb{Z}\{v_1, v_2\}$ is torsion-free, where $\mathbb{Z}\{v_1, v_2\}$ is the subgroup of $\mathbb{Z}^n$ generated by $\{v_1, v_2\}$. We claim that $\mathbb{Z}\{v_1, v_2\}$ is equal to $\mathbb{Z}\{v_1, v_2, \ldots, v_q\}$. Since the rank of $A$ is 2, for each $i \geq 3$ there is an integer $0 \neq \eta_i$ such that $\eta_i v_i \in \mathbb{Z}\{v_1, v_2\}$. Hence the image $\overline{v}_i$ of $v_i$ in $M$ is a torsion element, so $\overline{v}_i$ must be zero, that is, $v_i \in \mathbb{Z}\{v_1, v_2\}$. This completes the proof of the claim. Therefore

$$\Delta_2(A^t) = |T(\mathbb{Z}^n/\mathbb{Z}\{v_1, \ldots, v_q\})| = |T(\mathbb{Z}^n/\mathbb{Z}\{v_1, v_2\})| = 1,$$

where $T(\mathbb{Z}^n/\mathbb{Z}\{v_1, \ldots, v_q\})$ is the torsion subgroup of $\mathbb{Z}^n/\mathbb{Z}\{v_1, \ldots, v_q\}$, the first equality follows from [59, Theorem 3.9] and [59, pp. 187-188]. Therefore $\Delta_2(A^t) = 1$ and the Smith normal form of $A^t$ is the identity, as required.

We continue with the induction process by assuming $r \geq 3$. Let $w_1, \ldots, w_{2d}$ be the columns of $A^t$ and let $V_i$ be the linear space generated by $w_1, \ldots, w_{2i}$. Notice that, because of Eq. (2.7), for each odd integer $k$ the sum of rows $k$ and $k+1$ of the matrix $A$ is equal to $\mathbf{1} = (1, \ldots, 1)$, i.e., $w_k + w_{k+1} = \mathbf{1}$ for $k$ odd. Thus if $k$ is odd and we remove columns $w_k$ and $w_{k+1}$ from $A^t$ we obtain a submatrix whose rank is greater than or equal to $r - 1$. Thus after permuting columns we may assume

$$\dim(V_i) = \begin{cases} i+1 & \text{if } \ 1 \leq i \leq r-1, \\ r & \text{if } \ r \leq i. \end{cases} \tag{2.8}$$

Let $J$ be the square-free monomial ideal defined by the rows of the matrix $[w_1, \ldots, w_{2(r-1)}]$, where $w_1, \ldots, w_{2(r-1)}$ are column vectors, and let $\mathcal{D}$ be the clutter associated to the edge ideal $J$. If $v_i = (v_{i,1}, \ldots, v_{i,n})$ for $1 \leq i \leq q$, then $J$ is generated by the monomials:

$$x_1^{v_{1,1}} x_2^{v_{1,2}} \cdots x_{2(r-1)}^{v_{1,2(r-1)}}, \ldots, x_1^{v_{i,1}} x_2^{v_{i,2}} \cdots x_{2(r-1)}^{v_{i,2(r-1)}}, \ldots, x_1^{v_{q,1}} x_2^{v_{q,2}} \cdots x_{2(r-1)}^{v_{q,2(r-1)}}.$$

Thus $J$ is a minor of $I$ because $J$ is obtained from $I$ by making $x_i = 1$ for $i > 2(r-1)$. Then a minimal set of generators of $J$ consists of monomials of degree $r-1$, $\alpha_0(\mathcal{D}) = 2$, and $\mathcal{D}$ has a perfect matching. Furthermore $\mathcal{D}$ satisfies the packing property because $J$ is a minor of $I$. If $[w_1, \ldots, w_{2(r-1)}]$ diagonalizes (over the integers) to the identity matrix, so does $A^t$, this follows using arguments similar to those used for the case $\text{rank}(A) = 2$. Therefore we may harmlessly assume $d = r-1$, $I = J$, and $\mathcal{C} = \mathcal{D}$. This means that the matrix $A$ has rank $d+1$, the maximum possible.

Let $B$ be the matrix $[w_1, \ldots, w_{2(d-1)}]$ and let $I'$ be the monomial ideal defined by the rows of $B$, that is $I'$ is obtained from $I$ making $x_{2d-1} = x_{2d} = 1$. The matrix $B$ has rank $r-1$. Hence by induction hypothesis $B$ diagonalizes to a matrix $[I_{r-1}, \mathbf{0}]$, where $I_{r-1}$ is the identity matrix of size $r-1$. Recall that $f_1, f_2$ is a perfect matching of $\mathcal{C}$. Then by permuting rows and columns we may assume that the matrix $A^t$ is written as:

$$\underbrace{\begin{array}{ccccc}
10 & 10 & 10 & \cdots & 10 \\
01 & 01 & 01 & \cdots & 01 \\
\bigcirc & \bigcirc & \bigcirc & \cdots & \bigcirc \\
\vdots & \vdots & & & \vdots \\
\bigcirc & \bigcirc & \bigcirc & \cdots & \bigcirc \\
\bigcirc & \bigcirc & \bigcirc & \cdots & \bigcirc \\
\vdots & \vdots & & & \vdots \\
\bigcirc & \bigcirc & \bigcirc & \cdots & \bigcirc
\end{array}}_{B}
\begin{array}{cc}
10 & \leftarrow \\
01 & \\
10 & \leftarrow \\
\vdots & \vdots \\
10 & \leftarrow \\
01 & \\
\vdots & \\
01 &
\end{array}$$

where either a pair 1 0 or 0 1 must occur in the places marked with a circle and such that the number of $1's$ in the last column is greater than or equal to the number of $1's$ in any other column. Consider the matrix $C$ obtained from $A^t$ by removing the rows whose penultimate entry is equal to 1 (these are marked above with an arrow) and removing the last column. Let $K$ be the monomial ideal defined by the rows of $C$, that is $K$ is obtained from $I$ by making $x_{2d-1} = 0$ and $x_{2d} = 1$. By the choice of the last column and because of Eq. (2.8) it is seen that $K$ has height two. Since $K$ is a minor of $I$ it has the König property. Consequently the matrix $A^t$ has one of the following two

forms:

$$
\begin{array}{cccccc}
10 & 10 & 10 & \cdots & 10 & 10 \\
01 & 01 & 01 & \cdots & 01 & 01 \\
\bigcirc & \bigcirc & \bigcirc & \cdots & \bigcirc & 10 \\
\vdots & \vdots & & & \vdots & \vdots \\
\bigcirc & \bigcirc & \bigcirc & \cdots & \bigcirc & 10 \\
10 & 10 & 10 & \cdots & 10 & 01 \\
\bigcirc & \bigcirc & \bigcirc & \cdots & \bigcirc & 01 \\
\vdots & \vdots & & & \vdots & \vdots \\
\bigcirc & \bigcirc & \bigcirc & \cdots & \bigcirc & 01 \\
\multicolumn{5}{c}{\underbrace{\hspace{6em}}_{B}} &
\end{array}
\qquad\qquad
\begin{array}{ccccccl}
10 & 10 & 10 & \cdots & 10 & 10 & \text{row } v_1 \\
01 & 01 & 01 & \cdots & 01 & 01 & \text{row } v_2 \\
\bigcirc & \bigcirc & \bigcirc & \cdots & \bigcirc & 10 & \\
\vdots & \vdots & & & \vdots & \vdots & \\
\bigcirc & \bigcirc & \bigcirc & \cdots & \bigcirc & 10 & \\
\bigcirc & \bigcirc & \bigcirc & \cdots & \bigcirc & 01 & \text{row } v_j \\
\bigcirc & \bigcirc & \bigcirc & \cdots & \bigcirc & 01 & \text{row } v_{j+1} \\
\vdots & \vdots & & & \vdots & \vdots & \\
\bigcirc & \bigcirc & \bigcirc & \cdots & \bigcirc & 01 & \\
\multicolumn{5}{c}{\underbrace{\hspace{6em}}_{B}} & &
\end{array}
$$

where in the second case one has $v_j + v_{j+1} = (1, 1, \ldots, 1, 0, 2)$. In the second case, using row operations, we may replace $v_2$ by $v_j + v_{j+1} - v_2$. In the first case by permuting rows $v_2$ and $v_j$ we may replace $v_2$ by $(1, 0, \ldots, 1, 0, 0, 1)$. Thus in both cases, using row operations, we get that $A^t$ can be brought to the form:

$$
\begin{array}{cccccc}
10 & 10 & 10 & \cdots & 10 & 10 \\
10 & 10 & 10 & \cdots & 10 & 01 \\
\bigcirc & \bigcirc & \bigcirc & \cdots & \bigcirc & 10 \\
\vdots & \vdots & & & \vdots & \vdots \\
\bigcirc & \bigcirc & \bigcirc & \cdots & \bigcirc & 10 \\
\bigcirc & \bigcirc & \bigcirc & \cdots & \bigcirc & 01 \\
\vdots & \vdots & & & \vdots & \vdots \\
\bigcirc & \bigcirc & \bigcirc & \cdots & \bigcirc & 01 \\
\multicolumn{5}{c}{\underbrace{\hspace{6em}}_{B_1}} &
\end{array}
$$

where $B_1$ has rank $r-1$ and diagonalizes to an identity. Therefore it is readily seen that this matrix is equivalent to

$$
\left[\begin{array}{cccccccc}
I_{r-1} & 0 & 0 & \cdots & 0 & 0 & 0 \\
\mathbf{0} & 0 & 0 & \cdots & 0 & 1 & -1 \\
\mathbf{0} & 0 & 0 & \cdots & 0 & a_1 & b_1 \\
\vdots & \vdots & & & \vdots & \vdots & \vdots \\
\mathbf{0} & 0 & 0 & \cdots & 0 & a_s & b_s
\end{array}\right]
$$

for some integers $a_1, b_1, \ldots, a_s, b_s$. Next for $1 \leq i \leq s$ we multiply the second row by $-a_i$ and add it to row $i+2$. Then we add the last two columns. Using these row and column operations this matrix can be brought to the form:

$$\begin{bmatrix} I_{r-1} & 0 & 0 & \cdots & 0 & 0 & 0 \\ \mathbf{0} & 0 & 0 & \cdots & 0 & 1 & 0 \\ \mathbf{0} & 0 & 0 & \cdots & 0 & 0 & c_1 \\ \vdots & \vdots & & & \vdots & \vdots & \vdots \\ \mathbf{0} & 0 & 0 & \cdots & 0 & 0 & c_s \end{bmatrix}$$

for some integers $c_1, \ldots, c_s$. To finish the proof observe that $c_i = 0$ for all $i$ because this matrix has rank $r$. Hence this matrix reduces to $[I_r, \mathbf{0}]$, i.e., $A^t$ is equivalent to the identity matrix $[I_r, \mathbf{0}]$, as required. □

Next we show that clutters with nice algebraic and combinatorial properties have nice combinatorial optimization properties.

**Corollary 2.3.2** *Let $\mathcal{C}$ be a $d$-uniform clutter with a perfect matching such that $\mathcal{C}$ has the packing property and $\alpha_0(\mathcal{C}) = 2$. If $K[Ft]$ is normal, then $\mathcal{C}$ has the max-flow min-cut property.*

**Proof.** Let $B$ be the matrix with column vectors $(v_1, 1), \ldots, (v_q, 1)$. By Theorem 2.3.1 we have $\Delta_r(B) = 1$, where $r = \text{rank}(A)$. Notice that $r = \text{rank}(B)$ because $\mathcal{C}$ is uniform. According to Proposition 2.2.5(ii) the condition $\Delta_r(B) = 1$ is equivalent to the equality $\overline{K[Ft]} = A(P)$. Thus by the normality of $K[Ft]$ we get

$$K[Ft] = \overline{K[Ft]} = A(P).$$

By the Lehman theorem the polyhedron $Q(A)$ is integral (see Theorem 2.2.12). Hence by Proposition 2.2.5(i) we get that $\mathcal{C}$ satisfies the max-flow min-cut property. □

The next result gives some support to Conjecture 2.2.14.

**Corollary 2.3.3** *Let $\mathcal{C}$ be a $d$-uniform clutter with a perfect matching such that $\mathcal{C}$ has the packing property and $\alpha_0(\mathcal{C}) = 2$. If $v_1, \ldots, v_q$ are linearly independent, then $\mathcal{C}$ has the max-flow min-cut property.*

**Proof.** It follows at once from Corollary 2.3.2 because $K[Ft]$ is normal. □

A $d$-uniform clutter $\mathcal{C}$ is called 2-*partitionable* if its vertex set $X$ has a partition $X_1, \ldots, X_d$ such that $X_i$ is a minimal vertex cover of $\mathcal{C}$ and $|X_i| = 2$ for $i = 1, \ldots, d$. For instance the clutters of Theorem 2.3.1 are 2-partitionable.

Another instance is the clutter of minimal vertex covers of an unmixed bipartite graph (see Corollary 2.3.6). A clutter is called *unmixed* if all its minimal vertex covers have the same number of elements. A specific illustration of a famous 2-partitionable clutter is given in Example 2.3.5. Our next result is about this family of clutters.

Recall that an edge ideal $I \subset R$ is called *normal* if $\overline{I^i} = I^i$ for $i \geq 1$, where

$$\overline{I^i} = (\{x^a \in R | \, \exists \, p \geq 1; (x^a)^p \in I^{pi}\}) \subset R$$

is the integral closure of $I^i$. Also recall that $I$ is called *minimally non-normal* if $I$ is not normal and all its proper minors are normal. An edge ideal $I$ is normal if and only if its Rees algebra $R[It]$ is normal (see [87, Theorem 3.3.18]).

The other main result of this section is:

**Theorem 2.3.4** *Let $\mathcal{C}$ be a $d$-uniform clutter with a partition $X_1, \ldots, X_d$ of $X$ such that $X_i = \{x_{2i-1}, x_{2i}\}$ is a minimal vertex cover of $\mathcal{C}$ for all $i$. Then*

(a) $\mathrm{rank}(A) \leq d + 1$.

(b) *If $C$ is a minimal vertex cover of $\mathcal{C}$, then* $2 \leq |C| \leq d$.

(c) *If $\mathcal{C}$ satisfies the König property and there is a minimal vertex cover $C$ of $\mathcal{C}$ with $|C| = d \geq 3$, then* $\mathrm{rank}(A) = d + 1$.

(d) *If $I = I(\mathcal{C})$ is minimally non-normal and $\mathcal{C}$ satisfies the packing property, then* $\mathrm{rank}(A) = d + 1$.

**Proof.** (a) For each odd integer $k$ the sum of rows $k$ and $k + 1$ of the matrix $A$ is equal to $\mathbf{1} = (1, \ldots, 1)$. Thus the rank of $A$ is bounded by $d + 1$.

(b) By the pigeon hole principle, any minimal vertex cover $C$ of the clutter $\mathcal{C}$ satisfies $2 \leq |C| \leq d$.

(c) First notice that $C$ contains exactly one element of each $X_j$ because $X_j \not\subset C$ and $|C| = d$. Thus we may assume

$$C = \{x_1, x_3, \ldots, x_{2d-1}\}.$$

Consider the monomial $x^\alpha = x_2 x_4 \cdots x_{2d}$ and notice that $x_k x^\alpha \in I$ for each $x_k \in C$ because the monomial $x_k x^\alpha$ is clearly in every minimal prime of $I$. Writing $x_k = x_{2i-1}$ with $1 \leq i \leq d$ and

$$x^{\alpha_i} = \frac{x_{2i-1} x^\alpha}{x_{2i}} = \frac{x_{2i-1}(x_2 x_4 \cdots x_{2d})}{x_{2i}}. \qquad (2.9)$$

We claim that $x^{\alpha_i} \in I$, and so we conclude that $x^{\alpha_i}$ is a minimal generator of $I$. Assume that $x^{\alpha_i}$ is not in $I$. Hence

$$C_i = \{x_{2(1)-1}, ..., x_{2(i-1)-1}, x_{2i}, x_{2(i+1)-1}, ..., x_{2(d)-1}\}$$

is a minimal vertex cover. Thus the set

$$f = \{x_{2(1)}, ..., x_{2(i-1)}, x_{2i-1}, x_{2(i+1)}, ..., x_{2(d)}\}$$

is not an edge because $f \cap C_i = \emptyset$. Since $C_i$ is minimal, there is an edge $e$ disjoint from $C_i - \{x_{2i}\}$, thus either $e = f$ or $e = \{x_2, ..., x_{2d}\}$ a contradiction because $C$ is a minimal vertex cover. Hence $x^{\alpha_i}$ is a minimal generator of $I$. Thus we may assume $x^{\alpha_i} = x^{v_i}$ for $i = 1, \ldots, d$. The vector $\mathbf{1}$ belongs to the linear space generated by $v_1, \ldots, v_q$ because $\mathcal{C}$ has the König property, where $q$ is the number of edges of $\mathcal{C}$. By part (a) the rank of $A$ is at most $d+1$. Thus it suffices to prove that $v_1, \ldots, v_d, \mathbf{1}$ are linearly independent. Assume that

$$\lambda_1 v_1 + \cdots + \lambda_d v_d + \lambda_{d+1}\mathbf{1} = 0 \tag{2.10}$$

for some scalars $\lambda_1, \ldots, \lambda_{d+1}$. From Eq. (2.9) we have

$$v_i = e_{2i-1} - e_{2i} + \sum_{j=1}^{d} e_{2j} \tag{2.11}$$

for $i = 1, \ldots, d$. Hence using Eq. (2.10) we conclude that

$$\sum_{i=1}^{d} \lambda_i e_{2i-1} + \lambda_{d+1} \sum_{i=1}^{d} e_{2i-1} = \sum_{i=1}^{d} (\lambda_i + \lambda_{d+1}) e_{2i-1} = 0.$$

Hence $\lambda_{d+1} = -\lambda_i$ for $i = 1, \ldots, d$. Using Eq. (2.10) once more we get

$$\lambda_{d+1}(v_1 + \cdots + v_d - \mathbf{1}) = 0.$$

As $v_1 + \cdots + v_d - \mathbf{1}$ cannot be zero by Eq. (2.11), we obtain that $\lambda_{d+1} = 0$. Consequently $\lambda_i = 0$ for all $i$, as required.

(d) Let $x^\alpha t^b$ be a minimal generator of $\overline{R[It]}$ not in $R[It]$ and let $m = x^\alpha$. Then $m \in \overline{I^b} \setminus I^b$. Using [42, Proposition 4.3], one has $\deg(m) = bd$. Consider the $R$-module $N = \overline{I^b}/I^b$. Notice that the image of $m$ in $N$ is nonzero. The set of associated primes of $N$, denoted by $\mathrm{Ass}(N)$, is the set of prime ideals $\mathfrak{p}$ of $R$ such that $N$ contains a submodule isomorphic to $R/\mathfrak{p}$. By the hypothesis that $I$ is minimally non-normal, we have that $\mathfrak{m} = (x_1, \ldots, x_n)$ is the only associated prime of $N$. Hence, using [64, Theorem 9, p. 50], we have the following expression for the radical of the annihilator of $N$

$$\mathrm{rad}(\mathrm{ann}(N)) = \bigcap_{\mathfrak{p} \in \mathrm{Ass}(N)} \mathfrak{p} = \mathfrak{m} = (x_1, \ldots, x_n).$$

Consequently $\mathfrak{m}^r \subset \mathrm{ann}(N) = \{x \in R | xN = 0\}$ for some $r > 0$. Thus for $i$ odd we can write

$$x_i^r x^\alpha = (x^{v_1})^{a_1} \cdots (x^{v_q})^{a_q} x^\delta,$$

where $a_1 + \cdots + a_q = b$ and $\deg(x^\delta) = r$. If we write $x^\delta = x_i^{s_1} x_{i+1}^{s_2} x^\gamma$ with $x_i, x_{i+1}$ not in the support of $x^\gamma$, making $x_j = 1$ for $j \notin \{i, i+1\}$, it is not hard to see that $r = s_1 + s_2$ and $\gamma = 0$. Thus we get an equation:

$$x_i^{s_2} x^\alpha = (x^{v_1})^{a_1} \cdots (x^{v_q})^{a_q} x_{i+1}^{s_2}$$

with $s_2 > 0$. Using a similar argument we obtain an equation:

$$x_{i+1}^{w_1} x^\alpha = (x^{v_1})^{b_1} \cdots (x^{v_q})^{b_q} x_i^{w_1}$$

with $w_1 > 0$. Therefore

$$x_{i+1}^{s_2+w_1} (x^{v_1})^{a_1} \cdots (x^{v_q})^{a_q} = x_i^{s_2+w_1} (x^{v_1})^{b_1} \cdots (x^{v_q})^{b_q}.$$

Consider the group $\mathbb{Z}^n/\mathbb{Z}\mathcal{A}$, where $\mathcal{A} = \{v_1, \ldots, v_q\}$. Since this group is torsion-free, we get $e_i - e_{i+1} \in \mathbb{Z}\mathcal{A}$ for $i$ odd. Finally to conclude that $\mathrm{rank}(\mathrm{A}) = d + 1$ notice that $\mathbf{1} \in \mathbb{Z}\mathcal{A}$. □

**Example 2.3.5** [16, p. 12] Let $\mathcal{C}$ be the uniform clutter whose edges are

$$\{x_1, x_4, x_5\}, \{x_1, x_3, x_6\}, \{x_2, x_4, x_6\}, \{x_2, x_3, x_5\}$$

and let $A$ be its incidence matrix. Using Normaliz [11], it is not hard to see that the minimal vertex covers of $\mathcal{C}$ are

$$\begin{array}{llll} X_1 = \{x_1, x_2\}, & X_2 = \{x_3, x_4\}, & X_3 = \{x_5, x_6\}, & \\ C_4 = \{x_1, x_4, x_5\}, & C_5 = \{x_1, x_3, x_6\}, & C_6 = \{x_2, x_4, x_6\}, & C_7 = \{x_2, x_3, x_5\}. \end{array}$$

This clutter satisfies the hypotheses of Theorem 2.3.4 with $d = 3$ and has minimal vertex covers of sizes 2 and 3. Moreover the rank of $A$ is 4, $\mathcal{C}$ does not satisfy the König property and $Q(A)$ is integral.

For use below recall that a graph $G$ is *Cohen-Macaulay* if $R/I(G)$ is a Cohen-Macaulay ring. Any Cohen-Macaulay graph is unmixed [87, p. 169].

**Corollary 2.3.6** *Let $A$ be the incidence matrix of the clutter $\mathcal{C}$ of minimal vertex covers of an unmixed bipartite graph $G$. Then* (i) *$\mathcal{C}$ is 2-partitionable.* (ii) $\mathrm{rank}(A) \leq \alpha_0(G) + 1$*.* (iii) $\mathrm{rank}(A) = \alpha_0(G) + 1$ *if $G$ is Cohen-Macaulay.*

**Proof.** We set $d = \alpha_0(G)$. Clearly $\mathcal{C}$ is a $d$-uniform clutter because all minimal vertex covers of $G$ have size $d$. By [71, Theorem 78.1] the polyhedron $Q(A)$ is integral because $G$ is a bipartite graph and bipartite graphs have integral set covering polyhedron [42, Proposition 4.27]. Any minimal vertex cover of $\mathcal{C}$ is of the form $\{x_i, x_j\}$ for some edge $\{x_i, x_j\}$ of $G$. Therefore by Proposition 2.1.2 the clutter $\mathcal{C}$ is 2-partitionable. Thus we have shown (i). Then using part (a) of Theorem 2.3.4 we obtain (ii). To prove (iii) it suffices to prove that there are minimal vertex covers $C_1, \ldots, C_{d+1}$ of $G$ whose characteristic vectors are linearly independent. This follows by induction on the number of vertices and using the fact that any Cohen-Macaulay bipartite graph has at least one vertex of degree 1 [87, Theorem 6.4.4]. $\Box$

## 2.4 Triangulations and max-flow min-cut

Let $R = K[x_1, \ldots, x_n]$ be a polynomial ring over a field $K$ and let $\mathcal{C}$ be a $d$-uniform clutter with vertex set $X = \{x_1, \ldots, x_n\}$ and edge ideal $I = I(\mathcal{C})$. In what follows $F = \{x^{v_1}, \ldots, x^{v_q}\}$ will denote the minimal set of generators of $I$ and $\mathcal{A}$ will denote the set $\{v_1, \ldots, v_q\}$. The incidence matrix of $\mathcal{C}$, i.e., the $n \times q$ matrix with column vectors $v_1, \ldots, v_q$ will be denoted by $A$. One can think of the columns of $A$ as points in $\mathbb{R}^n$. The set $\mathcal{A}$ is called a *point configuration*.

Consider a homomorphism of $K$-algebras:

$$S = K[t_1, \ldots, t_q] \stackrel{\varphi}{\longrightarrow} K[F] \quad (t_i \stackrel{\varphi}{\longrightarrow} x^{v_i}),$$

where $S$ is a polynomial ring. The kernel of $\varphi$, denoted by $P$, is called the *toric ideal* of $K[F]$. For the rest of this section we assume that $\prec$ is a fixed term order for the set of monomials of $S$. We denote the initial ideal of $P$ by $\mathrm{in}_\prec(P)$. If the choice of a term order is clear, we may write $\mathrm{in}(P)$ instead of $\mathrm{in}_\prec(P)$. There is a one to one correspondence between simplicial complexes with vertex set $X$ and squarefree monomial ideals of $R$ via the Stanley-Reisner theory [75]:

$$\Delta \longmapsto I_\Delta = (x^a |\, x^a \text{ is squarefree and } \operatorname{supp}(x^a) \notin \Delta).$$

Notice that $\operatorname{rad}(\operatorname{in}(P))$, the radical of $\operatorname{in}(P)$, is a squarefree monomial ideal. Let $\Delta$ be the simplicial complex whose Stanley-Reisner ideal is $\operatorname{rad}(\operatorname{in}(P))$.

Let $\omega = (\omega_i) \in \mathbb{N}^q$ be an integral weight vector. If

$$f = f(t_1, \dots, t_q) = \lambda_1 t^{a_1} + \cdots + \lambda_s t^{a_s}$$

is a polynomial with $\lambda_1, \dots, \lambda_s$ in $K$, we define $\operatorname{in}_\omega(f)$, the *initial form* of $f$ relative to $\omega$, as the sum of all terms $\lambda_i t^{a_i}$ such that $\langle \omega, a_i \rangle$ is maximal. The ideal generated by all initial forms is denoted by $\operatorname{in}_\omega(P)$.

**Proposition 2.4.1** [77, Proposition 1.11] *For every term order* $\prec$*,* $\operatorname{in}_\prec(P) = \operatorname{in}_\omega(P)$ *for some non-negative integer weight vector* $\omega \in \mathbb{N}^q$*.*

**Theorem 2.4.2** [77, Theorem 8.3] *If* $\operatorname{in}(P) = \operatorname{in}_\omega(P)$*, then*

$$\Delta = \{\sigma |\, \exists\, c \in \mathbb{R}^n \text{ such that } \langle v_i, c \rangle = \omega_i \text{ if } t_i \in \sigma \text{ and } \langle v_i, c \rangle < \omega_i \text{ if } t_i \notin \sigma\}.$$

Let $\omega = (\omega_i) \in \mathbb{N}^q$ be a vector that represents the initial ideal $\operatorname{in}(P) = \operatorname{in}_\prec(P)$ with respect to a term order $\prec$, that is, $\operatorname{in}(P) = \operatorname{in}_\omega(P)$. Consider the primary decomposition of $\operatorname{rad}(\operatorname{in}(P))$ as an intersection of face ideals:

$$\operatorname{rad}(\operatorname{in}(P)) = \mathfrak{p}_1 \cap \mathfrak{p}_2 \cap \cdots \cap \mathfrak{p}_r,$$

where $\mathfrak{p}_1, \dots, \mathfrak{p}_r$ are ideals of $K[t_1, \dots, t_q]$ generated by subsets of $\{t_1, \dots, t_q\}$. Recall that the facets of $\Delta$ are given by

$$F_i = \{t_j |\, t_j \notin \mathfrak{p}_i\},$$

or equivalently by $\mathcal{A}_i = \{v_j |\, t_j \notin \mathfrak{p}_i\}$ if one identifies $t_i$ with $v_i$. According to Theorem 2.4.2, the family of facets (resp. cones or polytopes) of $\Delta$:

$$\{\mathcal{A}_1, \dots, \mathcal{A}_r\} \quad (\text{resp. } \{\mathbb{R}_+\mathcal{A}_1, \dots, \mathbb{R}_+\mathcal{A}_r\} \text{ or } \{\operatorname{conv}(\mathcal{A}_1), \dots, \operatorname{conv}(\mathcal{A}_r)\}),$$

form a regular triangulation of $\mathcal{A}$ (resp. $\mathbb{R}_+\mathcal{A}$ or $\operatorname{conv}(\mathcal{A})$). This means that

$$\operatorname{conv}(\mathcal{A}_1), \dots, \operatorname{conv}(\mathcal{A}_r)$$

are obtained by projection onto the first $n$ coordinates of the lower facets of

$$Q' = \mathrm{conv}((v_1, \omega_1), \ldots, (v_q, \omega_q)).$$

The regular triangulation $\{\mathcal{A}_1, \ldots, \mathcal{A}_r\}$ is called *unimodular* if $\mathbb{Z}\mathcal{A}_i = \mathbb{Z}\mathcal{A}$ for all $i$, where $\mathbb{Z}\mathcal{A}$ denotes the subgroup of $\mathbb{Z}^n$ spanned by $\mathcal{A}$. A major result of Sturmfels [77, Corollary 8.9] shows that this triangulation is unimodular if and only if $\mathrm{in}(P)$ is square-free.

We are interested in the following:

**Conjecture 2.4.3** If $\mathcal{C}$ is a uniform clutter that satisfies the max-flow min-cut property, then the rational polyhedral cone $\mathbb{R}_+\{v_1, \ldots, v_q\}$ has a unimodular regular triangulation.

**Example 2.4.4** Let $u_1, \ldots, u_r$ be the characteristic vectors of the collection of bases $\mathcal{B}$ of a transversal matroid $\mathcal{M}$. By [5, Proposition 2.1 and Theorem 4.2], the toric ideal of the subring $K[x^{u_1}, \ldots, x^{u_r}]$ has a square-free quadratic Gröbner basis. Therefore the cone $\mathbb{R}_+\{u_1, \ldots, u_r\}$ or the polytope $\mathrm{conv}(u_1, \ldots, u_r)$ has a unimodular regular triangulation. This gives support to Conjecture 2.4.3 because the clutter $\mathcal{B}$ has the max-flow min-cut property [42].

**Definition 2.4.5** A matrix $A$ with entries in $\{0, 1\}$ is called *balanced* if $A$ has no square submatrix of odd size with exactly two 1's in each row and column.

Recall that an integral matrix $A$ is $t$-*unimodular* if all the nonzero $r \times r$ sub-determinants of $A$ have absolute value equal to $t$, where $r$ is the rank of $A$. If $t = 1$ the matrix is called *unimodular*. If $A$ is $t$-unimodular, then any regular triangulation of $\mathbb{R}_+\{v_1, \ldots, v_q\}$ is unimodular [77], see [88, Proposition 5.20] for a very short proof of this fact.

The next result gives an interesting class of uniform clutters, coming from combinatorial optimization, that satisfy Conjecture 2.4.3. This result is surprising because not all balanced matrices are $t$-unimodular.

**Theorem 2.4.6** *Let $A$ be a balanced matrix with distinct column vectors $v_1, \ldots, v_q$. If $|v_i| = d$ for all $i$, then any regular triangulation of the cone $\mathbb{R}_+\{v_1, \ldots, v_q\}$ is unimodular.*

**Proof.** Let $\mathcal{A} = \{v_1, \ldots, v_q\}$ and let $\mathcal{A}_1, \ldots, \mathcal{A}_m$ be the elements of a regular triangulation of $\mathbb{R}_+\mathcal{A}$. Then $\dim \mathbb{R}_+\mathcal{A}_i = \dim \mathbb{R}_+\mathcal{A}$ and $\mathcal{A}_i$ is linearly independent for all $i$. Consider the clutter $\mathcal{C}$ whose edge ideal is $I = (x^{v_1}, \ldots, x^{v_q})$ and the subclutter $\mathcal{C}_i$ of $\mathcal{C}$ whose edges correspond to the vectors in $\mathcal{A}_i$. Let $A_i$ be the incidence matrix of $\mathcal{C}_i$. Since $A_i$ is a balanced matrix, using [71, Corollary 83.1a(iv), p. 1441], we get that the subclutter $\mathcal{C}_i$ has the max-flow min-cut property. Hence by Theorem 2.2.6 one has $\Delta_r(A_i) = 1$, where $r$ is the rank of $A_i$. Thus the invariant factors of $A_i$ are all equal to 1 (see [59, Theorem 3.9]). Therefore by the fundamental structure theorem of finitely generated abelian groups (see [59, p. 187]) the group $\mathbb{Z}^n/\mathbb{Z}\mathcal{A}_i$ is torsion-free for all $i$. Notice that $\dim \mathbb{R}_+\mathcal{A} = \operatorname{rank} \mathbb{Z}\mathcal{A}$ and $\dim \mathbb{R}_+\mathcal{A}_i = \operatorname{rank} \mathbb{Z}\mathcal{A}_i$ for all $i$. Since $r$ is equal to $\dim \mathbb{R}_+\mathcal{A}$, it follows rapidly that the quotient group $\mathbb{Z}\mathcal{A}/\mathbb{Z}\mathcal{A}_i$ is torsion-free and has rank 0 for all $i$. Consequently $\mathbb{Z}\mathcal{A} = \mathbb{Z}\mathcal{A}_i$ for all $i$, i.e., the triangulation is unimodular. $\Box$

If we do not require that $|v_i| = d$ for all $i$, this result is false even if $K[F]$ is homogeneous, i.e., even if there is $x_0 \in \mathbb{R}^n$ such that $\langle v_i, x_0 \rangle = 1$ for all $i$:

**Example 2.4.7** Consider the following matrix

$$A = \begin{pmatrix}
1 & 0 & 0 & 0 & 0 & 0 & 0 & 0 & 0 & 0 & 0 & 1 & 0 \\
0 & 1 & 0 & 0 & 0 & 0 & 0 & 0 & 0 & 0 & 0 & 1 & 0 \\
0 & 0 & 1 & 0 & 0 & 0 & 0 & 0 & 0 & 0 & 1 & 0 & 0 \\
0 & 0 & 0 & 1 & 0 & 0 & 0 & 0 & 0 & 0 & 1 & 0 & 0 \\
0 & 0 & 0 & 0 & 1 & 0 & 0 & 0 & 0 & 1 & 0 & 0 & 0 \\
0 & 0 & 0 & 0 & 0 & 1 & 0 & 0 & 0 & 1 & 0 & 0 & 0 \\
0 & 0 & 0 & 0 & 0 & 0 & 1 & 0 & 0 & 1 & 0 & 0 & 1 \\
0 & 0 & 0 & 0 & 0 & 0 & 0 & 1 & 0 & 0 & 1 & 0 & 1 \\
0 & 0 & 0 & 0 & 0 & 0 & 0 & 0 & 1 & 0 & 0 & 1 & 1 \\
0 & 0 & 0 & 0 & 0 & 0 & 0 & 0 & 0 & 1 & 1 & 1 & 1
\end{pmatrix}$$

Let $v_1, \ldots, v_{13}$ be the columns of $A$. It is not hard to see that the matrix $A$ is balanced. Using *Macaulay*2 [46] it is seen that the regular triangulation $\Delta$ of $\mathbb{R}_+\{v_1, \ldots, v_{13}\}$ determined by using the GRevLex order, on the polynomial ring $K[t_1, \ldots, t_{13}]$, has a simplex, namely $\{v_1, \ldots, v_6, v_{10}, \ldots, v_{13}\}$, which is not unimodular.

# Chapter 3

# Cohen-Macaulay Clutters with Combinatorial Optimization Properties and Parallelizations of Normal Edge Ideals

Let $\mathcal{C}$ be a uniform clutter and let $I = I(\mathcal{C})$ be its edge ideal. We prove that if $\mathcal{C}$ satisfies the packing property (resp. max-flow min-cut property), then there is a uniform Cohen-Macaulay clutter $\mathcal{C}_1$ satisfying the packing property (resp. max-flow min-cut property) such that $\mathcal{C}$ is a minor of $\mathcal{C}_1$. For arbitrary edge ideals of clutters we prove that the normality property is closed under parallelizations. Then we show some applications to edge ideals and clutters which are related to a conjecture of Conforti and Cornuéjols and to max-flow min-cut problems.

Let $R = K[x_1, \ldots, x_n]$ be a polynomial ring over a field $K$ and let $I$ be an ideal of $R$ minimally generated by a finite set $F = \{x^{v_1}, \ldots, x^{v_q}\}$ of square-free monomials. As usual we use the notation $x^a := x_1^{a_1} \cdots x_n^{a_n}$, where $a = (a_1, \ldots, a_n)$ is in $\mathbb{N}^n$. The *support* of a monomial $x^a$ is given by $\mathrm{supp}(x^a) = \{x_i \,|\, a_i > 0\}$. For technical reasons we shall assume that each variable $x_i$ occurs in at least one monomial of $F$.

A *clutter* with finite vertex set $X$ is a family of subsets of $X$, called edges, none of which is included in another. The set of vertices of a clutter $\mathcal{C}$ is denoted by $V(\mathcal{C})$ and the set of edges of $\mathcal{C}$ is denoted by $E(\mathcal{C})$. A clutter is called *d-uniform* if all its edges have exactly $d$ vertices. We associate to the ideal $I$ a *clutter* $\mathcal{C}$ by taking the set of indeterminates $X = \{x_1, \ldots, x_n\}$ as vertex set and $E = \{S_1, \ldots, S_q\}$ as edge set, where $S_k$ is the support of $x^{v_k}$. The vector $v_k$ is called the *characteristic vector* of $S_k$. The assignment $I \mapsto \mathcal{C}$ gives a natural one to one correspondence between the family of square-free monomial ideals and the family of clutters. The ideal $I$ is called the *edge ideal* of $\mathcal{C}$. To stress the relationship between $I$ and $\mathcal{C}$ we will use the notation $I = I(\mathcal{C})$. Edge ideals of graphs were introduced and studied in [73, 85]. Edge ideals of clutters also correspond to simplicial complexes via the Stanley-Reisner correspondence [75] and to facet ideals [31, 94]. The Cohen-Macaulay property of edge ideals has been recently studied in [12, 32, 49, 66, 81] using a combinatorial approach based on the notions of shellability, linear quotients, unmixedness, acyclicity and transitivity of digraphs, and the König property.

The aim of this chapter is to study the behavior, under certain operations, of various algebraic and combinatorial optimization properties of edge ideals and clutters such as the Cohen-Macaulay property, the normality, the torsion freeness, the packing and the max-flow min-cut properties. The study of edge ideals from the combinatorial optimization point of view was initiated in [7, 79] and continued in [19, 29, 40, 42, 45, 91], see also [53] . The Cohen-Macaulay and normality properties are two of the most interesting properties an edge ideal can have, see [10, 32, 75, 87] and [58, 84] respectively.

Recall that the *integral closure* of $I^i$, denoted by $\overline{I^i}$, is the ideal of $R$ given by

$$\overline{I^i} = (\{x^a \in R \mid \exists\, p \geq 1; (x^a)^p \in I^{pi}\}).$$

An ideal $I$ is called *normal* if $I^i = \overline{I^i}$ for all $i$. A clutter obtained from $\mathcal{C}$ by a sequence of deletions and duplications of vertices is called a *parallelization* of $\mathcal{C}$ and a clutter obtained from $\mathcal{C}$ by a sequence of deletions and contractions of vertices is called a *minor* of $\mathcal{C}$, see Section 3.1. It is known that the normality of $I(\mathcal{C})$ is closed under minors [29]. One of our main results shows that the normality of $I(\mathcal{C})$ is closed under parallelizations:

**Theorem 3.1.3** *Let $\mathcal{C}$ be a clutter and let $\mathcal{C}'$ be a parallelization of $\mathcal{C}$. If $I(\mathcal{C})$ is normal, then $I(\mathcal{C}')$ is normal.*

The ideal $I = I(\mathcal{C})$ is called *normally torsion free* if $I^i = I^{(i)}$ for all $i$, where $I^{(i)}$ is the *ith* symbolic power of $I$. As an application we prove that if $I(\mathcal{C})$ is normally torsion free and $\mathcal{C}'$ is a parallelization of $\mathcal{C}$, then $I(\mathcal{C}')$ is normally torsion free (Corollary 3.1.12). Let $A$ be the incidence matrix of $\mathcal{C}$, i.e., $A$ is the matrix with column vectors $v_1, \ldots, v_q$. A clutter $\mathcal{C}$ satisfies the *max-flow* *min-cut* (MFMC) property if both sides of the LP-duality equation

$$\min\{\langle w,x\rangle|\, x\geq 0; xA\geq \mathbf{1}\} = \max\{\langle y,\mathbf{1}\rangle|\, y\geq 0; Ay\leq w\}$$

have integral optimum solutions $x$ and $y$ for each non-negative integral vector $w$. A remarkable result of [45] (cf. [42, Theorem 4.6]) shows that $I(\mathcal{C})$ is normally torsion free if and only if $\mathcal{C}$ has the max-flow min-cut property. This fact makes a strong connection between commutative algebra and combinatorial optimization. It is known [71, Chapter 79] that a clutter $\mathcal{C}$ satisfies the max-flow min-cut property if and only if all parallelizations of the clutter $\mathcal{C}$ satisfy the König property (see Definition 3.1.7). As another application we give a proof of this fact using that the integrality of the polyhedron $\{x|\, x\geq 0; xA\geq \mathbf{1}\}$ is closed under parallelizations and minors and using that the normality of $I(\mathcal{C})$ is preserved under parallelizations and minors (Corollary 4.1.3).

A clutter $\mathcal{C}$ satisfies the *packing property* (PP for short) if all minors of $\mathcal{C}$ satisfy the König property. We say that a clutter $\mathcal{C}$ is *Cohen-Macaulay* if $R/I(\mathcal{C})$ is a Cohen-Macaulay ring, see [65]. The other main result of this chapter is:

**Theorem 3.2.3** *Let $\mathcal{C}$ be a $d$-uniform clutter on the vertex set $X$. Let*

$$Y=\{y_{ij}|\, 1\leq i\leq n;\ 1\leq j\leq d-1\}$$

*be a set of new variables, and let $\mathcal{C}'$ be the clutter with vertex set $V(\mathcal{C}')=X\cup Y$ and edge set*

$$E(\mathcal{C}')=E(\mathcal{C})\cup\{\{x_1,y_{11},\ldots,y_{1(d-1)}\},\ldots,\{x_n,y_{n1},\ldots,y_{n(d-1)}\}\}.$$

*Then the edge ideal $I(\mathcal{C}')$ is Cohen-Macaulay. If $\mathcal{C}$ satisfies* PP (*resp. max-flow min-cut*)*, then $\mathcal{C}'$ satisfies* PP (*resp. max-flow min-cut*)*.*

It is well known that if $\mathcal{C}$ satisfies the max-flow min-cut property, then $\mathcal{C}$ satisfies the packing property [16] (see Corollary 2.2.13). Conforti and Cornuéjols [15] conjecture that the converse is also true. Theorem 3.2.3 is interesting because it says that for uniform clutters it suffices to prove the conjecture for Cohen-Macaulay clutters, which have a rich structure. The Conforti-Cornuéjols conjecture has been studied in [22, 42, 45] using an algebraic approach based on certain algebraic properties of blowup algebras.

## 3.1 Normality and parallelizations

Let $\mathcal{C}$ be a clutter and let $I = I(\mathcal{C}) = (x^{v_1}, \ldots, x^{v_q})$ be its edge ideal. The incidence matrix of $\mathcal{C}$ whose column vectors are $v_1, \ldots, v_q$ will be denoted by $A$. Recall that the *Rees algebra* of $I$ is given by:

$$R[It] := R \oplus It \oplus \cdots \oplus I^i t^i \oplus \cdots \subset R[t],$$

where $t$ is a new variable. The Rees algebra of $I$ can be written as

$$R[It] \;=\; K[\{x^a t^b \,|\, (a,b) \in \mathbb{N}\mathcal{A}'\}]$$

where $\mathcal{A}' = \{(v_1, 1), \ldots, (v_q, 1), e_1, \ldots, e_n\}$, $e_i$ is the *ith* unit vector of $\mathbb{R}^{n+1}$ and $\mathbb{N}\mathcal{A}'$ is the subsemigroup of $\mathbb{N}^{n+1}$ spanned by $\mathcal{A}'$. According to [87] the integral closure of $R[It]$ in its field of fractions can be expressed as

$$\begin{aligned}\overline{R[It]} &= K[\{x^a t^b \,|\, (a,b) \in \mathbb{Z}\mathcal{A}' \cap \mathbb{R}_+\mathcal{A}'\}] \\ &= R \oplus \overline{I}t \oplus \overline{I^2}t^2 \oplus \cdots \oplus \overline{I^i}t^i \oplus \cdots,\end{aligned}$$

where $\overline{I^i}$ is the integral closure of $I^i$, $\mathbb{R}_+\mathcal{A}'$ is the cone spanned by $\mathcal{A}'$, and $\mathbb{Z}\mathcal{A}'$ is the subgroup spanned by $\mathcal{A}'$. Notice that $\mathbb{Z}\mathcal{A}' = \mathbb{Z}^{n+1}$. Hence $R[It]$ is normal if and only if any of the following two equivalent conditions hold:

(a) $\mathbb{N}\mathcal{A}' = \mathbb{Z}^{n+1} \cap \mathbb{R}_+\mathcal{A}'$, i.e., $\mathcal{A}'$ is an integral Hilbert basis of $\mathbb{R}_+\mathcal{A}'$.

(b) $I^i = \overline{I^i}$ for all $i \geq 1$.

If the second condition holds we say that $I$ is a *normal* ideal.

Let $\mathcal{C}$ be a clutter on the vertex set $X = \{x_1, \ldots, x_n\}$ and let $x_i \in X$. Then *duplicating* $x_i$ means extending $X$ by a new vertex $x_i'$ and replacing $E(\mathcal{C})$ by

$$E(\mathcal{C}) \cup \{(e \setminus \{x_i\}) \cup \{x_i'\} \,|\, x_i \in e \in E(\mathcal{C})\}.$$

The *deletion* of $x_i$, denoted by $\mathcal{C} \setminus \{x_i\}$, is the clutter formed from $\mathcal{C}$ by deleting the vertex $x_i$ and all edges containing $x_i$. A clutter obtained from $\mathcal{C}$ by a sequence of deletions and duplications of vertices is called a *parallelization*. If $w = (w_i)$ is a vector in $\mathbb{N}^n$, we denote by $\mathcal{C}^w$ the clutter obtained from $\mathcal{C}$ by deleting any vertex $x_i$ with $w_i = 0$ and duplicating $w_i - 1$ times any vertex $x_i$ if $w_i \geq 1$.

**Example 3.1.1** Let $G$ be the graph whose only edge is $\{x_1, x_2\}$ and let $w = (3, 3)$. Then $G^w = \mathcal{K}_{3,3}$ is the complete bipartite graph with bipartition $V_1 = \{x_1, x_1^2, x_1^3\}$ and $V_2 = \{x_2, x_2^2, x_2^3\}$.

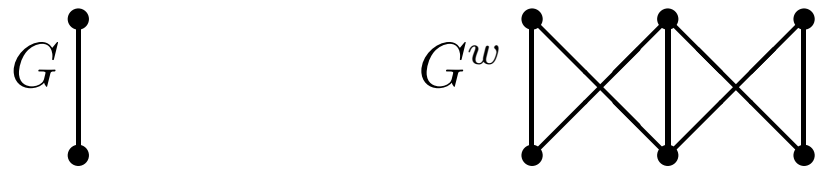

**Definition 3.1.2** Let $X' = \{x_{i_1}, \ldots, x_{i_r}, x_{j_1}, \ldots, x_{j_s}\}$ be a subset of $X$. A *minor* of $I$ is a proper ideal $I'$ of $R' = K[X \setminus X']$ obtained from $I$ by making $x_{i_k} = 0$ and $x_{j_\ell} = 1$ for all $k, \ell$. The ideal $I$ is considered to be a minor of itself. A *minor* of $\mathcal{C}$ is a clutter $\mathcal{C}'$ that corresponds to a minor $(0) \subsetneq I' \subsetneq R'$.

Notice that $\mathcal{C}'$ is obtained from $I'$ by considering the unique set of square-free monomials of $R'$ that minimally generate $I'$. If $I'$ is the ideal obtained from $I$ by making $x_i = 0$, then $I' = I(\mathcal{C} \setminus \{x_i\})$, i.e., making a variable equal to zero corresponds to a deletion.

It is known that the normality of $I(\mathcal{C})$ is closed under minors [29]. A main result of this section shows that the normality of $I(\mathcal{C})$ is closed under parallelizations.

**Theorem 3.1.3** *Let $\mathcal{C}$ be a clutter and let $\mathcal{C}'$ be a parallelization of $\mathcal{C}$. If $I(\mathcal{C})$ is normal, then $I(\mathcal{C}')$ is normal.*

**Proof.** From [29] we obtain that if $I(\mathcal{C})$ is normal and $\mathcal{C}'$ is a minor of $\mathcal{C}$, then $I(\mathcal{C}')$ is also normal. Thus we need only show that the normality of $I(\mathcal{C})$ is preserved when we duplicate a vertex of $\mathcal{C}$. Let $V(\mathcal{C}) = \{x_2, \ldots, x_n\}$ be the vertex set of $\mathcal{C}$ and let $\mathcal{C}'$ be the duplication of the vertex $x_2$. We denote the duplication of $x_2$ by $x_1$. We may assume that

$$I = I(\mathcal{C}) = (x_2 x^{w_1}, \ldots, x_2 x^{w_r}, x^{w_{r+1}}, \ldots, x^{w_q}),$$

where $x^{w_i} \in K[x_3, \ldots, x_n]$ for all $i$. We must show that the ideal

$$I(\mathcal{C}') = I + (x_1 x^{w_1}, \ldots, x_1 x^{w_r})$$

is normal. Consider the sets

$$\begin{aligned}
\mathcal{A} &= \{e_2, \ldots, e_n, (0,1,w_1,1), \ldots, (0,1,w_r,1), (0,0,w_{r+1},1), \ldots, (0,0,w_q,1)\},\\
\mathcal{A}' &= \mathcal{A} \cup \{e_1, (1,0,w_1,1), \ldots, (1,0,w_r,1)\}.
\end{aligned}$$

By hypothesis $\mathbb{Z}^{n+1} \cap \mathbb{R}_+\mathcal{A} = \mathbb{N}\mathcal{A}$. We must prove that $\mathbb{Z}^{n+1} \cap \mathbb{R}_+\mathcal{A}' = \mathbb{N}\mathcal{A}'$. It suffices to show that the left hand side is contained in the right hand side because the other inclusion always holds. Take an integral vector $(a,b,c,d)$ in $\mathbb{R}_+\mathcal{A}'$, where $a,b,d \in \mathbb{Z}$ and $c \in \mathbb{Z}^{n-2}$. Then

$$(a,b,c,d) = \sum_{i=1}^{r} \alpha_i(0,1,w_i,1) + \sum_{i=r+1}^{q} \alpha_i(0,0,w_i,1) + \sum_{i=1}^{r} \beta_i(1,0,w_i,1) + \sum_{i=1}^{n} \gamma_i e_i$$

for some $\alpha_i$, $\beta_i$, $\gamma_i$ in $\mathbb{R}_+$. Comparing entries one has

$$\begin{aligned} a &= \beta_1 + \cdots + \beta_r + \gamma_1, \\ b &= \alpha_1 + \cdots + \alpha_r + \gamma_2, \\ c &= \sum_{i=1}^{r}(\alpha_i + \beta_i)w_i + \sum_{i=r+1}^{q} \alpha_i w_i + \sum_{i=3}^{n} \gamma_i e_i, \\ d &= \sum_{i=1}^{r}(\alpha_i + \beta_i) + \sum_{i=r+1}^{q} \alpha_i. \end{aligned}$$

Consequently we obtain the equality

$$(0,a+b,c,d) = \sum_{i=1}^{r}(\alpha_i+\beta_i)(0,1,w_i,1) + \sum_{i=r+1}^{q} \alpha_i(0,0,w_i,1) + (\gamma_1+\gamma_2)e_2 + \sum_{i=3}^{n} \gamma_i e_i,$$

that is, the vector $(0,a+b,c,d)$ is in $\mathbb{Z}^{n+1} \cap \mathbb{R}_+\mathcal{A} = \mathbb{N}\mathcal{A}$. Thus there are $\lambda_i$, $\mu_i$ in $\mathbb{N}$ such that.

$$(0,a+b,c,d) = \sum_{i=1}^{r} \mu_i(0,1,w_i,1) + \sum_{i=r+1}^{q} \mu_i(0,0,w_i,1) + \sum_{i=2}^{n} \lambda_i e_i.$$

Comparing entries we obtain the equalities

$$\begin{aligned} a+b &= \mu_1 + \cdots + \mu_r + \lambda_2, \\ c &= \mu_1 w_1 + \cdots + \mu_q w_q + \lambda_3 e_3 + \cdots + \lambda_n e_n, \\ d &= \mu_1 + \cdots + \mu_q. \end{aligned}$$

Case (I): $b \leq \sum_{i=1}^{r} \mu_i$. If $b < \mu_1$, we set $b = \mu_1'$, $\mu_1' < \mu_1$, and define $\mu_1'' = \mu_1 - \mu_1'$. Otherwise pick $s \geq 2$ such that

$$\mu_1 + \cdots + \mu_{s-1} \leq b \leq \mu_1 + \cdots + \mu_s$$

Then $b = \mu_1 + \cdots + \mu_{s-1} + \mu_s'$, where $\mu_s' \leq \mu_s$. Set $\mu_s'' = \mu_s - \mu_s'$. Notice that

$$\begin{aligned} a+b &= \mu_1 + \cdots + \mu_r + \lambda_2 = a + \mu_1 + \cdots + \mu_{s-1} + \mu_s', \\ a &= \mu_s + \cdots + \mu_r + \lambda_2 - \mu_s' = \mu_{s+1} + \cdots + \mu_r + \mu_s'' + \lambda_2. \end{aligned}$$

Then

$$\begin{aligned} (a,b,c,d) &= \sum_{i=1}^{s-1} \mu_i(0,1,w_i,1) + \mu_s'(0,1,w_s,1) + \sum_{i=r+1}^{q} \mu_i(0,0,w_i,1) \\ &\quad + \mu_s''(1,0,w_s,1) + \sum_{i=s+1}^{r} \mu_i(1,0,w_i,1) + \lambda_2 e_1 + \sum_{i=3}^{n} \lambda_i e_i, \end{aligned}$$

that is, $(a,b,c,d) \in \mathbb{N}\mathcal{A}'$.

Case (II): $b > \sum_{i=1}^r \mu_i$. Then $b = \sum_{i=1}^r \mu_i + \lambda_2'$. Since

$$a + b = \mu_1 + \cdots + \mu_r + \lambda_2 = a + \mu_1 + \cdots + \mu_r + \lambda_2'$$

we get $a = \lambda_2 - \lambda_2'$. In particular $\lambda_2 \geq \lambda_2'$. Then

$$(a,b,c,d) = \sum_{i=1}^{r} \mu_i(0,1,w_i,1) + \sum_{i=r+1}^{q} \mu_i(0,0,w_i,1) + ae_1 + \lambda_2' e_2 + \sum_{i=3}^{n} \lambda_i e_i$$

that is, $(a,b,c,d) \in \mathbb{N}\mathcal{A}'$. □

Our next goal is to present some applications of this result, but first we need to prove a couple of lemmas and we need to recall some notions and results.

**Definition 3.1.4** A subset $C \subset X$ is a *minimal vertex cover* of the clutter $\mathcal{C}$ if: (i) every edge of $\mathcal{C}$ contains at least one vertex of $C$, and (ii) there is no proper subset of $C$ with the first property. If $C$ satisfies condition (i) only, then $C$ is called a *vertex cover* of $\mathcal{C}$.

**Definition 3.1.5** Let $A$ be the incidence matrix of $\mathcal{C}$. The clutter $\mathcal{C}$ satisfies the *max-flow min-cut* (MFMC) property if both sides of the LP-duality equation

$$\min\{\langle w,x\rangle \,|\, x \geq 0; xA \geq \mathbf{1}\} = \max\{\langle y,\mathbf{1}\rangle \,|\, y \geq 0; Ay \leq w\} \tag{3.1}$$

have integral optimum solutions $x$ and $y$ for each non-negative integral vector $w$.

Let $A$ be the *incidence matrix* of $\mathcal{C}$ whose column vectors are $v_1, \ldots, v_q$. The *set covering polyhedron* of $\mathcal{C}$ is given by:

$$Q(A) = \{x \in \mathbb{R}^n |\, x \geq 0;\, xA \geq \mathbf{1}\},$$

where $\mathbf{1} = (1, \ldots, 1)$. This polyhedron was studied in [42, 45] to characterize the max-flow min-cut property of $\mathcal{C}$ and to study certain algebraic properties of blowup algebras. A clutter $\mathcal{C}$ is said to be *ideal* if $Q(A)$ is an *integral polyhedron*, i.e., it has only integral vertices. The integral vertices of $Q(A)$ are precisely the characteristic vectors of the minimal vertex covers of $\mathcal{C}$ [42, Proposition 2.2].

**Theorem 3.1.6** ([29, 42, 45, 57]) *The following are equivalent*

(i) $\mathrm{gr}_I(R) = R[It]/IR[It]$ *is reduced, i.e.,* $\mathrm{gr}_I(R)$ *has no non-zero nilpotent elements.*

(ii) $R[It]$ *is normal and* $Q(A)$ *is an integral polyhedron.*

(iii) $I^i = I^{(i)}$ *for* $i \geq 1$*, where* $I^{(i)}$ *is the ith symbolic power of* $I$.

(iv) $\mathcal{C}$ *has the max-flow min-cut property.*

If condition (iii) is satisfied we say that $I$ is *normally torsion free*. A set of edges of the clutter $\mathcal{C}$ is *independent* or *stable* if no two of them have a common vertex. We denote the smallest number of vertices in any minimal vertex cover of $\mathcal{C}$ by $\alpha_0(\mathcal{C})$ and the maximum number of independent edges of $\mathcal{C}$ by $\beta_1(\mathcal{C})$. These numbers are related to min-max problems because they satisfy:

$$\begin{aligned}\alpha_0(\mathcal{C}) &\geq \min\{\langle \mathbf{1}, x\rangle |\, x \geq 0; xA \geq \mathbf{1}\} \\ &= \max\{\langle y, \mathbf{1}\rangle |\, y \geq 0; Ay \leq \mathbf{1}\} \geq \beta_1(\mathcal{C}).\end{aligned}$$

Notice that $\alpha_0(\mathcal{C}) = \beta_1(\mathcal{C})$ if and only if both sides of the equality have integral optimum solutions. These two numbers can be interpreted in terms of invariants of $I$. By [42] the height of the ideal $I$, denoted by $\mathrm{ht}(I)$, is equal to the *vertex covering number* $\alpha_0(\mathcal{C})$ and the *edge independence number* $\beta_1(\mathcal{C})$ is equal to the maximum $r$ such that there exists a regular sequence of $r$ monomials inside $I$.

**Definition 3.1.7** If $\alpha_0(\mathcal{C}) = \beta_1(\mathcal{C})$ we say that the clutter $\mathcal{C}$ (or the ideal $I$) has the *König property*.

**Definition 3.1.8** The clutter $\mathcal{C}$ (or the ideal $I$) satisfy the *packing property* (PP for short) if all its minors satisfy the König property, i.e., $\alpha_0(\mathcal{C}') = \beta_1(\mathcal{C}')$ for any minor $\mathcal{C}'$ of $\mathcal{C}$.

**Theorem 3.1.9** (A. Lehman; see [16, Theorem 1.8]) *If $\mathcal{C}$ has the packing property, then $Q(A)$ is integral.*

**Corollary 3.1.10** ([16]) *If the clutter $\mathcal{C}$ has the max-flow min-cut property, then $\mathcal{C}$ has the packing property.*

**Proof.** Assume that the clutter $\mathcal{C}$ has the max-flow min-cut property. This property is closed under taking minors. Thus it suffices to prove that $\mathcal{C}$ has the König property. We denote the incidence matrix of $\mathcal{C}$ by $A$. By hypothesis the LP-duality equation

$$\min\{\langle \mathbf{1}, x\rangle \,|\, x \geq 0; xA \geq \mathbf{1}\} = \max\{\langle y, \mathbf{1}\rangle \,|\, y \geq 0; Ay \leq \mathbf{1}\}$$

has optimum integral solutions $x$, $y$. To complete the proof notice that the left hand side of this equality is $\alpha_0(\mathcal{C})$ and the right hand side is $\beta_1(\mathcal{C})$. $\Box$

Conforti and Cornuéjols conjecture that the converse is also true:

**Conjecture 3.1.11** ([15]) If the clutter $\mathcal{C}$ has the packing property, then $\mathcal{C}$ has the max-flow min-cut property.

To the best of our knowledge this conjecture is open, see [16, Conjecture 1.6].

We are now ready to present the previously mentioned applications of Theorem 3.1.3.

**Corollary 3.1.12** *Let $\mathcal{C}$ be a clutter and let $\mathcal{C}'$ be a parallelization of $\mathcal{C}$. If $I(\mathcal{C})$ is normally torsion free, then $I(\mathcal{C}')$ is normally torsion free.*

**Proof.** Let $A$ and $A'$ be the incidence matrices of $\mathcal{C}$ and $\mathcal{C}'$ respectively. By Theorem 3.1.6 the ideal $I(\mathcal{C})$ is normal and $Q(A)$ is integral. From Theorem 3.1.3 the ideal $I(\mathcal{C}')$ is normal, and since the integrality of $Q(A)$ is closed under minors and parallelizations (see [42] and [71]) we get and $Q(A')$ is again integral. Thus applying Theorem 3.1.6 once more we get that $I(\mathcal{C}')$ is normally torsion free. $\Box$

**Corollary 3.1.13** *Let $\mathcal{C}$ be a clutter and let $\mathcal{C}'$ be a parallelization of $\mathcal{C}$. If $\mathcal{C}$ has the max-flow min-cut property, then $\mathcal{C}'$ has the König property. In particular $\mathcal{C}^w$ has the König property for all $w \in \mathbb{N}^n$.*

**Proof.** By Corollary 3.1.12 the clutter $\mathcal{C}'$ has the max-flow min-cut property. Thus applying Corollary 3.1.10 we obtain that $\mathcal{C}'$ has the König property. $\Box$

**Lemma 3.1.14** *Let $\mathcal{C}$ be a clutter and let $A$ be its incidence matrix. If $w = (w_i)$ is a vector in $\mathbb{N}^n$, then*

$$\beta_1(\mathcal{C}^w) \leq \max\{\langle y, \mathbf{1}\rangle \,|\, y \in \mathbb{N}^q;\ Ay \leq w\}.$$

**Proof.** We may assume that $w = (w_1, \ldots, w_m, 0, \ldots, 0)$, where $w_i \geq 1$ for $i = 1, \ldots, m$. Recall that for each $i$ the vertex $x_i$ is duplicated $w_i - 1$ times. We denote the duplications of $x_i$ by $x_i^2, \ldots, x_i^{w_i}$ and set $x_i^1 = x_i$. Thus the vertex set of $\mathcal{C}^w$ is equal to

$$V(\mathcal{C}^w) = \{x_1^1, \ldots, x_1^{w_1}, \ldots, x_i^1, \ldots, x_i^{w_i}, \ldots, x_m^1, \ldots, x_m^{w_m}\}.$$

There are $f_1, \ldots, f_{\beta_1}$ independent edges of $\mathcal{C}^w$, where $\beta_1 = \beta_1(\mathcal{C}^w)$. Each $f_i$ has the form

$$f_k = \{x_{k_1}^{j_{k_1}}, x_{k_2}^{j_{k_2}}, \ldots, x_{k_r}^{j_{k_r}}\} \qquad (1 \leq k_1 < \cdots < k_r \leq m;\ 1 \leq j_{k_i} \leq w_{k_i}).$$

We set $g_k = \{x_{k_1}^1, x_{k_2}^1, \ldots, x_{k_r}^1\} = \{x_{k_1}, x_{k_2}, \ldots, x_{k_r}\}$. By definition of $\mathcal{C}^w$ we get that $g_k \in E(\mathcal{C})$ for all $k$. We may re-order the $f_i$ so that

$$\underbrace{g_1 = g_2 = \cdots = g_{s_1}}_{s_1}, \underbrace{g_{s_1+1} = \cdots = g_{s_2}}_{s_2 - s_1}, \ldots, \underbrace{g_{s_{r-1}+1} = \cdots = g_{s_r}}_{s_r - s_{r-1}}$$

and $g_{s_1}, \ldots, g_{s_r}$ distinct, where $s_r = \beta_1$. Let $v_i$ be the characteristic vector of $g_{s_i}$. Set $y = s_1 e_1 + (s_2 - s_1)e_2 + \cdots + (s_r - s_{r-1})e_r$. We may assume that the incidence matrix $A$ of $\mathcal{C}$ has column vector $v_1, \ldots, v_q$. Then $y$ satisfies $\langle y, \mathbf{1}\rangle = \beta_1$. For each $k_i$ the number of variables of the form $x_{k_i}^\ell$ that occur in $f_1, \ldots, f_{\beta_1}$ is at most $w_{k_i}$ because the $f_i$ are pairwise disjoint. Hence for each $k_i$ the number of times that the variable $x_{k_i}^1$ occurs in $g_1, \ldots, g_{\beta_1}$ is at most $w_{k_i}$. Then

$$Ay = s_1 v_1 + (s_2 - s_1)v_2 + \cdots + (s_r - s_{r-1})v_r \leq w.$$

Therefore we obtain the required inequality. $\Box$

Let $\mathcal{C}$ be a clutter. For use below we denote the set of minimal vertex covers of $\mathcal{C}$ by $\Upsilon(\mathcal{C})$.

**Lemma 3.1.15** *Let $\mathcal{C}$ be a clutter and let $A$ be its incidence matrix. If $w \in \mathbb{N}^n$, then*

$$\min\left\{\sum_{x_i \in C} w_i \,\middle|\, C \in \Upsilon(\mathcal{C})\right\} = \alpha_0(\mathcal{C}^w).$$

**Proof.** We may assume that $w = (w_1, \ldots, w_m, w_{m+1}, \ldots, w_{m_1}, 0, \ldots, 0)$, where $w_i \geq 2$ for $i = 1, \ldots, m$, $w_i = 1$ for $i = m+1, \ldots, m_1$, and $w_i = 0$ for $i > m_1$. Thus for $i = 1, \ldots, m$ the vertex $x_i$ is duplicated $w_i - 1$ times. We denote the duplications of $x_i$ by $x_i^2, \ldots, x_i^{w_i}$ and set $x_i^1 = x_i$.

We first prove that the left hand side is less or equal than the right hand side. Let $C^w$ be a minimal vertex cover of $\mathcal{C}^w$ with $\alpha_0$ elements, where $\alpha_0 = \alpha_0(\mathcal{C}^w)$. We may assume that $C^w \cap \{x_1, \ldots, x_{m_1}\} = \{x_1, \ldots, x_s\}$. Note that $x_i^1, \ldots, x_i^{w_i}$ are in $C^w$ for $i = 1, \ldots, s$. Indeed since $C^w$ is a minimal vertex cover of $\mathcal{C}^w$, there exists an edge $e$ of $\mathcal{C}^w$ such that $e \cap C^w = \{x_i^1\}$. Then $(e \setminus \{x_i^1\}) \cup \{x_i^j\}$ is an edge of $C^w$ for $j = 1, \ldots, w_i$. Consequently $x_i^j \in C^w$ for $j = 1, \ldots, w_i$. Hence

$$w_1 + \cdots + w_s \leq |C^w| = \alpha_0. \tag{3.2}$$

On the other hand the set $C' = \{x_1, \ldots, x_s\} \cup \{x_{m_1+1}, \ldots, x_n\}$ is a vertex cover of $\mathcal{C}$. Let $D$ be a minimal vertex cover of $\mathcal{C}$ contained in $C'$. Let $e_D$ denote the characteristic vector of $D$. Then, since $w_i = 0$ for $i > m_1$, using Eq. (3.2) we get

$$\langle w, e_D \rangle = \sum_{x_i \in D} w_i = \sum_{x_i \in D \cap \{x_1, \ldots, x_s\}} w_i \quad \leq \sum_{x_i \in \{x_1, \ldots, x_s\}} w_i \;\leq\; \alpha_0.$$

This completes the proof of the asserted inequality.

Next we show that the right hand side of the inequality is less or equal than the left hand side. Let $C$ be a minimal vertex cover of $\mathcal{C}$. Note that the set

$$C' = \cup_{x_i \in C} \{x_i^1, \ldots, x_i^{w_i}\}$$

is a vertex cover of $\mathcal{C}^w$. Indeed any edge $e^w$ of $\mathcal{C}^w$ has the form $e^w = \{x_{i_1}^{j_1}, \ldots, x_{i_r}^{j_r}\}$ for some edge $e = \{x_{i_1}, \ldots, x_{i_r}\}$ of $\mathcal{C}$ and since $e$ is covered by $C$, we have that $e^w$ is covered by $C'$. Hence $\alpha_0(\mathcal{C}^w) \leq |C'| = \sum_{x_i \in C} w_i$. As $C$ was an arbitrary vertex cover of $\mathcal{C}$ we get the asserted inequality. $\square$

**Corollary 3.1.16** [71, Chapter 79] *Let $\mathcal{C}$ be a clutter. Then $\mathcal{C}$ satisfies the max-flow min-cut property if and only if $\beta_1(\mathcal{C}^w) = \alpha_0(\mathcal{C}^w)$ for all $w \in \mathbb{N}^n$.*

**Proof.** If $\mathcal{C}$ has the max-flow min-cut property, then $C^w$ has the König property by Corollary 3.1.13. Conversely if $\mathcal{C}^w$ has the König property for all $w \in \mathbb{N}^n$, then by Lemmas 3.1.14 and 3.1.15 both sides of the LP-duality equation

$$\min\{\langle w, x\rangle \,|\, x \geq 0; xA \geq \mathbf{1}\} = \max\{\langle y, \mathbf{1}\rangle \,|\, y \geq 0; Ay \leq w\}$$

have integral optimum solutions $x$ and $y$ for each non-negative integral vector $w$, i.e., $\mathcal{C}$ has the max-flow min-cut property. □

## 3.2 Cohen-Macaulay ideals with MFMC

One of the aims here is to show how to construct Cohen-Macaulay clutters satisfying max-flow min-cut, PP, and normality properties. Let $\mathcal{C}$ be a uniform clutter. A main result of this section proves that if $\mathcal{C}$ satisfies PP (resp. max-flow min-cut), then there is a uniform Cohen-Macaulay clutter $\mathcal{C}_1$ satisfying PP (resp. max-flow min-cut) such that $\mathcal{C}$ is a minor of $\mathcal{C}_1$. In particular for uniform clutters we prove that it suffices to show Conjecture 3.1.11 for Cohen-Macaulay clutters (see Corollary 3.2.4).

Let $R = K[x_1, \ldots, x_n]$ be a polynomial ring over a field $K$ and let $\mathcal{C}$ be a clutter on the vertex set $X$. As usual, in what follows, we denote the edge ideal of $\mathcal{C}$ by $I = I(\mathcal{C})$. Recall that $\mathfrak{p}$ is a minimal prime of $I = I(\mathcal{C})$ if and only if $\mathfrak{p} = (C)$ for some minimal vertex cover $C$ of $\mathcal{C}$ [87, Proposition 6.1.16]. Thus the primary decomposition of the edge ideal of $\mathcal{C}$ is given by

$$I(\mathcal{C}) = (C_1) \cap (C_2) \cap \cdots \cap (C_p),$$

where $C_1, \ldots, C_p$ are the minimal vertex covers of $\mathcal{C}$. In particular observe that the *height* of $I(\mathcal{C})$, denoted by $\operatorname{ht} I(\mathcal{C})$, equals the number of vertices in a minimum vertex cover of $\mathcal{C}$. Also notice that the associated primes of $I(\mathcal{C})$ are precisely the minimal primes of $I(\mathcal{C})$.

**Proposition 3.2.1** *Let $R[z_1, \ldots, z_\ell]$ be a polynomial ring over $R$. If $I$ is a normal ideal of $R$, then $J = (I, x_1 z_1 \cdots z_\ell)$ is a normal ideal of $R[z_1, \ldots, z_\ell]$.*

**Proof.** By induction on $p$ we will show $\overline{J^p} = J^p$ for all $p \geq 1$. If $p = 1$, then $\overline{J} = J$ because $J$ is square-free. Assume $\overline{J^i} = J^i$ for $i < p$ and $p \geq 2$. Let $y$ be a monomial in $\overline{J^p}$, then $y^m \in J^{pm}$, for some $m > 0$. Since $\overline{J^p} \subset \overline{J^{p-1}} = J^{p-1}$ we can write

$$y = z_1^{t_1} \cdots z_\ell^{t_\ell}(x_1 z_1 \cdots z_\ell)^r M f_1 \cdots f_{p-r-1},$$

where $M$ is a monomial with $z_i \notin \mathrm{supp}(M)$ for all $i$ and the $f_i$'s are monomials in $J$ with $z_i \notin \mathrm{supp}(f_j)$ for all $i, j$. We set $h = M f_1 \cdots f_{p-r-1}$. It suffices to show that $y \in J^p$. Since $y^m \in J^{pm}$ we have

$$y^m = z_1^{mt_1} \cdots z_\ell^{mt_\ell}(x_1 z_1 \cdots z_\ell)^{rm} h^m = N(x_1 z_1 \cdots z_\ell)^s g_1 \cdots g_{mp-s}, \tag{3.3}$$

where $N$ is a monomial, $z_i \notin \mathrm{supp}(g_j)$ for all $i, j$, and the $g_i$'s are monomials in $J$. We distinguish two cases:

Case (a): Assume $t_i = 0$ for some $i$, then $s \leq rm$ because $z_i^{rm}$ is the maximum power of $z_i$ that divides $y^m$. Making $z_j = 1$ for $j = 1, \ldots, \ell$ in Eq. (3.3) we get

$$x_1^{rm-s} h^m = N' g_1 \cdots g_{mp-s}.$$

Thus $h^m \in I^{(mp-s)-(rm-s)} = I^{m(p-r)}$. Therefore we get $h \in \overline{I^{p-r}} = I^{p-r}$ and $y = z_1^{t_1} \cdots z_\ell^{t_\ell}(x_1 z_1 \cdots z_\ell)^r h \in J^p$.

Case (b): If $t_i > 0$ for all $i$, we may assume $x_1 \notin \mathrm{supp}(M)$, otherwise $y \in J^p$. We may also assume $x_1 \notin \mathrm{supp}(f_i)$ for all $i$, otherwise it is not hard to see that we are back in case (a). Notice that $s \leq rm$, because $x_1 \notin \mathrm{supp}(h)$. From Eq. (3.3) it follows that $h \in \overline{I^{p-r}} = I^{p-r}$ and $y = z_1^{t_1} \cdots z_\ell^{t_\ell}(x_1 z_1 \cdots z_\ell)^r h \in J^p$. □

**Lemma 3.2.2** *Let $R[z_1, \ldots, z_\ell]$ be a polynomial ring over $R$ and let $I_1$ be the ideal obtained from $I$ by making $x_1 = 0$. Then:* (a) *if $I$ and $I_1$ satisfy the König property, then the ideal $J = (I, x_1 z_1 \cdots z_\ell)$ satisfies the König property, and* (b) *if $I$ satisfies* PP*, then $J$ satisfies* PP*.*

**Proof.** (a): If $\mathrm{ht}(I) = \mathrm{ht}(J)$, then $J$ satisfies König because $I$ does. Assume that $g = \mathrm{ht}(I) < \mathrm{ht}(J)$. Then $\mathrm{ht}(J) = g + 1$. Notice that every associated prime ideal of $I$ of height $g$ cannot contain $x_1$. We claim that $\mathrm{ht}(I_1) = g$. If $r = \mathrm{ht}(I_1) < g$, pick a minimal prime $\mathfrak{p}$ of $I_1$ of height $r$. Then $\mathfrak{p} + (x_1)$ is a prime ideal of height at most $g$ containing both $I$ and $x_1$, a contradiction. This proves the claim. Since $I_1$ satisfies König, there are $g$ independent monomials in $I_1$. Hence $h_1, \ldots, h_g, x_1 z_1 \cdots z_\ell$ are $g + 1$ independent monomials in $J$, as required. Part (b) follows readily from part (a). □

**Theorem 3.2.3** *Let $\mathcal{C}$ be a $d$-uniform clutter on the vertex set $X$. Let*

$$Y = \{y_{ij} | \, 1 \leq i \leq n; \; 1 \leq j \leq d-1\}$$

*be a set of new variables, and let $\mathcal{C}'$ be the clutter with vertex set $V(\mathcal{C}') = X \cup Y$ and edge set*

$$E(\mathcal{C}') = E(\mathcal{C}) \cup \{\{x_1, y_{11}, \ldots, y_{1(d-1)}\}, \ldots, \{x_n, y_{n1}, \ldots, y_{n(d-1)}\}\}.$$

*Then the edge ideal $I(\mathcal{C}')$ is Cohen-Macaulay. If $\mathcal{C}$ satisfies* PP (*resp. max-flow min-cut*), *then $\mathcal{C}'$ satisfies* PP (*resp. max-flow min-cut*).

**Proof.** Set $S = K[X \cup Y]$. First we prove that $I' = I(\mathcal{C}')$ is Cohen-Macaulay. Consider the ideal $I_0 = I + (x_1^d, \ldots, x_n^d)$. The idea is to show that, up to a relabeling of the variables, $S/I_0$ is a deformation of $S/I'$ by a regular sequence. Let $J$ be a copy of $I$ obtained from $I$ by making $x_i = y_{i1}$ for $i = 1, \ldots, n$ and let $I_1 = J + (y_{11}x_1^{d-1}, \ldots, y_{n1}x_n^{d-1})$. By [36] the sequence

$$\underline{h_1} = \{x_1 - y_{11}, \ldots, x_n - y_{n1}\}$$

is a regular sequence on $S/I_1$ and $(S/I_1)/(\underline{h_1}) = S/I_0$. The ring $S/I_0$ is Cohen-Macaulay because $R/I_0$ has dimension zero, hence $S/I_1$ is Cohen-Macaulay. Next consider the ideal $I_2 = J + (y_{11}y_{12}x_1^{d-2}, \ldots, y_{n1}y_{n2}x_n^{d-2})$. Again using [36] we get that the sequence

$$\underline{h_2} = \{x_1 - y_{12}, \ldots, x_n - y_{n2}\}$$

is a regular sequence on $S/I_2$ and $(S/I_2)/(\underline{h_2}) = S/I_1$. Thus $S/I_2$ is Cohen-Macaulay. A repeated use of this construction yields that $S/I_{d-1}$ is Cohen-Macaulay and $S/I' \simeq S/I_{d-1}$. The Cohen-Macaulayness of $I'$ can also be shown using [32, Theorem 8.2] or [66, Theorem 2.16]. If $\mathcal{C}$ satisfies PP, then from Lemma 3.2.2(b) it follows that $\mathcal{C}'$ satisfies PP. Assume that $\mathcal{C}$ satisfies MFMC. By Proposition 3.2.1 $S[I't]$ is normal. Since $\mathcal{C}'$ satisfies PP, by Lehman's theorem we get that $Q(A')$ is integral, where $A'$ is the incidence matrix of $\mathcal{C}'$. Therefore using Theorem 3.1.6 we conclude that $\mathcal{C}'$ has MFMC. □

Recall that a clutter $\mathcal{C}$ is called *Cohen-Macaulay* (CM for short) if $R/I(\mathcal{C})$ is a Cohen-Macaulay ring. Since $\mathcal{C}$ is a minor of $\mathcal{C}'$ we obtain:

**Corollary 3.2.4** *Let $\mathcal{C}$ be a uniform clutter. If $\mathcal{C}$ satisfies* PP (*resp. max-flow min-cut*), *then there is a uniform Cohen-Macaulay clutter $\mathcal{C}_1$ satisfying* PP (*resp. max-flow min-cut*) *such that $\mathcal{C}$ is a minor of $\mathcal{C}_1$.*

This result is interesting because it says that for uniform clutters it suffices to prove Conjecture 3.1.11 for Cohen-Macaulay clutters, which have a rich structure.

# Chapter 4

# Edge Ideals of Clique Clutters of Comparability Graphs

Let $(P, \prec)$ be a finite poset and let $G$ be its comparability graph. If $\mathrm{cl}(G)$ is the clutter of maximal cliques of $G$, we prove that $\mathrm{cl}(G)$ satisfies the max-flow min-cut property and that its edge ideal is normally torsion free. We prove that edge ideals of complete admissible uniform clutters are normally torsion free.

Let $(P, \prec)$ be a *partially ordered set* (*poset* for short) on the finite vertex set $X = \{x_1, \ldots, x_n\}$ and let $G$ be its *comparability graph*. Recall that the vertex set of $G$ is $X$ and the edge set of $G$ is the set of all unordered pairs $\{x_i, x_j\}$ such that $x_i$ and $x_j$ are comparable. A *clique* of $G$ is a subset of the set of vertices that induces a complete subgraph.

Let $\mathcal{C}$ be a *clutter* with finite vertex set $X$, that is, $\mathcal{C}$ is a family of subsets of $X$, called edges, none of which is included in another.

The set of vertices and edges of $\mathcal{C}$ are denoted by $V(\mathcal{C})$ and $E(\mathcal{C})$ respectively. The *incidence matrix* of $\mathcal{C}$ is the vertex-edge matrix whose columns are the characteristic vectors of the edges of $\mathcal{C}$. Let $R = K[x_1, \ldots, x_n]$ be a polynomial ring over a field $K$. The *edge ideal* of $\mathcal{C}$, denoted by $I(\mathcal{C})$, is the ideal of $R$ generated by all monomials $\prod_{x_i \in e} x_i$ such that $e \in E(\mathcal{C})$. The map $\mathcal{C} \mapsto I(\mathcal{C})$ gives a one to one correspondence between the family of clutters and the family of square-free monomial ideals. Edge ideals of clutters also correspond to simplicial complexes via the Stanley-Reisner correspondence [75].

The *clique clutter* of $G$, denoted by $\mathrm{cl}(G)$, is the clutter with vertex set $X$ whose edges are exactly the maximal cliques of $G$ with respect to inclusion.

One of our main algebraic results shows that the edge ideal $I = I(\mathrm{cl}(G))$ of $\mathrm{cl}(G)$ is normally torsion free (see Theorem 4.2.2), i.e., $I^i = I^{(i)}$ for $i \geq 1$, where $I^{(i)}$ is the *ith* symbolic power of $I$ (see Section 4.2). To prove this result we first show that the clique clutter of $G$ has the max-flow min-cut property (see Theorem 4.1.7). Then we use a remarkable result of [45] showing that an edge ideal $I(\mathcal{C})$, of a clutter $\mathcal{C}$, is normally torsion free if and only if $\mathcal{C}$ has the max-flow min-cut property. This fact makes a strong connection between commutative algebra and combinatorial optimization. Some other interesting links between these two areas can be found in [29, 42, 47, 48, 53, 81, 91] and in the references there. There are some other nice characterizations of the normally torsion free property that can be found in [57].

Let $f_1, \ldots, f_q$ be the maximal cliques of $G$ and let $v_k = \sum_{x_i \in f_k} e_i$ be the characteristic vector of $f_k$ for $k = 1, \ldots, q$, where $e_i$ is the *ith* unit vector of $\mathbb{R}^n$. The matrix $A$ with column vectors $v_1, \ldots, v_q$ is called the *vertex-clique matrix* of $G$ or the *incidence matrix* of $\mathrm{cl}(G)$. A *colouring* of the vertices of a graph is an assignment of colours to the vertices of the graph in such a way that adjacent vertices have distinct colours. The *chromatic number* of a graph is the minimal number of colours in a colouring. A graph is called *perfect* if for every induced subgraph $H$, the chromatic number of $H$ equals the size of the largest complete subgraph of $H$. It is well known that comparability graphs are perfect [71]. Thus $G$ is a perfect graph or equivalently the polytope

$$P(A) = \{x \mid x \geq 0;\ xA \leq \mathbf{1}\}$$

is integral, i.e., it has only integral vertices. Here $\mathbf{1}$ denotes the vector with all its entries equal to 1. We complement this fact by observing that the *set covering polyhedron*

$$Q(A) = \{x \mid x \geq 0;\ xA \geq \mathbf{1}\}$$

is also integral (see Corollary 4.1.9). Comparability graphs are interesting objects of study [3]. They have been nicely characterized in [39]. By [2] the clique clutter of any induced subgraph of $G$ satisfies the König property (see Definition 4.1.2). In the terminology of [8] comparability graphs are clique-perfect.

Let $I$ be a monomial ideal. Recall that the *integral closure* of $I^i$, denoted by $\overline{I^i}$, is the ideal of $R$ given by

$$\overline{I^i} = (\{x^a \in R \mid \exists\, p \in \mathbb{N} \setminus \{0\};\ (x^a)^p \in I^{pi}\}),$$

see for instance [87, Proposition 7.3.3]. As usual, we will use $x^a$ as an abbreviation for $x_1^{a_1} \cdots x_n^{a_n}$, where $a = (a_i)$ is a vector in $\mathbb{N}^n$. The ideal $I$ is called *normal* if $I^i = \overline{I^i}$ for all $i$. If $\overline{I} = I$, the ideal $I$ is called *integrally closed*. The normality property is one of most interesting properties an ideal can have, see [58, 84].

A subset $C \subset X$ is called a *minimal vertex cover* of the clutter $\mathcal{C}$ if: (i) every edge of $\mathcal{C}$ contains at least one vertex of $C$, and (ii) there is no proper subset of $C$ with the first property. Let $D_1, \ldots, D_s$ be the minimal vertex covers of $\overline{G}$, where $\overline{G}$ is the complement of $G$. The ideal of *vertex covers* of $\overline{G}$ is the square-free monomial ideal

$$I_c(\overline{G}) = (x^{u_1}, \ldots, x^{u_s}) \subset R,$$

where $x^{u_k} = \prod_{x_i \in D_k} x_i$. Note that $I_c(\overline{G})$ and $I(\mathrm{cl}(G))$ are dual of each other in the sense that $u_i + v_i = \mathbf{1}$ for $i = 1, \ldots, q$. In [91] it is shown that $I_c(\overline{G})$ is a normal ideal. We complement this fact by showing that the edge ideal $I(\mathrm{cl}(G))$ of $\mathrm{cl}(G)$ is a normal ideal (see Theorem 4.2.2). This is surprising because in general this duality does not preserve normality. As an application we prove that edge ideals of complete admissible uniform clutters are normally torsion free (see Theorem 4.2.3).

Along the chapter we introduce most of the notions that are relevant for our purposes. Our main references for combinatorial optimization and commutative algebra are [10, 71, 75, 84, 87]. In these references the reader will find the undefined terminology and notation that we use in what follows.

## 4.1 Maximal cliques of comparability graphs

In this section we introduce the max-flow min-cut property and prove our main combinatorial result, that is, we prove that the clique clutter of a comparability graph satisfies the max-flow min-cut property.

**Definition 4.1.1** Let $\mathcal{C}$ be a clutter and let $A$ be its incidence matrix. The clutter $\mathcal{C}$ satisfies the *max-flow min-cut* property if both sides of the LP-duality equation

$$\min\{\langle w, x\rangle \,|\, x \geq 0; xA \geq \mathbf{1}\} = \max\{\langle y, \mathbf{1}\rangle \,|\, y \geq 0; Ay \leq w\} \tag{4.1}$$

have integral optimum solutions $x$ and $y$ for each non-negative integral vector $w$.

Let $\mathcal{C}$ be a clutter. A set of edges of $\mathcal{C}$ is *independent* or *stable* if no two of them have a common vertex. We denote the smallest number of vertices in any minimal vertex cover of $\mathcal{C}$ by $\alpha_0(\mathcal{C})$ and the maximum number of independent edges of $\mathcal{C}$ by $\beta_1(\mathcal{C})$. These two numbers satisfy $\beta_1(\mathcal{C}) \leq \alpha_0(\mathcal{C})$.

**Definition 4.1.2** If $\beta_1(\mathcal{C}) = \alpha_0(\mathcal{C})$ we say that $\mathcal{C}$ has the *König property*.

Let $\mathcal{C}$ be a clutter on the vertex set $X = \{x_1, \ldots, x_n\}$ and let $x_i \in X$. Then *duplicating* $x_i$ means extending $X$ by a new vertex $x_i'$ and replacing $E(\mathcal{C})$ by

$$E(\mathcal{C}) \cup \{(e \setminus \{x_i\}) \cup \{x_i'\} | \, x_i \in e \in E(\mathcal{C})\}.$$

The *deletion* of $x_i$, denoted by $\mathcal{C} \setminus \{x_i\}$, is the clutter formed from $\mathcal{C}$ by deleting the vertex $x_i$ and all edges containing $x_i$. A clutter obtained from $\mathcal{C}$ by a sequence of deletions and duplications of vertices is called a *parallelization*. If $w = (w_i)$ is a vector in $\mathbb{N}^n$, we denote by $\mathcal{C}^w$ the clutter obtained from $\mathcal{C}$ by deleting any vertex $x_i$ with $w_i = 0$ and duplicating $w_i - 1$ times any vertex $x_i$ if $w_i \geq 1$.

The notion of parallelization can be used to give a characterization of the max-flow min-cut property which is suitable to study the clique clutter of the comparability graph of a poset.

**Theorem 4.1.3** [71, Chapter 79] *Let $\mathcal{C}$ be a clutter. Then $\mathcal{C}$ satisfies the max-flow min-cut property if and only if $\beta_1(\mathcal{C}^w) = \alpha_0(\mathcal{C}^w)$ for all $w \in \mathbb{N}^n$.*

**Lemma 4.1.4** *Let* $\mathrm{cl}(G)$ *be the clutter of maximal cliques of a graph $G$. If $G^1$ (resp.* $\mathrm{cl}(G)^1$*) is the graph (resp. clutter) obtained from $G$ (resp.* $\mathrm{cl}(G)$*) by duplicating the vertex $x_1$, then* $\mathrm{cl}(G)^1 = \mathrm{cl}(G^1)$.

**Proof.** Let $y_1$ be the duplication of $x_1$. Set $\mathcal{C} = \mathrm{cl}(G)$. First we prove that $E(\mathcal{C}^1) \subset E(\mathrm{cl}(G^1))$. Take $e \in E(\mathcal{C}^1)$. Case (i): Assume $y_1 \notin e$. Then $e \in E(\mathcal{C})$. Clearly $e$ is a clique of $G^1$. If $e \notin E(\mathrm{cl}(G^1))$, then $e$ can be extended to a maximal clique of $G^1$. Hence $e \cup \{y_1\}$ must be a clique of $G^1$. Note that $x_1 \notin e$ because $\{x_1, y_1\}$ is not an edge of $G^1$. Then $e \cup \{x_1\}$ is a clique of $G$, a contradiction. Thus $e$ is in $E(\mathrm{cl}(G^1))$. Case (ii): Assume $y_1 \in e$. Then there is $f \in E(\mathrm{cl}(G))$, with $x_1 \in f$, such that $e = (f \setminus \{x_1\}) \cup \{y_1\}$. Since $\{x, x_1\} \in E(G)$ for any $x$ in $f \setminus \{x_1\}$, one has that $\{x, y_1\} \in E(G^1)$ for any $x$ in $f \setminus \{x_1\}$. Then $e$ is a clique of $G^1$. If $e$ is not a maximal clique of $G^1$, there is $x \notin e$ which is adjacent in $G$ to any vertex of $f \setminus \{x_1\}$ and $x$ is adjacent to $y_1$ in $G^1$. In particular $x \neq x_1$. Then $x$ is adjacent in $G$ to $x_1$ and consequently $x$ is adjacent in $G$ to any vertex of $f$, a contradiction because $f$ is a maximal clique of $G$. Thus $e$ is in $\mathrm{cl}(G^1)$. Next we prove the inclusion $E(\mathrm{cl}(G^1)) \subset E(\mathcal{C}^1)$. Take $e \in E(\mathrm{cl}(G^1))$, i.e., $e$ is a maximal clique of $G^1$. Case (i): Assume $y_1 \notin e$. Then $e$ is a maximal clique of $G$, and so an edge of $\mathcal{C}^1$. Case (ii): Assume $y_1 \in e$. Then $e \setminus \{y_1\}$ is a clique of $G$ and $\{x, y_1\} \in E(G^1)$ for any $x$ in $e \setminus \{y_1\}$. Then $\{x, x_1\}$ is in $E(G)$ for any $x$ in $e \setminus \{y_1\}$. Hence $f = (e \setminus \{y_1\}) \cup \{x_1\}$ is a clique of $G$. Note that $f$ is a maximal clique of $G$. Indeed if $f$ is not a maximal clique of $G$, there is $x \in V(G) \setminus f$ which is adjacent in $G$ to every vertex of $e \setminus \{y_1\}$ and to $x_1$. Thus $x$ is adjacent to $y_1$ in $G^1$ and to every vertex in $e \setminus \{y_1\}$, i.e., $e \cup \{x\}$ is a clique of $G^1$, a contradiction. Thus $f \in \mathrm{cl}(G)$. Since $e = (f \setminus \{x_1\}) \cup \{y_1\}$ we obtain that $e \in E(\mathcal{C}^1)$. □

Unfortunately we do not have an analogous version of Lemma 4.1.4 valid for a deletion. In other words, if $G$ is a graph, the equality $\mathrm{cl}(G)^w = \mathrm{cl}(G^w)$, with $w$ an integral vector, fails in general (see Remark 4.1.5).

**Remark 4.1.5** Let $G$ be a graph. Let $G^1 = G \setminus \{x_1\}$ (resp. $\mathrm{cl}(G)^1 = \mathrm{cl}(G) \setminus \{x_1\}$) be the graph (resp. clutter) obtained from $G$ (resp. $\mathrm{cl}(G)$) by deleting the vertex $x_1$. The equality $\mathrm{cl}(G)^1 = \mathrm{cl}(G^1)$ fails in general. For instance if $G$ is a cycle of length three, then $E(\mathrm{cl}(G)^1) = \emptyset$ and $\mathrm{cl}(G^1)$ has exactly one edge.

Let $\mathcal{D}$ be a *digraph*, that is, $\mathcal{D}$ consists of a finite set $V(\mathcal{D})$ of vertices and a set $E(\mathcal{D})$ of ordered pairs of distinct vertices called edges. Let $A$, $B$ be two sets of vertices of $\mathcal{D}$. For use below recall that a (directed) path of $\mathcal{D}$ is called an *$A$–$B$ path* if it runs from a vertex in $A$ to a vertex in $B$. A set $C$ of vertices is called an *$A$–$B$ disconnecting* set if $C$ intersects each $A$–$B$ path. For convenience we recall the following classical result.

**Theorem 4.1.6** (Menger's theorem, see [71, Theorem 9.1]) *Let $\mathcal{D}$ be a digraph and let $A$, $B$ be two subsets of $V(\mathcal{D})$. Then the maximum number of vertex-disjoint $A$–$B$ paths is equal to the minimum size of an $A$–$B$ disconnecting vertex set.*

We come to the main result of this section.

**Theorem 4.1.7** *Let $(P, \prec)$ be a finite poset on the vertex set $X = \{x_1, \ldots, x_n\}$ and let $G$ be its comparability graph. If $\mathcal{C} = \mathrm{cl}(G)$ is the clutter of maximal cliques of $G$, then $\mathcal{C}$ satisfies the max-flow min-cut property.*

**Proof.** We can regard $P$ as a transitive digraph without cycles of length two with vertex set $X$ and edge set $E(P)$, i.e., the edges of $P$ are ordered pairs $(a,b)$ of distinct vertices with $a \prec b$ such that:

(i) $(a,b) \in E(P)$ and $(b,c) \in E(P) \Rightarrow (a,c) \in E(P)$ and
(ii) $(a,b) \in E(P) \Rightarrow (b,a) \notin E(P)$.

Note that because of the transitivity condition $P$ is in fact an acyclic digraph, that is, it has no directed cycles. Let $x_1$ be a vertex of $P$ and let $y_1$ be a new vertex. Consider the digraph $P^1$ with vertex set $X^1 = X \cup \{y_1\}$ and edge set

$$E(P^1) = E(P) \cup \{(y_1,x)|\,(x_1,x) \in E(P)\} \cup \{(x,y_1)|\,(x,x_1) \in E(P)\}.$$

The digraph $P^1$ is transitive. Indeed let $(a,b)$ and $(b,c)$ be two edges of $P^1$. If $y_1 \notin \{a,b,c\}$, then $(a,c) \in E(P) \subset E(P^1)$ because $P$ is transitive. If $y_1 = a$, then $(x_1,b)$ and $(b,c)$ are in $E(P)$. Hence $(x_1,c) \in E(P)$ and $(y_1,c) \in E(P^1)$. The cases $y_1 = b$ and $y_1 = c$ are treated similarly. Thus $P^1$ defines a poset $(P^1, \prec^1)$. The comparability graph $H$ of $P^1$ is precisely the graph $G^1$ obtained from $G$ by duplicating the vertex $x_1$ by the vertex $y_1$. To see this note that $\{x,y\}$ is an edge of $G^1$ if and only if $\{x,y\}$ is an edge of $G$ or $y = y_1$ and $\{x,x_1\}$ is an edge of $G$. Thus $\{x,y\}$ is an edge of $G^1$ if and only if $x$ is related to $y$ in $P$ or $y = y_1$ and $x$ is related to $y$ in $P^1$, i.e., $\{x,y\}$ is an edge of $G^1$ if and only if $\{x,y\}$ is an edge of $H$. From Lemma 4.1.4 we get that $\mathrm{cl}(G)^1 = \mathrm{cl}(G^1)$, where $\mathrm{cl}(G)^1$ is the clutter obtained from $\mathrm{cl}(G)$ by duplicating the vertex $x_1$ by the vertex $y_1$. Altogether we obtain that the clutter $\mathrm{cl}(G)^1$ is the clique clutter of the comparability graph $G^1$ of the poset $P^1$.

By Theorem 4.1.3 it suffices to prove that $\mathrm{cl}(G)^w$ has the König property for all $w \in \mathbb{N}^n$. Since duplications commute with deletions, by permuting vertices, we may assume that $w = (w_1,\ldots,w_r,0,\ldots,0)$, where $w_i \geq 1$ for $i = 1,\ldots,r$. Consider the clutter $\mathcal{C}_1$ obtained from $\mathrm{cl}(G)$ by duplicating $w_i - 1$ times the vertex $x_i$ for $i = 1,\ldots,r$. We denote the vertex set of $\mathcal{C}_1$ by $X_1$. By successively applying the fact that $\mathrm{cl}(G)^1 = \mathrm{cl}(G^1)$, we conclude that there is a poset $P_1$ with comparability graph $G_1$ and vertex set $X_1$ such that $\mathcal{C}_1 = \mathrm{cl}(G_1)$. As before we regard $P_1$ as a transitive acyclic digraph.

Let $A$ and $B$ be the set of minimal and maximal elements of the poset $P_1$, i.e., the elements of $A$ and $B$ are the sources and sinks of $P_1$ respectively. We set $S = \{x_{r+1},\ldots,x_n\}$. Consider the digraph $\mathcal{D}$ whose vertex set is $V(\mathcal{D}) = X_1 \setminus S$ and whose edge set is defined as follows. A pair $(x,y)$ in $V(\mathcal{D}) \times V(\mathcal{D})$ is in $E(\mathcal{D})$ if and only if $(x,y) \in E(P_1)$ and there is no vertex $z$ in $X_1$ with $x \prec z \prec y$. Notice that $\mathcal{D}$ is a sub-digraph of $P_1$ which is not necessarily the digraph of a poset. We set $A_1 = A \setminus S$ and $B_1 = B \setminus S$. Note that $\mathcal{C}^w = \mathcal{C}_1 \setminus S$, the clutter obtained from $\mathcal{C}_1$ by removing all vertices of $S$ and all edges sharing a vertex with $S$. If every edge of $\mathcal{C}_1$ intersects $S$, then $E(\mathcal{C}^w) = \emptyset$ and there is nothing to prove. Thus we may assume that there is a maximal clique $K$ of $G_1$ disjoint form $S$. Note that by the maximality of $K$ and by the transitivity of $P_1$ we get that $K$ contains at least one source and one sink of $P_1$, i.e., $A_1 \neq \emptyset$ and $B_1 \neq \emptyset$ (see argument below).

The maximal cliques of $G_1$ not containing any vertex of $S$ correspond exactly to the $A_1$–$B_1$ paths of $\mathcal{D}$. Indeed let $c = \{v_1, \ldots, v_s\}$ be a maximal clique of $G_1$ disjoint from $S$. Consider the sub-poset $P_c$ of $P_1$ induced by $c$. Note that $P_c$ is a tournament, i.e., $P_c$ is an oriented graph (no-cycles of length two) such that any two vertices of $P_c$ are comparable. By [3, Theorem 1.4.5] any tournament has a Hamiltonian path, i.e., a spanning oriented path. Therefore we may assume that

$$v_1 \prec v_2 \prec \cdots \prec v_{s-1} \prec v_s$$

By the maximality of $c$ we get that $v_1$ is a source of $P_1$, $v_s$ is a sink of $P_1$, and $(v_i, v_{i+1})$ is an edge of $\mathcal{D}$ for $i = 1, \ldots, s-1$. Thus $c$ is an $A_1$–$B_1$ path of $\mathcal{D}$, as required. Conversely let $c = \{v_1, \ldots, v_s\}$ be an $A_1$–$B_1$ path of $\mathcal{D}$. Clearly $c$ is a clique of $P_1$ because $P_1$ is a poset. Assume that $c$ is not a maximal clique of $G_1$. Then there is a vertex $v \in X_1 \setminus c$ such that $v$ is related to every vertex of $c$. Since $v_1, v_s$ are a source and a sink of $P_1$ respectively we get $v_1 \prec v \prec v_s$. We claim that $v_i \prec v$ for $i = 1, \ldots, s$. By induction assume that $v_i \prec v$ for some $1 \leq i < s$. If $v \prec v_{i+1}$, then $v_i \prec v \prec v_{i+1}$, a contradiction to the fact that $(v_i, v_{i+1})$ is an edge of $\mathcal{D}$. Thus $v_{i+1} \prec v$. Making $i = s$ we get that $v_s \prec v$, a contradiction. This proves that $c$ is a maximal clique of $G_1$. Therefore, since the maximal cliques of $G_1$ not containing any vertex in $S$ are exactly the edges of $\mathcal{C}^w = \mathcal{C}_1 \setminus S$, by Menger's theorem (see Theorem 4.1.6) we obtain that $\mathcal{C}^w$ satisfies the König property. □

**Example 4.1.8** Consider the posets given by the transitive digraphs:

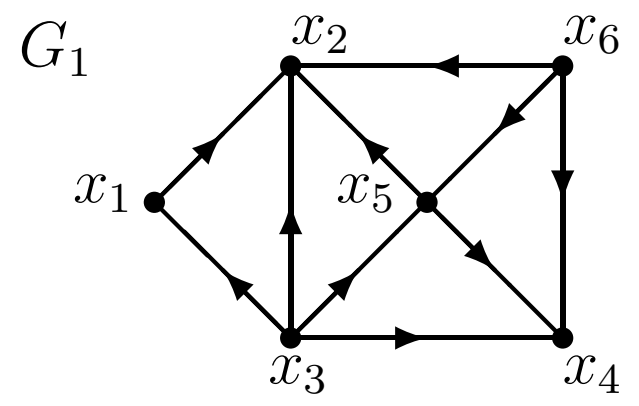

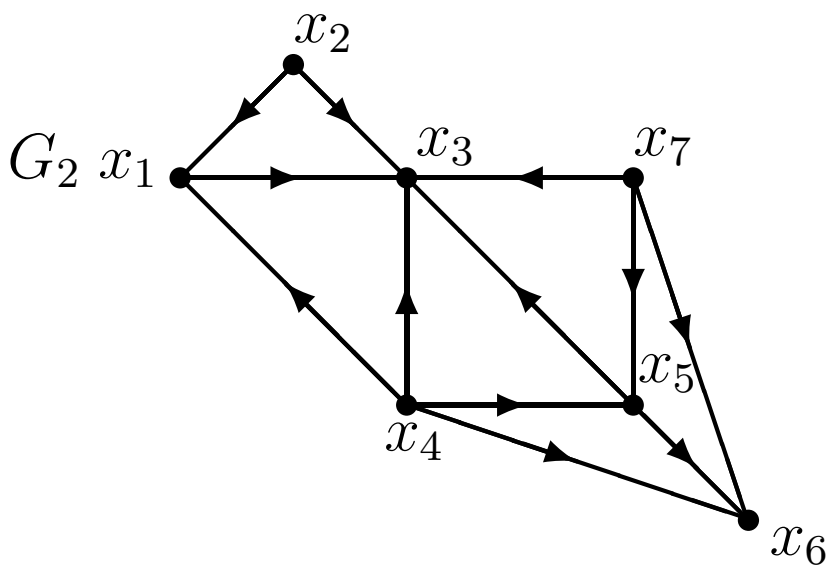

Then $I(\text{cl}(G_1)) = (x_1x_2x_3, x_2x_3x_5, x_2x_5x_6, x_4x_5x_6, x_3x_4x_5)$ and $I(\text{cl}(G_2)) = (x_1x_2x_3, x_1x_3x_4, x_3x_4x_5, x_3x_5x_7, x_4x_5x_6, x_5x_6x_7)$

It is well known that if $G$ is a comparability graph and $A$ is the vertex-clique matrix of $G$, then $G$ is perfect and the polytope $P(A) = \{x|\, x \geq 0;\, xA \leq \mathbf{1}\}$ is integral (see [16, 71]). The next result complement this fact.

**Corollary 4.1.9** *Let $G$ be a comparability graph and let $A$ be the vertex-clique matrix of $G$. Then the polyhedron $Q(A) = \{x|\, x \geq 0;\, xA \geq \mathbf{1}\}$ is integral.*

**Proof.** By Theorem 4.1.7 the clique clutter cl($G$) has the max-flow min-cut property. Thus the system $xA \geq \mathbf{1}$; $x \geq 0$ is totally dual integral, i.e., the maximum in Eq. (4.1) has an integral optimum solution $y$ for each integral vector $w$ with finite maximum. Hence $Q(A)$ has only integral vertices by [71, Theorem 5.22]. □

**Corollary 4.1.10** *Let $P = (X, \prec)$ be a poset on the vertex $X$ and let $G$ be its comparability graph. Then the maximum number of disjoint maximal cliques of $G$ is equal to the minimum size of a set intersecting all maximal cliques of $G$.*

**Proof.** By Theorem 4.1.7, the clique clutter of $G$ satisfies the max-flow min-cut property. Then by Theorem 4.1.3, the clique clutter of $G$ satisfies the König property. □

**Proposition 4.1.11** *Let $P = (X, \prec)$ be a poset on the vertex $X$. Then the maximum size of a chain of $P$ is equal to the minimum number of disjoint anti-chains into which $X$ can be decomposed.*

**Proof.** Let $L$ be a chain of maximum size and let $m = |L|$. Clearly $m$ is less or equal to the minimum number of disjoint anti-chains into which $X$ can be decomposed, because if $X_1, ..., X_k$ is a decomposition of $X$ into disjoint anti-chains, then $|L \cap X_i| \leq i$ for all $i$, so $|L| \leq k$. Let $x \in X$ and let $h(x)$ be the size of the longest chain of $P$ with end vertex $x$. Since any two vertices $x, y$ with $h(x) = h(y)$ are incomparable, we get a decomposition

$$X = \{x|h(x) = 1\} \cup \{x|h(x) = 2\} \cup \cdots \cup \{x|h(x) = m\}$$

into at most $m$ disjoint anti-chains, so min $\leq$ max. $\Box$

**Corollary 4.1.12** *If $G$ is the comparability graph of a poset $P = (X, \prec)$, then $G$ is perfect.*

**Proof.** Let $\omega(G)$ and $\chi(G)$ be the clique number and the chromatic number of $G$ respectively. By Proposition 4.1.11, we have $\omega(G) = \chi(G)$. As the class of comparability graphs is closed under taking induced subgraphs, we get that $G$ is a perfect graph. $\Box$

Recall that the complement of a perfect graph is also perfect (Perfect Graph Theorem, see [13]).

**Theorem 4.1.13** (*Dilworth decomposition theorem*) *Let $P = (X, \prec)$ be a poset on $X$. Then the maximum size of an anti-chain is equal to the minimum number of disjoint chains into which $X$ can be decomposed.*

**Proof.** Let $G$ be the comparability graph of $P$. By Corollary 4.1.12, $G$ is perfect. So, the complement $G'$ of $G$ is also perfect. In particular $\omega(G') = \chi(G')$ and the result follows readily. $\Box$

## 4.2 Normally torsion freeness and normality

Let $\mathcal{C}$ be a clutter and let $I = I(\mathcal{C}) \subset R$ be its edge ideal. If $C_1, \ldots, C_s$ are the minimal vertex covers of $\mathcal{C}$, then the primary decomposition of $I$ is

$$I = \mathfrak{p}_1 \cap \mathfrak{p}_2 \cap \cdots \cap \mathfrak{p}_s, \tag{4.2}$$

where $\mathfrak{p}_i$ is the ideal of $R$ generated by $C_i$. The *ith symbolic power* of $I$, denoted by $I^{(i)}$, is given by $I^{(i)} = \mathfrak{p}_1{}^i \cap \cdots \cap \mathfrak{p}_s^i$.

**Theorem 4.2.1** ([45]) *Let $\mathcal{C}$ be a clutter, let $A$ be the incidence matrix of $\mathcal{C}$, and let $I = I(\mathcal{C})$ be its edge ideal. The following are equivalent*

(i) *$I$ is normal and $Q(A) = \{x | \, x \geq 0;\, xA \geq \mathbf{1}\}$ is an integral polyhedron.*

(ii) *$I^i = I^{(i)}$ for $i \geq 1$.*

(iii) *$\mathcal{C}$ has the max-flow min-cut property.*

**Theorem 4.2.2** *If $G$ is a comparability graph and $\mathrm{cl}(G)$ is its clique clutter, then the edge ideal $I = I(\mathrm{cl}(G))$ of $\mathrm{cl}(G)$ is normally torsion free and normal.*

**Proof.** It follows from Theorems 4.1.7 and 4.2.1. □

**Complete admissible uniform clutters** In this paragraph we introduce a family of clique clutters of comparability graphs. Let $d \geq 2$, $g \geq 2$ be two integers and let

$$X^1 = \{x_1^1, \ldots, x_g^1\},\, X^2 = \{x_1^2, \ldots, x_g^2\},\, \ldots, X^d = \{x_1^d, \ldots, x_g^d\}$$

be disjoint sets of variables. The clutter $\mathcal{C}$ with vertex set $X = X^1 \cup \cdots \cup X^d$ and edge set

$$E(\mathcal{C}) = \{\{x_{i_1}^1, x_{i_2}^2, \ldots, x_{i_d}^d\} | \, 1 \leq i_1 \leq i_2 \leq \cdots \leq i_d \leq g\}$$

is called a *complete admissible uniform clutter*. The edge ideal of this clutter was introduced and studied in [34]. It was also studied in [49]. This ideal has many good properties, for instance $I(\mathcal{C})$ and its Alexander dual are Cohen-Macaulay and have linear resolutions (see [34, Proposition 4.5, Lemma 4.6]). For a thorough study of Cohen-Macaulay admissible clutters see [49, 66].

**Theorem 4.2.3** *If $\mathcal{C}$ is a complete admissible uniform clutter, then its edge ideal $I(\mathcal{C})$ is normally torsion free and normal.*

**Proof.** Let $(P, \prec)$ be the poset with vertex set $X$ and partial order given by $x_k^\ell \prec x_p^m$ if and only if $1 \leq \ell < m \leq d$ and $1 \leq k \leq p \leq g$. We denote the comparability graph of $P$ by $G$. We claim that $E(\mathcal{C}) = E(\mathrm{cl}(G))$, where $\mathrm{cl}(G)$ is the clique clutter of $G$. Let $f = \{x_{i_1}^1, x_{i_2}^2, \ldots, x_{i_d}^d\}$ be an edge of $\mathcal{C}$, i.e., we have $1 \leq i_1 \leq i_2 \leq \cdots \leq i_d \leq g$. Clearly $f$ is a clique of $G$. If $f$ is not maximal, then there is a vertex $x_k^\ell$ not in $f$ which is adjacent in $G$ to every vertex of $f$. In particular $x_k^\ell$ must be comparable to $x_{i_\ell}^\ell$, which is impossible.

Thus $f$ is an edge of $\mathrm{cl}(G)$. Conversely let $f$ be an edge of $\mathrm{cl}(G)$. We can write $f = \{x_{i_1}^{k_1}, x_{i_2}^{k_2}, \ldots, x_{i_s}^{k_s}\}$, where $k_1 < \cdots < k_s$ and $i_1 \leq \cdots \leq i_s$. By the maximality of $f$ we get that $s = d$ and $k_i = i$ for $i = 1, \ldots, d$. Thus $f$ is an edge of $\mathcal{C}$. Hence by Theorem 4.2.2 we obtain that $I(\mathcal{C})$ is normally torsion free and normal. □

**Acknowledgments.** We thank Seth Sullivant for pointing out an alternative proof of Theorem 4.2.2 based on the fact that the edge ideal of a comparability graph is differentially perfect (see [78, Section 4]). One of the consequences of being differentially perfect is that if every maximal chain in the poset has the same fixed length $k$, then the edge ideal of the clique clutter of its comparability graph is normally torsion free. The general case, that is, the situation where not all maximal chains have the same length, can be reduced to the special case above.

# Chapter 5

# Symbolic Rees Algebras, Vertex Covers and Irreducible Representations of Rees Cones

The relation between facets of Rees cones of clutters and irreducible $b$-vertex covers is examined. Let $G$ be a simple graph and let $I_c(G)$ be its ideal of vertex covers. We give a graph theoretical description of the irreducible $b$-vertex covers of $G$, i.e., we describe the minimal generators of the symbolic Rees algebra of $I_c(G)$. As an application we recover an explicit description of the edge cone of a graph. Then we study the irreducible $b$-vertex covers of the blocker of $G$, i.e., we study the minimal generators of the Symbolic Rees algebra of the edge ideal of $G$. We give a graph theoretical description of the irreducible binary $b$-vertex covers of the blocker of $G$. It is shown that they are in one to one correspondence with the irreducible induced subgraphs of $G$. As a byproduct we obtain a method, using Hilbert bases, to obtain all irreducible induced subgraphs of $G$. Irreducible graphs are studied. We show how to build irreducible graphs and give a method to construct irreducible $b$-vertex covers of the blocker of $G$ with high degree relative to the number of vertices of $G$.

A *clutter* $\mathcal{C}$ with vertex set $X = \{x_1, \ldots, x_n\}$ is a family of subsets of $X$, called edges, none of which is included in another. The set of vertices and edges of $\mathcal{C}$ are denoted by $V(\mathcal{C})$ and $E(\mathcal{C})$ respectively. A basic example of a clutter is a graph. Let $R = K[x_1, \ldots, x_n]$ be a polynomial ring over a field $K$. The *edge ideal* of $\mathcal{C}$, denoted by $I(\mathcal{C})$, is the ideal of $R$ generated by all monomials $\prod_{x_i \in e} x_i$ such that $e \in E(\mathcal{C})$. The assignment $\mathcal{C} \mapsto I(\mathcal{C})$ establishes a natural one to one correspondence between the family of clutters and the family of square-free monomial ideals. Let $\mathcal{C}$ be a clutter and let $F = \{x^{v_1}, \ldots, x^{v_q}\}$ be the minimal set of generators of its edge ideal $I = I(\mathcal{C})$. As usual we use $x^a$ as an abbreviation for $x_1^{a_1} \cdots x_n^{a_n}$, where $a = (a_1, \ldots, a_n) \in \mathbb{N}^n$. The $n \times q$ matrix with column vectors $v_1, \ldots, v_q$ will be denoted by $A$, it is called the *incidence matrix* of $\mathcal{C}$.

The *blowup algebra* studied here is the *symbolic Rees algebra*:

$$R_s(I) = R \oplus I^{(1)}t \oplus \cdots \oplus I^{(i)}t^i \oplus \cdots \subset R[t],$$

where $t$ is a new variable and $I^{(i)}$ is the $ith$ symbolic power of $I$. Closely related to $R_s(I)$ is the *Rees algebra* of $I$:

$$R[It] := R \oplus It \oplus \cdots \oplus I^i t^i \oplus \cdots \subset R[t].$$

The study of symbolic powers of edge ideals was initiated in [73] and further elaborated on in [1, 29, 42, 45, 51, 78, 91]. By a result of Lyubeznik [63], $R_s(I)$ is a $K$-algebra of finite type. In general the minimal set of generators of $R_s(I)$ as a $K$-algebra is very hard to describe in terms of $\mathcal{C}$ (see [1]). There are two exceptional cases. If the clutter $\mathcal{C}$ has the max-flow min-cut property, then by a result of [45] we have $I^i = I^{(i)}$ for all $i \geq 1$, i.e., $R_s(I) = R[It]$. If $G$ is a perfect graph, then the minimal generators of $R_s(I(G))$ are in one to one correspondence with the cliques (complete subgraphs) of $G$ [91], see also [78]. We shall be interested in studying the minimal set of generators of $R_s(I)$ using polyhedral geometry. Let $G$ be a graph and let $I_c(G)$ be the Alexander dual of $I(G)$, see definition below. Some of the main results of this chapter are graph theoretical descriptions of the minimal generators of $R_s(I(G))$ and $R_s(I_c(G))$. In Sections 5.1 and 5.2 we show that both algebras encode combinatorial information of the graph which can be decoded using integral Hilbert bases.

The *Rees cone* of $I$, denoted by $\mathbb{R}_+(I)$, is the polyhedral cone consisting of the non-negative linear combinations of the set

$$\mathcal{A}' = \{e_1, \ldots, e_n, (v_1, 1), \ldots, (v_q, 1)\} \subset \mathbb{R}^{n+1},$$

where $e_i$ is the $ith$ unit vector.

A subset $C \subset X$ is a *minimal vertex cover* of the clutter $\mathcal{C}$ if: (i) every edge of $\mathcal{C}$ contains at least one vertex of $C$, and (ii) there is no proper subset of $C$ with the first property. If $C$ satisfies condition (i) only, then $C$ is called a *vertex cover* of $\mathcal{C}$. Let $\mathfrak{p}_1,\ldots,\mathfrak{p}_s$ be the minimal primes of the edge ideal $I=I(\mathcal{C})$ and let

$$C_k=\{x_i|\, x_i\in\mathfrak{p}_k\}\qquad (k=1,\ldots,s)$$

be the corresponding minimal vertex covers of $\mathcal{C}$, see [87, Proposition 6.1.16]. Recall that the primary decomposition of the edge ideal of $\mathcal{C}$ is given by

$$I(\mathcal{C})=(C_1)\cap(C_2)\cap\cdots\cap(C_s).$$

In particular observe that the height of $I(\mathcal{C})$ equals the number of vertices in a minimum vertex cover of $\mathcal{C}$. The *ith* symbolic power of $I$ is given by

$$I^{(i)}=S^{-1}I^i\cap R\ \text{ for } i\geq 1,$$

where $S=R\setminus\cup_{k=1}^{s}\mathfrak{p}_i$ and $S^{-1}I^i$ is the localization of $I^i$ at $S$. In our situation the *ith* symbolic power of $I$ has a simple expression: $I^{(i)}=\mathfrak{p}_1^i\cap\cdots\cap\mathfrak{p}_s^i$, see [87]. The Rees cone of $I$ is a finitely generated rational cone of dimension $n+1$. Hence by the finite basis theorem [93, Theorem 4.11] there is a unique irreducible representation

$$\mathbb{R}_+(I)=H_{e_1}^+\cap H_{e_2}^+\cap\cdots\cap H_{e_{n+1}}^+\cap H_{\ell_1}^+\cap H_{\ell_2}^+\cap\cdots\cap H_{\ell_r}^+ \tag{5.1}$$

such that each $\ell_k$ is in $\mathbb{Z}^{n+1}$, the non-zero entries of each $\ell_k$ are relatively prime, and none of the closed halfspaces $H_{e_1}^+,\ldots,H_{e_{n+1}}^+,H_{\ell_1}^+,\ldots,H_{\ell_r}^+$ can be omitted from the intersection. Here $H_a^+$ denotes the closed halfspace $H_a^+=\{x|\,\langle x,a\rangle\geq 0\}$ and $H_a$ stands for the hyperplane through the origin with normal vector $a$, where $\langle\ ,\ \rangle$ denotes the standard inner product. The *facets* (i.e., the proper faces of maximum dimension or equivalently the faces of dimension $n$) of the Rees cone are exactly:

$$F_1=H_{e_1}\cap\mathbb{R}_+(I),\ldots,F_{n+1}=H_{e_{n+1}}\cap\mathbb{R}_+(I),H_{\ell_1}\cap\mathbb{R}_+(I),\ldots,H_{\ell_r}\cap\mathbb{R}_+(I).$$

According to [29, Lemma 3.1] we may always assume that $\ell_k=-e_{n+1}+\sum_{x_i\in C_k}e_i$ for $1\leq k\leq s$, i.e., each minimal vertex cover of $\mathcal{C}$ determines a facet of the Rees cone and every facet of the Rees cone satisfying $\langle\ell_k,e_{n+1}\rangle=-1$ must be of the form $\ell_k=-e_{n+1}+\sum_{x_i\in C_k}e_i$ for some minimal vertex cover $C_k$ of $\mathcal{C}$. This is quite interesting because this is saying that the Rees cone of $I(\mathcal{C})$ is a carrier of combinatorial information of the clutter $\mathcal{C}$. Thus we can extract the primary decomposition of $I(\mathcal{C})$ from the irreducible representation of $\mathbb{R}_+(I(\mathcal{C}))$.

Rees cones have been used to study algebraic and combinatorial properties of blowup algebras of square-free monomial ideals and clutters [29, 42]. Blowup algebras are interesting objects of study in algebra and geometry [82].

The ideal of *vertex covers* of $\mathcal{C}$ is the square-free monomial ideal

$$I_c(\mathcal{C}) = (x^{u_1}, \ldots, x^{u_s}) \subset R,$$

where $x^{u_k} = \prod_{x_i \in C_k} x_i$. Often the ideal $I_c(\mathcal{C})$ is called the *Alexander dual* of $I(\mathcal{C})$. The clutter $\Upsilon(\mathcal{C})$ associated to $I_c(\mathcal{C})$ is called the *blocker* of $\mathcal{C}$, see [16]. Notice that the edges of $\Upsilon(\mathcal{C})$ are precisely the minimal vertex covers of $\mathcal{C}$. If $G$ is a graph, then $R_s(I_c(G))$ is generated as a $K$-algebra by elements of degree in $t$ at most two [51, Theorem 5.1]. One of the main result of Section 5.1 is a graph theoretical description of the minimal generators of $R_s(I_c(G))$ (see Theorem 5.1.9). As an application we recover an explicit description given in [80] of the edge cone of a graph (Corollary 5.1.10).

The symbolic Rees algebra of the ideal $I_c(\mathcal{C})$ can be interpreted in terms of "$k$-vertex covers" [51] as we now explain. Let $a = (a_1, \ldots, a_n) \neq 0$ be a vector in $\mathbb{N}^n$ and let $b \in \mathbb{N}$. We say that $a$ is a *$b$-vertex cover* of $I$ (or $\mathcal{C}$) if $\langle v_i, a\rangle \geq b$ for $i = 1, \ldots, q$. Often we will call a $b$-vertex cover simply a *$b$-cover*. This notion plays a role in combinatorial optimization [71, Chapter 77, p. 1378] and algebraic combinatorics [51, 52].

The *algebra of covers* of $I$ (or $\mathcal{C}$), denoted by $R_c(I)$, is the $K$-subalgebra of $K[t]$ generated by all monomials $x^a t^b$ such that $a$ is a $b$-cover of $I$. We say that a $b$-cover $a$ of $I$ is *reducible* if there exists an $i$-cover $c$ and a $j$-cover $d$ of $I$ such that $a = c + d$ and $b = i + j$. If $a$ is not reducible, we call *$a$ irreducible*. The irreducible 0 and 1 covers of $\mathcal{C}$ are the unit vector $e_1, \ldots, e_n$ and the incidence vectors $u_1, \ldots, u_s$ of the minimal vertex covers of $\mathcal{C}$, respectively. The minimal generators of $R_c(I)$ as a $K$-algebra correspond to the irreducible covers of $I$. Notice the following dual descriptions:

$$\begin{aligned} I^{(b)} &= (\{x^a \mid \langle (a,b), \ell_i\rangle \geq 0 \text{ for } i = 1, \ldots, s\}) \\ &= (\{x^a \mid \langle a, u_i\rangle \geq b \text{ for } i = 1, \ldots, s\}), \\ J^{(b)} &= (\{x^a \mid \langle a, v_i\rangle \geq b \text{ for } i = 1, \ldots, q\}), \end{aligned}$$

where $J = I_c(\mathcal{C})$. Hence $R_c(I) = R_s(J)$ and $R_c(J) = R_s(I)$.

In general each $\ell_i$ occurring in Eq. (5.1) determines a minimal generator of $R_s(I_c(\mathcal{C}))$. Indeed if we write $\ell_i = (a_i, -d_i)$, where $a_i \in \mathbb{N}^n$, $d_i \in \mathbb{N}$, then $a_i$ is an irreducible $d_i$-cover of $I$ (Lemma 5.1.3). Let $F_{n+1}$ be the facet of $\mathbb{R}_+(I)$ determined by the hyperplane $H_{e_{n+1}}$. Thus we have a map $\psi$:

$$\begin{array}{rcl} \{\text{Facets of } \mathbb{R}_+(I(\mathcal{C}))\} \setminus \{F_{n+1}\} & \xrightarrow{\psi} & R_s(I_c(\mathcal{C})) \\ H_{\ell_k} \cap \mathbb{R}_+(I) & \xrightarrow{\psi} & x^{a_k} t^{d_k}, \text{ where } \ell_k = (a_k, -d_k) \\ H_{e_i} \cap \mathbb{R}_+(I) & \xrightarrow{\psi} & x_i \end{array}$$

whose image provides a good approximation for the minimal set of generators of $R_s(I_c(\mathcal{C}))$ as a $K$-algebra. Likewise the facets of $\mathbb{R}_+(I_c(\mathcal{C}))$ give an approximation for the minimal set of generators of $R_s(I(\mathcal{C}))$. In Example 5.1.6 we show a connected graph $G$ for which the image of the map $\psi$ does not generate $R_s(I_c(G))$. For balanced clutters, i.e., for clutters without odd cycles, the image of the map $\psi$ generates $R_s(I_c(\mathcal{C}))$. This follows from [42, Propositions 4.10 and 4.11]. In particular the image of the map $\psi$ generate $R_s(I_c(\mathcal{C}))$ when $\mathcal{C}$ is a bipartite graph. It would be interesting to characterize when the irreducible representation of the Rees cone determine the irreducible covers.

The *Simis cone* of $I$ is the rational polyhedral cone:

$$\text{Cn}(I) = H_{e_1}^+ \cap \cdots \cap H_{e_{n+1}}^+ \cap H_{(u_1,-1)}^+ \cap \cdots \cap H_{(u_s,-1)}^+,$$

Simis cones were introduced in [29] to study symbolic Rees algebras of square-free monomial ideals. If $\mathcal{H}$ is an integral Hilbert basis of $\text{Cn}(I)$, then $R_s(I(\mathcal{C}))$ equals $K[\mathbb{N}\mathcal{H}]$, the semigroup ring of $\mathbb{N}\mathcal{H}$ (see [29, Theorem 3.5]). This result is interesting because it allows us to compute the minimal generators of $R_s(I(\mathcal{C}))$ using Hilbert bases. The program *Normaliz* [11] is suitable for computing Hilbert bases. There is a description of $\mathcal{H}$ valid for perfect graphs [91].

If $G$ is a perfect graph, the irreducible $b$-covers of $\Upsilon(G)$ correspond to cliques of $G$ [91] (cf. Corollary 5.2.5). In this case, setting $\mathcal{C} = \Upsilon(G)$, it turns out that the image of $\psi$ generates $R_s(I_c(\Upsilon(G)))$. Notice that $I_c(\Upsilon(G))$ is equal to $I(G)$.

In Section 5.2 we introduce and study the concept of an irreducible graph. A $b$-cover $a = (a_1, \ldots, a_n)$ is called *binary* if $a_i \in \{0, 1\}$ for all $i$. We present a graph theoretical description of the irreducible binary $b$-vertex covers of the blocker of $G$ (see Theorem 5.2.7). It is shown that they are in one to one correspondence with the irreducible induced subgraphs of $G$. As a byproduct we obtain a method, using Hilbert bases, to obtain all irreducible induced subgraphs of $G$ (see Corollary 5.2.10). We give a simple procedure to build irreducible graphs (Proposition 5.2.17) and give a method to construct irreducible $b$-vertex covers of the blocker of $G$ with high degree relative to the number of vertices of $G$ (see Corollaries 5.2.22 and 5.2.23).

Along the chapter we introduce most of the notions that are relevant for our purposes. For unexplained terminology and notation we refer to [20, 60] and [65, 82].

## 5.1 Blowup algebras of ideals of vertex covers

Let $G$ be a simple graph with vertex set $X = \{x_1, \dots, x_n\}$. In what follows we shall always assume that $G$ has no isolated vertices. Here we will give a graph theoretical description of the irreducible $b$-covers of $G$, i.e., we will describe the symbolic Rees algebra of $I_c(G)$.

Let $S$ be a set of vertices of $G$. The *neighbor set* of $S$, denoted by $N_G(S)$, is the set of vertices of $G$ that are adjacent with at least one vertex of $S$. The set $S$ is called *independent* if no two vertices of $S$ are adjacent. The empty set is regarded as an independent set whose incidence vector is the zero vector. Notice the following duality: $S$ is a maximal independent set of $G$ (with respect to inclusion) if and only if $X \setminus S$ is a minimal vertex cover of $G$.

**Lemma 5.1.1** *If* $a = (a_i) \in \mathbb{N}^n$ *is an irreducible* $k$*-cover of* $G$*, then* $0 \leq k \leq 2$ *and* $0 \leq a_i \leq 2$ *for* $i = 1, \dots, n$.

**Proof.** Recall that $a$ is a $k$-cover of $G$ if and only if $a_i + a_j \geq k$ for each edge $\{x_i, x_j\}$ of $G$. If $k = 0$ or $k = 1$, then by the irreducibility of $a$ it is seen that either $a = e_i$ for some $i$ or $a = e_{i_1} + \cdots + e_{i_r}$ for some minimal vertex cover $\{x_{i_1}, \dots, x_{i_r}\}$ of $G$. Thus we may assume that $k \geq 2$.

Case (I): $a_i \geq 1$ for all $i$. Clearly $\mathbf{1} = (1, \dots, 1)$ is a 2-cover. If $a - \mathbf{1} \neq 0$, then $a - \mathbf{1}$ is a $k - 2$ cover and $a = \mathbf{1} + (a - \mathbf{1})$, a contradiction to $a$ being an irreducible $k$-cover. Hence $a = \mathbf{1}$. Pick any edge $\{x_i, x_j\}$ of $G$. Since $a$ is a $k$-cover, we get $2 = a_i + a_j \geq k$ and $k$ must be equal to 2.

Case (II): $a_i = 0$ for some $i$. We may assume $a_i = 0$ for $1 \leq i \leq r$ and $a_i \geq 1$ for $i > r$. Notice that the set $S = \{x_1, \dots, x_r\}$ is independent because if $\{x_i, x_j\}$ is an edge and $1 \leq i < j \leq r$, then $0 = a_i + a_j \geq k$, a contradiction. Consider the neighbor set $N_G(S)$ of $S$. We may assume that $N_G(S) = \{x_{r+1}, \dots, x_s\}$. Observe that $a_i \geq k \geq 2$ for $i = r + 1, \dots, s$, because $a$ is a $k$-cover. Write

$$a = (0, \dots, 0, a_{r+1} - 2, \dots, a_s - 2, a_{s+1} - 1, \dots, a_n - 1) + (\underbrace{0, \dots, 0}_{r}, \underbrace{2, \dots, 2}_{s-r}, \underbrace{1, \dots, 1}_{n-s}) = c + d.$$

Clearly $d$ is a 2-cover. If $c \neq 0$, using that $a_i \geq k \geq 2$ for $r+1 \leq i \leq s$ and $a_i \geq 1$ for $i > s$ it is not hard to see that $c$ is a $(k-2)$-cover. This gives a contradiction, because $a = c + d$. Hence $c = 0$. Therefore $a_i = 2$ for $r < i \leq s$, $a_i = 1$ for $i > s$, and $k = 2$. □

The next result complements the fact that the symbolic Rees algebra of $I_c(G)$ is generated by monomials of degree in $t$ at most two [51, Theorem 5.1].

**Corollary 5.1.2** *$R_s(I_c(G))$ is generated as a $K$-algebra by monomials of degree in $t$ at most two and total degree at most $2n$.*

**Proof.** Let $x^a t^k$ be a minimal generator of $R_s(I_c(G))$ as a $K$-algebra. Then $a = (a_1, \ldots, a_n)$ is an irreducible $k$-cover of $G$. By Lemma 5.1.1 we obtain that $0 \leq k \leq 2$ and $0 \leq a_i \leq 2$ for all $i$. If $k = 0$ or $k = 1$, we get that the degree of $x^a t^k$ is at most $n$ because when $k = 0$ or 1 one has $a = e_i$ or $a = \sum_{x_i \in C_k} e_i$ for some minimal vertex cover $C_k$ of $G$, respectively. If $k = 2$, by the proof of Lemma 5.1.1 either $a = \mathbf{1}$ or $a_i = 0$ for some $i$. Thus $\deg(x^a) \leq 2(n-1)$. □

Let $I = I(G)$ be the edge ideal of $G$. For use below consider the vectors $\ell_1, \ldots, \ell_r$ that occur in the irreducible representation of $\mathbb{R}_+(I)$ given in Eq. (5.1).

**Lemma 5.1.3** *Let $\mathcal{C}$ be a clutter and let $I = I(\mathcal{C})$ be its edge ideal. If $\ell_k = (a_k, -d_k)$ is any of the vectors that occur in Eq.* (5.1)*, where $a_k \in \mathbb{N}^n$, $d_k \in \mathbb{N}$, then $a_k$ is an irreducible $d_k$-cover of $\mathcal{C}$.*

**Proof.** We proceed by contradiction assume there is a $d_k'$-cover $a_k'$ and a $d_k''$-cover $a_k''$ such that $a_k = a_k' + a_k''$ and $d_k = d_k' + d_k''$. Set $F' = H_{(a_k', -d_k')} \cap \mathbb{R}_+(I)$ and $F'' = H_{(a_k'', -d_k'')} \cap \mathbb{R}_+(I)$. Clearly $F', F''$ are proper faces of $\mathbb{R}_+(I)$ and $F = \mathbb{R}_+(I) \cap H_{\ell_k} = F' \cap F''$. Recall that any proper face of $\mathbb{R}_+(I)$ is the intersection of those facets that contain it (see [93, Theorem 3.2.1(vii)]). Applying this fact to $F'$ and $F''$ it is seen that $F' \subset F$ or $F'' \subset F$, i.e., $F = F'$ or $F = F''$. We may assume $F = F'$. Hence $H_{(a_k', -d_k')} = H_{\ell_k}$. Taking orthogonal complements we get that $(a_k', -d_k') = \lambda(a_k, -d_k)$ for some $\lambda \in \mathbb{Q}_+$, because the orthogonal complement of $H_{\ell_k}$ is one dimensional and it is generated by $\ell_k$. Since the non-zero entries of $\ell_k$ are relatively prime, we may assume that $\lambda \in \mathbb{N}$. Thus $d_k' = \lambda d_k \geq d_k \geq d_k'$ and $\lambda$ must be equal to 1. Hence $a_k = a_k'$ and $a_k''$ must be zero, a contradiction. □

**Corollary 5.1.4** *If $\ell_i = (\ell_{i1}, \ldots, \ell_{in}, -\ell_{i(n+1)})$, then $0 \leq \ell_{ij} \leq 2$ for $j = 1, \ldots, n$ and $1 \leq \ell_{i(n+1)} \leq 2$.*

**Proof.** It suffices to observe that $(\ell_{i1},\dots,\ell_{in})$ is an irreducible $\ell_{i(n+1)}$-cover of $G$ and to apply Lemma 5.1.1 and Lemma 5.1.3 $\Box$

**Lemma 5.1.5** $a=(1,\dots,1)$ *is an irreducible* 2*-cover of* $G$ *if and only if* $G$ *is non bipartite.*

**Proof.** $\Rightarrow$) We proceed by contradiction assuming that $G$ is a bipartite graph. Then $G$ has a bipartition $(V_1,V_2)$. Set $a'=\sum_{x_i\in V_1}e_i$ and $a''=\sum_{x_i\in V_2}e_i$. Since $V_1$ and $V_2$ are minimal vertex covers of $G$, we can decompose $a$ as $a=a'+a''$, where $a'$ and $a''$ are 1-covers, which is impossible.

$\Leftarrow$) Notice that $a$ cannot be the sum of a 0-cover and a 2-cover. Indeed if $a=a'+a''$, where $a'$ is a 0-cover and $a''$ is a 1-cover, then $a''$ has an entry $a_i$ equal to zero. Pick an edge $\{x_i,x_j\}$ incident with $x_i$, then $\langle a'',e_i+e_j\rangle\le 1$, a contradiction. Thus we may assume that $a=c+d$, where $c,d$ are 1-covers. Let $C_r$ be an odd cycle of $G$ of length $r$. Notice that any vertex cover of $C_r$ must contain a pair of adjacent vertices because $r$ is odd. Clearly a vertex cover of $G$ is also a vertex cover of the subgraph $C_r$. Hence the vertex covers of $G$ corresponding to $c$ and $d$ must contain a pair of adjacent vertices, a contradiction because $c$ and $d$ are complementary vectors and the complement of a vertex cover is an independent set. $\Box$

**Example 5.1.6** Consider the following graph $G$:

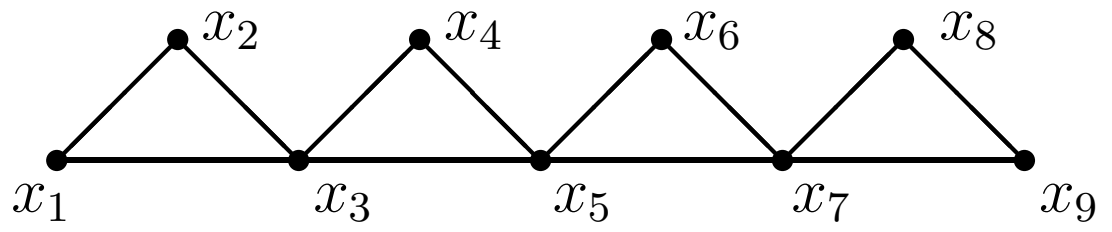

Using *Normaliz* [11] it is seen that the vector $a=(1,1,2,0,2,1,1,1,1)$ is an irreducible 2-cover of $G$ such that the supporting hyperplane $H_{(a,-2)}$ does not define a facet of the Rees cone of $I(G)$. Thus, in general, the image of $\psi$ described at the beginning of the chapter does not determine $R_s(I_c(G))$. We may use Lemma 5.1.5 to construct non-connected graphs with this property.

**Definition 5.1.7** Let $A$ be the incidence matrix of a clutter $\mathcal{C}$. A clutter $\mathcal{C}$ has a cycle of length $r$ if there is a square sub-matrix of $A$ of size $r\ge 3$ with exactly two 1's in each row and column. A clutter without odd cycles is called *balanced.*

**Proposition 5.1.8** ([42, Proposition 4.11]) *If* $\mathcal{C}$ *is a balanced clutter, then*

$$R_s(I_c(\mathcal{C}))=R[I_c(\mathcal{C})t].$$

This result was first shown for bipartite graphs in [40, Corollary 2.6] and later generalized to balanced clutters [42] using an algebro-combinatorial description of clutters with the max-flow min-cut property [45].

Let $S$ be a set of vertices of a graph $G$, the *induced subgraph* $\langle S\rangle$ is the maximal subgraph of $G$ with vertex set $S$. The next result has used in [35] to show that any associated prime of $I_c(G)^2$ is generated by the vertices of an odd hole of $G$.

**Theorem 5.1.9** *Let $0\neq a=(a_i)\in\mathbb{N}^n$ and let $\Upsilon(G)$ be the family of minimal vertex covers of $G$.*

(i) *If $G$ is bipartite, then $a$ is an irreducible $b$-cover of $G$ if and only if $b=0$ and $a=e_i$ for some $1\leq i\leq n$ or $b=1$ and $a=\sum_{x_i\in C}e_i$ for some $C\in\Upsilon(G)$.*

(ii) *If $G$ is non-bipartite, then $a$ is an irreducible $b$-cover if and only if $a$ has one of the following forms*:

(a) (0-*covers*) $b=0$ *and* $a=e_i$ *for some* $1\leq i\leq n$,

(b) (1-*covers*) $b=1$ *and* $a=\sum_{x_i\in C}e_i$ *for some* $C\in\Upsilon(G)$,

(c) (2-*covers*) $b=2$ *and* $a=(1,\ldots,1)$,

(d) (2-*covers*) $b=2$ *and up to permutation of vertices*

$$a=(\underbrace{0,\ldots,0}_{|A|},\underbrace{2,\ldots,2}_{|N_G(A)|},1,\ldots,1)$$

*for some independent set of vertices $A\neq\emptyset$ of $G$ such that*

($\mathrm{d}_1$) *$N_G(A)$ is not a vertex cover of $G$ and $V\neq A\cup N_G(A)$,*

($\mathrm{d}_2$) *the induced subgraph $\langle V\setminus(A\cup N_G(A))\rangle$ has no isolated vertices and is not bipartite.*

**Proof.** (i) $\Rightarrow$) Since $G$ is a bipartite graph, by Proposition 5.1.8, we have the equality $R_s(J)=R[Jt]$, where $J=I_c(G)$ is the ideal of vertex covers of $G$. Thus the minimal set of generator of $R_s(J)$ as a $K$-algebra is the set

$$\{x_1,\ldots,x_n,x^{u_1}t,\ldots,x^{u_r}t\},$$

where $u_1,\ldots,u_r$ are the incidence vectors of the minimal vertex covers of $G$. By hypothesis $a$ is an irreducible $b$-cover of $G$, i.e., $x^at^b$ is a minimal generator of $R_s(I_c(\mathcal{C}))$. Therefore either $a = e_i$ for some $i$ and $b = 0$ or $a = u_i$ for some $i$ and $b = 1$. The converse follows readily and is valid for any graph or clutter.

(ii) $\Rightarrow$) By Lemma 5.1.1 we have $0 \leq b \leq 2$ and $0 \leq a_i \leq 2$ for all $i$. If $b = 0$ or $b = 1$, then clearly $a$ has the form indicated in (a) or (b) respectively.

Assume $b = 2$. If $a_i \geq 1$ for all $i$, then $a_i = 1$ for all $i$, otherwise if $a_i = 2$ for some $i$, then $a - e_i$ is a 2-cover and $a = e_i + (a - e_i)$, a contradiction. Hence $a = \mathbf{1}$. Thus we may assume that $a$ has the form

$$a = (0, \ldots, 0, 2, \ldots, 2, 1, \ldots, 1).$$

We set $A = \{x_i |\, a_i = 0\} \neq \emptyset$, $B = \{x_i |\, a_i = 2\}$, and $C = V \setminus (A \cup B)$. Observe that $A$ is an independent set because $a$ is a 2-cover and $B = N_G(A)$ because $a$ is irreducible. Hence it is seen that conditions ($\mathrm{d}_1$) and ($\mathrm{d}_2$) are satisfied. Using Lemma 5.1.5, the proof of the converse is direct. □

**Edge cones of graphs** Let $G$ be a connected graph and let $\mathcal{A} = \{v_1, \ldots, v_q\}$ be the set of all vectors $e_i + e_j$ such that $\{x_i, x_j\}$ is an edge of $G$. The *edge cone* of $G$, denoted by $\mathbb{R}_+\mathcal{A}$, is defined as the cone generated by $\mathcal{A}$. Below we give an explicit combinatorial description of the edge cone.

Let $A$ be an *independent set* of vertices of $G$. The supporting hyperplane of the edge cone of $G$ defined by

$$\sum_{x_i \in N_G(A)} x_i - \sum_{x_i \in A} x_i = 0$$

will be denoted by $H_A$.

Edge cones and their representations by closed halfspaces are a useful tool to study $a$-invariants of edge subrings [86, 79]. The following result is a prototype of these representations. As an application we give a direct proof of the next result using Rees cones.

**Corollary 5.1.10** ([80, Corollary 2.8]) *A vector* $a = (a_1, \ldots, a_n) \in \mathbb{R}^n$ *is in* $\mathbb{R}_+\mathcal{A}$ *if and only if* $a$ *satisfies the following system of linear inequalities*

$$\begin{array}{rcl} a_i & \geq & 0, \quad i = 1, \ldots, n; \\ \sum_{x_i \in N_G(A)} a_i - \sum_{x_i \in A} a_i & \geq & 0, \quad \textit{for all independent sets } A \subset V(G). \end{array}$$

**Proof.** Set $\mathcal{B} = \{(v_1, 1), \ldots, (v_q, 1)\}$ and $I = I(G)$. Notice the equality

$$\mathbb{R}_+(I) \cap \mathbb{R}\mathcal{B} = \mathbb{R}_+\mathcal{B}, \tag{5.2}$$

where $\mathbb{R}\mathcal{B}$ is $\mathbb{R}$-vector space spanned by $\mathcal{B}$. Consider the irreducible representation of $\mathbb{R}_+(I)$ given in Eq. (5.1) and write $\ell_i = (a_i, -d_i)$, where $0 \neq a_i \in \mathbb{N}^n$, $0 \neq d_i \in \mathbb{N}$. Next we show the equality:

$$\mathbb{R}_+\mathcal{A} = \mathbb{R}\mathcal{A} \cap \mathbb{R}^n_+ \cap H^+_{(2a_1/d_1 - \mathbf{1})} \cap \cdots \cap H^+_{(2a_r/d_r - \mathbf{1})}, \tag{5.3}$$

where $\mathbf{1} = (1, \ldots, 1)$. Take $\alpha \in \mathbb{R}_+\mathcal{A}$. Clearly $\alpha \in \mathbb{R}\mathcal{A} \cap \mathbb{R}^n_+$. We can write

$$\alpha = \lambda_1 v_1 + \cdots + \lambda_q v_q \;\Rightarrow\; |\alpha| = 2(\lambda_1 + \cdots + \lambda_q) = 2b.$$

Thus $(\alpha, b) = \lambda_1(v_1, 1) + \cdots + \lambda_q(v_q, 1)$, i.e., $(\alpha, b) \in \mathbb{R}_+\mathcal{B}$. Hence from Eq. (5.2) we get $(\alpha, b) \in \mathbb{R}_+(I)$ and

$$\langle (\alpha, b), (a_i, -d_i) \rangle \geq 0 \;\Rightarrow\; \langle \alpha, a_i \rangle \geq b d_i = (|\alpha|/2) d_i = |\alpha|(d_i/2).$$

Writing $\alpha = (\alpha_1, \ldots, \alpha_n)$ and $a_i = (a_{i1}, \ldots, a_{in})$, the last inequality gives:

$$\alpha_1 a_{i_1} + \cdots + \alpha_n a_{in} \geq (\alpha_1 + \cdots + \alpha_n)(d_i/2) \;\Rightarrow\; \langle \alpha, a_i - (d_i/2)\mathbf{1} \rangle \geq 0.$$

Then $\langle \alpha, 2a_i/d_i - \mathbf{1} \rangle \geq 0$ and $\alpha \in H^+_{(2a_i/d_i - \mathbf{1})}$ for all $i$, as required. This proves that $\mathbb{R}_+\mathcal{A}$ is contained in the right hand side of Eq. (5.3). The other inclusion follows similarly. Now by Lemma 5.1.3 we obtain that $a_i$ is an irreducible $d_i$-cover of $G$. Therefore, using the explicit description of the irreducible $b$-covers of $G$ given in Theorem 5.1.9, we get the equality

$$\mathbb{R}_+\mathcal{A} = \left( \bigcap_{A \in \mathcal{F}} H^+_A \right) \bigcap \left( \bigcap_{i=1}^{n} H^+_{e_i} \right),$$

where $\mathcal{F}$ is the collection of all the independent sets of vertices of $G$. From this equality the assertion follows at once. □

The edge cone of $G$ encodes information about both the Hilbert function of the edge subring $K[G]$ (see [79]) and the graph $G$ itself. As a simple illustration, we recover the following version of the marriage theorem for bipartite graphs, see [6]. Recall that a pairing by an independent set of edges of all the vertices of a graph $G$ is called a *perfect matching* or a 1-*factor*.

**Corollary 5.1.11** *If $G$ is a bipartite connected graph, then $G$ has a perfect matching if and only if $|A| \leq |N_G(A)|$ for every independent set of vertices $A$ of $G$.*

**Proof.** Notice that the graph $G$ has a perfect matching if and only if the vector $\beta = (1, 1, \ldots, 1)$ is in $\mathbb{N}\mathcal{A}$. By [79, Lemma 2.9] we have the equality $\mathbb{Z}^n \cap \mathbb{R}_+\mathcal{A} = \mathbb{N}\mathcal{A}$. Hence $\beta$ is in $\mathbb{N}\mathcal{A}$ if and only if $\beta \in \mathbb{R}_+\mathcal{A}$. Thus the result follows from Corollary 5.1.10. □

## 5.2 Symbolic Rees algebras of edge ideals

Let $G$ be a graph with vertex set $X = \{x_1, \ldots, x_n\}$ and let $I = I(G)$ be its edge ideal. As before we denote the clutter of minimal vertex covers of $G$ by $\Upsilon(G)$. The clutter $\Upsilon(G)$ is called the *blocker* of $G$. Recall that the symbolic Rees algebra of $I(G)$ is given by

$$R_s(I(G)) = K[x^a t^b |\, a \text{ is an irreducible } b\text{-cover of } \Upsilon(G)] \tag{5.4}$$

and the set $\mathcal{B} = \{x^a t^b |\, a \text{ is an irreducible } b\text{-cover of } \Upsilon(G)\}$ is the minimal set of generators of $R_s(I(G))$ as a $K$-algebra. The main purpose of this section is to study the symbolic Rees algebra of $I(G)$ and to explain the difficulties in finding a combinatorial representation for the minimal set of generators of this algebra.

**Lemma 5.2.1** *Let* $0 \neq a = (a_1, \ldots, a_m, 0, \ldots, 0) \in \mathbb{N}^n$ *and let* $a' = (a_1, \ldots, a_m)$. *If* $0 \neq b \in \mathbb{N}$, *then* $a$ *is an irreducible* $b$*-cover of* $\Upsilon(G)$ *if and only if* $a'$ *is an irreducible* $b$*-cover of* $\Upsilon(\langle S \rangle)$, *where* $S = \{x_1, \ldots, x_m\}$.

**Proof.** It suffices to prove that $a$ is a $b$-cover of the blocker of $G$ if and only if $a'$ is a $b$-cover of the blocker of $\langle S \rangle$.

$\Rightarrow$) The induced subgraph $\langle S \rangle$ is not a discrete graph. Take a minimal vertex cover $C'$ of $\langle S \rangle$. Set $C = C' \cup (V(G) \setminus S)$. Since $C$ is a vertex cover of $G$ such that $C \setminus \{x_i\}$ is not a vertex cover of $G$ for every $x_i \in C'$, there is a minimal vertex cover $C_\ell$ of $G$ such that $C' \subset C_\ell \subset C$ and $C' = C_\ell \cap S$. Notice that

$$\sum_{x_i \in C'} a_i = \sum_{x_i \in C_\ell \cap S} a_i = \langle a, u_\ell \rangle \geq b,$$

where $u_\ell$ is the incidence vector of $C_\ell$. Hence $\sum_{x_i \in C'} a_i \geq b$, as required.

$\Leftarrow$) Take a minimal vertex cover $C_\ell$ of $G$. Then $C_\ell \cap S$ contains a minimal vertex cover $C'_\ell$ of $\langle S \rangle$. Let $u_\ell$ (resp. $u'_\ell$) be the incidence vector of $C_\ell$ (resp. $C'_\ell$). Notice that

$$\langle a, u_\ell \rangle = \sum_{x_i \in C_\ell \cap S} a_i \geq \sum_{x_i \in C'_\ell} a_i = \langle a', u'_\ell \rangle \geq b.$$

Hence $\langle a, u_\ell \rangle \geq b$, as required. $\Box$

We denote a complete subgraph of $G$ with $r$ vertices by $\mathcal{K}_r$. Let $v$ be a vertex of $G$, the *neighbor set* $N_G(v)$ of $v$ is the set of vertices of $G$ adjacent to $v$.

**Lemma 5.2.2** *Let $G$ be a graph and let $a = (a_1, \ldots, a_n)$ be an irreducible $b$-cover of $\Upsilon(G)$ such that $a_i \geq 1$ for all $i$. If $\langle N_G(x_n)\rangle = \mathcal{K}_r$, then $a_i = 1$ for all $i$, $b = r$, $n = r + 1$, and $G = \mathcal{K}_n$.*

**Proof.** We may assume that $N_G(x_n) = \{x_1, \ldots, x_r\}$. We set

$$c = e_1 + \cdots + e_r + e_n; \quad d = (a_1 - 1, \ldots, a_r - 1, a_{r+1}, \ldots, a_{n-1}, a_n - 1).$$

Notice that $\langle x_1, \ldots, x_r, x_n\rangle = \mathcal{K}_{r+1}$. Thus $c$ is an $r$-cover of $\Upsilon(G)$ because any minimal vertex cover of $G$ must intersect all edges of $\mathcal{K}_{r+1}$. By the irreducibility of $a$, there exists a minimal vertex cover $C_\ell$ of $G$ such that $\sum_{x_i \in C_\ell} a_i = b$. Clearly we have $b \geq g \geq r$, where $g$ is the height of $I(G)$. Let $C_k$ be an arbitrary minimal vertex cover of $G$. Since $C_k$ contains exactly $r$ vertices of $\mathcal{K}_{r+1}$, from the inequality $\sum_{x_i \in C_k} a_i \geq b$ we get $\sum_{x_i \in C_k} d_i \geq b - r$, where $d_1, \ldots, d_n$ are the entries of $d$. Therefore $d = 0$; otherwise if $d \neq 0$, then $d$ is a $(b - r)$-cover of $\Upsilon(G)$ and $a = c + d$, a contradiction to the irreducibility of $a$. It follows that $g = r$, $n = r + 1$, $a_i = 1$ for $1 \leq i \leq r$, $a_n = 1$, and $G = \mathcal{K}_n$. □

*Notation* We regard $\mathcal{K}_0$ as the empty set with zero elements. A sum over an empty set is defined to be 0.

**Proposition 5.2.3** *Let $G$ be a graph and let $J = I_c(G)$ be its ideal of vertex covers. Then the set*

$$F = \{(a_i) \in \mathbb{R}^{n+1} |\ \textstyle\sum_{x_i \in \mathcal{K}_r} a_i = (r - 1)a_{n+1}\} \cap \mathbb{R}_+(J)$$

*is a facet of the Rees cone $\mathbb{R}_+(J)$.*

**Proof.** If $\mathcal{K}_r = \emptyset$, then $r = 0$ and $F = H_{e_{n+1}} \cap \mathbb{R}_+(J)$, which is clearly a facet because $e_1, \ldots, e_n \in F$. If $r = 1$, then $F = H_{e_i} \cap \mathbb{R}_+(J)$ for some $1 \leq i \leq n$, which is a facet because $e_j \in F$ for $j \notin \{i, n + 1\}$ and there is at least one minimal vertex cover of $G$ not containing $x_i$. We may assume that $X' = \{x_1, \ldots, x_r\}$ is the vertex set of $\mathcal{K}_r$ and $r \geq 2$. For each $1 \leq i \leq r$ there is a minimal vertex cover $C_i$ of $G$ not containing $x_i$. Notice that $C_i$ contains $X' \setminus \{x_i\}$. Let $u_i$ be the incidence vector of $C_i$. Since the rank of $u_1, \ldots, u_r$ is $r$, it follows that the set

$$\{(u_1, 1), \ldots, (u_r, 1), e_{r+1}, \ldots, e_n\}$$

is a linearly independent set contained in $F$, i.e., $\dim(F) = n$. Hence $F$ is a facet of $\mathbb{R}_+(J)$ because the hyperplane that defines $F$ is a supporting hyperplane. □

**Proposition 5.2.4** *Let $G$ be a graph and let $0 \neq a = (a_1, \ldots, a_n) \in \mathbb{N}^n$. If*

(a) $a_i \in \{0, 1\}$ *for all $i$, and*

(b) $\langle \{x_i | \, a_i > 0\} \rangle = \mathcal{K}_{r+1}$,

*then $a$ is an irreducible $r$-cover of $\Upsilon(G)$.*

**Proof.** By Proposition 5.2.3, the closed halfspace $H^+_{(a,-r)}$ occurs in the irreducible representation of the Rees cone $\mathbb{R}_+(J)$, where $J = I_c(G)$. Hence $a$ is an irreducible $r$-cover by Lemma 5.1.3. □

A *clique* of a graph $G$ is a set of vertices that induces a complete subgraph. We will also call a complete subgraph of $G$ a clique. Symbolic Rees algebras are related to perfect graphs as is seen below. We refer to [16, 20, 71] and the references there for the theory of perfect graphs.

*Notation* The *support* of $x^a = x_1^{a_1} \cdots x_n^{a_n}$ is $\mathrm{supp}(x^a) = \{x_i \, | \, a_i > 0\}$.

**Corollary 5.2.5** ([91]) *If $G$ is a graph, then*

$$K[x^a t^r | \, x^a \ \textit{square-free} \ ; \langle \mathrm{supp}(x^a) \rangle = \mathcal{K}_{r+1}; \, 0 \leq r < n] \subset R_s(I(G))$$

*with equality if and only if $G$ is a perfect graph.*

**Proof.** The inclusion follows from Proposition 5.2.4. If $G$ is a perfect graph, then by [91, Corollary 3.3] the equality holds. Conversely if the equality holds, then by Lemma 5.1.3 and Proposition 5.2.3 we have

$$\mathbb{R}_+(I_c(G)) = \left\{(a_i) \in \mathbb{R}^{n+1} | \, \textstyle\sum_{x_i \in \mathcal{K}_r} a_i \geq (r-1)a_{n+1}; \, \forall \, \mathcal{K}_r \subset G \right\}. \tag{5.5}$$

Hence a direct application of [91, Proposition 2.2] gives that $G$ is perfect. □

The *vertex covering number* of $G$, denoted by $\alpha_0(G)$, is the number of vertices in a minimum vertex cover of $G$ (the cardinality of any smallest vertex cover in $G$). Notice that $\alpha_0(G)$ equals the height of $I(G)$. If $H$ is a discrete graph, i.e., all the vertices of $H$ are isolated, we set $I(H) = 0$ and $\alpha_0(H) = 0$.

**Lemma 5.2.6** *Let $G$ be a graph. If $a = e_1 + \cdots + e_r$ is an irreducible $b$-cover of $\Upsilon(G)$, then $b = \alpha_0(H)$, where $H = \langle x_1, \ldots, x_r \rangle$.*

**Proof.** The case $b = 0$ is clear. Assume $b \geq 1$. Let $C_1, \ldots, C_s$ be the minimal vertex covers of $G$ and let $u_1, \ldots, u_s$ be their incidence vectors. Notice that $\langle a, u_i \rangle = b$ for some $i$. Indeed if $\langle a, u_i \rangle > b$ for all $i$, then $a - e_1$ is a $b$-cover of $\Upsilon(G)$ and $a = (a - e_1) + e_1$, a contradiction. Hence

$$b = \langle a, u_i \rangle = |\{x_1, \ldots, x_r\} \cap C_i| \geq \alpha_0(H).$$

This proves that $b \geq \alpha_0(H)$. Notice that $H$ is not a discrete graph. Then we can pick a minimal vertex cover $A$ of $H$ such that $|A| = \alpha_0(H)$. The set

$$C = A \cup (V(G) \setminus \{x_1, \ldots, x_r\})$$

is a vertex cover of $G$. Hence there is a minimal vertex cover $C_\ell$ of $G$ such that $A \subset C_\ell \subset C$. Observe that $C_\ell \cap \{x_1, \ldots, x_r\} = A$. Thus we get $\langle a, u_\ell \rangle = |A| \geq b$, i.e., $\alpha_0(H) \geq b$. Altogether we have $b = \alpha_0(H)$. □

**Theorem 5.2.7** *Let $G$ be a graph and let $a = (1, \ldots, 1)$. Then $a$ is a reducible $\alpha_0(G)$-cover of $\Upsilon(G)$ if and only if there are $H_1$ and $H_2$ induced subgraphs of $G$ such that*

(i) *$V(G)$ is the disjoint union of $V(H_1)$ and $V(H_2)$, and*

(ii) *$\alpha_0(G) = \alpha_0(H_1) + \alpha_0(H_2)$.*

**Proof.** $\Rightarrow$) We may assume that $a_1 = e_1 + \cdots + e_r$, $a_2 = a - a_1$, $a_i$ is a $b_i$-cover of $\Upsilon(G)$, $b_i \geq 1$ for $i = 1, 2$, and $\alpha_0(G) = b_1 + b_2$. Consider the subgraphs $H_1 = \langle x_1, \ldots, x_r \rangle$ and $H_2 = \langle x_{r+1}, \ldots, x_n \rangle$. Let $A$ be a minimal vertex cover of $H_1$ with $\alpha_0(H_1)$ vertices. Since $C = A \cup (V(G) \setminus \{x_1, \ldots, x_r\})$ is a vertex cover $G$, there is a minimal vertex cover $C_k$ of $G$ such that $A \subset C_k \subset C$. Hence

$$|A| = |C_k \cap \{x_1, \ldots, x_r\}| = \langle a_1, u_k \rangle \geq b_1,$$

and $\alpha_0(H_1) \geq b_1$. Using a similar argument we get that $\alpha_0(H_2) \geq b_2$. If $C_\ell$ is a minimal vertex cover of $G$ with $\alpha_0(G)$ vertices, then $C_\ell \cap V(H_i)$ is a vertex cover of $H_i$. Therefore

$$b_1 + b_2 = \alpha_0(G) = |C_\ell| = \sum_{i=1}^{2} |C_\ell \cap V(H_i)| \geq \sum_{i=1}^{2} \alpha_0(H_i) \geq b_1 + b_2,$$

and consequently $\alpha_0(G) = \alpha_0(H_1) + \alpha_0(H_2)$.

$\Leftarrow$) We may assume $V(H_1) = \{x_1, \ldots, x_r\}$ and $V(H_2) = V(G) \setminus V(H_1)$. Set $a_1 = e_1 + \cdots + e_r$ and $a_2 = a - a_1$. For any minimal vertex cover $C_k$ of $G$, we have that $C_k \cap V(H_i)$ is a vertex cover of $H_i$. Hence

$$\langle a_1, u_k \rangle = |C_k \cap \{x_1, \ldots, x_r\}| \geq \alpha_0(H_1),$$

where $u_k$ is the incidence vector of $C_k$. Consequently $a_1$ is an $\alpha_0(H_1)$-cover of $\Upsilon(G)$. Similarly we obtain that $a_2$ is an $\alpha_0(H_2)$-cover of $\Upsilon(G)$. Therefore $a$ is a reducible $\alpha_0(G)$-cover of $\Upsilon(G)$. □

**Definition 5.2.8** A graph satisfying conditions (i) and (ii) is called a *reducible* graph. If $G$ is not reducible, it is called *irreducible.*

These notions appear in [27]. As far as we know there is no structure theorem for irreducible graphs. Examples of irreducible graphs include complete graphs, odd cycles, and complements of odd cycles. Below we give a method, using Hilbert bases, to obtain all irreducible induced subgraphs of $G$.

By [60, Lemma 5.4] there exists a finite set $\mathcal{H} \subset \mathbb{N}^{n+1}$ such that

(a) $\mathrm{Cn}(I(G)) = \mathbb{R}_+\mathcal{H}$, and

(b) $\mathbb{Z}^{n+1} \cap \mathbb{R}_+\mathcal{H} = \mathbb{N}\mathcal{H}$,

where $\mathbb{N}\mathcal{H}$ is the additive subsemigroup of $\mathbb{N}^{n+1}$ generated by $\mathcal{H}$.

**Definition 5.2.9** The set $\mathcal{H}$ is called an integral *Hilbert basis* of $\mathrm{Cn}(I(G))$.

If we require $\mathcal{H}$ to be minimal (with respect inclusion), then $\mathcal{H}$ is unique [69].

**Corollary 5.2.10** *Let $G$ be a graph and let $\alpha = (a_1, \ldots, a_n, b) \in \{0,1\}^n \times \mathbb{N}$. Then $\alpha$ is an element of the minimal integral Hilbert basis of* $\mathrm{Cn}(I(G))$ *if and only if the induced subgraph $H = \langle \{x_i | \, a_i = 1\} \rangle$ is irreducible with $b = \alpha_0(H)$.*

**Proof.** The map $(a_1, \ldots, a_n, b) \mapsto x_1^{a_1} \cdots x_n^{a_n} t^b$ establishes a one to one correspondence between the minimal integral Hilbert basis of $\mathrm{Cn}(I(G))$ and the minimal generators of $R_s(I(G))$ as a $K$-algebra. Thus the result follows from Lemma 5.2.1 and Theorem 5.2.7. □

**Example 5.2.11** (E. Reyes) Consider the graph $G$ shown below. Let $I$ be the edge ideal of $G$ and let $\mathcal{H}$ be the Hilbert basis of $\mathrm{Cn}(I)$. Using the map $(a_1,\ldots,a_n,b) \mapsto x_1^{a_1}\cdots x_n^{a_n}t^b$, together with *Normaliz* [11], it is seen that $G$ has exactly 49 irreducible subgraphs. Since $\alpha_0(G)=6$ and the vector $(1,\ldots,1,6)$ is not in $\mathcal{H}$ we obtain that $G$ is a reducible graph.

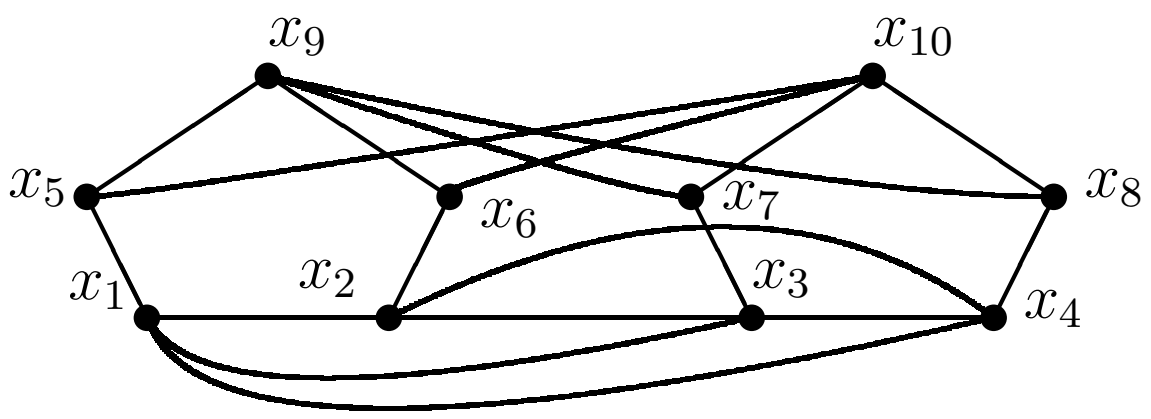

Fig. 4. Graph $G$

The following parallelization $G^{(1,\ldots,1,2)}$ of $G$ is irreducible because $(1,\ldots,1,2)\in\mathcal{H}$.

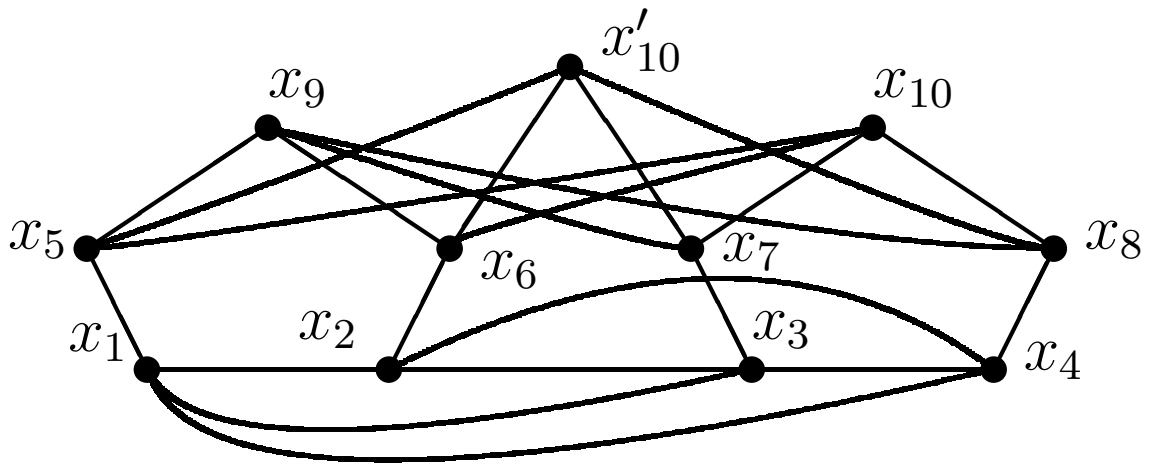

Fig. 5. Graph $G^{(1,\ldots,1,2)}$

The next result shows that irreducible graphs occur naturally in the theory of perfect graphs.

**Proposition 5.2.12** *A graph $G$ is perfect if and only if the only irreducible induced subgraphs of $G$ are the complete subgraphs.*

**Proof.** $\Rightarrow$) Let $H$ be an irreducible induced subgraph of $G$. We may assume that $V(H)=\{x_1,\ldots,x_r\}$. Set $a'=(1,\ldots,1)\in\mathbb{N}^r$ and $a=(a',0\ldots,0)\in\mathbb{N}^n$. By Theorem 5.2.7, $a'$ is an irreducible $\alpha_0(H)$ cover of $\Upsilon(H)$. Then by Lemma 5.2.1, $a'$ is an irreducible $\alpha_0(H)$ cover of $\Upsilon(G)$. Since $x_1\cdots x_r t^{\alpha_0(H)}$ is a minimal generator of $R_s(I(G))$, using Corollary 5.2.5 we obtain that $\alpha_0(H)=r-1$ and that $H$ is a complete subgraph of $G$ on $r$ vertices.

$\Leftarrow$) In [13] it is shown that $G$ is a perfect graph if and only if no induced subgraph of $G$ is an odd cycle of length at least five or the complement of one. Since odd cycles and and their complements are irreducible subgraphs. It follows that $G$ is perfect. $\Box$

**Definition 5.2.13** A graph $G$ is called *vertex critical* if $\alpha_0(G\setminus\{x_i\}) < \alpha_0(G)$ for all $x_i \in V(G)$.

**Remark 5.2.14** If $x_i$ is any vertex of a graph $G$ and $\alpha_0(G\setminus\{x_i\}) < \alpha_0(G)$, then $\alpha_0(G\setminus\{x_i\}) = \alpha_0(G) - 1$

**Lemma 5.2.15** *If $G$ is an irreducible graph, then $G$ is a connected graph and $\alpha_0(G\setminus\{x_i\}) = \alpha_0(G) - 1$ for $i = 1,\ldots,n$.*

**Proof.** Let $G_1,\ldots,G_r$ be the connected components of $G$. Since $\alpha_0(G)$ is equal to $\sum_i \alpha_0(G_i)$, we get $r = 1$. Thus $G$ is connected. To complete the proof it suffices to prove that $\alpha_0(G\setminus\{x_i\}) < \alpha_0(G)$ for all $i$ (see Remark 5.2.14). If $\alpha_0(G\setminus\{x_i\}) = \alpha_0(G)$, then $G = H_1 \cup H_2$, where $H_1 = G\setminus\{x_i\}$ and $V(H_2) = \{x_i\}$, a contradiction. $\Box$

**Definition 5.2.16** *The cone $C(G)$, over a graph $G$, is obtained by adding a new vertex $v$ to $G$ and joining every vertex of $G$ to $v$.*

The next result can be used to build irreducible graphs. In particular it follows that cones over irreducible graphs are irreducible.

**Proposition 5.2.17** *Let $G$ be a graph with $n$ vertices and let $H$ be a graph obtained from $G$ by adding a new vertex $v$ and some new edges joining $v$ with $V(G)$. If $a = (1,\ldots,1) \in \mathbb{N}^n$ is an irreducible $\alpha_0(G)$-cover of $\Upsilon(G)$ such that $\alpha_0(H) = \alpha_0(G) + 1$, then $a' = (a, 1)$ is an irreducible $\alpha_0(H)$-cover of $\Upsilon(H)$.*

**Proof.** Clearly $a'$ is an $\alpha_0(H)$-cover of $\Upsilon(H)$. Assume that $a' = a_1' + a_2'$, where $0 \neq a_i'$ is a $b_i'$-cover of $\Upsilon(H)$ and $b_1' + b_2' = \alpha_0(H)$. We may assume that $a_1' = (1,\ldots,1,0,\ldots,0)$ and $a_2' = (0,\ldots,0,1,\ldots,1)$. Let $a_i$ be the vector in $\mathbb{N}^n$ obtained from $a_i'$ by removing its last entry. Set $v = x_{n+1}$. Take a minimal vertex cover $C_k$ of $G$ and consider $C_k' = C_k \cup \{x_{n+1}\}$. Let $u_k'$ (resp. $u_k$) be the incidence vector of $C_k'$ (resp. $C_k$). Then

$$\langle a_1, u_k\rangle = \langle a_1', u_k'\rangle \geq b_1' \text{ and } \langle a_2, u_k\rangle + 1 = \langle a_2', u_k'\rangle \geq b_2',$$

consequently $a_1$ is a $b'_1$-cover of $\Upsilon(G)$. If $b'_2 = 0$, then $a_1$ is an $\alpha_0(H)$-cover of $\Upsilon(G)$, a contradiction; because if $u$ is the incidence vector of a minimal vertex cover of $G$ with $\alpha_0(G)$ elements, then we would obtain $\alpha_0(G) \geq \langle u, a_1 \rangle \geq \alpha_0(H)$, which is impossible. Thus $b'_2 \geq 1$, and $a_2$ is a $(b'_2 - 1)$-cover of $\Upsilon(G)$ if $a_2 \neq 0$. Hence $a_2 = 0$, because $a = a_1 + a_2$ and $a$ is irreducible. This means that $a'_2 = e_{n+1}$ is a $b'_2$-cover of $\Upsilon(H)$, a contradiction. Therefore $a'$ is an irreducible $\alpha_0(H)$-cover of $\Upsilon(H)$, as required. □

**Definition 5.2.18** A graph $G$ is called *edge critical* if $\alpha_0(G \setminus e) < \alpha_0(G)$ for all $e \in E(G)$.

**Proposition 5.2.19** *If $G$ is a connected edge critical graph, then $G$ is irreducible.*

**Proof.** Assume that $G$ is reducible. Then there are induced subgraphs $H_1$, $H_2$ of $G$ such that $V(H_1)$, $V(H_2)$ is a partition of $V(G)$ and $\alpha_0(G) = \alpha_0(H_1) + \alpha_0(H_2)$. Since $G$ is connected there is an edge $e = \{x_i, x_j\}$ with $x_i$ a vertex of $H_1$ and $x_j$ a vertex of $H_2$. Pick a minimal vertex cover $C$ of $G \setminus e$ with $\alpha_0(G) - 1$ vertices. As $E(H_i)$ is a subset of $E(G \setminus e) = E(G) \setminus \{e\}$ for $i = 1, 2$, we get that $C$ covers all edges of $H_i$ for $i = 1, 2$. Hence $C$ must have at least $\alpha_0(G)$ elements, a contradiction. □

**Corollary 5.2.20** *The following implications hold for any connected graph*:

$$\textit{edge critical} \implies \textit{irreducible} \implies \textit{vertex critical.}$$

**Finding generators of symbolic Rees algebras using cones** In [1] Bahiano showed that if $H$ is the graph obtained by taking a cone over a pentagon, then

$$R_s(I(H)) = R[I(H)t][x_1 \cdots x_5 t^3, x_1 \cdots x_6 t^4, x_1 \cdots x_5 x_6^2 t^5].$$

This simple example shows that taking a cone over an irreducible graph tends to increase the degree in $t$ of the generators of the symbolic Rees algebra. Other examples using this "cone process" have been shown in [51, Example 5.5].

Let $G$ be a graph with vertex set $V(G) = \{x_1, \ldots, x_n\}$. The aim here is to give a general procedure—based on the irreducible representation of the Rees cone of $I_c(G)$—to construct generators of $R_s(I(H))$ of high degree in $t$, where $H$ is a graph constructed from $G$ by recursively taking cones over graphs already constructed.

By the finite basis theorem [93, Theorem 4.11] there is a unique irreducible representation

$$\mathbb{R}_+(I_c(G)) = H_{e_1}^+ \cap H_{e_2}^+ \cap \cdots \cap H_{e_{n+1}}^+ \cap H_{\alpha_1}^+ \cap H_{\alpha_2}^+ \cap \cdots \cap H_{\alpha_p}^+ \tag{5.6}$$

such that each $\alpha_k$ is in $\mathbb{Z}^{n+1}$, the non-zero entries of each $\alpha_k$ are relatively prime, and none of the closed halfspaces $H_{e_1}^+, \ldots, H_{e_{n+1}}^+, H_{\alpha_1}^+, \ldots, H_{\alpha_p}^+$ can be omitted from the intersection. For use below we assume that $\alpha$ is any of the vectors $\alpha_1, \ldots, \alpha_p$ that occur in the irreducible representation. Thus we can write $\alpha = (a_1, \ldots, a_n, -b)$ for some $a_i \in \mathbb{N}$ and for some $b \in \mathbb{N}$.

**Lemma 5.2.21** *Let $H$ be the cone over $G$. If*

$$\beta = (a_1, \ldots, a_n, (\textstyle\sum_{i=1}^n a_i) - b, -\sum_{i=1}^n a_i) = (\beta_1, \ldots, \beta_{n+1}, -\beta_{n+2})$$

*and $a_i \geq 1$ for all $i$, then $F = H_\beta \cap \mathbb{R}_+(I_c(H))$ is a facet of $\mathbb{R}_+(I_c(H))$.*

**Proof.** First we prove that $\mathbb{R}_+(I_c(H)) \subset H_\beta^+$, i.e., $H_\beta$ is a supporting hyperplane of the Rees cone. By Lemma 5.1.3, $(a_1, \ldots, a_n)$ is an irreducible $b$-cover of $\Upsilon(G)$. Hence there is $C \in \Upsilon(G)$ such that $\sum_{x_i \in C} a_i = b$. Therefore $\beta_{n+1}$ is greater or equal than 1. This proves that $e_1, \ldots, e_{n+1}$ are in $H_\beta^+$. Let $C$ be any minimal vertex cover of $H$ and let $u = \sum_{x_i \in C} e_i$ be its characteristic vector. Case (i): If $x_{n+1} \notin C$, then $C = \{x_1, \ldots, x_n\}$ and

$$\sum_{x_i \in C} \beta_i = \sum_{i=1}^n a_i = \beta_{n+2},$$

that is, $(u, 1) \in H_\beta^+$. Case (ii): If $x_{n+1} \in C$, then $C_1 = C \setminus \{x_{n+1}\}$ is a minimal vertex cover of $G$. Hence

$$\sum_{x_i \in C} \beta_i = \sum_{x_i \in C_1} \beta_i + \beta_{n+1} \geq b + \beta_{n+1} = \beta_{n+2},$$

that is, $(u, 1) \in H_\beta^+$. Therefore $\mathbb{R}_+(I_c(H)) \subset H_\beta^+$. To prove that $F$ is a facet we must show it has dimension $n+1$ because the dimension of $\mathbb{R}_+(I_c(H))$ is $n+2$. We denote the characteristic vector of a minimal vertex cover $C_k$ of $G$ by $u_k$. By hypothesis there are minimal vertex covers $C_1, \ldots, C_n$ of $G$ such that the vectors $(u_1, 1), \ldots, (u_n, 1)$ are linearly independent and

$$\langle (a, -b), (u_k, 1) \rangle = 0 \iff \langle a, u_k \rangle = b, \tag{5.7}$$

for $k = 1, \ldots, n$. Therefore

$$\langle(\beta_1, \ldots, \beta_{n+1}), (u_k, 1)\rangle = \beta_{n+2} \quad \text{and} \quad \langle(\beta_1, \ldots, \beta_{n+1}), (1, \ldots, 1, 0)\rangle = \beta_{n+2},$$

i.e., the set $\mathcal{B} = \{(u_1, 1), \ldots, (u_n, 1), (1, \ldots, 1, 0\}$ is contained in $H_\beta$. Since

$$C_1 \cup \{x_{n+1}\}, \ldots, C_n \cup \{x_{n+1}\}, \{x_1, \ldots, x_n\}$$

are minimal vertex covers of $H$, the set $\mathcal{B}$ is also contained in $\mathbb{R}_+(I_c(H))$ and consequently in $F$. Thus its suffices to prove that $\mathcal{B}$ is linearly independent. If $(1, \ldots, 1, 0)$ is a linear combination of $(u_1, 1), \ldots, (u_n, 1)$, then we can write

$$(1, \ldots, 1) = \lambda_1 u_1 + \cdots + \lambda_n u_n$$

for some scalars $\lambda_1, \ldots, \lambda_n$ such that $\sum_{i=1}^n \lambda_i = 0$. Hence from Eq. (5.7) we get

$$|a| = \langle(1, \ldots, 1), a\rangle = \lambda_1\langle u_1, a\rangle + \cdots + \lambda_n\langle u_n, a\rangle = (\lambda_1 + \cdots + \lambda_n)b = 0,$$

a contradiction. □

**Corollary 5.2.22** *If $a_i \geq 1$ for all $i$, then $x_1^{\beta_1} \cdots x_{n+1}^{\beta_{n+1}} t^{\beta_{n+2}}$ is a minimal generator of $R_s(I(H))$.*

**Proof.** By Lemma 5.2.21, $F = H_\beta \cap \mathbb{R}_+(I_c(H))$ is a facet of $\mathbb{R}_+(I_c(H))$. Therefore using Lemma 5.1.3, the vector $(\beta_1, \ldots, \beta_{n+1})$ is an irreducible $\beta_{n+2}$-cover of $\Upsilon(H)$, i.e., $x_1^{\beta_1} \cdots x_{n+1}^{\beta_{n+1}} t^{\beta_{n+2}}$ is a minimal generator of $R_s(I(H))$. □

**Corollary 5.2.23** *Let $G_0 = G$ and let $G_r$ be the cone over $G_{r-1}$ for $r \geq 1$. If $\alpha = (1, \ldots, 1, -g)$, then*

$$(\underbrace{1, \ldots, 1}_{n}, \underbrace{n-g, \ldots, n-g}_{r})$$

*is an irreducible $n + (r-1)(n-g)$ cover of $G_r$. In particular $R_s(I(G_r))$ has a generator of degree in $t$ equal to $n + (r-1)(n-g)$.*

As a very particular example of our construction consider:

**Example 5.2.24** Let $G = C_s$ be an odd cycle of length $s = 2k+1$. Note that $\alpha_0(C_s) = (s+1)/2 = k+1$. Then by Corollary 5.2.23

$$x_1 \cdots x_s x_{s+1}^k \cdots x_{s+r}^k t^{rk+k+1}$$

is a minimal generator of $R_s(I(G_r))$. This illustrates that the degree in $t$ of the minimal generators of $R_s(I(G_r))$ is much larger than the number of vertices of the graph $G_r$ [51].

# Index of Notation

# Index